\documentclass[10pt]{amsart}
\usepackage{latexsym,amsmath,amssymb,graphicx}
\usepackage[all,web]{xy}
\usepackage{maplestd2e}

\DeclareFontFamily{OT1}{rsfs10}{}
\DeclareFontShape{OT1}{rsfs10}{m}{n}{ <-> rsfs10 }{}
\DeclareMathAlphabet{\mathscript}{OT1}{rsfs10}{m}{n}

\DeclareMathOperator{\Spec}{Spec}   
\DeclareMathOperator{\Ext}{Ext}     
\DeclareMathOperator{\ext}{\mathcal{E}\it{x}\hskip1pt \it{t}\hskip1pt} 
\DeclareMathOperator{\hm}{\mathcal{H}\it{o}\hskip1pt \it{m}\hskip1pt} 
\DeclareMathOperator{\Aut}{Aut}     
\DeclareMathOperator{\rk}{rk}       
\DeclareMathOperator{\Sing}{Sing}   
\DeclareMathOperator{\cone}{Cn}     
\DeclareMathOperator{\crk}{crk}     
\DeclareMathOperator{\Hess}{Hess}    

\title[Milnor and Tyurina numbers of hypersurface singularities]{MAPLE
subroutines for computing Milnor and Tyurina numbers of
hypersurface singularities with application to Arnol'd adjacencies.}

\author[M. Rossi and L.Terracini]{Michele Rossi and Lea Terracini}

\address{Dipartimento di Matematica, Universit\`a di Torino,
via Carlo Alberto 10, 10123 Torino} \email{michele.rossi@unito.it,
lea.terracini@unito.it}

\thanks{This work has been developed despite the effects of the Italian law 133/08 (http://groups.google.it/group/scienceaction ).
        This law drastically reduces public funds to public Italian universities, which is particularly dangerous for scientific free research,
        and it will prevent young researchers from getting a position, either temporary or tenured, in Italy.
        The authors are protesting against this law to obtain its repeal.}

\def\p2{\mathbb{P}^2}
\def\p3{\mathbb{P}^3}
\def\p4{\mathbb{P}^4}

\def\rk{\operatorname{rk}}

\def\GL{\operatorname{GL}}

\def\C{\mathbb{C}}
\def\R{\mathbb{R}}

\def\Q{\mathbb{Q}}
\def\N{\mathbb{N}}
\def\T{\mathbb{T}}

\def\Cx{\mathbb{C}[\mathbf{x}]}
\def\Cxx{\mathbb{C}\{\mathbf{x}\}}

\theoremstyle{plain}
\newtheorem{theorem}{Theorem}[section]
\newtheorem{proposition}[theorem]{Proposition}
\newtheorem{thm-def}[theorem]{Theorem--Definition}
\newtheorem{corollary}[theorem]{Corollary}

\newtheorem{lemma}[theorem]{Lemma}

\theoremstyle{remark}
\newtheorem{remark}[theorem]{Remark}

\newtheorem{example}[theorem]{Example}

\newtheorem*{caveat}{Caveat Lector}

\theoremstyle{definition}
\newtheorem{definition}[theorem]{Definition}

\newtheorem*{step I}{Step I}
\newtheorem*{step II}{Step II}
\newtheorem*{step III}{Step III}
\newtheorem*{step IV}{Step IV}

\newtheorem*{acknowledgements}{Acknowledgements}

\newcommand{\oneline}{\vskip12pt}
\newcommand{\halfline}{\vskip6pt}
\begin{document}

 \DefineParaStyle{Maple Heading 4}
 \DefineParaStyle{Maple Heading 2}
 \DefineParaStyle{Maple Text Output}
 \DefineParaStyle{Maple Bullet Item}
 \DefineParaStyle{Maple Warning}
 \DefineParaStyle{Maple Error}
 \DefineParaStyle{Maple Dash Item}
 \DefineParaStyle{Maple Heading 3}
 \DefineParaStyle{Maple Heading 1}
 \DefineParaStyle{Maple Title}
 \DefineParaStyle{Maple Normal}
 \DefineCharStyle{Maple 2D Input}
 \DefineCharStyle{Maple Maple Input}
 \DefineCharStyle{Maple 2D Output}
 \DefineCharStyle{Maple 2D Math}
 \DefineCharStyle{Maple Hyperlink}

\begin{abstract}
In the present paper MAPLE subroutines computing Milnor and
Tyurina numbers of an isolated algebraic hypersurface singularity
are presented and described in detail. They represents examples,
and perhaps the first ones, of a MAPLE implementation of
\emph{local monomial ordering}.

\noindent As an application, the last section is devoted to writing down equations of algebraic stratifications of Kuranishi spaces of simple Arnol'd singularities: they geometrically represents, by means of inclusions of algebraic subsets, the partial ordering on classes of simple singularities induced by the \emph{adjacency} relation.
\end{abstract}

\maketitle

\tableofcontents

Two basic invariants of a complex analytic complete intersection
singularity $p\in \overline{U}=\mathbf{f}^{-1}(0)$ are its
\emph{Milnor number} $\mu(p)$ and its \emph{Tyurina number}
$\tau(p)$. The former essentially "counts the number" of vanishing
cycles in the intermediate cohomology of a nearby smoothing
$U_t=\mathbf{f}^{-1}(t)$ of $\overline{U}$, which actually turns
out to be the \emph{multiplicity of $p$ as a critical point} of
the map $\mathbf{f}$. The latter "counts the dimension" of the
base space of a versal deformation of $p\in \overline{U}$, which
actually turn out to be the \emph{multiplicity of $p$ as a
singularity} of the complex space $\overline{U}$. After the
Looijenga--Steenbrink Theorem \cite{LS} it is a well known fact
that $\tau(p)\leq \mu(p)$.

\halfline The purpose of the present paper is to present MAPLE
subroutines \cite{Maple} allowing to compute those invariants in
the case of an isolated algebraic hypersurface singularity
(i.h.s.). In fact, from the computational point of view, their
calculation could be very intricate to end up and a computer
employment may be needed in most situations. Let us underline that
the actual originality of our procedure is not so much based on
the effective computation of those invariants as on its
implementation with a very interdisciplinary and worldwide
diffused and known math software like MAPLE is. In fact computer
algebra packages computing those singularities' invariants already
exists (one for all is SINGULAR
\cite{Singular}). But our hope is that routines here presented
could be useful to all those who are interesting, for any reason,
in a concrete evaluation of those invariants without being so
 much motivated for learning an entire computer algebra
package. Once implemented, the present routine is so easy that
even an undergraduate student may use it!

In other words we believe that the present routine could be an
interesting and, as far as we know, the first example of a MAPLE
implementation of \emph{local monomial ordering} (l.m.o.). In fact
term orders, i.e. usual monomial ordering (recently called
\emph{global} in contrast with the word \emph{local}), are already
implemented in MAPLE with the command \texttt{MonomialOrder} in
the \texttt{Groebner} package. Actually user-defined term orders
are also allowed and, in particular, l.m.o.'s can be easily
defined as the \emph{opposite} of a standard g.m.o. (pure
lexicographic, graded lexicographic, reverse ... etc.). The
problem is that Buchberger S-procedure may not end up determining
a \emph{normal form} since a l.m.o. is not a well-order on the
contrary of g.m.o.'s. Then the present MAPLE routines start with an implementation of the Mora algorithm
for determining a
\emph{weak normal form} and then monomial bases of Milnor and
Tyurina ideals in the complex ring of convergent power series
$\C\{x_1,\ldots ,x_{n+1}\}$ (see \cite{Mora}, \cite{Greuel-Pfister} Algorithm 1.7.6, \cite{Decker-Lossen} Algorithm 9.22).
As a consequence a monomial basis of the Kuranishi space, parameterizing small versal
deformations of the given i.h.s., is obtained, allowing to concretely write
down these deformations even for more intricate cases. Actually we get procedures which are able to perform calculation in $(\C[\lambda])[\mathbf{x}]$, where $\lambda=(\lambda_1,\ldots,\lambda_r)$ is an $r$-tuple of parameters e.g. coordinates of the Kuranishi space.

\halfline An interesting application of this latter feature is that of \emph{writing down equations} of algebraic stratifications in Kuranishi spaces of Arnol'd simple singularities, giving an \emph{explicit} geometric interpretation, by means of inclusions of algebraic subsets, of \emph{Arnol'd's adjacency partial order relation} over classes of
simple singularities (see Section \ref{Arnold-ex}). This section ends up with an explicit list of the most specialized 1-parameter deformations of simple singularities realizing adjacencies between distinct classes of simple singularities (see \ref{speciali}).

\begin{acknowledgements} We are greatly indebt with G.~M.~Greuel
who timely pointed out to us serious mistakes in the first version
of the this note. His concise and sharp remarks considerably
helped us to improve our routine and the final product. Authors
would also like to thank A.~Albano for enlightening conversations.
\end{acknowledgements}

\section{Milnor number of an isolated hypersurface
singularity}

From the topological point of view a \emph{good representative}
$\overline{U}$ of an \emph{isolated hypersurface singularity
(i.h.s.)} is the zero locus of a holomorphic map
\begin{equation}\label{polinomio}
    f:\C^{n+1}\longrightarrow\C\quad,\quad n\geq 0
\end{equation}
admitting an isolated critical point in $0\in\C^4$.

\noindent Set $U_T:=f^{-1}(T)$ where $T$ is a small enough
neighborhood of $0\in \C$. Then we can assume that \emph{$f$ is a
submersion over $U_T\setminus\{0\}$}: therefore
$\overline{U}=f^{-1}(0)=U_0$ and $U_t$ is a \emph{local smoothing}
of $\overline{U}$, for $0\neq t\in T$.

\begin{definition} Let $X\subset\C^m$ be a subset. The following subset of $\C^m$
\[
\cone _0 (X) := \left\{ tx\ |\ \forall t\in [0,1]\subset\R\ ,\
\forall x\in X \right\}
\]
will be called the \emph{cone projecting $X$}.
\end{definition}

\begin{theorem}[Local topology of a isolated hypersurface
singularity, \cite{Milnor68} Theorem 2.10, Theorem
5.2]\label{i.h.s.topologia} Let $D_{\varepsilon}$ denote the
closed ball of radius $\varepsilon>0$, centered in $0\in\C^{n+1}$,
whose boundary is the $2n+1$--dimensional sphere
$S_{\varepsilon}^{2n+1}$. Then, for $\varepsilon$ small enough,
the intersection
\[
\overline{B}_{\varepsilon}:=\overline{U}\cap D_{\varepsilon}
\]
is homeomorphic to the cone $\cone _0 (K)$ projecting
$K:=\overline{U}\cap S_{\varepsilon}^{2n+1}$, which is called the
\emph{knot} or \emph{link of the singularity $0\in \overline{U}$}.
\end{theorem}

\begin{thm-def}[Local homology type of the smoothing \cite{Milnor68},
Theorems 5.11, 6.5, 7.2]\label{milnor} Set $\widetilde{U}:=U_{t}$
for some $0\neq t\in T$. Then, for $\varepsilon$ small enough, the
intersection
\[
\widetilde{B}_{\varepsilon}:=\widetilde{U} \cap D_{\varepsilon}
\]
(called the \emph{Milnor fibre} of $f$) has the homology type of a
bouquet of $n$--dimensional spheres. In particular the $n$--th
Betti number $b_n(\widetilde{B}_{\varepsilon})$ (called the
\emph{Milnor number} $m_p$ of $p$) coincides with the
\emph{multiplicity} of the critical point $0\in \C^{n+1}$ of $f$
as a solution to the following collection of equations
\[
\frac{\partial f}{\partial x_1}=\frac{\partial f}{\partial
x_2}=\cdots =\frac{\partial f}{\partial x_{n+1}}=0\ .
\]
\end{thm-def}

\subsection{Milnor number from the algebraic point of view}
\label{m-t ihs}

Theorem \ref{milnor} allows the following algebraic interpretation
of the Milnor number.

Let $\mathcal{O}_0$ be the local ring of germs of holomorphic
function of $\C^{n+1}$ at the origin. By definition of holomorphic
function and the identity principle we have that $\mathcal{O}_0$
\emph{is isomorphic to the ring of convergent power series}
$\C\{x_1,\ldots,x_{n+1}\}$. A \emph{germ of hypersurface
singularity} is defined as the Stein complex space
\begin{equation}\label{germe}
    U_0:=\Spec(\mathcal{O}_{f,0})
\end{equation}
where $\mathcal{O}_{f,0}:=\mathcal{O}_0/(f)$ and $f$ is the germ
represented by the series expansion of the holomorphic function
(\ref{polinomio}).

\begin{definition}[Milnor number of an i.h.s, see e.g.
\cite{Looijenga}]\label{milnor nmb} The \emph{Milnor number} of
the hypersurface singularity $0\in U_0$ is defined as \emph{the
multiplicity of the critical point $0\in\C^{n+1}$ of $f$ as a
solution of the system of partials of $f$} (\cite{Milnor68} \S 7)
which is
\begin{equation}\label{milnor alg.}
  \mu_f(0)  = \dim_{\C}\left(\mathcal{O}_0/J_f\right) =
    \dim_{\C}\left(\C\{x_1,\ldots,x_{n+1}\}/J_f\right)
\end{equation}
where $\dim_{\C}$ means ``dimension as a $\C$--vector space" and
$J_f$ is the jacobian ideal $J_f:=\left(\frac{\partial f}{\partial
x_1},\ldots,\frac{\partial f}{\partial x_{n+1}}\right)$. For
shortness we will denote the Milnor number (\ref{milnor alg.}) by
$\mu(0)$ whenever $f$ is clear.
\end{definition}

\section{Tyurina number of an isolated hypersurface singularity}

\subsection{Deformations of complex spaces}\label{deformazioniCY}
Let $\mathcal{X}\stackrel{x} {\longrightarrow}B$ be a \emph{flat},
surjective map of complex spaces such that $B$ is connected and
there exists a special point $0\in B$ whose fibre $X=x^{-1}(0)$
may be singular. Then $\mathcal{X}$ is called \emph{a deformation
family of $X$}. If the fibre $X_b=x^{-1}(b)$ is smooth, for some
$b\in B$, then $X_b$ is called \emph{a smoothing of $X$}.

\noindent Let $\Omega_{X}$ be the sheaf of holomorphic
differential forms on $X$ and consider the
\emph{Lichtenbaum--Schlessinger cotangent sheaves}
\cite{Lichtenbaum-Schlessinger} of $X$, $\Theta^i_{X} = \ext
^i\left(\Omega_{X},\mathcal{O}_{X}\right)$. Then $\Theta^0_X = \hm
\left(\Omega_X,\mathcal{O}_X\right)=: \Theta_X$ is the ``tangent"
sheaf of $X$ and $\Theta^i_X$ is supported over $\Sing(X)$, for
any $i>0$. Consider the associated \emph{local} and \emph{global
deformation objects}
\[
    T^i_X:= H^0(X,\Theta^i_X)\quad,\quad
    \T^i_X:=\Ext ^i \left(\Omega^1_X,\mathcal{O}_X\right)\ ,\
    i=0,1,2.
\]
Then by the \emph{local to global spectral
    sequence} relating the global $\Ext$ and sheaf $\ext$
    (see \cite{Grothendieck57} and \cite{Godement} II, 7.3.3) we get
    \begin{equation*}
    \xymatrix{E^{p,q}_2=H^p\left(X,\Theta_X^q\right)
               \ar@{=>}[r] & \mathbb{T}_X^{p+q}}
    \end{equation*}
giving that
\begin{eqnarray}
    &&\T^0_X \cong T^0_X \cong H^0(X,\Theta_X)\ , \label{T0} \\
    \label{T-lisci}&&\text{if $X$ is smooth then}\quad \T^i_X\cong H^i(X,\Theta_X)\ ,  \label{T1'}\\
    &&\text{if $X$ is Stein then}\quad T^i_X\cong\T^i_X\ .
    \label{Ti-Ti}
\end{eqnarray}
Recall that $\mathcal{X}\stackrel{x}{\longrightarrow}B$ is called
a \emph{versal} deformation family of $X$ if for any deformation
family $(\mathcal{Y},X)\stackrel{y}{\longrightarrow}
    (C,o)$ of $X$ there exists a map of pointed complex spaces $h:(U,o)\rightarrow
    (B,0)$, defined on a neighborhood $o\in U\subset C$,
    such that $\mathcal{Y}|_U$ is the \emph{pull--back} of
    $\mathcal{X}$ by $h$ i.e.
    \[
        \xymatrix{&\mathcal{Y}|_U=U\times_B \mathcal{X}\ar[r]\ar[d]^-y&\mathcal{X}\ar[d]^-x\\
                   C& U\ar@{_{(}->}[l]\ar[r]^-h&B}
    \]

\begin{theorem}[Douady--Grauert--Palamodov \cite{Douady74},
\cite{Grauert74}, \cite{Palamodov72} and \cite{Palamodov} Theorems
5.4 and 5.6]\label{DGP teorema} Every compact complex space $X$
has an effective versal deformation
$\mathcal{X}\stackrel{x}{\longrightarrow}B$ which is a proper map
and a versal deformation of each of its fibers. Moreover the germ
of analytic space $(B,0)$ (the \emph{Kuranishi space of $X$}) is
isomorphic to the germ of analytic space $(q^{-1}(0), 0)$, where
$q:\T^1_X\rightarrow \T^2_X$ is a suitable holomorphic map (the
\emph{obstruction map}) such that $q(0)=0$.
\end{theorem}

In particular if $q\equiv 0$ (e.g. when $\T^2_X=0$) then $(B,0)$
turns out to be isomorphic to the germ of a neighborhood of the
origin in $\T^1_X$.

\subsection{Deformations of an i.h.s.} Let us consider the germ of
i.h.s. $U_0:=\Spec(\mathcal{O}_{f,0})$ as defined in
(\ref{germe}).

\begin{definition}[Tyurina number of an i.h.s.]\label{Tyurina}
The \emph{Tyurina number} of the i.h.s. $0\in U_0$ is
\[
    \tau_f(0) := \dim_{\C}\T^1_{U_0} \stackrel{\text{(\ref{Ti-Ti})}}=
    \dim_{\C}T^1_{U_0} = h^0(U_0,\Theta^1_{U_0})
\]
often denoted simply by $\tau(0)$ whenever $f$ is clear. Since
$U_0$ is Stein, the obstruction map $q$ in Theorem \ref{DGP
teorema} is trivial and the Tyurina number $\tau(0)$ turns out to
give \emph{the dimension of the Kuranishi space of $U_0$.}
\end{definition}

\begin{proposition}[see e.g. \cite{Stevens}]\label{T^1}
If $0\in U_0=\Spec(\mathcal{O}_{f,0})$ is the
germ of an i.h.s. then $\T^1_{U_0}\cong \mathcal{O}_{f,0}/J_f$ and
\begin{equation}\label{Tyurina alg.}
    \tau_f(0)=\dim_{\C}\left(\C\{x_1,\ldots,x_{n+1}\}/I_f\right)\ .
\end{equation}
where $I_f:=(f)+J_f$. Then, recalling (\ref{milnor alg.}),
$\tau(0)\leq\mu(0)$. In particular $\tau(0)$ gives the
multiplicity of $0$ as a singular point of the complex space germ
$U_0$.
\end{proposition}

\section{Milnor and Tyurina numbers of a polynomial}

Let us consider the polynomial algebra $\Cx:=\C[x_1,...,x_{n+1}]$
and let $I\subset\Cx$ be an ideal. If we consider the natural
inclusion $\Cx\subset\Cxx :=\C\{x_1,\ldots ,x_{n+1}\}$ then we get
\begin{equation}\label{disuguaglianza}
    \dim _{\C}\left(\Cxx / I\cdot\Cxx\right)\leq \dim _{\C}\left(\Cx / I\right)
\end{equation}
since the algebra $\Cxx$ contains also non constant units.

\begin{example}\label{esempio}
Let us consider the ideal $I=(x+x^2)\subset \C[x]$. Then
$\C[x]/I=\langle 1,x \rangle_{\C}$ while
$\C\{x\}/I\cdot\C\{x\}=\langle 1 \rangle_{\C}$, since $x\in I\cdot
\C\{x\}$ being associated with the generator $x+x^2$ by the unit
$1+x\in\C\{x\}$. This means that
\begin{equation}\label{disuguaglianza-stretta}
    1=\dim _{\C}\left(\C\{x\} / I\cdot\C\{x\}\right) < \dim _{\C}\left(\C[x] /
    I\right)=2\ .
\end{equation}
In particular the first equality means that, recalling Definition
\ref{milnor nmb} and formula (\ref{milnor alg.}), the Milnor
number of $0\in\C$, as a critical point of the polynomial map
$f(x)=({1\over 2}+{1\over 3}x)x^2$, turns out to be $\mu (0) =1$.
\end{example}

Let us then set the following

\begin{definition}[Milnor and Tyurina numbers of a polynomial]
\label{polyMT}
Given a polynomial $f\in\Cx$ the following dimension
\begin{equation}\label{Gmilnor}
    \mu(f):=\dim _{\C}\left(\Cx / J_f\right)
\end{equation}
is called \emph{the Milnor number of the polynomial $f$}.
Analogously the dimension
\begin{equation}\label{Gtyurina}
    \tau(f):=\dim _{\C}\left(\Cx / I_f\right)
\end{equation}
where $I_f:=(f)+J_f$, is called \emph{the Tyurina number of the
polynomial $f$}. Then clearly $\tau(f)\leq \mu(f)$.
\end{definition}

The inequality (\ref{disuguaglianza}) then gives that $\mu_f
(p)\leq\mu(f)$ and $\tau_f(p)\leq\tau(f)$, for any point
$p\in\C^{n+1}$ and any polynomial $f\in\Cx$. Moreover Milnor and
Tyurina numbers of points and polynomials are related by the
following

\begin{proposition}[see e.g. \cite{Greuel-Pfister} \S A.9]\label{Loc-to-Glob}
For any $f\in\C[x_1,\ldots ,x_{n+1}]$
\begin{eqnarray}\label{loc-to-glob}
  \mu(f) &=& \sum_{p\in\C^{n+1}} \mu_f(p) \\
\nonumber
  \tau(f) &=& \sum_{p\in\C^{n+1}} \tau_f(p)
\end{eqnarray}
\end{proposition}

\begin{remark}
Observe that sums on the right terms of (\ref{loc-to-glob}) are
actually finite and well--defined since $\mu(p)\neq 0$ if and only
if $p$ is a critical point of the polynomial map $f$ and
$\tau(p)\neq 0$ if and only if $p$ is a singular point of the
$n$--dimensional algebraic hypersurface
$f^{-1}(0)\subset\C^{n+1}$.
\end{remark}

\begin{example}[Example \ref{esempio} continued] Critical points of the polynomial
map $f(x)=({1\over 2}+{1\over 3}x)x^2$ are given by $0$ and $-1$.
The translation $x\mapsto x-1$ transforms $f$ into the polynomial
$g(x)={1\over 3}({1\over 2}+x)(x-1)^2$ and
\[
    \mu_f(-1)=\mu_g(0)=\dim_{\C}\left(\C\{x\}/(x^2-x)\right)=1\ .
\]
The inequality (\ref{disuguaglianza-stretta}) then gives
$2=\mu(f)=\mu_f(0)+\mu_f(-1)$, according with the first equality in
(\ref{loc-to-glob}). Moreover $0$ is the unique singular point of
the 0--dimensional hypersurface $f^{-1}(0)=\{-3/2,0\}\subset\C$
and
\begin{eqnarray*}
  \tau(f) &=& \dim_{\C}\left(\C[x]\left/\left(\left({1\over 2}+{1\over
    3}x\right)x^2,x+x^2\right)\right)\right.=1 \\
  \tau_f(0) &=& \dim_{\C}\left(\C\{x\}\left/\left(\left({1\over 2}+{1\over
    3}x\right)x^2,x+x^2\right)\right)\right.=1
\end{eqnarray*}
according with the second equality in (\ref{loc-to-glob}).
\end{example}

\subsection{Milnor and Tyurina numbers of weighted homogeneous polynomials}
Let us recall that a polynomial $f\in\C[x_1,\dots ,x_{n+1}]$ is
called \emph{weighted homogeneous (w.h.p.)} or
\emph{quasi--homogeneous} if there exist  $n+1$ \emph{positive}
rational numbers $\mathbf{w}=(w_1,\ldots ,w_{n+1})\in \Q^{n+1}$
such that
\[
    \sum_{i=1}^{n+1} w_i\alpha_i =1
\]
for any monomial
$\mathbf{x}^{\alpha}:=\prod_{i=1}^{n+1}x_i^{\alpha_i}$ appearing in
$f$; $\mathbf{w}$ is then called \emph{the vector of (rational)
weights of} $f$ and the \emph{generalized Euler formula}
\begin{equation}\label{Eulero}
    f=\sum_{i=1}^{n+1}w_i x_i \frac{\partial f}{\partial x_i}
\end{equation}
follows immediately for any w.h.p. $f\in\Cx$ admitting the same
vector of weights. For any $i=1,\ldots,n+1$, let
$(p_i,q_i)\in\N^2$ be the unique ordered couple of positive
coprime integers such that $p_i/q_i$ is the reduced fraction
representing the positive rational number $w_i$. Calling $d$ the
least common factor of denominators $q_1,\ldots,q_{n+1}$, the
positive integers $d_i:=dw_i$ satisfy the following \emph{weighted
homogeneity relation}
\begin{equation}\label{omogeneità}
    \forall\lambda\in\C\quad
    f\left(\lambda^{d_1}x_1,\ldots,\lambda^{d_{n+1}}x_{n+1}\right) = \lambda^d\
    f(x_1,\dots,x_{n+1})\ .
\end{equation}
For this reason $\mathbf{d}=(d_1,\ldots,d_{n+1})$ is called
\emph{the vector of integer weights of $f$} and $d=|\mathbf{d}|$
is called the \emph{degree of $f$}.

\begin{proposition}\label{MT-whp} Given a polynomial $f\in\Cx$ with a finite number of critical points, the following
assertions are equivalent:
\begin{itemize}
    \item[{\rm(\emph{a})}] $f$ is a w.h.p.,
    \item[{\rm(\emph{b})}] $\tau(f)=\mu(f)$,
    \item[{\rm(\emph{c})}] $\forall p\in\C^{n+1}\quad \tau(p)=\mu(p)$.
\end{itemize}
In particular \emph{(\emph{c})} implies that the set of critical
points of $f$ coincides with the set of singular points of
$f^{-1}(0)$. Moreover \emph{(\emph{a})} means actually that the
origin $0\in\C^{n+1}$ is the unique possible critical point of $f$
and then the unique possible singular point of $f^{-1}(0)$.
Therefore \emph{(\emph{b})} and \emph{((\emph{\emph{c}})} give
that
\begin{equation}\label{uguaglianze}
    \mu(0)=\mu(f)=\tau(f)=\tau(0)\ .
\end{equation}
\end{proposition}

\begin{proof} $(a)\Rightarrow (c)$. The generalized Euler formula
(\ref{Eulero}) implies that $I_f = J_f$ and (\emph{c}) follows
immediately by (\ref{milnor alg.}) and (\ref{Tyurina alg.}). In
particular if $\mathbf{x}=(x_1,\ldots,x_{n+1})$ is a critical
point of $f$, then (\ref{Eulero}) and (\ref{omogeneità}) give that
$(\lambda^{d_1}x_1,\ldots,\lambda^{d_{n+1}}x_{n+1})$ is a critical
point of $f$ for any complex number $\lambda$. Then $\mathbf{x}=0$
since $f$ admits at most a finite number of critical points,
meaning that $f$ admits at most the origin as a critical point.

$(c)\Rightarrow (b)$. This follows immediately by Proposition
\ref{Loc-to-Glob}.

$(b)\Rightarrow (a)$.
Since $J_f\subseteq I_f$ then we
get the natural surjective map of $\C$--algebras
\[
    \xymatrix{\Cx/J_f\ar@{->>}[r]&\Cx/I_f}
\]
The hypothesis $\tau(f)=\mu(f)$ implies that it is also injective
 which suffices to show that $J_f=I_f$.
 Hence $f\in J_f$ and a
famous result by K.~Saito \cite{Saito} allows to conclude that $f$
is a w.h.p..
\end{proof}

\begin{definition}[Weighted homogeneous singularity -
w.h.s.]\label{whs} A $n$--dimensional i.h.s. $0\in
U_0=\Spec\mathcal{O}_{f,0}$ is called \emph{weighted homogeneous}
(or \emph{quasi-homogeneous}) if there exists a w.h.p.
$F\in\C[x_1,\ldots,x_{n+1}]$ such that $U_0\cong \Spec
\mathcal{O}_{F,0}$, as germs of complex spaces.
\end{definition}

\begin{remark}\label{contact=right} Definition \ref{whs} is equivalent to require that
there exists an automorphism $\phi^*$ of $\mathcal{O}_0=\Cxx$,
induced by a biholomorphic local coordinates change
$(\C^{n+1},0)\stackrel{\phi}{\rightarrow}(\C^{n+1},0)$, such that
$\phi^*(f):=f\circ\phi=F$ (see \cite{Greuel-etal} Lemma 2.13).
\end{remark}

\begin{proposition}[Characterization of a w.h.s]\label{whs-caratterizzazione} $0\in
U_0=\Spec\mathcal{O}_{f,0}$ is a w.h.s. if and only if
$\tau_f(0)=\mu_f(0)$.
\end{proposition}

\begin{proof} The statement follows immediately by Proposition
\ref{MT-whp}, keeping in mind Remark \ref{contact=right} and
observing that the jacobian ideals $J_f$ and $J_F$ can be obtained
each other by multiplying the jacobian matrix of the
coordinate change $\phi$, which is clearly invertible in a
neighborhood of $0$.
\end{proof}

\subsection{An example: weighted homogeneous cDV singularities}
An example of an isolated hypersurface singularity is given by a
\emph{compound Du Val} (cDV) singularity which is a 3--fold point
$p$ such that, for a hyperplane section $H$ through $p$, $p\in H$
is a Du Val surface singularity i.e. an A--D--E singular point
(see \cite{Reid80}, \S 0 and \S 2, and \cite{BPvdV84}, chapter
III). Then a cDV point $p$ is a germ of hypersurface singularity
$0\in {U}_0:=\Spec (\mathcal{O}_{f,0})$, where $f$ is the
polynomial
    \begin{equation}\label{cDV}
    f(x,y,z,t) := x^2 + q(y,z) + t\ g(x,y,z,t)\ ,
    \end{equation}
    such that $g(x,y,z,t)$ is a generic element of the maximal
    ideal
    $\mathfrak{m}_0:=(x,y,z,t)\subset\C[x,y,z,t]$ and
    \begin{eqnarray}\label{DV}
    A_n &:& q(y,z):=y^2 + z^{n+1}\quad\text{for}\ n\geq 1 \\
    \nonumber
    D_n &:& q(y,z):=y^2z + z^{n-1}\quad\text{for}\ n\geq 4 \\
    \nonumber
    E_6 &:& q(y,z):=y^3 + z^4   \\
    \nonumber
    E_7 &:& q(y,z):=y^3 + yz^3   \\
    \nonumber
    E_8 &:& q(y,z):=y^3 + z^5
    \end{eqnarray}
    In particular if
    \begin{equation}\label{Arnold}
        g(x,y,z,t)=t
    \end{equation}
    then $f=0$ in (\ref{cDV}) \emph{is said to define an Arnol'd
    simple (threefold) singularity} (\cite{Arnol'd}, \cite{Arnol'd&c} \S 15 and in particular \cite{Arnol'd-ST} \S I.2.3) denoted by $A_n,D_n,E_{6,7,8}$, respectively.

    \noindent The index ($n,6,7,8$) turns out to be the
    Milnor number of the surface Du Val singularity
    $0\in U_0\cap\{t=0\}$  or equivalently its Tyurina number,
    since a Du Val singular point always admits a weighted
    homogeneous local equation. When a cDV point is
    defined by a weighted homogeneous polynomial $f$, a classical result of J.~Milnor and
    P.~Orlik allows to compare this index with its Milnor (and then Tyurina) number.
    In particular we get
\begin{eqnarray}\label{weights}
  w(x) &=& 1/2 \\
\nonumber
  w(y) &=& \left\{\begin{array}{cc}
    1/2 & \text{if $p$ is $cA_n$,} \\
    (n-2)/(2n-2) & \text{if $p$ is $cD_n$,} \\
    1/3 & \text{if $p$ is $cE_{6,7,8}$.} \\
  \end{array}\right. \\
\nonumber
  w(z) &=& \left\{\begin{array}{cc}
    1/(n+1) & \text{if $p$ is $cA_n$,} \\
    1/(n-1) & \text{if $p$ is $cD_n$,} \\
    1/4 & \text{if $p$ is $cE_6$,} \\
    2/9 & \text{if $p$ is $cE_7$,} \\
    1/5 & \text{if $p$ is $cE_8$.} \\
  \end{array}\right.
\end{eqnarray}

\begin{theorem}[Milnor--Orlik \cite{Milnor-Orlik}, Thm. 1]\label{Milnor-Orlik}
Let $f(x_1,\ldots,x_{n+1})$ be a w.h.p., with rational weights
$w_1,\ldots,w_{n+1}$, admitting an isolated critical point at the
origin. Then the Milnor number of the origin is given by
\[
    \mu (0) =
    \left(w_1^{-1}-1\right)\left(w_2^{-1}-1\right)\cdots\left(w_{n+1}^{-1}-1\right)\ .
\]
\end{theorem}

By putting weights (\ref{weights}) in the previous Milnor--Orlik
formula we get the following

\begin{corollary}\label{Milnor w.h.}
Le $0\in U_0$ be a w.h. cDV point of index $n$. Then
\[
    \tau(0)=\mu(0)= n \left(w(t)^{-1}-1\right)\ .
\]
In particular for Arnol'd simple singularities we get
$n=\tau(0)=\mu(0)$, as can also be directly checked by the
definition.
\end{corollary}

\subsection{The algebraic computation via Gr\"{o}bner basis}
Let us consider an ideal $I\subset\Cx$ and
let $L(I)$ denote the ideal generated by \emph{leading monomials},
with respect to a fixed \emph{monomial order}, of elements in $I$.
It is a well known fact (see e.g. \cite{Cox} \S 5.3) that the
following are isomorphisms of $\C$--vector spaces
\begin{equation}\label{groebner}
    \xymatrix@1{\C[x_1,\ldots,x_{n+1}]/I\ \ar[r]^-{\cong}& \ \C[x_1,\ldots,x_{n+1}]/L(I)\ \ar[r]^-{\cong}& \ \left\langle M\setminus L(I)\right\rangle_{\C}}
\end{equation}
where $M$ is the set (actually a multiplicative monoid) of all monomials $\mathbf{x}^{\alpha}$.
Consider
a polynomial $f\in\Cx$. By Definition \ref{polyMT} and isomorphisms (\ref{groebner}), the computation
of $\mu(f)$ and $\tau(f)$ reduces to calculate
\[
    \dim_{\C}\left(\C[x_1,\ldots,x_{n+1}]/L(I)\right)=\left| M\setminus L(I)\right|
\]
where $L(I)$ is the ideal generated by \emph{leading monomials},
with respect to a fixed monomial order, of elements in the ideal
$I\subset\C[x_1,\ldots,x_{n+1}]$ which is either the jacobian
ideal $J_f$ or the ideal $I_f=(f)+J_f$, respectively. The point is
then \emph{determining a Gr\"{o}bner basis of $I$ w.r.t. the fixed
monomial order}, which can be realized e.g. by the
\texttt{Groebner Package} of MAPLE.

\begin{remark} The MAPLE computation of Milnor and Tyurina numbers of polynomials
is realized by procedures \texttt{PolyMilnor} and
\texttt{PolyTyurina}, whose concrete description is postponed to
section \ref{maple}. Their usefulness is clarified by Proposition
\ref{MT-whp} and in particular by equations (\ref{uguaglianze}).
In fact in the case of a w.h.p. $f$ admitting an isolated critical
point in $0\in\C^{n+1}$ there is no need of working with power
series (and then with \emph{local monomial orders}) to determine
$\mu(0)$ and $\tau(0)$. Since the Buchberger algorithm implemented
with MAPLE turns out to be more efficient when running by usual
term orders, the use of \emph{global} procedures
\texttt{PolyMilnor} and \texttt{PolyTyurina} has to be preferred
to the use of their \emph{local} counterparts \texttt{Milnor} and
\texttt{Tyurina}, when possible.
\end{remark}

\section{Monomial Ordering}

First of all let us recall what is usually meant by a monomial
order. For more details the interested reader in remanded e.g. to
\cite{Cox} \S 2.2, \cite{Greuel-Pfister} \S 9 and
\cite{Decker-Lossen} \S 1.

Let $M$ be the multiplicative monoide of monomials
$\mathbf{x}^{\alpha}=\prod_{i=1}^{n+1}x_i^{\alpha_i}$\ : clearly
$\log : M\stackrel{\cong}{\rightarrow} \N^{n+1}$. A \emph{(global)
monomial order on $\Cx$} is a \emph{total order relation} $\leq$
on $M$ which is
\begin{itemize}
    \item[(i)] \emph{multiplicative} i.e. $\forall \alpha,\beta,\gamma\in\N^{n+1}
    \quad \mathbf{x}^{\alpha}\leq\mathbf{x}^{\beta}\Rightarrow
    \mathbf{x}^{\alpha}\cdot\mathbf{x}^{\gamma}\leq\mathbf{x}^{\beta}\cdot\mathbf{x}^{\gamma}$\
    ,
    \item[(ii)] \emph{a well-ordering} i.e. every nonempty subset of
    $M$ has a \emph{smallest element}.
\end{itemize}
Since $\Cx$ is a noetherian ring a multiplicative total order on
$M$ is a well-ordering if and only if
\begin{itemize}
    \item[(ii')] $\forall i=1,\ldots,n+1\quad 1 < x_i$\ .
\end{itemize}

\begin{definition}[Local and global monomial orders, \cite{Greuel-Pfister} Definition 1.2.4]\label{lmo}
In the following a m.o. on $\Cx$ will denote simply a total order
relation $\leq$ on $M$ which is multiplicative i.e. satisfying
(i). A m.o. will be called \emph{global (g.m.o.)} if also (ii), or
equivalently (ii'), is satisfied. Moreover a m.o. will be called
\emph{local (l.m.o.)} if
\begin{itemize}
    \item[(ii'')] $\forall i=1,\ldots,n+1\quad x_i < 1$\ .
\end{itemize}
\end{definition}

\subsection{Localizations in $\Cx$ and rings implemented by monomial orders}

Given a m.o. $\leq$ on $M$ consider the following subset of $\Cx$
\begin{equation*}
    S:=\left\{f\in\Cx\ |\ L(f)\in\C\setminus\{0\}\right\}
\end{equation*}
where $L(f)$ is the \emph{leading monomial} of $f$ w.r.t. $\leq$\
. Since $S$ is a multiplicative subset of $\Cx$ we can consider
the localization $S^{-1}\Cx$. Then (ii') and (ii'') give
immediately the following

\begin{proposition}\label{anelli}
\begin{eqnarray*}
  S^{-1}\Cx &=& \Cx\quad\Leftrightarrow\quad \text{$\leq$ is a g.m.o.} \\
  S^{-1}\Cx &=& \Cx_{(\mathbf{x})}\quad\Leftrightarrow\quad \text{$\leq$ is a l.m.o.}
\end{eqnarray*}
where $\Cx_{(\mathbf{x})}$ is the localization of $\Cx$ at the
maximal ideal $(\mathbf{x})\subset\Cx$.
\end{proposition}

By Taylor power series expansion of locally holomorphic functions,
there is a natural inclusion $\Cx_{(\mathbf{x})}\subset\Cxx$
giving the following commutative diagram, for every ideal
$I\subset\Cx_{(\mathbf{x})}$:
\begin{equation*}
    \xymatrix{\Cx_{(\mathbf{x})}\ar@{^{(}->}[r]\ar@{->>}[d]&\Cxx\ar@{->>}[d]\\
               \left.\Cx_{(\mathbf{x})}\right/I\ar@{^{(}->}[r]&\left.\Cxx\right/I\cdot\Cxx }
\end{equation*}

\begin{proposition}[for a proof see e.g. \cite{Decker-Lossen} Proposition 9.4]
If $\dim_{\C}\left(\left.\Cx_{(\mathbf{x})}\right/I\right)$ is
finite then the inclusion
$\left.\Cx_{(\mathbf{x})}\right/I\subset\left.\Cxx\right/I\cdot\Cxx$
is an isomorphism of $\C$--algebras. In particular both the
underlying vector spaces have the same dimension.
\end{proposition}

\begin{theorem}[\cite{Decker-Lossen} Theorem 9.29] Let $\leq$ be a
m.o. on $\Cx$ and $I\subset S^{-1}\Cx$ be an ideal. Then
\[
    \dim_{\C}\left(S^{-1}\Cx/I\right)=\dim_{\C}\left(\Cx/L(I)\right)
\]
where $L(I)$ is the ideal generated by \emph{leading monomials},
with respect to $\leq$\ , of polynomials in $I$. In particular if
it is finite then $M\setminus L(I)$ represents a basis
of the vector space $S^{-1}\Cx$.
\end{theorem}

\begin{corollary} Let $f\in\Cx$ admit an isolated critical point
at $0\in\C^{n+1}$, $\leq$ be a m.o. on $\Cx$ and $I=J_f$ (resp.
$I=I_f$). Then
\[
    \dim_{\C}\left(\Cx/L(I)\right)=\left\{\begin{array}{cc}
      \mu(f)\ (\text{resp.}\ \tau(f)) & \text{if $\leq$ is global} \\
      \mu_f(0)\ (\text{resp.}\ \tau_f(0)) & \text{if $\leq$ is local} \\
    \end{array}\right.
\]
\end{corollary}

\begin{remark} The point is
then determining a \emph{standard basis of $I$} w.r.t. the fixed
m.o. $\leq$\ . If the latter is a \emph{global} one, a standard
basis is a usual Gr\"{o}bner basis which is obtained by applying
the Buchberger algorithm. In MAPLE this is implemented by the
\texttt{Groebner Package}.

\noindent On the other hand, if $\leq$ is a \emph{local} m.o. the
Buchberger algorithm does no more work. In fact $\leq$ is no more
a well-ordering and the division algorithm employed by the
Buchberger algorithm for determining \emph{normal forms} of
\emph{S-polynomials} may do not terminate. This problem can be dodged
by means of a \emph{weak normal form algorithm} (\emph{weakNF}),
firstly due to F.~Mora \cite{Mora}, and of a \emph{standard basis
algorithm} (\emph{SB}) which replaces the Buchberger algorithm.
This is precisely what has been implemented in SINGULAR since 1990
(see \cite{Greuel-Pfister} \S 1.7 and references thereof). The aim
of the following section \ref{maple} is to present a MAPLE
implementation of \emph{weakNF} and \emph{SB} algorithms to yield
a procedure computing Milnor and Tyurina numbers of points.
\end{remark}

\section{MAPLE subroutines in detail}\label{maple}

The present section is devoted to present and describe in detail
MAPLE subroutines computing Milnor and Tyurina numbers of critical
points of a polynomial $f$. They are available as MAPLE 12 files
\texttt{.mw} at \cite{Maple}. They are composed by several
procedures, the most important of which are the following:
\begin{itemize}
    \item \texttt{weakNF}\quad which is the MAPLE implementation of the
    weak normal form Algorithm 1.7.6 in \cite{Greuel-Pfister}:
    \begin{eqnarray*}
      &\text{Input}& :\quad  f\in\Cx\ ,\ \text{$G$ a finite subset in $\Cx$}\ ,\\
            &&         \quad   \text{\texttt{variables} $\subseteq\textbf{x}=\{x_1,\dots,x_{n+1}\}$}\ , \\
            &&          \quad \text{$STO =$ a MAPLE type  \texttt{ShortMonomialOrder} expression} \\
      &\text{Output}&:\quad
      h\in\left(\C[\mathbf{x}\setminus{\text{\texttt{variables}}}]
      \right)[\text{\texttt{variables}}]\\
      && \quad       \text{which is a polynomial weak normal form of $f$ w.r.t.
      $G$}\ ;
    \end{eqnarray*}
    \item \texttt{SB}\quad which is the MAPLE implementation of the
    standard basis Algorithm 1.7.1 in \cite{Greuel-Pfister}:
    \begin{eqnarray*}
      &\text{Input}& :\quad \text{$G$ a finite subset in $\Cx$, \texttt{weakNF}}\ , \\
      &\text{Output}&:\quad \text{a standard basis $S$ of $I:=(G)\subset\Cx$}\ ;
    \end{eqnarray*}
    \item \texttt{Milnor}\quad which is the procedure computing
    the Milnor number of an isolated critical point:
    \begin{eqnarray*}
      &\text{Input}& :\quad F\in\Cx\ ,\ \text{\texttt{SB}}\ , \\
      &&\quad \text{Optional: \texttt{variables}, by default assigned by
      \texttt{indets(F)}}\\
      &&\quad \text{Optional: a l.m.o., by default assigned by
      \texttt{tdeg\_\ min(variables)}}\\
       &\text{Output}& :\quad \text{A standard basis $G$ of}\
       J_F\subset\left(\C[\mathbf{x}\setminus{\text{\texttt{variables}}}]
      \right)[\text{\texttt{variables}}],\ L(G),\\
      &&\quad\quad  M\setminus L(J_F)\ ,\ \mu_F(0)\ ;
    \end{eqnarray*}
    \item \texttt{PolyMilnor}\quad which is the procedure computing
    the Milnor number of a polynomial:
    \begin{eqnarray*}
      &\text{Input}& :\quad F\in\Cx\ ,\ \text{\texttt{Basis} in the \texttt{Groebner} package}\ , \\
      &&\quad \text{Optional: \texttt{variables}, by default assigned by
      \texttt{indets(F)}}\\
      &&\quad \text{Optional: a g.m.o., by default assigned by
      \texttt{tdeg(variables)}}\\
       &\text{Output}& :\quad \text{A Gr\"{o}bner basis $G$ of}\
       J_F\subset\left(\C[\mathbf{x}\setminus{\text{\texttt{variables}}}]
      \right)[\text{\texttt{variables}}],\ L(G),\\
      &&\quad\quad  M\setminus L(J_F)\ ,\ \mu(F)\ ;
    \end{eqnarray*}
    \item \texttt{Tyurina}\quad which is the procedure computing
    the Tyurina number of an isolated singular point:
    \begin{eqnarray*}
      &\text{Input}& :\quad F\in\Cx\ ,\ \text{\texttt{SB}}\ , \\
      &&\quad \text{Optional: \texttt{variables}, by default assigned by
      \texttt{indets(F)}}\\
      &&\quad \text{Optional: a l.m.o., by default assigned by
      \texttt{tdeg\_\ min(variables)}}\\
       &\text{Output}& :\quad \text{A standard basis $G$ of}\
       I_F\subset\left(\C[\mathbf{x}\setminus{\text{\texttt{variables}}}]
      \right)[\text{\texttt{variables}}],\ L(G),\\
      &&\quad\quad  M\setminus L(I_F)\ ,\ \tau_F(0)\ ;
    \end{eqnarray*}
    \item \texttt{PolyTyurina}\quad which is the procedure computing
    the Tyurina number of a polynomial:
    \begin{eqnarray*}
      &\text{Input}& :\quad F\in\Cx\ ,\ \text{\texttt{Basis} in the \texttt{Groebner} package}\ , \\
      &&\quad \text{Optional: \texttt{variables}, by default assigned by
      \texttt{indets(F)}}\\
      &&\quad \text{Optional: a g.m.o., by default assigned by
      \texttt{tdeg(variables)}}\\
       &\text{Output}& :\quad \text{A Gr\"{o}bner basis $G$ of}\
       I_F\subset\left(\C[\mathbf{x}\setminus{\text{\texttt{variables}}}]
      \right)[\text{\texttt{variables}}],\ L(G),\\
      &&\quad\quad  M\setminus L(I_F)\ ,\ \tau(F)\ .
    \end{eqnarray*}
\end{itemize}

\subsection{The subroutines}\label{subroutine}

 Preambles to introduce the useful MAPLE packages:

\begin{maplegroup}
\begin{mapleinput}
\mapleinline{active}{1d}{with(Groebner):
}{}
\end{mapleinput}
\end{maplegroup}
\begin{maplegroup}
\begin{mapleinput}
\mapleinline{active}{1d}{with(PolynomialIdeals):
}{}
\end{mapleinput}
\end{maplegroup}
\begin{maplegroup}
\begin{mapleinput}
\mapleinline{active}{1d}{with(Ore_algebra): }{}
\end{mapleinput}
\end{maplegroup}

\halfline

\noindent A first control procedure :

\halfline

\begin{maplegroup}
\begin{mapleinput}
\mapleinline{active}{1d}{localorglobal := proc (STO, variables)
A:= poly_algebra(op(variables)); TP := MonomialOrder(A, STO);
nuu:= 1; muu := 1;
for i to nops(variables) do if TestOrder(1,variables[i], TP)
then nuu := 0 else muu := 0 end if end do;
if nuu = 1 then Lo else if muu = 1 then Gl else Mi
end if end if end proc:}{}
\end{mapleinput}
\end{maplegroup}

\halfline

\noindent Actually \texttt{localorglobal} is not an essential
procedure in the present routine: its meaning is simply that of
giving a feedback about what kind of m.o. the user is employing
and stopping the procedure running with a wrong term order: in
fact \texttt{Milnor} and \texttt{Tyurina} may give wrong output
when running with a g.m.o.; on the other hand \texttt{PolyMilnor}
and \texttt{PolyTyurina} may not terminate when running with a
l.m.o..

\subsubsection{Implementing local monomial orders}
The following three procedures give the core of the \emph{MAPLE
implementation of l.m.o.'s} for determining standard basis of
ideals in $\Cx_{(\mathbf{x})}$. The first procedure introduces the
\emph{ecart} concept which is the main ingredient in the Mora
algorithm \texttt{weakNF}. It is defined following
\cite{Greuel-Pfister} Definition 1.7.5:

\halfline

\begin{maplegroup}
\begin{mapleinput}
\mapleinline{active}{1d}{ecart := proc (f, variables, STO)
degree(f, variables)-degree(LeadingMonomial(f, STO), variables)
end proc:}{}
\end{mapleinput}
\end{maplegroup}

\halfline

\noindent Then we give the Mora algorithm for determining a weak
normal form of a polynomial $f\in\Cx$ w.r.t. a finite subset of
polynomials $G\subset\Cx$ (see \cite{Greuel-Pfister} Algorithm
1.7.6 and \cite{Decker-Lossen} Algorithm 9.22):

\halfline

\begin{maplegroup}
\begin{mapleinput}
\mapleinline{active}{1d}{weakNF := proc (f, G, variables, STO)
h := f;
TT := G;
TTh := \{\}; for i to nops(TT) do
if divide(LeadingMonomial(h, STO),LeadingMonomial(TT[i], STO))
then TTh := \{TT[i], op(TTh)\}
end if end do;
while (h <> 0 and TTh <> \{\}) do
L := [op(TTh)];
L1 := sort(L, proc (t1, t2) options operator, arrow;
ecart(t1,variables,STO) <= ecart(t2,variables,STO) end proc);
g := L1[1];
if ecart(h, variables, STO) < ecart(g, variables, STO)
then TT := \{h, op(TT)\} end if;
h := SPolynomial(h, g, STO);
TTh := \{\};
for i to nops(TT) do
if divide(LeadingMonomial(h, STO), LeadingMonomial(TT[i],STO))
then TTh := \{TT[i], op(TTh)\} end if end do end do;
h end proc
}{}
\end{mapleinput}
\end{maplegroup}

\halfline

\noindent At last the standard basis procedure giving the analogue
of Buchberger algorithm with l.m.o.'s:

\halfline

\begin{maplegroup}
\begin{mapleinput}
\mapleinline{active}{1d}{SB := proc (G, variables, STO)
S := G;
P:= \{seq(seq(\{G[i], G[j]\},j=i+1..nops(G)),i=1..nops(G))\};
while P <> \{\} do
P1 := P[1];
P := `minus`(P,\{P1\});
h := weakNF(SPolynomial(P1[1],P1[2],STO), S, variables, STO);
if h <> 0 then P := \{seq(\{h, S[i]\}, i = 1 .. nops(S)),op(P)\};
S := \{h, op(S)\} end if end do;
S end proc:}{}
\end{mapleinput}
\end{maplegroup}

\halfline

\subsubsection{Computing Milnor and Tyurina numbers}\label{MTproc}
We are now in a position to introduce the procedures computing
Milnor and Tyurina numbers of critical points of a polynomial
$F\in\Cx$. They actually give some more output. Precisely, saying
$I=J_F$ (resp. $I=I_F$) \texttt{Milnor} (resp. \texttt{Tyurina})
returns:
\begin{itemize}
    \item a standard basis $G$ of the ideal
    $I\cdot\Cx_{(\mathbf{x})}\subset\Cx_{(\mathbf{x})}$,
    \item the leading monomial basis $L(G)$ w.r.t. a fixed m.o.,
    \item $M\setminus L(I)$ representing ($\mod I\cdot\Cx$) a basis of the quotient vector space
     $\Cx_{(\mathbf{x})}\left/I\cdot\Cx_{(\mathbf{x})}\right.$
    \item and its dimension over $\C$, which is $\mu_F(0)$ (resp. $\tau_F(0)$).
\end{itemize}
Each of the previous output may be listed separately by means of
the sub-procedures \texttt{MilnorGroebnerBasis, MilnorLT,
MilnorBasis} and \texttt{MilnorNumber} (respectively,
\texttt{TyurinaGroebnerBasis, TyurinaLT, TyurinaBasis} and
\texttt{TyurinaNumber}).

\halfline

\noindent The \texttt{Milnor} procedure:

\halfline

\begin{maplegroup}
\begin{mapleinput}
\mapleinline{active}{1d}{Milnor := proc (F, variables::set:=indets(F),
U::anything:=tdeg_min(op(variables)))
if type(U,ShortMonomialOrder) = false then error
"invalid input: your ShortMonomialOrder is not well defined"
 else if localorglobal(U, variables) <> Lo then error
"invalid input: your Short Monomial Order should be LOCAL"
else r := nops(variables);
J := [seq(diff(F, variables[i]), i = 1..r)];
G := SB(J,variables, U);
Ini := [seq(LeadingMonomial(G[i], U), i = 1..nops(G))];
massimo:=r*max(seq(degree(Ini[i],variables),i =1..nops(Ini)));
if IsZeroDimensional(<(op(Ini))>,\{op(variables)\}) = false then
error "the given critical point is not isolated"
else
L := sort([op(subs(uq = 1, convert(map(expand,
series(1/(product(1-variables[k]*uq, k = 1 .. r)),
uq, massimo+1)), polynom)))],
proc (t1, t2) options operator, arrow;
TestOrder(t1, t2, U) end proc);
M := [];
for j to nops(L) do up := 1;
for k to nops(Ini) do
if divide(L[j], Ini[k]) then up := 0
end if end do;
if up = 1 then M := [op(M), L[j]]
end if end do;
[G,Ini, M, nops(M)]
end if end if end if end proc:}{}
\end{mapleinput}
\end{maplegroup}

\halfline

\noindent The associated sub--procedures:

\begin{maplegroup}
\begin{mapleinput}
\mapleinline{active}{1d}{MilnorGroebnerBasis:=proc(F,variables::set:=indets(F),
T::anything:=tdeg_min(op(variables)))
Milnor(F,variables,T)[1]
end proc: }{}
\end{mapleinput}
\end{maplegroup}

\halfline

\begin{maplegroup}
\begin{mapleinput}
\mapleinline{active}{1d}{MilnorLT:=proc(F,variables::set:=indets(F),
T::anything:=tdeg_min(op(variables)))
Milnor(F, variables, T)[2]
end proc:}{}
\end{mapleinput}
\end{maplegroup}

\halfline

\begin{maplegroup}
\begin{mapleinput}
\mapleinline{active}{1d}{MilnorBasis:=proc(F,variables::set:=indets(F),
T::anything:=tdeg_min(op(variables)))
Milnor(F, variables,T)[3]
end proc: }{}
\end{mapleinput}
\end{maplegroup}

\halfline

\begin{maplegroup}
\begin{mapleinput}
\mapleinline{active}{1d}{MilnorNumber:=proc(F,variables::set:=indets(F),
T::anything:=tdeg_min(op(variables)))
Milnor(F,variables,T)[4] end proc: }{}
\end{mapleinput}
\end{maplegroup}

\halfline

\noindent The Tyurina procedure:

\halfline

\begin{maplegroup}
\begin{mapleinput}
\mapleinline{active}{1d}{Tyurina := proc(F,variables::set:=indets(F),
U::anything := tdeg_min(op(variables)))
if type(U,ShortMonomialOrder)=false then error
"invalid input: your ShortMonomialOrder is not well defined":
else if localorglobal(U, variables)<>Lo then error
"invalid input: your Short Monomial Order should be LOCAL":
else r := nops(variables):
K := [F, seq(diff(F, variables[i]), i =1 .. r)]:
H := SB(K, variables,U);
Ini := [seq(LeadingMonomial(H[i], U), i = 1..nops(H))];
massimo:=r*max(seq(degree(Ini[i]),i = 1..nops(Ini)));
if IsZeroDimensional(<(op(Ini))>,\{op(variables)\}) = false then
error "the given singular point is not isolated"
else
L := sort([op(subs(uq = 1, convert(map(expand,
series(1/(product(1-variables[k]*uq, k = 1 .. r)),
uq, massimo+1)), polynom)))],
proc (t1, t2) options operator, arrow;
TestOrder(t1, t2, U) end proc);
M := [];
for j to nops(L) do up :=1;
for k to nops(Ini) do
if divide(L[j], Ini[k]) then up := 0
end if end do;
if up = 1 then M := [op(M), L[j]]
end if end do;
[H,Ini,M,nops(M)] end if end if end if end proc: }{}
\end{mapleinput}
\end{maplegroup}

\halfline

\noindent The associated sub--procedures:

\halfline

\begin{maplegroup}
\begin{mapleinput}
\mapleinline{active}{1d}{TyurinaGroebnerBasis := proc(F,variables::set := indets(F),
T::anything :=tdeg_min(op(variables)))
Tyurina(F,variables,T)[1]
end proc: }{}
\end{mapleinput}
\end{maplegroup}

\halfline

\begin{maplegroup}
\begin{mapleinput}
\mapleinline{active}{1d}{TyurinaLT := proc (F,variables::set :=indets(F),
T::anything := tdeg_min(op(variables)))
Tyurina(F,variables,T)[2]
end proc: }{}
\end{mapleinput}
\end{maplegroup}

\halfline

\begin{maplegroup}
\begin{mapleinput}
\mapleinline{active}{1d}{TyurinaBasis := proc (F, variables::set:= indets(F),
T::anything := tdeg_min(op(variables)))
Tyurina(F,variables,T)[3]
end proc: }{}
\end{mapleinput}
\end{maplegroup}

\halfline

\begin{maplegroup}
\begin{mapleinput}
\mapleinline{active}{1d}{TyurinaNumber := proc (F, variables::set:= indets(F),
T::anything := tdeg_min(op(variables)))
Tyurina(F,variables, T)[4]
end proc:}{}
\end{mapleinput}
\end{maplegroup}

\halfline

Ultimately the following procedures allows to compute Milnor and
Tyurina numbers $\mu(F)$ and $\tau(F)$ of a polynomial $F\in\Cx$.
They give the same output of \texttt{Milnor} and \texttt{Tyurina}
but for an ideal $I\subset\Cx$, since they works with a g.m.o.. As
before each output may be listed separately by means of analogous
sub--procedures.

\halfline

\noindent The procedure computing $\mu(F)$:

\halfline

\begin{maplegroup}
\begin{mapleinput}
\mapleinline{active}{1d}{PolyMilnor := proc (F,variables::set:=indets(F),
U::anything := tdeg(op(variables)))
if type(U, ShortMonomialOrder) = false then error
"invalid input: your ShortMonomialOrder is not well defined":
else if localorglobal(U, variables) <> Gl then error
"invalid input: your Short Monomial Order should be GLOBAL":
else r := nops(variables);
J := [seq(diff(F, variables[i]), i = 1 .. r)];
G := Basis(J, U);
Ini := [seq(LeadingMonomial(G[i], U), i = 1 .. nops(G))];
massimo:=r*max(seq(degree(Ini[i]),i=1..nops(Ini)));
if IsZeroDimensional(<(op(Ini))>, \{op(variables)\}) = false then
error "there are non isolated critical points"
else
L := sort([op(subs(uq = 1, convert(map(expand,
series(1/(product(1-variables[k]*uq, k = 1 .. r)),
uq, massimo+1)), polynom)))],
proc (t1, t2) options operator, arrow;
TestOrder(t1, t2, U) end proc);
M := [];
for j to nops(L) do up := 1;
for k to nops(Ini) do
if divide(L[j], Ini[k]) then up := 0
end if end do;
if up = 1 then M := [op(M), L[j]]
end if end do;
[G, Ini, M, nops(M)]
end if end if end if end proc:}{}
\end{mapleinput}
\end{maplegroup}

\halfline

\noindent The associated sub--procedures:

\halfline

\begin{maplegroup}
\begin{mapleinput}
\mapleinline{active}{1d}{PolyMilnorGroebnerBasis := proc (F,variables::set := indets(F),
T::anything :=tdeg(op(variables)))
PolyMilnor(F,variables,T)[1]
end proc: }{}
\end{mapleinput}
\end{maplegroup}

\halfline

\begin{maplegroup}
\begin{mapleinput}
\mapleinline{active}{1d}{PolyMilnorLT := proc (F, variables::set:= indets(F),
T::anything := tdeg(op(variables)))
PolyMilnor(F, variables, T)[2]
end proc:}{}
\end{mapleinput}
\end{maplegroup}

\halfline

\begin{maplegroup}
\begin{mapleinput}
\mapleinline{active}{1d}{PolyMilnorBasis := proc (F,variables::set:= indets(F),
T::anything :=tdeg(op(variables)))
PolyMilnor(F,variables,T)[3]
end proc: }{}
\end{mapleinput}
\end{maplegroup}

\halfline

\begin{maplegroup}
\begin{mapleinput}
\mapleinline{active}{1d}{PolyMilnorNumber:=proc(F,variables::set:= indets(F),
T::anything :=tdeg(op(variables)))
PolyMilnor(F,variables,T)[4]
end proc: }{}
\end{mapleinput}
\end{maplegroup}

\halfline

\noindent The procedure computing $\tau(F)$:

\halfline

\begin{maplegroup}
\begin{mapleinput}
\mapleinline{active}{1d}{PolyTyurina := proc (F, variables::set :=indets(F),
U::anything := tdeg(op(variables)))
if type(U,ShortMonomialOrder) = false then error
"invalid input: your ShortMonomialOrder is not well defined":
else if localorglobal(U, variables) <> Gl then error
"invalid input: your Short Monomial Order should be GLOBAL"
else r := nops(variables);
K := [F,seq(diff(F, variables[i]), i = 1 .. r)];
H := Basis(K, U);
Ini := [seq(LeadingMonomial(H[i], U), i = 1 .. nops(H))];
massimo:=r*max(seq(degree(Ini[i]),i=1..nops(Ini)));
if IsZeroDimensional(<(op(Ini))>, \{op(variables)\}) = false then
error "there are non isolated singular points"
else
L := sort([op(subs(uq = 1, convert(map(expand,
series(1/(product(1-variables[k]*uq, k = 1 .. r)),
uq, massimo+1)), polynom)))],
proc (t1, t2) options operator, arrow;
TestOrder(t1, t2, U) end proc);
M := [];
for j to nops(L) do up := 1;
for k to nops(Ini) do
if divide(L[j], Ini[k]) then up := 0
end if end do;
if up = 1 then M := [op(M), L[j]]
end if end do;
[H, Ini, M, nops(M)]
end if end if end if end proc: }{}
\end{mapleinput}
\end{maplegroup}

\halfline

\noindent The associated sub--procedures:

\halfline

\begin{maplegroup}
\begin{mapleinput}
\mapleinline{active}{1d}{PolyTyurinaGroebnerBasis := proc(F,variables::set := indets(F),
T::anything := tdeg(op(variables)))
PolyTyurina(F,variables,T)[1]
end proc: }{}
\end{mapleinput}
\end{maplegroup}

\halfline

\begin{maplegroup}
\begin{mapleinput}
\mapleinline{active}{1d}{PolyTyurinaLT := proc (F,variables::set:= indets(F),
T::anything := tdeg(op(variables)))
PolyTyurina(F,variables,T)[2]
end proc: }{}
\end{mapleinput}
\end{maplegroup}

\halfline

\begin{maplegroup}
\begin{mapleinput}
\mapleinline{active}{1d}{PolyTyurinaBasis := proc (F,variables,
T::anything := tdeg(op(variables)))
PolyTyurina(F,variables,T)[3]
end proc: }{}
\end{mapleinput}
\end{maplegroup}

\halfline

\begin{maplegroup}
\begin{mapleinput}
\mapleinline{active}{1d}{PolyTyurinaNumber := proc (F,variables::set := indets(F),
T::anything := tdeg(op(variables)))
PolyTyurina(F, variables, T)[4]
end proc: }{}
\end{mapleinput}
\end{maplegroup}

\subsection{Some user friendly examples}\label{user_friendly}

\noindent Once implemented the subroutine \ref{subroutine} needs
quite simple and minimal commands to work. As a first example let
us start, for comparison, by a problem already studied by using SINGULAR in
\cite{Greuel-etal} Example 2.7.2(2).

\begin{example}\label{Greuel-ex} Let us study critical points of $F(x,y):=x^5+y^5+x^2y^2$ and
singularities of $F^{-1}(0)$. Hence we have to type:

\halfline

\begin{maplegroup}
\begin{mapleinput}
\mapleinline{active}{1d}{F:=x\symbol{94}5+y\symbol{94}5+x\symbol{94}2*y\symbol{94}2:
}{}
\end{mapleinput}
\end{maplegroup}

\halfline

\noindent Let us find, at first, the critical points of $F$, by solving the
algebraic system of partial derivatives:

\halfline

\begin{maplegroup}
\begin{mapleinput}
\mapleinline{active}{1d}{solve(\{diff(F, x), diff(F, y)\}, [x,y]);
}{}
\end{mapleinput}
\mapleresult
\begin{maplelatex}
\mapleinline{inert}{2d}{[[x = 0, y = 0], [x = 0, y = 0], [x = 0, y
= 0], [x = 0, y = 0], [x = 0, y = 0], }
{\[\displaystyle
[[x=0,y=0],[x=0,y=0],[x=0,y=0],[x=0,y=0],[x=0,y=0],\]}
\mapleinline{inert}{2d}{ [x = 0, y = 0], [x = 0, y = 0], [x = 0, y
= 0], [x = -2/5, y = -2/5], }
{\[\displaystyle
[x=0,y=0],[x=0,y=0],[x=0,y=0],[x=-2/5,y=-2/5],\]}
\mapleinline{inert}{2d}{[x = 2/5-(2/5)*RootOf(_Z^4-_Z^3+_Z^2-_Z+1,
label = _L2)^3}
{\[\displaystyle [x=2/5-2/5\, \left( {\it RootOf}
\left( {{\it \_Z}}^{4}-{{\it \_Z}}^{3}+{{\it \_Z}}^{2}-{\it
\_Z}+1,{\it label}= {\it \_L2} \right)  \right) ^{3}\]}
\mapleinline{inert}{2d}{+(2/5)*RootOf(_Z^4-_Z^3+_Z^2-_Z+1, label =
_L2)^2}
{\[\displaystyle \quad + 2/5\, \left( {\it RootOf} \left( {{\it
\_Z}}^{4}- {{\it \_Z}}^{3}+{{\it \_Z}}^{2}-{\it \_Z}+1,{\it
label}={\it \_L2} \right) \right) ^{2}\]}
\mapleinline{inert}{2d}{-(2/5)*RootOf(_Z^4-_Z^3+_Z^2-_Z+1, label =
_L2),}
{\[\displaystyle \quad -2/5\,{\it RootOf} \left( {{\it \_Z}}^{4}
-{{\it \_Z}}^{3}+{{\it \_Z}}^{2}-{\it \_Z}+1,{\it label}={\it
\_L2} \right) ,\]}
\mapleinline{inert}{2d}{ y =
(2/5)*RootOf(_Z^4-_Z^3+_Z^2-_Z+1, label = _L2)]]}
{\[\displaystyle
\ y=2/5\,{\it RootOf} \left( {{\it \_Z}}^{4}-{{\it \_Z}}^{3}+{{\it
\_Z}}^{2}-{\it \_Z}+1,{\it label}={\it \_L2} \right) ]]\]}
\end{maplelatex}
\end{maplegroup}

\halfline

\noindent Then $F$ admits 6 critical points: the repetition of the
solution in the origin means that this point is a multiple
solution. Then we have to expect $\mu(0)>0$.

\noindent Singular points of $F^{-1}(0)$ are given by:

\halfline

\begin{maplegroup}
\begin{mapleinput}
\mapleinline{active}{1d}{solve(\{F, diff(F, x), diff(F, y)\}, [x,y]) }{}
\end{mapleinput}
\mapleresult
\begin{maplelatex}
\mapleinline{inert}{2d}{[[x = 0, y = 0], [x = 0, y = 0], [x = 0, y = 0], [x = 0, y = 0], [x = 0, y = 0],}
{\[\displaystyle [[x=0,y=0],[x=0,y=0],[x=0,y=0],[x=0,y=0],[x=0,y=0],\]}
\mapleinline{inert}{2d}{[x = 0, y = 0], [x = 0, y = 0], [x = 0, y = 0], [x = 0, y = 0], [x = 0, y = 0],}
{\[\displaystyle [x=0,y=0],[x=0,y=0],[x=0,y=0],[x=0,y=0],[x=0,y=0],\]}
\mapleinline{inert}{2d}{[x = 0, y = 0], [x = 0, y = 0], [x = 0, y = 0], [x = 0, y = 0], [x = 0, y = 0],}
{\[\displaystyle [x=0,y=0],[x=0,y=0],[x=0,y=0],[x=0,y=0],[x=0,y=0],\]}
\mapleinline{inert}{2d}{[x = 0, y = 0], [x = 0, y = 0], [x = 0, y = 0], [x = 0, y = 0], [x = 0, y = 0],}
{\[\displaystyle [x=0,y=0],[x=0,y=0],[x=0,y=0],[x=0,y=0],[x=0,y=0],\]}
\mapleinline{inert}{2d}{[x = 0, y = 0], [x = 0, y = 0], [x = 0, y = 0], [x = 0, y = 0], [x = 0, y = 0],}
{\[\displaystyle [x=0,y=0],[x=0,y=0],[x=0,y=0],[x=0,y=0],[x=0,y=0],\]}
\mapleinline{inert}{2d}{[x = 0, y = 0], [x = 0, y = 0], [x = 0, y = 0], [x = 0, y = 0], [x = 0, y = 0],}
{\[\displaystyle [x=0,y=0],[x=0,y=0],[x=0,y=0],[x=0,y=0],[x=0,y=0],\]}
\mapleinline{inert}{2d}{[x = 0, y = 0], [x = 0, y = 0], [x = 0, y = 0], [x = 0, y = 0], [x = 0, y = 0],}
{\[\displaystyle [x=0,y=0],[x=0,y=0],[x=0,y=0],[x=0,y=0],[x=0,y=0],\]}
\mapleinline{inert}{2d}{[x = 0, y = 0]]}
{\[\displaystyle [x=0,y=0]]\]}
\end{maplelatex}
\end{maplegroup}

\halfline

\noindent Then $F^{-1}(0)$ has a unique singular point in the
origin. Therefore, by the second formula in (\ref{loc-to-glob}),
$\tau(F)=\tau(0)$ meaning that \texttt{TyurinaNumber} and
\texttt{PolyTyurinaNumber} give the same number. To get this
number one simply have to type either

\halfline

\begin{maplegroup}
\begin{mapleinput}
\mapleinline{active}{1d}{TyurinaNumber(F) }{}
\end{mapleinput}
\mapleresult
\begin{maplelatex}
\mapleinline{inert}{2d}{10}{\[\displaystyle 10\]}
\end{maplelatex}
\end{maplegroup}

\halfline

\noindent or

\halfline

\begin{maplegroup}
\begin{mapleinput}
\mapleinline{active}{1d}{PolyTyurinaNumber(F) }{}
\end{mapleinput}
\mapleresult
\begin{maplelatex}
\mapleinline{inert}{2d}{10}{\[\displaystyle 10\]}
\end{maplelatex}
\end{maplegroup}

\halfline

\begin{itemize}
    \item \emph{In general, if $F^{-1}(0)$ admits a unique singular point
    the procedure} \newline
    \texttt{PolyTyurina} \emph{has to be preferred, since it turns
    out to be more efficient.}
\end{itemize}
To compute Milnor numbers let us start by $\mu(F)$, by typing

\halfline

\begin{maplegroup}
\begin{mapleinput}
\mapleinline{active}{1d}{PolyMilnorNumber(F) }{}
\end{mapleinput}
\mapleresult
\begin{maplelatex}
\mapleinline{inert}{2d}{16}{\[\displaystyle 16\]}
\end{maplelatex}
\end{maplegroup}

\halfline

\noindent Since $\tau(0)\leq\mu(0)$ and $F$ admits 6 critical
points, by (\ref{loc-to-glob}) we have to expect $10\leq\mu(0)\leq
11$. Furthermore the 5 critical points different from the origin can be exchanged each other under the action of the order 5 cyclic group
$$\left\langle\left(
                            \begin{array}{cc}
                              -\epsilon^3 & 0 \\
                              0 & \epsilon^2 \\
                            \end{array}
                          \right)^i\ |\ \epsilon^5 + 1=0\ ,\ 1\leq i\leq 5
\right\rangle \subset \GL (2,\C) $$
which is also a subgroup of $\Aut (F)$. Then they
cannot assume Milnor number greater than 1, giving $\mu(0)=11$. In
fact:

\halfline

\begin{maplegroup}
\begin{mapleinput}
\mapleinline{active}{1d}{MilnorNumber(F) }{}
\end{mapleinput}
\mapleresult
\begin{maplelatex}
\mapleinline{inert}{2d}{11}{\[\displaystyle 11\]}
\end{maplelatex}
\end{maplegroup}

\halfline

\noindent Then in this case, the use of \texttt{solve},
\texttt{PolyTyurina} and \texttt{PolyMilnor}, may avoid to employ
\texttt{Tyurina} and \texttt{Milnor} \emph{which turn out to be in
general less efficient procedures.}
\end{example}

\begin{example}[w.h. polynomials] What observed at the end of the
previous Example \ref{Greuel-ex} \emph{is obviously true for a
w.h.p.}, after Proposition \ref{MT-whp}. Let us in fact consider
the $E_6$ 3--dimensional singularity $0\in F^{-1}(0)$ where
\[
    F(x,y,z,t)= x^2+y^3+z^4+t^2\ .
\]
By Corollary \ref{Milnor w.h.}, since $w(t)=1/2$, one has to
expect $\tau(0)=\mu(0)=6$. This fact can be checked by the quicker procedure \texttt{PolyMilnor}:

\halfline

\begin{maplegroup}
\begin{mapleinput}
\mapleinline{active}{1d}{F:=x\symbol{94}2+y\symbol{94}3+z\symbol{94}4+t\symbol{94}2:
}{}
\end{mapleinput}
\end{maplegroup}
\begin{maplegroup}
\begin{mapleinput}
\mapleinline{active}{1d}{PolyMilnorNumber(F);}{}
\end{mapleinput}
\mapleresult
\begin{maplelatex}
\mapleinline{inert}{2d}{6}{\[\displaystyle 6\]}
\end{maplelatex}
\end{maplegroup}
\end{example}

\begin{example}[Non--isolated singularities] All the procedures
presented in \ref{MTproc} stop, giving an error message, if the
considered polynomial admits \emph{non--isolated singularities}.
Consider, in fact, $F= x^2z^2+y^2z^2+ x^2y^2$ admitting the union
of the three coordinate axes as the locus of critical (and
singular, since $F$ is homogeneous) points as can easily checked
by typing:

\halfline

\begin{maplegroup}
\begin{mapleinput}
\mapleinline{active}{1d}{F:=x\symbol{94}2*z\symbol{94}2+y\symbol{94}2*z\symbol{94}2+x\symbol{94}2*y\symbol{94}2: }{}
\end{mapleinput}
\end{maplegroup}
\begin{maplegroup}
\begin{mapleinput}
\mapleinline{active}{1d}{solve(\{diff(F, x), diff(F, y), diff(F,z)\}, [x, y, z]); }{}
\end{mapleinput}
\mapleresult
\begin{maplelatex}
\mapleinline{inert}{2d}{[[x = 0, y = y, z = 0], [x = 0, y = y, z =
0], [x = x, y = 0, z = 0], [x = 0, y = 0, z = z]]}{\[\displaystyle
[[x=0,y=y,z=0],[x=0,y=y,z=0],[x=x,y=0,z=0],[x=0,y=0,z=z]]\]}
\end{maplelatex}
\end{maplegroup}

\halfline

\noindent Then:

\halfline

\begin{maplegroup}
\begin{mapleinput}
\mapleinline{active}{1d}{PolyMilnorNumber(F); }{}
\end{mapleinput}
\mapleresult
\begin{Maple Error}{
Error, (in PolyMilnor) there are non isolated critical points}\end{Maple
Error}

\end{maplegroup}

\halfline

\noindent The user will obtain the similar error messages by running
any of the other procedures.
\end{example}

\begin{remark}[Be careful with variables!]\label{rem_variabili} The second input of any
procedure in \ref{MTproc} is the set \texttt{variables} of
variables one wants to work with. It is an \emph{optional input}
meaning that by default \texttt{variables} is assumed to be the
set \texttt{indets(F)} of variables appearing in the polynomial
$F$. This means that if the user is interesting in consider the
cylinder $F^{-1}(0)$ where $F:\C^3\rightarrow\C$ is the polynomial
map $F(x,y,z)= y^2-x(x-1)(x-2)$ then he have to type:

\halfline

\begin{maplegroup}
\begin{mapleinput}
\mapleinline{active}{1d}{F := y\symbol{94}2-x*(x-1)*(x-2): }{}
\end{mapleinput}
\end{maplegroup}
\begin{maplegroup}
\begin{mapleinput}
\mapleinline{active}{1d}{PolyMilnor(F, \{x, y, z\}); }{}
\end{mapleinput}
\mapleresult
\begin{Maple Error}{
Error, (in PolyMilnor) there are non isolated critical points}\end{Maple
Error}
\end{maplegroup}

\noindent which is right since $F$ do not admit isolated critical
points as can be checked by:

\halfline

\begin{maplegroup}
\begin{mapleinput}
\mapleinline{active}{2d}{}{$$}
\mapleinline{active}{1d}{solve(\{diff(F, x), diff(F, y), diff(F,z)\}, [x,y,z]); }{}
\end{mapleinput}
\mapleresult
\begin{maplelatex}
\mapleinline{inert}{2d}{[[x = RootOf(3*_Z^2-6*_Z+2), y = 0, z =
z]]}{\[\displaystyle [[x={\it RootOf} \left( 3\,{{\it
\_Z}}^{2}-6\,{\it \_Z}+2 \right) ,y=0,z=z]]\]}
\end{maplelatex}
\end{maplegroup}

\halfline

\noindent By the way, if the set of variables is not specified
then by default it is assumed to be $\{x,y\}$, meaning that $F$ is
considered as a polynomial map from $\C^2$ to $\C$. In this case
$F$ admits only the two isolated critical points

\halfline

\begin{maplegroup}
\begin{mapleinput}
\mapleinline{active}{1d}{solve({diff(F, x), diff(F, y)}, [x,y])}{}
\end{mapleinput}
\mapleresult
\begin{maplelatex}
\mapleinline{inert}{2d}{[[x = RootOf(3*_Z^2-6*_Z+2, label = _L1), y = 0]]}
{\[\displaystyle [[x={\it RootOf} \left( 3\,{{\it \_Z}}^{2}-6\,{\it \_Z}+2,{\it label}={\it \_L1}
 \right) ,y=0]]\]}
\end{maplelatex}
\end{maplegroup}

\halfline

\noindent In fact

\halfline

\begin{maplegroup}
\begin{mapleinput}
\mapleinline{active}{1d}{PolyMilnor(F); }{}
\end{mapleinput}
\mapleresult
\begin{maplelatex}
\mapleinline{inert}{2d}{[[y, 3*x^2-6*x+2], [y, x^2], [1, x],
2]}{\[\displaystyle [[y,3\,{x}^{2}-6\,x+2],[y,{x}^{2}],[1,x],2]\]}
\end{maplelatex}
\end{maplegroup}

\halfline

\noindent Read the output as follows: the first output is the
Gr\"{o}bner basis $G$ of $J_F$, the second output is the list of
leading monomials of elements in $G$, the third output is $M\setminus L(J_F)$ whose cardinality is precisely the fourth output. Since the problem is
symmetric w.r.t. the $y$--axis, $\mu(F)=2$ implies that each
critical point has Milnor number 1.

\noindent Observe that the origin is not a critical point of $F$
both as a polynomial map from $\C^3$ and from $\C^2$. In fact

\halfline

\begin{maplegroup}
\begin{mapleinput}
\mapleinline{active}{1d}{Milnor(F);}{}
\end{mapleinput}
\mapleresult
\begin{maplelatex}
\mapleinline{inert}{2d}{[[-(x-1)*(x-2)-x*(x-2)-x*(x-1), 2*y], [1,
y], [], 0]}{\[\displaystyle [[- \left( x-1 \right)  \left( x-2
\right) -x \left( x-2 \right) -x \left( x-1 \right)
,2\,y],[1,y],[],0]\]}
\end{maplelatex}
\end{maplegroup}
\begin{maplegroup}
\begin{mapleinput}
\mapleinline{active}{1d}{Milnor(F, \{x, y, z\}) }{}
\end{mapleinput}
\mapleresult
\begin{maplelatex}
\mapleinline{inert}{2d}{[[0, -(x-1)*(x-2)-x*(x-2)-x*(x-1), 2*y],
[1, 1, y], [], 0]}{\[\displaystyle [[0,- \left( x-1 \right)
\left( x-2 \right) -x \left( x-2 \right) -x \left( x-1 \right)
,2\,y],[1,1,y],[],0]\]}
\end{maplelatex}
\end{maplegroup}

\halfline

\noindent In this case the first output is a standard bases of
$J_F$ whose leading monomials, w.r.t. a fixed l.m.o.,  give the
second output.

\noindent At last let us observe that the zero locus $F^{-1}(0)$
is smooth both as a subset of $\C^2$ and of $\C^3$, in fact

\halfline

\begin{maplegroup}
\begin{mapleinput}
\mapleinline{active}{1d}{PolyTyurina(F);}{}
\end{mapleinput}
\mapleresult
\begin{maplelatex}
\mapleinline{inert}{2d}{[[1], [1], [], 0]}{\[\displaystyle
[[1],[1],[],0]\]}
\end{maplelatex}
\end{maplegroup}
\begin{maplegroup}
\begin{mapleinput}
\mapleinline{active}{1d}{PolyTyurina(F, \{x, y, z\}); }{}
\end{mapleinput}
\mapleresult
\begin{maplelatex}
\mapleinline{inert}{2d}{[[1], [1], [], 0]}{\[\displaystyle
[[1],[1],[],0]\]}
\end{maplelatex}
\end{maplegroup}
\end{remark}

\subsection{Optional input: some more subtle utilities}\label{optional}

All the procedure described in \ref{MTproc} require three input:
the polynomial $F$ and two further optional input, precisely
\begin{itemize}
    \item a set of variables, by default set as the variables \texttt{indets(F)}
    appearing in $F$,
    \item a monomial order, by default set either as \texttt{tdeg(variables)}, which is the
    MAPLE command for the graduated reverse lexicographic g.m.o., or as \texttt{tdeg\_\
    min(variables)}, which is the l.m.o. defined as the
    \texttt{tdeg} opposite: the former is clearly introduced in
    \texttt{PolyMilnor} and \texttt{PolyTyurina} and the latter in
    \texttt{Milnor} and \texttt{Tyurina}.
\end{itemize}
Introducing different choices may show interesting possibilities
of our subroutine.

\subsubsection{Monomial ordering}
It is a well known fact that
the graduated reverse lexicographic g.m.o. is in general the more
efficient monomial ordering for Buchberger algorithm: this is the
reason for the default choices in \texttt{PolyMilnor} and
\texttt{PolyTyurina}. Anyway, if needed, these procedures may run
with many further g.m.o.: e.g. if, for any reason, the user will
prefer to run \texttt{PolyTyurina} w.r.t. the pure lexicographic
g.m.o. he will have to type, in the following case of a
deformation of a threefold $E_8$ singularity:

\halfline

\begin{maplegroup}
\begin{mapleinput}
\mapleinline{active}{1d}{F:=x\symbol{94}2+y\symbol{94}3+z\symbol{94}5+t\symbol{94}2+y*z+y\symbol{94}2+z\symbol{94}2+y*z\symbol{94}2+z\symbol{94}3+y*z\symbol{94}3+z\symbol{94}4;
}{}
\end{mapleinput}
\mapleresult
\begin{maplelatex}
\mapleinline{inert}{2d}{F :=
x^2+y^3+z^5+t^2+y*z^2+z^3+y*z^3+z^4}{\[\displaystyle F\, :=
\,{x}^{2}+{y}^{3}+{z}^{5}+{t}^{2}+y{z}^{2}+{z}^{3}+y{z}^{3}+{z}^{4}\]}
\end{maplelatex}
\end{maplegroup}
\begin{maplegroup}
\begin{mapleinput}
\mapleinline{active}{1d}{PolyMilnor(F); }{}
\end{mapleinput}
\mapleresult
\begin{maplelatex}
\mapleinline{inert}{2d}{[[x, t, 3*y*z^2-15*y^2*z+4*z^2+2*y*z+3*y^2, 3*y^2+z^2+z^3,}
{\[\displaystyle [[x,t,3\,y{z}^{2}-15\,{y}^{2}z+4\,{z}^{2}+2\,yz+3\,{y}^{2},3\,{y}^{2}+{z}^{2}+{z}^{3},\]}
\mapleinline{inert}{2d}{3375*y^4-1677*y^2*z-639*y^3+224*z^2+124*y*z+114*y^2,}
{\[\displaystyle 3375\,{y}^{4}-1677\,{y}^{2}z-639\,{y}^{3}+224\,{z}^{2}+124\,yz+114\,{y}^{2},\]}
\mapleinline{inert}{2d}{225*y^3*z-324*y^2*z-18*y^3+88*z^2+38*y*z+93*y^2],}
{\[\displaystyle 225\,{y}^{3}z-324\,{y}^{2}z-18\,{y}^{3}+88\,{z}^{2}+38\,yz+93\,{y}^{2}], \]}
\mapleinline{inert}{2d}{[x, t, y*z^2, z^3, y^4, y^3*z], [1, y, z, y^2, y*z, z^2, y^3, y^2*z], 8]}
{\[\displaystyle [x,t,y{z}^{2},{z}^{3},{y}^{4},{y}^{3}z],[1,y,z,{y}^{2},yz,{z}^{2},{y}^{3},{y}^{2}z],8]\]}
\end{maplelatex}
\end{maplegroup}
\begin{maplegroup}
\begin{mapleinput}
\mapleinline{active}{1d}{PolyMilnor(F,plex(x,y,z,t)); }{}
\end{mapleinput}
\mapleresult
\begin{maplelatex}
\mapleinline{inert}{2d}{[[t, 31*z^3+88*z^4+159*z^5+129*z^6+75*z^7,}
{\[\displaystyle [[t,31\,{z}^{3}+88\,{z}^{4}+159\,{z}^{5}+129\,{z}^{6}+75\,{z}^{7},\]}
\mapleinline{inert}{2d}{128*z^3+319*z^4+237*z^5+225*z^6+69*z^2+46*y*z, 3*y^2+z^2+z^3, x],}
{\[\displaystyle 128\,{z}^{3}+319\,{z}^{4}+237\,{z}^{5}+225\,{z}^{6}+69\,{z}^{2}+46\,yz,3\,{y}^{2}+{z}^{2}+{z}^{3},x],\]}
\mapleinline{inert}{2d}{[t, z^7, y*z, y^2, x], [1, z, z^2, z^3, z^4, z^5, z^6, y], 8]}
{\[\displaystyle [t,{z}^{7},yz,{y}^{2},x],[1,z,{z}^{2},{z}^{3},{z}^{4},{z}^{5},{z}^{6},y],8]\]}
\end{maplelatex}
\end{maplegroup}

\halfline

\noindent Observe how different are the two Gr\"{o}bner bases and
consequently the leading monomial bases and associated
bases of quotient vector spaces. The user may also verify how much
slower is \texttt{plex} w.r.t. the default \texttt{tdeg} by
running by himself the routines.

\noindent In particular, running \texttt{Tyurina} w.r.t. different
l.m.o.'s gives \emph{different monomial basis of the Kuranishi
space}:

\halfline

\begin{maplegroup}
\begin{mapleinput}
\mapleinline{active}{1d}{TyurinaBasis(F); }{}
\end{mapleinput}
\mapleresult
\begin{maplelatex}
\mapleinline{inert}{2d}{[z^2, z, y, 1]}
{\[\displaystyle [{z}^{2},z,y,1]\]}
\end{maplelatex}
\end{maplegroup}
\begin{maplegroup}
\begin{mapleinput}
\mapleinline{active}{1d}{TyurinaBasis(F,plex_min(x,y,z,t)); }{}
\end{mapleinput}
\mapleresult
\begin{maplelatex}
\mapleinline{inert}{2d}{[y^2, y, z, 1]}
{\[\displaystyle [{y}^{2},y,z,1]\]}
\end{maplelatex}
\end{maplegroup}

\halfline

\begin{remark} For what concerns efficiency of l.m.o.'s in \texttt{Milnor} and \texttt{Tyurina} it turns out that sometimes \texttt{plex\_min} is more efficient than \texttt{tdeg\_min}, as observed in the last Section \ref{Arnold-ex} when proving Theorems \ref{E7} and \ref{E8}. But we do not know if this is a general fact, then we keep \texttt{tdeg\_min} as the default l.m.o. both in \texttt{Milnor} and \texttt{Tyurina}, for coherence with the default choice of \texttt{tdeg} for their global counterparts.
\end{remark}

\subsubsection{Variables}\label{variabili} The default choice for the optional input
\texttt{variables} as the set of those variables appearing in the
given polynomial $F$ (\texttt{indets(F)}) has been thought to make
our routine more user friendly. Anyway this choice may hide some
important subtleties, as already pointed out in Remark
\ref{rem_variabili} in the case \texttt{variables} has been chosen
as a greater set of variables w.r.t. \texttt{indets(F)}. Here we
want to underline a significant potentiality of our routine when
\texttt{variables} is chosen to be a \emph{strictly smaller subset
of} \texttt{indets(F)}.

\noindent Let us set $F=x^3+x^4+xy^2$. Then we get:

\halfline

\begin{maplegroup}
\begin{mapleinput}
\mapleinline{active}{1d}{F :=x\symbol{94}3+y\symbol{94}4+x*y\symbol{94}2: }{}
\end{mapleinput}
\end{maplegroup}
\begin{maplegroup}
\begin{mapleinput}
\mapleinline{active}{1d}{Milnor(F); }{}
\end{mapleinput}
\mapleresult
\begin{maplelatex}
\mapleinline{inert}{2d}{[{6*x^3-4*y^4, 3*x^2+y^2, 4*y^3+2*x*y},
[y^2, x*y, x^3], [x^2, x, y, 1], 4]}{\[\displaystyle [ \left\{
6\,{x}^{3}-4\,{y}^{4},3\,{x}^{2}+{y}^{2},4\,{y}^{3}+2\,xy \right\}
,[{y}^{2},xy,{x}^{3}],[{x}^{2},x,y,1],4]\]}
\end{maplelatex}
\end{maplegroup}
\begin{maplegroup}
\begin{mapleinput}
\mapleinline{active}{1d}{Tyurina(F); }{}
\end{mapleinput}
\mapleresult
\begin{maplelatex}
\mapleinline{inert}{2d}{[{x^3+y^4+x*y^2, 3*x^2+y^2, 4*y^3+2*x*y,
2*x^3-y^4}, [y^2, x*y, x*y^2, x^3], [x^2, x, y, 1],
4]}{\[\displaystyle [ \left\{
{x}^{3}+{y}^{4}+x{y}^{2},3\,{x}^{2}+{y}^{2},4\,{y}^{3}+2\,xy,2\,{x}^{3}-{y}^{4}
\right\} ,[{y}^{2},xy,x{y}^{2},{x}^{3}],[{x}^{2},x,y,1],4]\]}
\end{maplelatex}
\end{maplegroup}

\halfline

\noindent Then $T^1\cong\langle x^2,x,y,1\rangle_{\C}$ and
\begin{equation}\label{1-parameter}
    F_t:= F + t x^2 = x^3+y^4+xy^2 + tx^2
\end{equation}
is a non--trivial 1--parameter small deformation of $F$ such that,
for any fixed $t\in\C$, $F_t$ has a critical point in $0\in\C^2$
which is also a singular point of the plane curve $F_t^{-1}(0)$.
We are interested in studying Milnor and Tyurina numbers of this
singularity \emph{for any $t\in\C$}\ . This can be performed by a
careful use of the variables input. Let us first of all observe
that if no optional input are added then we get

\halfline

\begin{maplegroup}
\begin{mapleinput}
\mapleinline{active}{1d}{Ft := F+t*x\symbol{94}2; }{}
\end{mapleinput}
\mapleresult
\begin{maplelatex}
\mapleinline{inert}{2d}{Ft := x^3+y^4+x*y^2+t*x^2}{\[\displaystyle
{\it Ft}\, := \,{x}^{3}+{y}^{4}+x{y}^{2}+t{x}^{2}\]}
\end{maplelatex}
\end{maplegroup}
\begin{maplegroup}
\begin{mapleinput}
\mapleinline{active}{1d}{Milnor(Ft); }{}
\end{mapleinput}
\mapleresult
\begin{Maple Error}{
Error, (in Milnor) the given critical point is not
isolated}\end{Maple Error}
\end{maplegroup}

\halfline

\noindent In fact, by default \texttt{Milnor} considers $F_t$ as a
polynomial map defined over $\C^3(x,y,t)$. By forcing
\texttt{Milnor} to work with variables $\{x,y\}$ only, then $F_t$
is considered as a polynomial map defined over $\C^2$ with
coefficient ring $\C[t]$, i.e. $F_t\in\left(\C[t]\right)[x,y]$,
giving:

\halfline

\begin{maplegroup}
\begin{mapleinput}
\mapleinline{active}{1d}{Milnor(Ft, \{x, y\}) }{}
\end{mapleinput}
\mapleresult
\begin{maplelatex}
\mapleinline{inert}{2d}{[{3*x^2+y^2+2*t*x, (1-4*t)*y^3+3*y*x^2,
4*y^3+2*x*y}, [x, y^3, x*y], [y^2, y, 1], 3]}{\[\displaystyle [
\left\{ 3\,{x}^{2}+{y}^{2}+2\,tx, \left( 1-4\,t \right)
{y}^{3}+3\,y{x}^{2},4\,{y}^{3}+2\,xy \right\}
,[x,{y}^{3},xy],[{y}^{2},y,1],3]\]}
\end{maplelatex}
\end{maplegroup}
\begin{maplegroup}
\begin{mapleinput}
\mapleinline{active}{1d}{Tyurina(Ft, \{x, y\}) }{}
\end{mapleinput}
\mapleresult
\begin{maplelatex}
\mapleinline{inert}{2d}{[{3*x^2+y^2+2*t*x, x^3+y^4+x*y^2+t*x^2,
(1-4*t)*y^3+3*y*x^2, 4*y^3+2*x*y},}
{\[\displaystyle [ \left\{ 3\,{x}^{2}+{y}^{2}+2\,tx,{x}^{3}+{y}^{4}
+x{y}^{2}+t{x}^{2}, \left( 1-4\,t \right) {y}^{3}+3\,y{x}^{2},4\,{y}^{3}+2\,xy \right\} ,
\]}
\mapleinline{inert}{2d}{ [x, y^3, x^2, x*y], [y^2, y, 1], 3]}
{\[\displaystyle [x,{y}^{3},{x}^{2},xy],[{y}^{2},y,1],3] \]}
\end{maplelatex}
\end{maplegroup}

\halfline

\noindent This means that, \emph{for generic \footnote{$t$ is
treated as a variable without any evaluation.} $t$},
 $\tau_{F_t}(0)=3=\mu_{F_t}(0)$.  Moreover \emph{by looking at the leading coefficients of
the given standard basis of $J_F$ we get all the relations
defining non--generic values for $t$}, precisely:

\halfline

\begin{maplegroup}
\begin{mapleinput}
\mapleinline{active}{1d}{MB := MilnorGroebnerBasis(Ft, \{x, y\})
}{}
\end{mapleinput}
\mapleresult
\begin{maplelatex}
\mapleinline{inert}{2d}{MB := {3*x^2+y^2+2*t*x,
(1-4*t)*y^3+3*y*x^2, 4*y^3+2*x*y}}{\[\displaystyle {\it MB}\, :=
\, \left\{ 3\,{x}^{2}+{y}^{2}+2\,tx, \left( 1-4\,t \right)
{y}^{3}+3\,y{x}^{2},4\,{y}^{3}+2\,xy \right\} \]}
\end{maplelatex}
\end{maplegroup}
\begin{maplegroup}
\begin{mapleinput}
\mapleinline{active}{1d}{for i from 1 to nops(MB) do
LeadingTerm(MB[i],tdeg_min(x,y)) end do; }{}
\end{mapleinput}
\mapleresult
\begin{maplelatex}
\mapleinline{inert}{2d}{2*t, x}{\[\displaystyle 2\,t,\,x\]}
\end{maplelatex}
\mapleresult
\begin{maplelatex}
\mapleinline{inert}{2d}{1-4*t, y^3}{\[\displaystyle
1-4\,t,\,{y}^{3}\]}
\end{maplelatex}
\mapleresult
\begin{maplelatex}
\mapleinline{inert}{2d}{2, x*y}{\[\displaystyle 2,\,xy\]}
\end{maplelatex}
\end{maplegroup}

\halfline

\noindent For $t=0$ we do not have any deformation of $F$, giving
$\tau_{F_0}(0)=4=\mu_{F_0}$. But the further relation $1-4t=0$
gives:

\halfline

\begin{maplegroup}
\begin{mapleinput}
\mapleinline{active}{1d}{t := 1/4: Milnor(Ft, \{x, y\}) }{}
\end{mapleinput}
\mapleresult
\begin{maplelatex}
\mapleinline{inert}{2d}{[{3*x^2+y^2+(1/2)*x, 36*y^3*x^2+12*y^5,
4*y^3+2*x*y}, [x, y^5, x*y], [y^4, y^3, y^2, y, 1],
5]}{\[\displaystyle [ \left\{
3\,{x}^{2}+{y}^{2}+1/2\,x,36\,{y}^{3}{x}^{2}+12\,{y}^{5},4\,{y}^{3}+2\,xy
\right\} ,[x,{y}^{5},xy],[{y}^{4},{y}^{3},{y}^{2},y,1],5]\]}
\end{maplelatex}
\end{maplegroup}
\begin{maplegroup}
\begin{mapleinput}
\mapleinline{active}{1d}{Tyurina(Ft, \{x, y\}) }{}
\end{mapleinput}
\mapleresult
\begin{maplelatex}
\mapleinline{inert}{2d}{}
{\[\displaystyle \left[ \left\{ {x}^{3}+{y}^{4}+x{y}^{2}+1/4\,{x}^{2},{\frac {3}{64}}\,{y}^{4}{x}^{2}+{\frac {1}{64}}\,{y}^{6}\\
\mbox{}-{\frac {27}{256}}\,{x}^{6}-{\frac
{9}{256}}\,{x}^{4}{y}^{2},4\,{y}^{3}+2\,xy,\right.\right.\]}
\mapleinline{inert}{2d}{}
{\[\displaystyle \left.\left.\quad 36\,{y}^{3}{x}^{2}+12\,{y}^{5},3\,{x}^{2}+{y}^{2}+1/2\,x \right\}
,[x,{y}^{6},{x}^{2},{y}^{5},xy],[{y}^{4},{y}^{3},{y}^{2},y,1],5\right]\]}
\end{maplelatex}
\end{maplegroup}

\halfline

\noindent Therefore
\[
    \tau_{F_t}(0)=\mu_{F_t}(0)=\left\{\begin{array}{cc}
      5 & \text{for $t=1/4$}\ , \\
      4 & \text{for $t=0$}\ , \\
      3 & \text{otherwise}\ . \\
    \end{array}\right.
\]
In particular Proposition \ref{whs-caratterizzazione} implies
that, \emph{for any $t$, $0\in\Spec\mathcal{O}_{F_t,0}$ is a w.h.
singularity}, in spite of the fact that $F_t$ is never a w.h.
polynomial.

\subsection{An efficiency remark: global to local subroutines}

After numerous applications of the previous routines the reader
will convince himself that the \texttt{SB} procedure turns out to
be less efficient than the Buchberger algorithm as implemented in
MAPLE. As a consequence our routines can be arranged in the
following decreasing sequence of efficiency:
\[
\text{\texttt{PolyMilnor}} > \text{\texttt{PolyTyurina}} >
\text{\texttt{Milnor}} > \text{\texttt{Tyurina}}\ .
\]
A slight improvement of \texttt{Milnor} and \texttt{Tyurina}
efficiency can be obtained by applying the \texttt{SB} algorithm
to a Gr\"{o}bner basis of $J_f$ and $I_f$ rather than to their
original generators. What is obtained is a sort of ``pasting" of
global and local routines, precisely:
\begin{itemize}
    \item \texttt{MILNOR}\quad which is a procedure computing
     Milnor numbers of both a polynomial $f$ and of a critical point of $f$:
    \begin{eqnarray*}
      &\text{Input}& :\quad F\in\Cx\ ,\ \text{\texttt{Basis} in the \texttt{Groebner} package}\ ,
      \ \text{\texttt{SB}}\ , \\
      &&\quad \text{Optional: \texttt{variables}, by default assigned by
      \texttt{indets(F)}}\\
      &&\quad \text{Optional: a couple [l.m.o., g.m.o.], by default assigned
      by}\\
      &&\hskip2truecm\text{ [\texttt{tdeg\_min(variables), tdeg(variables)}]}\\
       &\text{Output}& :\  \text{a Gr\"{o}bner basis $G$ of}\
       J_F\subset\left(\C[\mathbf{x}\setminus{\text{\texttt{variables}}}]
      \right)[\text{\texttt{variables}}],\ L_{Gl}(G),\\
      &&\quad\quad M\setminus L_{Gl}(J_F)  \ ,\ \mu(F),\\
      &&\quad \text{a standard basis $B$ of}\
       J_F\subset\left(\C[\mathbf{x}\setminus{\text{\texttt{variables}}}]
      \right)[\text{\texttt{variables}}],\ L_{Lo}(B),\\
      &&\quad \quad M\setminus L_{Lo}(J_F)  \ ,\ \mu_F(0);
    \end{eqnarray*}
    \item \texttt{TYURINA}\quad which is a the procedure computing
    the Tyurina numbers of both a polynomial $f$ and of a critical point of $f$:
    \begin{eqnarray*}
      &\text{Input}& :\quad F\in\Cx\ ,\ \text{\texttt{Basis} in the \texttt{Groebner} package}\ ,
      \ \text{\texttt{SB}}\ , \\
      &&\quad \text{Optional: \texttt{variables}, by default assigned by
      \texttt{indets(F)}}\\
      &&\quad \text{Optional: a couple [l.m.o., g.m.o.], by default assigned
      by}\\
      &&\hskip2truecm\text{ [\texttt{tdeg\_min(variables), tdeg(variables)}]}\\
       &\text{Output}& :\ \text{a Gr\"{o}bner basis $G$ of}\
       J_F\subset\left(\C[\mathbf{x}\setminus{\text{\texttt{variables}}}]
      \right)[\text{\texttt{variables}}],\ L_{Gl}(G),\\
      &&\quad  \quad M\setminus L_{Gl}(I_F)  \ ,\ \tau(F),\\
      &&\quad \text{a standard basis $B$ of}\
       I_F\subset\left(\C[\mathbf{x}\setminus{\text{\texttt{variables}}}]
      \right)[\text{\texttt{variables}}],\ L_{Lo}(B),\\
      &&\quad \quad M\setminus L_{Lo}(I_F)  \ ,\ \tau_F(0);
    \end{eqnarray*}
\end{itemize}
where $L_{Lo}$ and $L_{Gl}$ mean \emph{leading monomials w.r.t the
given local and global m.o., respectively}.

\subsubsection{The MAPLE details}

Let us start with the procedure computing Milnor numbers:

\oneline

\begin{maplegroup}
\begin{mapleinput}
\mapleinline{active}{1d}{MILNOR := proc (F, variables::set := indets(F),
SMOS::anything := [tdeg_min(op(variables)),tdeg(op(variables))])
U := SMOS[1]; V := SMOS[2];
if type(U, ShortMonomialOrder) = false then error
"invalid input: the first ShortMonomialOrder is not well defined"
else if localorglobal(U, variables) <> Lo then error
"invalid input: the first ShortMonomialOrder should be LOCAL"
else if type(V, ShortMonomialOrder) = false then error
"invalid input: the second ShortMonomialOrder is not well defined"
else if localorglobal(V, variables) <> Gl then error
"invalid input: the second ShortMonomialOrder should be GLOBAL"  }{}
\end{mapleinput}
\begin{mapleinput}
\mapleinline{active}{1d}{else r := nops(variables);
J:= [seq(diff(F, variables[i]), i = 1 .. r)];
J := Basis(J, V);
IniJ := [seq(LeadingMonomial(J[i], V), i = 1 .. nops(J))];
massimo:=r*max(seq(degree(IniJ[i],variables),i = 1..nops(IniJ)));
if IsZeroDimensional(<(op(IniJ))>,\{op(variables)\}) = false then error
"the given critical point is not isolated"  }{}
\end{mapleinput}
\begin{mapleinput}
\mapleinline{active}{1d}{else
LUU := [op(subs(vq= 1, convert(map(expand,
series(1/(product(1-variables[k]*vq, k =1 .. r)),
vq, massimo+1)), polynom)))];
L := sort(LUU, proc (t1,t2) options operator, arrow;
TestOrder(t1, t2, V) end proc);    }{}
\end{mapleinput}
\begin{mapleinput}
\mapleinline{active}{1d}{N :=[];
for j to nops(L) do up := 1;
for k to nops(IniJ) do
if divide(L[j], IniJ[k]) then up := 0
end if end do;
if up = 1 then N:= [op(N), L[j]]
end if end do;  }{}
\end{mapleinput}
\begin{mapleinput}
\mapleinline{active}{1d}{G := SB(J, variables, U);
Ini := [seq(LeadingMonomial(G[i], U), i = 1 .. nops(G))];
massimo:=r*max(seq(degree(Ini[i],variables),i = 1..nops(Ini)));
LUU := [op(subs(vq = 1, convert(map(expand,
series(1/(product(1-variables[h]*vq, h = 1 .. r)), vq,
massimo+1)), polynom)))];
L := sort(LUU, proc (t1, t2) options operator, arrow;
TestOrder(t1, t2, U) end proc); }{}
\end{mapleinput}
\begin{mapleinput}
\mapleinline{active}{1d}{M := [];
for j to nops(L) do up := 1;
for k to nops(Ini) do
if divide(L[j], Ini[k]) then up := 0
end if end do;
if up = 1 then M := [op(M), L[j]]
end if end do;
[J, IniJ, N, nops(N), G, Ini, M, nops(M)]
end if end if end if end if end if end proc: }{}
\end{mapleinput}
\end{maplegroup}

\halfline

\noindent The associated sub--procedures:

\halfline

\begin{maplegroup}
\begin{mapleinput}
\mapleinline{active}{1d}{PolyMILNORGroebnerBasis := proc (F,variables::set:=indets(F),
SMOS::anything:=[tdeg_min(op(variables)),tdeg(op(variables))])
MILNOR(F,variables,SMOS)[1] end proc: }{}
\end{mapleinput}
\end{maplegroup}

\halfline

\begin{maplegroup}
\begin{mapleinput}
\mapleinline{active}{1d}{PolyMILNORLT := proc (F, variables::set:= indets(F),
SMOS::anything:=[tdeg_min(op(variables)),tdeg(op(variables))])
MILNOR(F, variables, SMOS)[2] end proc: }{}
\end{mapleinput}
\end{maplegroup}

\halfline

\begin{maplegroup}
\begin{mapleinput}
\mapleinline{active}{1d}{PolyMILNORBasis := proc (F,variables::set := indets(F),
SMOS::anything:=[tdeg_min(op(variables)),tdeg(op(variables))])
MILNOR(F,variables,SMOS)[3] end proc: }{}
\end{mapleinput}
\end{maplegroup}

\halfline

\begin{maplegroup}
\begin{mapleinput}
\mapleinline{active}{1d}{PolyMILNORNumber:=proc (F, variables::set:= indets(F),
SMOS::anything:=[tdeg_min(op(variables)),tdeg(op(variables))])
MILNOR(F,variables,SMOS)[4] end proc: }{}
\end{mapleinput}
\end{maplegroup}

\halfline

\begin{maplegroup}
\begin{mapleinput}
\mapleinline{active}{1d}{MILNORGroebnerBasis := proc(F,variables::set := indets(F),
SMOS::anything:=[tdeg_min(op(variables)),tdeg(op(variables))])
MILNOR(F,variables,SMOS)[5] end proc: }{}
\end{mapleinput}
\end{maplegroup}

\halfline

\begin{maplegroup}
\begin{mapleinput}
\mapleinline{active}{1d}{MILNORLT := proc (F, variables::set := indets(F),
SMOS::anything:=[tdeg_min(op(variables)),tdeg(op(variables))])
MILNOR(F, variables, SMOS)[6] end proc: }{}
\end{mapleinput}
\end{maplegroup}

\halfline

\begin{maplegroup}
\begin{mapleinput}
\mapleinline{active}{1d}{MILNORBasis := proc (F,variables::set :=indets(F),
SMOS::anything:=[tdeg_min(op(variables)),tdeg(op(variables))])
MILNOR(F, variables,SMOS)[7] end proc: }{}
\end{mapleinput}
\end{maplegroup}

\halfline

\begin{maplegroup}
\begin{mapleinput}
\mapleinline{active}{1d}{MILNORNumber:=proc (F,variables::set :=indets(F),
SMOS::anything:=[tdeg_min(op(variables)),tdeg(op(variables))])
MILNOR(F,variables,SMOS)[8] end proc: }{}
\end{mapleinput}
\end{maplegroup}

\halfline

\noindent Then the routine computing Tyurina numbers:

\halfline

\begin{maplegroup}
\begin{mapleinput}
\mapleinline{active}{1d}{TYURINA := proc (F, variables::set := indets(F),
SMOS::anything := [tdeg_min(op(variables)),tdeg(op(variables))])
U := SMOS[1]; V := SMOS[2];
if type(U, ShortMonomialOrder) = false then error
"invalid input: the first ShortMonomialOrder is not well defined"
else if localorglobal(U, variables) <> Lo then error
"invalid input: the first ShortMonomialOrder should be LOCAL"
else if type(V,ShortMonomialOrder) = false then error
"invalid input: the second ShortMonomialOrder is not well defined"
else if localorglobal(V, variables) <> Gl then error
"invalid input: the second ShortMonomialOrder should be GLOBAL"}{}
\end{mapleinput}
\begin{mapleinput}
\mapleinline{active}{1d}{else r :=nops(variables);
J := [F,seq(diff(F, variables[i]), i = 1 .. r)];
J := Basis(J, V);
IniJ := [seq(LeadingMonomial(J[i], V), i = 1 .. nops(J))];
massimo:=r*max(seq(degree(IniJ[i],variables),i = 1..nops(IniJ)));
if IsZeroDimensional(<(op(IniJ))>,\{op(variables)\})= false then error
"the given singular point is not isolated" }{}
\end{mapleinput}
\begin{mapleinput}
\mapleinline{active}{1d}{else LUU := [op(subs(vq = 1, convert(map(expand,
series(1/(product(1-variables[k]*vq, k = 1 .. r)),
vq, massimo+1)), polynom)))];
L := sort(LUU, proc (t1, t2) options operator, arrow;
TestOrder(t1, t2, V) end proc);}{}
\end{mapleinput}
\begin{mapleinput}
\mapleinline{active}{1d}{N := [];
for j to nops(L) do up := 1;
for k to nops(IniJ) do
if divide(L[j], IniJ[k]) then up := 0
end if end do;
if up = 1 then N := [op(N),L[j]]
end if end do;}{}
\end{mapleinput}
\begin{mapleinput}
\mapleinline{active}{1d}{G := SB(J, variables, U);
Ini :=[seq(LeadingMonomial(G[i], U), i = 1 .. nops(G))];
massimo:=r*max(seq(degree(Ini[i],variables), i = 1..nops(Ini)));
LUU := [op(subs(vq = 1, convert(map(expand,
series(1/(product(1-variables[h]*vq, h = 1 .. r)),
vq, massimo+1)), polynom)))];
L := sort(LUU, proc (t1, t2) options operator, arrow;
TestOrder(t1, t2, U) end proc); }{}
\end{mapleinput}
\begin{mapleinput}
\mapleinline{active}{1d}{M := [];
for j to nops(L) do up := 1;
for k to nops(Ini) do
if divide(L[j], Ini[k]) then up :=
0 end if end do;
if up = 1 then M := [op(M), L[j]]
end if end do;
[J, IniJ, N, nops(N), G, Ini, M, nops(M)]
end if end if end if end if end if end proc:}{}
\end{mapleinput}
\end{maplegroup}

\halfline

\noindent The associated sub--procedures:

\halfline

\begin{maplegroup}
\begin{mapleinput}
\mapleinline{active}{1d}{PolyTYURINAGroebnerBasis := proc(F,variables::set:=indets(F),
SMOS::anything:=[tdeg_min(op(variables)),tdeg(op(variables))])
TYURINA(F,variables,SMOS)[1] end proc: }{}
\end{mapleinput}
\end{maplegroup}

\halfline

\begin{maplegroup}
\begin{mapleinput}
\mapleinline{active}{1d}{PolyTYURINALT := proc (F, variables::set:= indets(F),
SMOS::anything:=[tdeg_min(op(variables)),tdeg(op(variables))])
TYURINA(F, variables, SMOS)[2] end proc: }{}
\end{mapleinput}
\end{maplegroup}

\halfline

\begin{maplegroup}
\begin{mapleinput}
\mapleinline{active}{1d}{PolyTYURINABasis := proc (F, variables::set := indets(F),
SMOS::anything:=[tdeg_min(op(variables)),tdeg(op(variables))])
TYURINA(F,variables,SMOS)[3] end proc: }{}
\end{mapleinput}
\end{maplegroup}

\halfline

\begin{maplegroup}
\begin{mapleinput}
\mapleinline{active}{1d}{PolyTYURINANumber:=proc (F,variables::set := indets(F),
SMOS::anything:=[tdeg_min(op(variables)),tdeg(op(variables))])
TYURINA(F,variables,SMOS)[4] end proc: }{}
\end{mapleinput}
\end{maplegroup}

\halfline

\halfline

\begin{maplegroup}
\begin{mapleinput}
\mapleinline{active}{1d}{TYURINAGroebnerBasis := proc(F,variables::set := indets(F),
SMOS::anything:=[tdeg_min(op(variables)),tdeg(op(variables))])
TYURINA(F,variables,SMOS)[5] end proc: }{}
\end{mapleinput}
\end{maplegroup}

\halfline

\begin{maplegroup}
\begin{mapleinput}
\mapleinline{active}{1d}{TYURINALT := proc (F, variables::set :=indets(F),
SMOS::anything:=[tdeg_min(op(variables)),tdeg(op(variables))])
TYURINA(F, variables, SMOS)[6] end proc: }{}
\end{mapleinput}
\end{maplegroup}

\halfline

\begin{maplegroup}
\begin{mapleinput}
\mapleinline{active}{1d}{TYURINABasis := proc (F,variables::set :=indets(F),
SMOS::anything:=[tdeg_min(op(variables)),tdeg(op(variables))])
TYURINA(F, variables,SMOS)[7] end proc: }{}
\end{mapleinput}
\end{maplegroup}

\halfline

\begin{maplegroup}
\begin{mapleinput}
\mapleinline{active}{1d}{TYURINANumber:=proc (F,variables::set :=indets(F),
SMOS::anything:=[tdeg_min(op(variables)),tdeg(op(variables))])
TYURINAnor(F,variables,SMOS)[8] end proc:}{}
\end{mapleinput}
\end{maplegroup}

\halfline

\subsubsection{An example of application.} Let us consider the
same polynomial $F(x,y)=x^3+y^4+xy^2$ given in \ref{variabili}.
Procedures \texttt{MILNOR} and \texttt{TYURINA} give \emph{all}
the information we could get by applying all the introduced
routines, precisely:

\halfline

\begin{maplegroup}
\begin{mapleinput}
\mapleinline{active}{1d}{MILNOR(F) }{}
\end{mapleinput}
\mapleresult
\begin{maplelatex}
\mapleinline{inert}{2d}{}
{\[\displaystyle [[3\,{x}^{2}+{y}^{2},2\,{y}^{3}+xy],[{x}^{2},{y}^{3}],[1,y,x,{y}^{2},xy,x{y}^{2}],6,\]}
\mapleinline{inert}{2d}{}
{\[\displaystyle \left\{ -2\,{y}^{4}+3\,{x}^{3},2\,{y}^{3}+xy,3\,{x}^{2}+{y}^{2} \right\},
[{y}^{2},xy,{x}^{3}],[{x}^{2},x,y,1],4]\]}
\end{maplelatex}
\end{maplegroup}
\begin{maplegroup}
\begin{mapleinput}
\mapleinline{active}{1d}{TYURINA(F) }{}
\end{mapleinput}
\mapleresult
\begin{maplelatex}
\mapleinline{inert}{2d}{}
{\[\displaystyle [[xy,3\,{x}^{2}+{y}^{2},{y}^{3}],[xy,{x}^{2},{y}^{3}],[1,y,x,{y}^{2}],4, \]}
\mapleinline{inert}{2d}{}
{\[\displaystyle \left\{ xy,3\,{x}^{3},3\,{x}^{2}+{y}^{2},{y}^{3} \right\} ,
[{y}^{2},{x}^{3}{y}^{3},xy],[{x}^{2},x,y,1],4]\]}
\end{maplelatex}
\end{maplegroup}

\halfline

\noindent In particular it turns out that $F$ admits some further
critical point which is not a singular point of $F^{-1}(0)$. By

\halfline

\begin{maplegroup}
\begin{mapleinput}
\mapleinline{active}{1d}{solve(\{diff(F,x),diff(F,y)\},[x,y]); }{}
\end{mapleinput}
\mapleresult
\begin{maplelatex}
\mapleinline{inert}{2d}{}
{\[\displaystyle [[x=0,y=0],[x=0,y=0],[x=-1/2\, \left( {\it RootOf} \left( 3\,{{\it \_Z}}^{2}
+1,{\it label}={\it \_L1} \right)  \right) ^{2},\]}
\mapleinline{inert}{2d}{}
{\[\displaystyle y=1/2\,{\it RootOf} \left( 3\,{{\it \_Z}}^{2}+1,{\it label}={\it \_L1} \right)],[x=0,y=0]]\]}
\end{maplelatex}
\end{maplegroup}

\halfline

\noindent it follows that there are precisely two further critical
points of $F$ having Milnor number 1 and Tyurina number 0. Let us
now type:

\halfline

\begin{maplegroup}
\begin{mapleinput}
\mapleinline{active}{1d}{TyB := TyurinaBasis(F) }{}
\end{mapleinput}
\mapleresult
\begin{maplelatex}
\mapleinline{inert}{2d}{TyB := [x^2, x, y, 1]}{\[\displaystyle
{\it TyB}\, := \,[{x}^{2},x,y,1]\]}
\end{maplelatex}
\end{maplegroup}
\begin{maplegroup}
\begin{mapleinput}
\mapleinline{active}{1d}{T := nops(TyB): }{}
\end{mapleinput}
\end{maplegroup}
\begin{maplegroup}
\begin{mapleinput}
\mapleinline{active}{1d}{F[Lambda] := F+sum(lambda[i]*TyB[T-i],i= 0 .. T-1); }{}
\end{mapleinput}
\mapleresult
\begin{maplelatex}
\mapleinline{inert}{2d}{F[Lambda] :=
x^3+y^4+x*y^2+lambda[0]+lambda[1]*y+lambda[2]*x+lambda[3]*x^2}{\[\displaystyle
F_{{\Lambda}}\, :=
\,{x}^{3}+{y}^{4}+x{y}^{2}+\lambda_{{0}}+\lambda_{{1}}y+\lambda_{{2}}x+\lambda_{{3}}{x}^{2}\]}
\end{maplelatex}
\end{maplegroup}
\begin{maplegroup}
\begin{mapleinput}
\mapleinline{active}{1d}{TYURINA(F[Lambda], \{x, y\}); }{}
\end{mapleinput}
\mapleresult
\begin{maplelatex}
\mapleinline{inert}{2d}{[[1], [1], [], 0, [1], [1], [],
0]}{\[\displaystyle [[1],[1],[],0,[1],[1],[],0]\]}
\end{maplelatex}
\end{maplegroup}
\begin{maplegroup}
\begin{mapleinput}
\mapleinline{active}{1d}{Tyurina(F[Lambda], \{x, y\}); }{}
\end{mapleinput}
\mapleresult
\begin{Maple Warning}{
Warning,  computation interrupted}\end{Maple Warning}
\end{maplegroup}

\halfline

\noindent We had to interrupt the calculation of \texttt{Tyurina}
since it wasn't able to produce any output after considerable
time, on the contrary of \texttt{TYURINA} which quickly produced
the given (trivial) output.

\begin{remark} What observed in \ref{user_friendly} and \ref{optional}, about optional inputs for \texttt{Milnor} and \texttt{Tyurina}, applies analogously for \texttt{MILNOR} and \texttt{TYURINA}.
\end{remark}

\section{Application: adjacencies of Arnol'd simple singularities}\label{Arnold-ex}

Let us consider the classes of Arnol'd simple singularities
$A_n,D_n,E_6,E_7,E_8$. Recall that a class of singularities $B$ is
said to be \emph{adjacent} to a class of singularities $A$ (notation
$A\leftarrow B$) if any singularity in $B$ can be deformed to a
singularity in $A$ by an arbitrarily small deformation (see \cite{Arnol'd-ST} \S I.2.7, \cite{Arnol'd&c} \S
15.0 and \cite{Looijenga} \S 7.C). Adjacency turns out to be a \emph{partial order relation} on the set of singularities' equivalence classes.

\noindent In the following we will employ the optional input on
variables, as observed in \ref{variabili}, to show explicit
equations of algebraic \emph{stratifications} of Kuranishi
spaces verifying the following Arnol'd adjacency diagram
\begin{equation}\label{adiacenza}
    \xymatrix{A_1&\ar[l]A_2&\ar[l]A_3&\ar[l]A_4&\ar[l]A_5&\ar[l]A_6&\ar[l]A_7&\cdots\ar[l]A_{n-1}\cdots\\
                &&\ar[u]D_4&\ar[u]\ar[l]D_5&\ar[u]\ar[l]D_6&\ar[u]\ar[l]D_7&\ar[l]\ar[u]D_8&\cdots\ar[l]\ar[u]D_n\cdots\\
                &&&\ar[ruu]\ar[u]E_6&\ar[ruu]\ar[u]\ar[l]E_7&\ar[ruu]\ar[u]\ar[l]E_8&&}
\end{equation}
(see \cite{Arnol'd}, \cite{Arnol'd&c} and in particular \cite{Arnol'd-ST} \S I.2.7). Such a stratification gives a geometric interpretation, by means of inclusions of algebraic subsets, of the partial order relation induced by adjacency.

Let us first of all observe that the Kuranishi space $T^1$ of an
Arnol'd simple $m$--fold singularity of type $A_n,D_n,E_6,E_7,E_8$
do not depend on its dimension $m$, since partials of quadratic terms
are linear generators of $I_f$ eliminating the associated variable
from $M\setminus L(I_f)$. Therefore the study of diagram (\ref{adiacenza}) can be
reduced to the case of Arnol'd simple \emph{curve} singularities
whose local equations are given in (\ref{DV}): this is actually guaranteed by the following Morse Splitting Lemma \ref{Morse}.

In the sequel we will need the following notation and results, essentially due to V.I.~Arnol'd \cite{Arnol'd&c}. We refer the interested reader to books \cite{Arnol'd-ST}, \cite{Arnol'd&c}, \cite{Looijenga} and \cite{Greuel-etal} for details and proofs.

\begin{definition}[Co-rank of a critical point] The \emph{co-rank of a critical point $p$} of a holomorphic function $f$ defined over an open subset of $\C^n$ is the number
\[
    \crk_f(p):= n - \rk \left(\Hess_f(p)\right)
\]
where $\Hess_f(p)$ is the Hessian matrix of $f$ in $p$ (for shortness the function $f$ will be omitted when clear from the context). In particular if $\crk(p)=0$ then $p$ is called a \emph{non--degenerate}, or \emph{Morse}, \emph{critical point}.
\end{definition}

\begin{lemma}[Morse Splitting Lemma, \cite{Arnol'd&c} \S11.1, \cite{Looijenga} (7.16) and \cite{Greuel-etal} Theorem I.2.47]\label{Morse} Assume that $f\in\mathfrak{m}^2\subset\Cxx$ and  $\crk_f(0)=n-k$.  Then $f$ is equivalent \footnote{In the sense of Definition \ref{whs} and Remark \ref{contact=right}.} to the following germ of singularity
\[
    \sum_{i=1}^k x_i^2+g(x_{k+1},\ldots,x_n)
\]
where $g\in\mathfrak{m}^3$ is uniquely determined (up to equivalence).
\end{lemma}

\begin{theorem}[\cite{Greuel-etal}, Theorems I.2.46, I.2.48, I.2.51, I.2.53]\label{classificazione}\hfill
\begin{enumerate}
  \item For $f\in\mathfrak{m}^2\subset\Cxx$ the following facts are equivalent:
  \begin{itemize}
    \item $\crk(0)=0$ i.e. $0$ is a non--degenerate critical point of $f$,
    \item $\mu(0)=1=\tau(1)$,
    \item $0\in f^{-1}(0)$ is, up to equivalence, a node i.e. a simple $A_1$ singularity.
  \end{itemize}
  \item For $f\in\mathfrak{m}^2\subset\Cxx$ the following facts are equivalent:
  \begin{itemize}
    \item $\crk(0)\leq 1$ and $\mu(0)=m$,
    \item $0\in f^{-1}(0)$ is equivalent to a simple $A_m$ singularity.
  \end{itemize}
  \item For $f\in\mathfrak{m}^3\subset\C\{y,z\}$ the following facts are equivalent:
  \begin{itemize}
    \item the 3--jet $f^{(3)}$ of $f$ factors into at least two distinct linear factors and $\mu(0)=m\geq 4$,
    \item $0\in f^{-1}(0)$ is equivalent to a simple $D_m$ singularity.
  \end{itemize}
  \item For $f\in\mathfrak{m}^3\subset\C\{y,z\}$ the following facts are equivalent:
  \begin{itemize}
    \item the 3--jet $f^{(3)}$ of $f$ has a unique liner factor (of multiplicity 3) and $\mu(f)\leq 8$,
    \item $0\in f^{-1}(0)$ is equivalent to a simple singularity of type $E_6$, $E_7$ or $E_8$ and $\mu(0)=6$, $7$ or $8$ respectively.
  \end{itemize}
\end{enumerate}
\end{theorem}

Let us now introduce a non-standard notation, useful to describe a nice geometric property of stratifications via algebraic subsets of a simple Arnol'd singularity's Kuranishi space, as explained in the following statements. Consider the following square of subset inclusions
\begin{equation}\label{square}
    \xymatrix{B\ar@{^{(}->}[r]&D\\
                A\ar@{^{(}->}[u]\ar@{^{(}->}[r]&C\ar@{^{(}->}[u]}
\end{equation}
Then $A\subseteq B\cap C$ , necessarily.

\begin{figure}[h]
  \includegraphics[width=5cm, angle=-90]{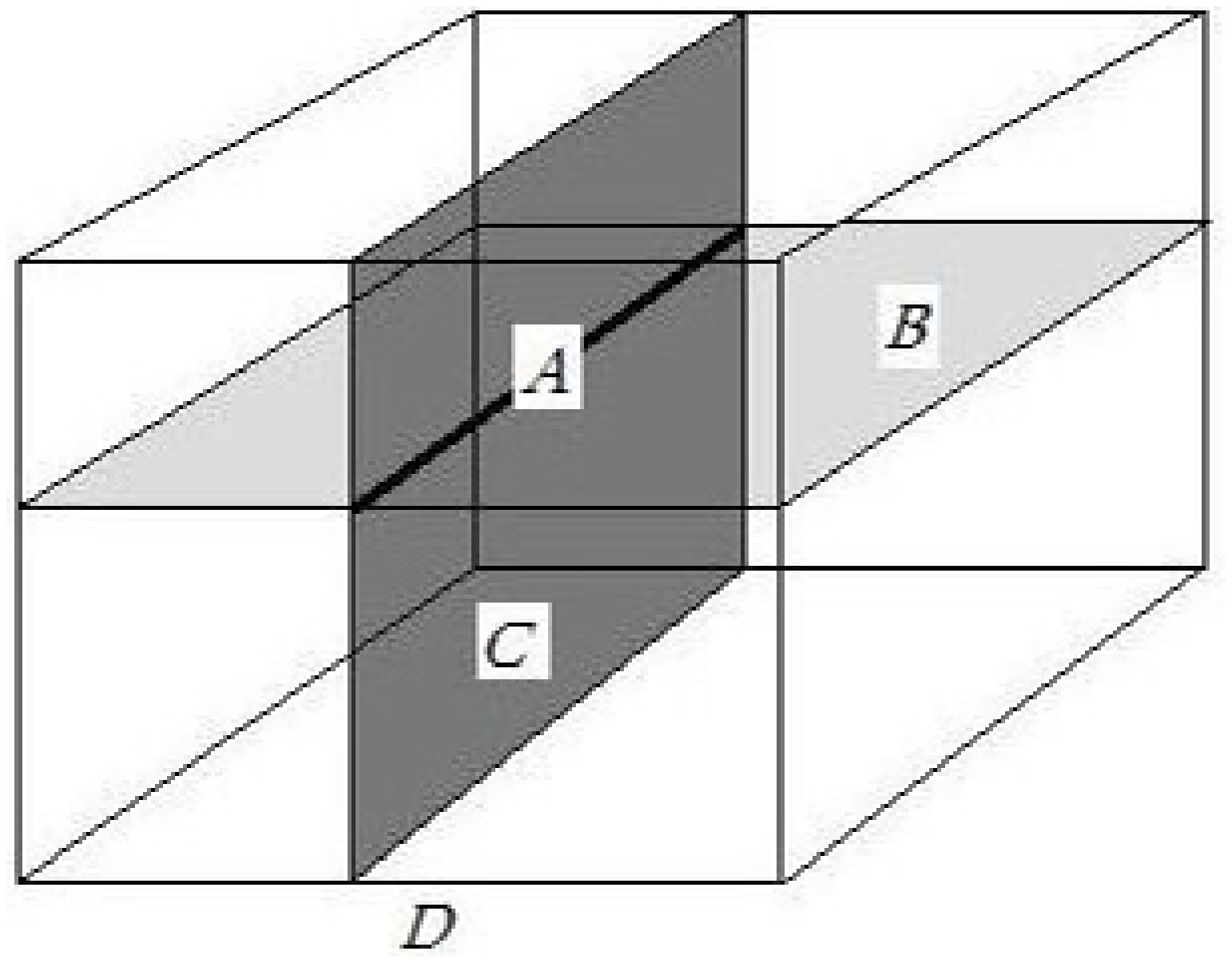}\\
  \caption{The c.i.p. for the inclusions' square (\ref{square})}\label{cip}
\end{figure}

\begin{figure}
  \includegraphics[width=5cm, angle=-90]{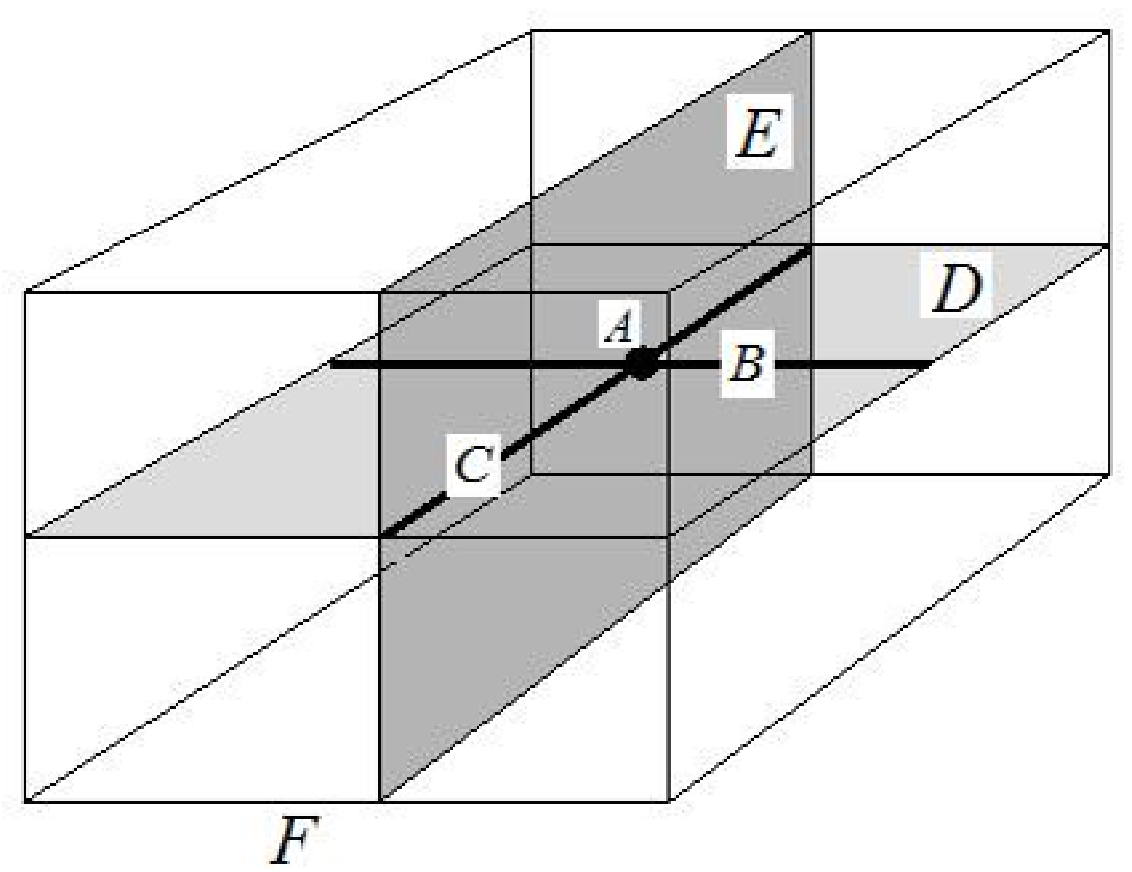}\\
  \caption{The sequence (\ref{2squares}) of two c.i.p. squares}\label{2cip}
\end{figure}

\begin{definition}[Complete Intersection Property - (c.i.p.)] A square of subset inclusions (\ref{square}) is said to admit the \emph{complete intersection property} if $$A=B\cap C\ .$$
For shortness we will say that (\ref{square}) \emph{is a c.i.p. square}. The geometric meaning of c.i.p. in (\ref{square}) is explained by Figure \ref{cip}, while Figure \ref{2cip} describes geometrically the following \emph{sequence of two c.i.p. squares}
\begin{equation}\label{2squares}
    \xymatrix{B\ar@{^{(}->}[r]&D\ar@{^{(}->}[r]&F\\
                A\ar@{^{(}->}[r]\ar@{^{(}->}[u]&C\ar@{^{(}->}[r]\ar@{^{(}->}[u]&E\ar@{^{(}->}[u]}
\end{equation}
meaning that $A=B\cap C= B\cap D\cap E$.
Moreover the following inclusions' diagram
\begin{equation}\label{square-union}
    \xymatrix{\ B\cup F\ \ar@{^{(}->}[dr]\ &E\ \ar@{_{(}->}[l]\ar@{^{(}->}[dr]&\\
        A\ar@{^{(}->}[u]\ar@{^{(}->}[dr]&\ D\ &\ G\ \ar@{_{(}->}[l]\\
        &C \ar@{^{(}->}[u]&}
\end{equation}
is called \emph{a union of c.i.p. squares} if $A=B\cap C$ and $E=F\cap G$. A particular case, occurring in the following, is when $C=G$\ : then diagram (\ref{square-union}) becomes the following one
\begin{equation}\label{square-union2}
    \xymatrix{&B\cup F\ar@{^{(}->}[r]&D\\
              &E\ar@{^{(}->}[r]\ar@{^{(}->}[u]&C=G\ar@{^{(}->}[u]\\
              A\ar@/^/@{^{(}->}[ruu] \ar@/_/@{^{(}->}[rru]&& }
\end{equation}
and we will say this diagram to represent \emph{a hinged union of c.i.p. squares}, whose hinge is the inclusion $C\hookrightarrow D$.

\noindent At last the following inclusions' diagram
\begin{equation}\label{reducible-square}
    \xymatrix{&B\ar@{^{(}->}[r]&D\\
              &E\ar@{^{(}->}[r]\ar@{^{(}->}[u]&C\ar@{^{(}->}[u]\\
              A\ar@/^/@{^{(}->}[ruu] \ar@/_/@{^{(}->}[rru]&& }
\end{equation}
is  called \emph{a reducible c.i.p. square} if $A\cup E = B\cap C$ i.e. if $$\xymatrix{B\ar@{^{(}->}[r]&D\\A\cup E\ar@{^{(}->}[r]\ar@{^{(}->}[u]&C\ar@{^{(}->}[u]}$$
is a c.i.p. square.

\end{definition}

\subsection{Outline of the following results}\label{outline}

Statements and proofs of the following Theorems \ref{An} , \ref{Dn} , \ref{E6} , \ref{E7} and \ref{E8} have the same structure we are going to outline here. Precisely their statements describe set theoretical stratifications by algebraic subsets of the Kuranishi space $T^1$ of simple hypersurface singularities $A_n, D_n$ and $E_n$ (the latter with $6\leq n\leq 8$). Their proofs go on by the following steps:
\begin{enumerate}
  \item look for the critical points of a generic small deformation of our initial simple singularity, by solving the polynomial system assigned by partial derivatives (Jacobian ideal generators): they turn out to be precisely $n$ points (with $6\leq n\leq 8$ in the $E_n$ case);
  \item imposing one of the previous critical points to be actually a singular point means defining a hypersurface $\mathcal{L}\subset T^1$ : we have then $n$ of such hypersurfaces, one for each critical point;
  \item any of those hypersurfaces is then stratified by nested algebraic subsets defined by a progressive vanishing of leading coefficients in Jacobian ideals' standard bases of more and more specialized deformations. More precisely, the general strategy is that of looking at the leading monomials ordered by the choice of a suitable l.m.o.: then imposing the vanishing of only the leading coefficient associated with the smallest leading monomial realizes ``horizontal adjacencies", while imposing the vanishing of all the leading coefficients gives ``vertical adjacencies", in diagrams (\ref{An-adiacenze}), (\ref{Dn-adiacenza}), (\ref{E6-adiacenza}), (\ref{E7-adiacenza}) and (\ref{E8-adiacenza}).
\end{enumerate}
The last step (3) is obtained by a systematic use of routines \texttt{Milnor}, \texttt{Tyurina}, \texttt{MILNOR} and \texttt{TYURINA} previously described, running with suitable l.m.o.'s defined on a strict subset of variables appearing in the polynomial equation of a generic small deformation, as explained in \ref{variabili}. Such a procedure allows to explicitly write down relations on deformation parameters and then equations of the algebraic stratifications.

\begin{caveat} The large use of MAPLE routines in step (3) makes difficult to follow and understand algebraic arguments without the support of MAPLE 11 or 12. Hence the reader is warmly advised of repeating by himself MAPLE commands as described in the following. Moreover since the following results are here reported as a significant application of the previously presented subroutines, their role, inputs and outputs will be explained and listed in detail; on the contrary the use of standard MAPLE procedures, like e.g. \texttt{eliminate} or \texttt{EliminationIdeal}, will be left to the reader who may helpfully consult the excellent on-line MAPLE help system.
\end{caveat}

\subsection{Simple singularities of $A_n$ type.}

\begin{theorem}\label{An} Let $T^1$ be the Kuranishi space of a simple
$N$--dimensional singular point $0\in f^{-1}(0)$ with
\[
    f(x_1,\ldots,x_{N+1}) = \sum_{i=1}^N x_i^2\ +\
    x_{N+1}^{n+1}\quad\text{(for $n\geq 1$)}\ .
\]
The subset of $T^1$ parameterizing small deformations of $0\in f^{-1}(0)$ to a simple node (i.e. an $A_1$ singularity) is the union of $n$ hypersurfaces. Moreover, calling $\mathcal{L}$ any of those hypersurfaces, there exists a stratification of nested algebraic subsets
\begin{equation}\label{An-inclusioni}
    \xymatrix@1{\mathcal{L}&{\ \mathcal{V}_2}\ar@{_{(}->}[l]&{\ \cdots}\ar@{_{(}->}[l]
    &{\ \mathcal{V}_2^m}\ar@{_{(}->}[l]&{\ \cdots}\ar@{_{(}->}[l]
    &{\ \mathcal{V}_2^n}\ar@{_{(}->}[l]}
\end{equation}
verifying the Arnol'd's adjacency diagram
\begin{equation}\label{An-adiacenze}
    \xymatrix@1{A_1& A_2\ar[l]& \cdots\ar[l] & A_m\ar[l] &\cdots\ar[l] &
    A_n\ar[l]}
\end{equation}
where
\begin{itemize}
    \item $\mathcal{L}$ is the hypersurface of $T^1$ defined by equation
(\ref{lamda-condizione:A_n}), keeping in mind
(\ref{A_n,z}),
    \item $\mathcal{V}_2^m:=\bigcap_{k=2}^m \mathcal{V}_k$ where
$\mathcal{V}_k$ are hypersurfaces of $\mathcal{L}$ defined by the
vanishing of variables $v_k$ introduced by
(\ref{lamda-condizione:A_n-2}).
\end{itemize}
\end{theorem}

\begin{proof} Let us follow the outline previously exposed in \ref{outline}.

\vskip 5pt\noindent (1) By the Morse Splitting Lemma \ref{Morse} we can reduce to the case $N=1$ with $f(y,z)= y^2 + z^{n+1}$
for $n\geq 1$. Then Proposition \ref{T^1} gives
\[T^1 \cong \C[y,z]/(y, z^n) \cong
\langle 1,z,\ldots, z^{n-1}\rangle_{\C}
\]
and,  given $\Lambda =
(\lambda_0,\ldots,\lambda_{n-1})\in T^1$, the associated
deformation of $U_0=\Spec(\mathcal{O}_{f,0})$ is
\[
    U_{\Lambda} = \{f_{\Lambda}(y,z) := f(y,z) + \sum_{i=0}^{n-1} \lambda_i
    z^i=0
    \}\ .
\]
A solution of the jacobian system of partial derivatives is then
given by $p_{\Lambda}=(0,z_{\Lambda})$ where $z_{\Lambda}$ is a zero of the following polynomial
\begin{equation}\label{A_n,z}
    (n+1)z^n + \sum_{i=1}^{n-1}
    i\lambda_i z^{i-1} \in  \C[\mathbf{\lambda}][z]\ .
\end{equation}
This means that $f_{\Lambda}$ admits precisely $n$ critical points.

\vskip 5pt
\noindent (2) Imposing $p_{\Lambda}\in U_{\Lambda}$, which is asking for one of the previous critical points to be actually a singular point of $U_{\Lambda}$, defines the
following hypersurface in $T^1$
\begin{equation}\label{lamda-condizione:A_n}
    p_{\Lambda}\in U_{\Lambda}\ \Longleftrightarrow\ \Lambda\in
    \mathcal{L}:=\{ z_{\Lambda}^{n+1} + \sum_{i=0}^{n-1} \lambda_i z_{\Lambda}^i =
    0\}\subset T^1\ .
\end{equation}
Notice that we get precisely $n$ such hypersurfaces, one for each critical point of $f_{\Lambda}$.
\begin{itemize}
    \item \emph{$\mathcal{L}$ is a hypersurface of the Kuranishi space
parameterizing small deformations of $U_0$ admitting a singular
point of type at least $A_1$ in the origin}.
\end{itemize}

\noindent After translating $z\mapsto z+z_{\Lambda}$, we get
\begin{eqnarray}\label{fLambda-An}
\nonumber
  f_{\Lambda}(y,z+z_{\Lambda}) = f(y,z) &+& \left(z_{\Lambda}^{n+1}
     + \sum_{i=0}^{n-1} \lambda_i z_{\Lambda}^i\right) \\
      &+& \left((n+1)z_{\Lambda}^n + \sum_{i=1}^{n-1}
    i\lambda_i z_{\Lambda}^{i-1}\right)z \\
\nonumber
    &+& \sum_{k=2}^{n-1} \left( {n+1\choose k} z_{\Lambda}^{n+1-k} +
    \sum_{i=k}^{n-1}{i\choose k} \lambda_i
    z_{\Lambda}^{i-k}\right)z^k \\
\nonumber
    &+& (n+1)z_{\Lambda}\ z^n\quad .
\end{eqnarray}
Then by (\ref{lamda-condizione:A_n}) and (\ref{A_n,z}) the origin
(i.e. $p_{\Lambda}\in U_{\Lambda}$) is at least a node: in fact by
running \texttt{Milnor} and \texttt{Tyurina} e.g. when $n=7$, we
get (set $L:=z_{\Lambda}$):

\halfline

\begin{maplegroup}
\begin{mapleinput}
\mapleinline{active}{1d}{F:=y\symbol{94}2 + z\symbol{94}8; }{}
\end{mapleinput}
\mapleresult
\begin{maplelatex}
\mapleinline{inert}{2d}{F := y^2+z^8}{\[\displaystyle F\, :=
\,{y}^{2}+{z}^{8}\]}
\end{maplelatex}
\end{maplegroup}
\begin{maplegroup}
\begin{mapleinput}
\mapleinline{active}{1d}{TyB := TyurinaBasis(F) }{}
\end{mapleinput}
\mapleresult
\begin{maplelatex}
\mapleinline{inert}{2d}{TyB := [z^6, z^5, z^4, z^3, z^2, z,
1]}{\[\displaystyle {\it TyB}\, :=
\,[{z}^{6},{z}^{5},{z}^{4},{z}^{3},{z}^{2},z,1]\]}
\end{maplelatex}
\end{maplegroup}
\begin{maplegroup}
\begin{mapleinput}
\mapleinline{active}{1d}{T := nops(TyB):
F[Lambda] := F+sum(lambda[i]*TyB[T-i], i = 0 .. T-1) }{}
\end{mapleinput}
\mapleresult
\begin{maplelatex}
\mapleinline{inert}{2d}{F[Lambda] := y^2+z^8+lambda[0]+lambda[1]*z+lambda[2]*z^2+lambda[3]*z^3+lambda[4]*z^4+lambda[5]*z^5+lambda[6]*z^6}{\[\displaystyle F_{{\Lambda}}\, := \,{y}^{2}+{z}^{8}+\lambda_{{0}}+\lambda_{{1}}z+\lambda_{{2}}{z}^{2}+\lambda_{{3}}{z}^{3}+\lambda_{{4}}{z}^{4}+\lambda_{{5}}{z}^{5}\\
\mbox{}+\lambda_{{6}}{z}^{6}\]}
\end{maplelatex}
\end{maplegroup}
\begin{maplegroup}
\begin{mapleinput}
\mapleinline{active}{1d}{solve(\{diff(F[Lambda], y), diff(F[Lambda], z)\}, [y, z]) }{}
\end{mapleinput}
\mapleresult
\begin{maplelatex}
\mapleinline{inert}{2d}{[[y = 0, z = RootOf(8*_Z^7+lambda[1]+2*lambda[2]*_Z+3*lambda[3]*_Z^2+4*lambda[4]*_Z^3+5*lambda[5]*_Z^4+6*lambda[6]*_Z^5)]]}{\[\displaystyle [[y=0,z={\it RootOf} \left( 8\,{{\it \_Z}}^{7}+\lambda_{{1}}+2\,\lambda_{{2}}{\it \_Z}\\
\mbox{}+3\,\lambda_{{3}}{{\it \_Z}}^{2}+4\,\lambda_{{4}}{{\it
\_Z}}^{3}+5\,\lambda_{{5}}{{\it \_Z}}^{4}+6\,\lambda_{{6}}{{\it
\_Z}}^{5} \right) ]]\]}
\end{maplelatex}
\end{maplegroup}

\halfline

\noindent Notice that the polynomial in \emph{RootOf} is precisely
(\ref{A_n,z}) for $n=7$.

\halfline

\begin{maplegroup}
\begin{mapleinput}
\mapleinline{active}{1d}{z := Z+L: }{}
\end{mapleinput}
\end{maplegroup}
\begin{maplegroup}
\begin{mapleinput}
\mapleinline{active}{1d}{collect(F[Lambda], [y, Z], 'distributed')
}{}
\end{mapleinput}
\mapleresult
\begin{maplelatex}
\mapleinline{inert}{2d}{}
{\[\displaystyle {y}^{2}+{Z}^{8}+8\,L{Z}^{7}+ \left( \lambda_{{6}}+28\,{L}^{2} \right) {Z}^{6}+
\left( 56\,{L}^{3}+\lambda_{{5}}+6\,\lambda_{{6}}L \right) {Z}^{5}+\]}
\mapleinline{inert}{2d}{}
{\[\displaystyle \left( 70\,{L}^{4}+5\,\lambda_{{5}}L+15\,\lambda_{{6}}{L}^{2}+\lambda_{{4}} \right) {Z}^{4}
+ \left( 10\,\lambda_{{5}}{L}^{2}+\lambda_{{3}}+56\,{L}^{5}+4\,\lambda_{{4}}L+20\,\lambda_{{6}}{L}^{3} \right) {Z}^{3}\]}
\mapleinline{inert}{2d}{}
{\[\displaystyle + \left( 10\,\lambda_{{5}}{L}^{3}+3\,\lambda_{{3}}L+6\,\lambda_{{4}}{L}^{2}+15\,
\lambda_{{6}}{L}^{4}+\lambda_{{2}}+28\,{L}^{6} \right) {Z}^{2}+\]}
\mapleinline{inert}{2d}{}
{\[\displaystyle \left( 4\,\lambda_{{4}}{L}^{3}+8\,{L}^{7}+5\,\lambda_{{5}}{L}^{4}+6\,\lambda_{{6}}{L}^{5}+2\,
\lambda_{{2}}L+3\,\lambda_{{3}}{L}^{2}+\lambda_{{1}} \right)Z\]}
\mapleinline{inert}{2d}{}
{\[\displaystyle +\lambda_{{2}}{L}^{2}+\lambda_{{5}}{L}^{5}+\lambda_{{6}}{L}^{6}+{L}^{8}+\lambda_{{0}}+
\lambda_{{3}}{L}^{3}+\lambda_{{4}}{L}^{4}+\lambda_{{1}}L\]}
\end{maplelatex}
\end{maplegroup}

\halfline

\noindent Compare (\ref{fLambda-An}) for n=7 with the previous
output expression.

\halfline

\begin{maplegroup}
\begin{mapleinput}
\mapleinline{active}{1d}{SD := [seq(Z\symbol{94}i, i = 2 .. 7)]
}{}
\end{mapleinput}
\mapleresult
\begin{maplelatex}
\mapleinline{inert}{2d}{SD := [Z^2, Z^3, Z^4, Z^5, Z^6,
Z^7]}{\[\displaystyle {\it SD}\, :=
\,[{Z}^{2},{Z}^{3},{Z}^{4},{Z}^{5},{Z}^{6},{Z}^{7}]\]}
\end{maplelatex}
\end{maplegroup}
\begin{maplegroup}
\begin{mapleinput}
\mapleinline{active}{1d}{for j from 2 to 7 do
FL[j]:=y\symbol{94}2+Z\symbol{94}8+sum(v[i]*SD[i-1], i = j .. 7) end do }{}
\end{mapleinput}
\mapleresult
\begin{maplelatex}
\mapleinline{inert}{2d}{FL[2]:=y^2+Z^8+v[2]*Z^2+v[3]*Z^3+v[4]*Z^4+v[5]*Z^5+v[6]*Z^6+v[7]*Z^7}{\[\displaystyle
{\it FL}_{{2}}\, :=
\,{y}^{2}+{Z}^{8}+v_{{2}}{Z}^{2}+v_{{3}}{Z}^{3}+v_{{4}}{Z}^{4}+v_{{5}}{Z}^{5}+v_{{6}}{Z}^{6}+v_{{7}}{Z}^{7}\]}
\end{maplelatex}
\mapleresult
\begin{maplelatex}
\mapleinline{inert}{2d}{FL[3] :=
y^2+Z^8+v[3]*Z^3+v[4]*Z^4+v[5]*Z^5+v[6]*Z^6+v[7]*Z^7}{\[\displaystyle
{\it FL}_{{3}}\, :=
\,{y}^{2}+{Z}^{8}+v_{{3}}{Z}^{3}+v_{{4}}{Z}^{4}+v_{{5}}{Z}^{5}+v_{{6}}{Z}^{6}+v_{{7}}{Z}^{7}\]}
\end{maplelatex}
\mapleresult
\begin{maplelatex}
\mapleinline{inert}{2d}{FL[4] :=
y^2+Z^8+v[4]*Z^4+v[5]*Z^5+v[6]*Z^6+v[7]*Z^7}{\[\displaystyle {\it
FL}_{{4}}\, :=
\,{y}^{2}+{Z}^{8}+v_{{4}}{Z}^{4}+v_{{5}}{Z}^{5}+v_{{6}}{Z}^{6}+v_{{7}}{Z}^{7}\]}
\end{maplelatex}
\mapleresult
\begin{maplelatex}
\mapleinline{inert}{2d}{FL[5] :=
y^2+Z^8+v[5]*Z^5+v[6]*Z^6+v[7]*Z^7}{\[\displaystyle {\it
FL}_{{5}}\, :=
\,{y}^{2}+{Z}^{8}+v_{{5}}{Z}^{5}+v_{{6}}{Z}^{6}+v_{{7}}{Z}^{7}\]}
\end{maplelatex}
\mapleresult
\begin{maplelatex}
\mapleinline{inert}{2d}{FL[6] :=
y^2+Z^8+v[6]*Z^6+v[7]*Z^7}{\[\displaystyle {\it FL}_{{6}}\, :=
\,{y}^{2}+{Z}^{8}+v_{{6}}{Z}^{6}+v_{{7}}{Z}^{7}\]}
\end{maplelatex}
\mapleresult
\begin{maplelatex}
\mapleinline{inert}{2d}{FL[7] := y^2+Z^8+v[7]*Z^7}{\[\displaystyle
{\it FL}_{{7}}\, := \,{y}^{2}+{Z}^{8}+v_{{7}}{Z}^{7}\]}
\end{maplelatex}
\end{maplegroup}
\begin{maplegroup}
\begin{mapleinput}
\mapleinline{active}{1d}{for j from 2 to 7 do Milnor(FL[9-j],\{y,Z\}) end do }{}
\end{mapleinput}
\mapleresult
\begin{maplelatex}
\mapleinline{inert}{2d}{[[2*y, 8*Z^7+7*v[7]*Z^6], [y, Z^6], [Z^5,
Z^4, Z^3, Z^2, Z, 1], 6]}{\[\displaystyle
[[2\,y,8\,{Z}^{7}+7\,v_{{7}}{Z}^{6}],[y,{Z}^{6}],[{Z}^{5},{Z}^{4},{Z}^{3},{Z}^{2},Z,1],6]\]}
\end{maplelatex}
\mapleresult
\begin{maplelatex}
\mapleinline{inert}{2d}{[[2*y, 8*Z^7+6*v[6]*Z^5+7*v[7]*Z^6], [y,
Z^5], [Z^4, Z^3, Z^2, Z, 1], 5]}{\[\displaystyle
[[2\,y,8\,{Z}^{7}+6\,v_{{6}}{Z}^{5}+7\,v_{{7}}{Z}^{6}],[y,{Z}^{5}],[{Z}^{4},{Z}^{3},{Z}^{2},Z,1],5]\]}
\end{maplelatex}
\mapleresult
\begin{maplelatex}
\mapleinline{inert}{2d}{[[2*y, 8*Z^7+5*v[5]*Z^4+6*v[6]*Z^5+7*v[7]*Z^6], [y, Z^4], [Z^3, Z^2, Z, 1], 4]}{\[\displaystyle [[2\,y,8\,{Z}^{7}+5\,v_{{5}}{Z}^{4}+6\,v_{{6}}{Z}^{5}+7\,v_{{7}}{Z}^{6}],[y,{Z}^{4}],[{Z}^{3},{Z}^{2\\
\mbox{}},Z,1],4]\]}
\end{maplelatex}
\mapleresult
\begin{maplelatex}
\mapleinline{inert}{2d}{[[2*y, 8*Z^7+4*v[4]*Z^3+5*v[5]*Z^4+6*v[6]*Z^5+7*v[7]*Z^6], [y, Z^3], [Z^2, Z, 1], 3]}{\[\displaystyle [[2\,y,8\,{Z}^{7}+4\,v_{{4}}{Z}^{3}+5\,v_{{5}}{Z}^{4}+6\,v_{{6}}{Z}^{5}+7\,v_{{7}}{Z}^{6}],[y,{Z}^{3\\
\mbox{}}],[{Z}^{2},Z,1],3]\]}
\end{maplelatex}
\mapleresult
\begin{maplelatex}
\mapleinline{inert}{2d}{[[2*y, 8*Z^7+3*v[3]*Z^2+4*v[4]*Z^3+5*v[5]*Z^4+6*v[6]*Z^5+7*v[7]*Z^6], [y, Z^2], [Z, 1], 2]}{\[\displaystyle [[2\,y,8\,{Z}^{7}+3\,v_{{3}}{Z}^{2}+4\,v_{{4}}{Z}^{3}+5\,v_{{5}}{Z}^{4}+6\,v_{{6}}{Z}^{5}+7\,v_{{7}}{Z}^{6}\\
\mbox{}],[y,{Z}^{2}],[Z,1],2]\]}
\end{maplelatex}
\mapleresult
\begin{maplelatex}
\mapleinline{inert}{2d}{[[2*y, 8*Z^7+2*v[2]*Z+3*v[3]*Z^2+4*v[4]*Z^3+5*v[5]*Z^4+6*v[6]*Z^5+7*v[7]*Z^6], [y, Z], [1], 1]}{\[\displaystyle [[2\,y,8\,{Z}^{7}+2\,v_{{2}}Z+3\,v_{{3}}{Z}^{2}+4\,v_{{4}}{Z}^{3}+5\,v_{{5}}{Z}^{4}+6\,v_{{6}}{Z}^{5}\\
\mbox{}+7\,v_{{7}}{Z}^{6}],[y,Z],[1],1]\]}
\end{maplelatex}
\end{maplegroup}
\begin{maplegroup}
\begin{mapleinput}
\mapleinline{active}{1d}{for j from 2 to 7 do Tyurina(FL[9-j],\{y, Z\}) end do }{}
\end{mapleinput}
\mapleresult
\begin{maplelatex}
\mapleinline{inert}{2d}{[[y^2+Z^8+v[7]*Z^7, 2*y, 8*Z^7+7*v[7]*Z^6], [y^2, y, Z^6], [Z^5, Z^4, Z^3, Z^2, Z, 1], 6]}{\[\displaystyle [[{y}^{2}+{Z}^{8}+v_{{7}}{Z}^{7},2\,y,8\,{Z}^{7}+7\,v_{{7}}{Z}^{6}],[{y}^{2},y,{Z}^{6}],[{Z}^{5},{Z}^{4},{Z}^{3},{Z}^{2\\
\mbox{}},Z,1],6]\]}
\end{maplelatex}
\mapleresult
\begin{maplelatex}
\mapleinline{inert}{2d}{[[y^2+Z^8+v[6]*Z^6+v[7]*Z^7, 2*y, 8*Z^7+6*v[6]*Z^5+7*v[7]*Z^6], [y^2, y, Z^5], [Z^4, Z^3, Z^2, Z, 1], 5]}{\[\displaystyle [[{y}^{2}+{Z}^{8}+v_{{6}}{Z}^{6}+v_{{7}}{Z}^{7},2\,y,8\,{Z}^{7}+6\,v_{{6}}{Z}^{5}+7\,v_{{7}}{Z}^{6}\\
\mbox{}],[{y}^{2},y,{Z}^{5}],[{Z}^{4},{Z}^{3},{Z}^{2},Z,1],5]\]}
\end{maplelatex}
\mapleresult
\begin{maplelatex}
\mapleinline{inert}{2d}{}
{\[\displaystyle [[{y}^{2}+{Z}^{8}+v_{{5}}{Z}^{5}+v_{{6}}{Z}^{6}+v_{{7}}{Z}^{7},2\,y,8\,{Z}^{7}+5\,v_{{5}}{Z}^{4}+6\,v_{{6}}{Z}^{5}
+7\,v_{{7}}{Z}^{6}],[{y}^{2},y,{Z}^{4}],\]}
\mapleinline{inert}{2d}{}
{\[\displaystyle [{Z}^{3},{Z}^{2},Z,1],4]\]}
\end{maplelatex}
\mapleresult
\begin{maplelatex}
\mapleinline{inert}{2d}{}
{\[\displaystyle [[{y}^{2}+{Z}^{8}+v_{{4}}{Z}^{4}+v_{{5}}{Z}^{5}+v_{{6}}{Z}^{6}+v_{{7}}{Z}^{7},2\,y,8\,{Z}^{7}+
4\,v_{{4}}{Z}^{3}\\
\mbox{}+5\,v_{{5}}{Z}^{4}+6\,v_{{6}}{Z}^{5}+7\,v_{{7}}{Z}^{6}],\]}
\mapleinline{inert}{2d}{}
{\[\displaystyle [{y}^{2},y,{Z}^{3}],[{Z}^{2},Z,1],3]\]}
\end{maplelatex}
\mapleresult
\begin{maplelatex}
\mapleinline{inert}{2d}{}
{\[\displaystyle [[{y}^{2}+{Z}^{8}+v_{{3}}{Z}^{3}+v_{{4}}{Z}^{4}+v_{{5}}{Z}^{5}+v_{{6}}{Z}^{6}+v_{{7}}{Z}^{7},
2\,y,8\,{Z}^{7}+3\,v_{{3}}{Z}^{2}+4\,v_{{4}}{Z}^{3}+5\,v_{{5}}{Z}^{4}+\]}
\mapleinline{inert}{2d}{}
{\[\displaystyle 6\,v_{{6}}{Z}^{5}+7\,v_{{7}}{Z}^{6}],[{y}^{2},y,{Z}^{2}],[Z,1],2]\]}
\end{maplelatex}
\mapleresult
\begin{Maple Warning}{
Warning,  computation interrupted}\end{Maple Warning}

\end{maplegroup}
\begin{maplegroup}
\begin{mapleinput}
\mapleinline{active}{1d}{for j from 2 to 7 do TYURINA(FL[9-j],\{y, Z\}) end do }{}
\end{mapleinput}
\mapleresult
\begin{maplelatex}
\mapleinline{inert}{2d}{[[y, Z^6], [y, Z^6], [1, Z, Z^2, Z^3, Z^4, Z^5], 6, [y, Z^6], [y, Z^6], [Z^5, Z^4, Z^3, Z^2, Z, 1], 6]}{\[\displaystyle [[y,{Z}^{6}],[y,{Z}^{6}],[1,Z,{Z}^{2},{Z}^{3},{Z}^{4},{Z}^{5}],6,[y,{Z}^{6}],[y,{Z}^{6}],[{Z}^{5},{Z}^{4},{Z}^{3},{Z}^{2},Z,1\\
\mbox{}],6]\]}
\end{maplelatex}
\mapleresult
\begin{maplelatex}
\mapleinline{inert}{2d}{[[y, Z^5], [y, Z^5], [1, Z, Z^2, Z^3, Z^4], 5, [y, Z^5], [y, Z^5], [Z^4, Z^3, Z^2, Z, 1], 5]}{\[\displaystyle [[y,{Z}^{5}],[y,{Z}^{5}],[1,Z,{Z}^{2},{Z}^{3},{Z}^{4}],5,[y,{Z}^{5}],[y,{Z}^{5}],[{Z}^{4},{Z}^{3},{Z}^{2},Z,1],5\\
\mbox{}]\]}
\end{maplelatex}
\mapleresult
\begin{maplelatex}
\mapleinline{inert}{2d}{[[y, Z^4], [y, Z^4], [1, Z, Z^2, Z^3], 4,
[y, Z^4], [y, Z^4], [Z^3, Z^2, Z, 1], 4]}{\[\displaystyle
[[y,{Z}^{4}],[y,{Z}^{4}],[1,Z,{Z}^{2},{Z}^{3}],4,[y,{Z}^{4}],[y,{Z}^{4}],[{Z}^{3},{Z}^{2},Z,1],4]\]}
\end{maplelatex}
\mapleresult
\begin{maplelatex}
\mapleinline{inert}{2d}{[[y, Z^3], [y, Z^3], [1, Z, Z^2], 3, [y,
Z^3], [y, Z^3], [Z^2, Z, 1], 3]}{\[\displaystyle
[[y,{Z}^{3}],[y,{Z}^{3}],[1,Z,{Z}^{2}],3,[y,{Z}^{3}],[y,{Z}^{3}],[{Z}^{2},Z,1],3]\]}
\end{maplelatex}
\mapleresult
\begin{maplelatex}
\mapleinline{inert}{2d}{[[y, Z^2], [y, Z^2], [1, Z], 2, [y, Z^2],
[y, Z^2], [Z, 1], 2]}{\[\displaystyle
[[y,{Z}^{2}],[y,{Z}^{2}],[1,Z],2,[y,{Z}^{2}],[y,{Z}^{2}],[Z,1],2]\]}
\end{maplelatex}
\mapleresult
\begin{maplelatex}
\mapleinline{inert}{2d}{[[Z, y], [Z, y], [1], 1, [Z, y], [Z, y],
[1], 1]}{\[\displaystyle [[Z,y],[Z,y],[1],1,[Z,y],[Z,y],[1],1]\]}
\end{maplelatex}
\end{maplegroup}
\begin{maplegroup}
\begin{mapleinput}
\mapleinline{active}{2d}{}{\[\]}
\end{mapleinput}
\end{maplegroup}

\halfline

\noindent We reported here both the output obtained with
\texttt{Tyurina} and \texttt{TYURINA} to give a further account of
differences between the two procedures: after some time, we had to
stop \texttt{Tyurina} during the calculation of the last output,
while \texttt{TYURINA} was able to quickly conclude the
calculation; in fact the standard bases of $I_{FL_j}$ obtained by
the latter procedure are noticeably simpler than those obtained by
the former procedure. In the following we will employ
\texttt{TYURINA} every time \texttt{Tyurina} will be inefficient.

\vskip 5pt
\noindent (3) First of all let us observe that, for any $j=2,\ldots,7$,
$\tau(0)=\mu(0)$ meaning that the origin deforms always as a w.h.
singularity, by Proposition \ref{whs-caratterizzazione}. Moreover
$\crk(0)\leq 1$, since the rank of the Hessian matrix is always at
least 1 for the contribution of $y^2$. Then \emph{the origin
deforms always as a simple $A_{\mu}$ singularity, where $\mu$ is
its Milnor number} by Theorem \ref{classificazione}(2). In
particular the standard basis of $J_{FL_2}$ is given by
$\{2y,8Z^7+2v_2Z+3v_3Z^2+4v_4Z^3+5v_5Z^4+6v_6Z^5+7v_7Z^6\}$ whose
leading coefficient w.r.t. the default l.m.o. is $2v_2$. This
means that
\begin{itemize}
    \item \emph{$v_2=0$ is the equation of the codimension 1
subvariety $\mathcal{V}_2\subset\mathcal{L}$ parameterizing small
deformation of $U_0$ admitting a singularity of type at least
$A_2$ in the origin}.
\end{itemize}
Analogously the standard basis of $J_{FL_3}$ is
$\{2y,8Z^7+3v_3Z^2+4v_4Z^3+5v_5Z^4+6v_6Z^5+7v_7Z^6\}$ whose
leading coefficient is $3v_3$ meaning that
\begin{itemize}
    \item \emph{$v_3=0$ is the equation
of the subvariety $\mathcal{V}_3\subset\mathcal{L}\subset T^1$
such that $\mathcal{V}_2\cap\mathcal{V}_3\subset \mathcal{L}$ is
the codimension 2 subvariety parameterizing small deformation of
$U_0$ admitting a singularity of type at least $A_3$ in the
origin},
\end{itemize}
and so on: in general, for $n\geq 3$, setting
\begin{eqnarray}\label{lamda-condizione:A_n-2}
    v_k &:=&{n+1\choose k} z_{\Lambda}^{n+1-k} +
    \sum_{i=k}^{n-1}{i\choose k} \lambda_i
    z_{\Lambda}^{i-k}  \quad,\quad k=2,\ldots,n-1 \\
 \nonumber
    v_n&:=&z_{\Lambda}
\end{eqnarray}
and defining codimension 1 subvarieties $\mathcal{V}_k:=\{v_k=0\}$
of $\mathcal{L}$, then
\begin{itemize}
    \item \emph{$p_{\Lambda}\in U_{\Lambda}$ turns out to
be a $A_m\ ( 2\leq m\leq n)$ simple hypersurface singularity if
and only if $\Lambda$ is the generic element of the codimension
$m$ subvariety $\mathcal{V}_2^m:=\bigcap_{k=2}^m
\mathcal{V}_k\subset T^1$}.
\end{itemize}
This gives precisely the nested stratification
(\ref{An-inclusioni}) verifying the top row (\ref{An-adiacenze})
in diagram (\ref{adiacenza}) of Arnol'd's adjacencies.
\end{proof}

\subsection{Simple singularities of $D_n$ type.}

\begin{theorem}\label{Dn} Let $T^1$ be the Kuranishi space of a simple
$N$--dimensional singular point $0\in f^{-1}(0)$ with
\[
    f(x_1,\ldots,x_{N+1}) = \sum_{i=1}^{N-1} x_i^2\ +\
    x_N^2\ x_{N+1}+ x_{N+1}^{n-1}\quad\text{(for $n\geq 4$)}
\]
The subset of $T^1$ parameterizing small deformations of $0\in f^{-1}(0)$ to a simple node is the union of $n$ hypersurfaces. Moreover, calling $\mathcal{L}$ any of those hypersurfaces, there exists a stratification of nested
algebraic subsets giving rise to the following sequence of inclusions and c.i.p. squares
\begin{equation}\label{Dn-inclusioni}
    \xymatrix{\mathcal{L}&\ar@{_{(}->}[l]\ \mathcal{W}_2&\ar@{_{(}->}[l]\
                \mathcal{V}_0^1\cup\mathcal{W}_2^3&
                \ar@{_{(}->}[l]\cdots \mathcal{W}_2^m\cdots&\ar@{_{(}->}[l]\ \mathcal{W}_2^{n-1}\\
                &&\ar@{_{(}->}[u]\
                \mathcal{V}_0^2&\ar@{_{(}->}[u]\ar@{_{(}->}[l]\cdots
                \mathcal{V}_0^{m-1}\cdots&\ar@{_{(}->}[u]\ar@{_{(}->}[l]\ \mathcal{V}_0^{n-2}=\{0\}}
\end{equation}
verifying the Arnol'd's adjacency diagram
\begin{equation}\label{Dn-adiacenza}
    \xymatrix{A_1&\ar[l]A_2&\ar[l]A_3&\cdots\ar[l]A_m\cdots&\ar[l]A_{n-1}\\
                &&\ar[u]D_4&\ar[u]\cdots\ar[l]D_m\cdots&\ar[u]\ar[l]D_n}
\end{equation}
where
\begin{itemize}
    \item $\mathcal{L}$ is the hypersurface of $T^1$ defined by equation
(\ref{lamda-condizione:D_n}), keeping in mind
(\ref{D_n,z,y}),
    \item $\mathcal{V}_0^m:=\bigcap_{k=0}^m \mathcal{V}_k$
and $\mathcal{W}_2^m:=\bigcap_{k=2}^m \mathcal{W}_k$ where
$\mathcal{V}_k, \mathcal{W}_k$ are hypersurfaces of $\mathcal{L}$
defined by equations (\ref{VW-Dn}), keeping in mind definitions
(\ref{vk-Dn}).
\end{itemize}
\end{theorem}

\begin{remark}\label{D6asFig2}
The fact that every square in diagram (\ref{Dn-inclusioni}) is c.i.p. means that
\begin{equation}\label{intersezioni}
    \mathcal{V}_0^{m-1} =
    \mathcal{V}_0^{m-2}\cap\mathcal{W}_0^{m}\quad \text{for any $4\leq
    m\leq n-1$}\ .
\end{equation}
In particular Figure 2 represents geometrically the stratification of $\mathcal{V}_0^1\cup\mathcal{W}_2^3\subset T^1$ in the case of a $D_6$ singularity, by setting
\[
    F=\mathcal{V}_0^1\cup\mathcal{W}_2^3\ ,\ E=\mathcal{V}_0^2\ ,\ D=\mathcal{W}_2^4\ ,\ C=\mathcal{V}_0^3\ ,\ B=\mathcal{W}_2^5\ ,\ A=\mathcal{V}_0^4=\{0\}\ .
\]
\end{remark}

\begin{proof} Following the outline \ref{outline}.

\vskip 5pt\noindent (1) By the Morse Splitting Lemma \ref{Morse}, our problem can be reduced to the case $N=1$ with
$f(y,z)= y^2z + z^{n-1}$ for $n\geq 4$. Let us start by typing

\halfline

\begin{maplegroup}
\begin{mapleinput}
\mapleinline{active}{1d}{for n from 4 to 10 do FD[n]:=y\symbol{94}2*z+z\symbol{94}(n-1) end do }{}
\end{mapleinput}
\mapleresult
\begin{maplelatex}
\mapleinline{inert}{2d}{FD[4] := y^2*z+z^3}{\[\displaystyle {\it
FD}_{{4}}\, := \,{y}^{2}z+{z}^{3}\]}
\end{maplelatex}
\mapleresult
\begin{maplelatex}
\mapleinline{inert}{2d}{FD[5] := y^2*z+z^4}{\[\displaystyle {\it
FD}_{{5}}\, := \,{y}^{2}z+{z}^{4}\]}
\end{maplelatex}
\mapleresult
\begin{maplelatex}
\mapleinline{inert}{2d}{FD[6] := y^2*z+z^5}{\[\displaystyle {\it
FD}_{{6}}\, := \,{y}^{2}z+{z}^{5}\]}
\end{maplelatex}
\mapleresult
\begin{maplelatex}
\mapleinline{inert}{2d}{FD[7] := y^2*z+z^6}{\[\displaystyle {\it
FD}_{{7}}\, := \,{y}^{2}z+{z}^{6}\]}
\end{maplelatex}
\mapleresult
\begin{maplelatex}
\mapleinline{inert}{2d}{FD[8] := y^2*z+z^7}{\[\displaystyle {\it
FD}_{{8}}\, := \,{y}^{2}z+{z}^{7}\]}
\end{maplelatex}
\mapleresult
\begin{maplelatex}
\mapleinline{inert}{2d}{FD[9] := y^2*z+z^8}{\[\displaystyle {\it
FD}_{{9}}\, := \,{y}^{2}z+{z}^{8}\]}
\end{maplelatex}
\mapleresult
\begin{maplelatex}
\mapleinline{inert}{2d}{FD[10] := y^2*z+z^9}{\[\displaystyle {\it
FD}_{{10}}\, := \,{y}^{2}z+{z}^{9}\]}
\end{maplelatex}
\end{maplegroup}
\begin{maplegroup}
\begin{mapleinput}
\mapleinline{active}{1d}{for n from 4 to 10 do TB[n]:=TyurinaBasis(FD[n]) end do}{}
\mapleinline{active}{2d}{}{$$}
\end{mapleinput}
\mapleresult
\begin{maplelatex}
\mapleinline{inert}{2d}{TB[4]:= [z^2, z, y, 1]}{\[\displaystyle
{\it TB}_{{4}}\, := \,[{z}^{2},z,y,1]\]}
\end{maplelatex}
\mapleresult
\begin{maplelatex}
\mapleinline{inert}{2d}{TB[5] := [z^3, z^2, z, y,
1]}{\[\displaystyle {\it TB}_{{5}}\, :=
\,[{z}^{3},{z}^{2},z,y,1]\]}
\end{maplelatex}
\mapleresult
\begin{maplelatex}
\mapleinline{inert}{2d}{TB[6] := [z^4, z^3, z^2, z, y,
1]}{\[\displaystyle {\it TB}_{{6}}\, :=
\,[{z}^{4},{z}^{3},{z}^{2},z,y,1]\]}
\end{maplelatex}
\mapleresult
\begin{maplelatex}
\mapleinline{inert}{2d}{TB[7] := [z^5, z^4, z^3, z^2, z, y,
1]}{\[\displaystyle {\it TB}_{{7}}\, :=
\,[{z}^{5},{z}^{4},{z}^{3},{z}^{2},z,y,1]\]}
\end{maplelatex}
\mapleresult
\begin{maplelatex}
\mapleinline{inert}{2d}{TB[8] := [z^6, z^5, z^4, z^3, z^2, z, y,
1]}{\[\displaystyle {\it TB}_{{8}}\, :=
\,[{z}^{6},{z}^{5},{z}^{4},{z}^{3},{z}^{2},z,y,1]\]}
\end{maplelatex}
\mapleresult
\begin{maplelatex}
\mapleinline{inert}{2d}{TB[9] := [z^7, z^6, z^5, z^4, z^3, z^2, z,
y, 1]}{\[\displaystyle {\it TB}_{{9}}\, :=
\,[{z}^{7},{z}^{6},{z}^{5},{z}^{4},{z}^{3},{z}^{2},z,y,1]\]}
\end{maplelatex}
\mapleresult
\begin{maplelatex}
\mapleinline{inert}{2d}{TB[10] := [z^8, z^7, z^6, z^5, z^4, z^3,
z^2, z, y, 1]}{\[\displaystyle {\it TB}_{{10}}\, :=
\,[{z}^{8},{z}^{7},{z}^{6},{z}^{5},{z}^{4},{z}^{3},{z}^{2},z,y,1]\]}
\end{maplelatex}
\end{maplegroup}

\halfline

\noindent In general, for $n\geq 4$, the Kuranishi space $T^1$
will be then given by
\[
T^1 = \langle 1,y,z,\ldots, z^{n-2}\rangle_{\C}\ .
\]
Given $\Lambda = (\lambda_0, \lambda,
\lambda_1,\ldots,\lambda_{n-2})\in T^1$, the associated small
deformation of $U_0$ is
\[
    U_{\Lambda} = \{f_{\Lambda}(y,z) := f(y,z) + \lambda y + \sum_{i=0}^{n-2} \lambda_i
    z^i =0
    \}\ .
\]
A solution of the jacobian system of partials is then given by a solution
$p_{\Lambda}=(y_{\Lambda},z_{\Lambda})$ of the following polynomial system in $\C[\lambda][y,z]$
\begin{equation}\label{D_n,z,y}
    \left\{\begin{array}{c}
      2yz + \lambda = 0 \\
      (n-1)z^{n-2} + y^2 + \sum_{i=1}^{n-2} i \lambda_i z^{i-1}= 0\\
    \end{array}\right.
\end{equation}
giving precisely $n$ critical points for $f_{\Lambda}$.

\vskip 5pt\noindent (2) Imposing that one of those critical points, say $p_{\Lambda}$, is actually a singular point of $U_{\Lambda}$ means to require that
\begin{equation}\label{lamda-condizione:D_n}
    p_{\Lambda}\in U_{\Lambda}\ \Longleftrightarrow\ \Lambda\in
    \mathcal{L}:=\{ y_{\Lambda}^2z_{\Lambda} + z_{\Lambda}^{n-1} +
    \lambda y_{\Lambda} + \sum_{i=0}^{n-2} \lambda_i z_{\Lambda}^i =
    0\}\subset T^1\
\end{equation}
where, as above, $\mathcal{L}$ is one of the $n$ hypersurfaces of $T^1$
parameterizing small deformations of $0\in U_0$ to nodes. After
translating $y\mapsto y+y_{\Lambda},z\mapsto z+z_{\Lambda}$, we
get
\begin{eqnarray}\label{fLambda-Dn}
\nonumber
  f_{\Lambda}(y+y_{\Lambda},z+z_{\Lambda}) = f &+&
    \left(y_{\Lambda}^2z_{\Lambda} + z_{\Lambda}^{n-1} +
    \lambda y_{\Lambda} + \sum_{i=0}^{n-2} \lambda_i z_{\Lambda}^i\right) \\
    &+& \left(2y_{\Lambda}z_{\Lambda} + \lambda\right)y \\
\nonumber
    &+& \left((n-1)z_{\Lambda}^{n-2} + y_{\Lambda}^2 + \sum_{i=1}^{n-2} i \lambda_i z_{\Lambda}^{i-1}\right)z \\
\nonumber
    &+& 2y_{\Lambda}\ y z\ +\  z_{\Lambda}\ y^2\\
\nonumber
    &+& \sum_{k=2}^{n-2} \left( {n-1\choose k} z_{\Lambda}^{n-1-k} +
    \sum_{i=k}^{n-2}{i\choose k} \lambda_i
    z_{\Lambda}^{i-k}\right)z^k
\end{eqnarray}
Then by (\ref{lamda-condizione:D_n}) and (\ref{D_n,z,y}) the
origin (i.e. $p_{\Lambda}\in U_{\Lambda}$) is at least a node.

\vskip 5pt\noindent (3) We want at first recover, by means of a suitable algebraic stratification of the Kuranishi space $T^1$, the upper row of the adjacency diagram (\ref{adiacenza}). At this purpose consider the case $n=6$:

\halfline

\begin{maplegroup}
\begin{mapleinput}
\mapleinline{active}{1d}{T := nops(TB[6]);}{\[\]}
\end{mapleinput}
\end{maplegroup}
\begin{maplegroup}
\begin{mapleinput}
\mapleinline{active}{1d}{F[Lambda] := FD[6]+lambda[0]+lambda*y+
sum(lambda[i-1]*TB[6][T-i], i = 2 .. T-1); }{\[\]}
\end{mapleinput}
\mapleresult
\begin{maplelatex}
\mapleinline{inert}{2d}{`F&Lambda;` := y^2*z+z^5+lambda[0]+lambda*y+lambda[1]*z+lambda[2]*z^2+lambda[3]*z^3+lambda[4]*z^4}{\[\displaystyle F\Lambda \, := \,{y}^{2}z+{z}^{5}+\lambda_{{0}}+\lambda\,y+\lambda_{{1}}z+\lambda_{{2}}{z}^{2}+\lambda_{{3}}{z}^{3}+\lambda_{{4}}{z}^{4}\]}
\end{maplelatex}
\end{maplegroup}

\halfline

\noindent Compare $F\Lambda$ with $f_{\Lambda}$ above.

\halfline
\begin{maplegroup}
\begin{mapleinput}
\mapleinline{active}{1d}{solve({diff(F[Lambda], y), diff(F[Lambda], z)}, [y, z]); }{\[\]}
\end{mapleinput}
\mapleresult
\begin{maplelatex}
\mapleinline{inert}{2d}{}
{\[\displaystyle [[y=-1/2\,{\frac {\lambda}{{\it RootOf} \left( {\lambda}^{2}+20\,{{\it \_Z}}^{6}+4\,\lambda_{{1}}{{\it \_Z}}^{2}+8\,\lambda_{{2}}{{\it \_Z}}^{3}+12\,\lambda_{{3}}{{\it \_Z}}^{4}+16\,\lambda_{{4}}{{\it \_Z}}^{5} \right) }},\]}
\mapleinline{inert}{2d}{}
{\[\displaystyle \quad z={\it RootOf} \left( {\lambda}^{2}+20\,{{\it \_Z}}^{6}+4\,\lambda_{{1}}{{\it \_Z}}^{2}+8\,\lambda_{{2}}{{\it \_Z}}^{3}\\
\mbox{}+12\,\lambda_{{3}}{{\it \_Z}}^{4}+16\,\lambda_{{4}}{{\it \_Z}}^{5} \right) ]]\]}
\end{maplelatex}
\end{maplegroup}

\halfline

\noindent Compare these solutions with $(y_{\Lambda},z_{\Lambda})$ in (\ref{D_n,z,y}). Set $K:=y_{\Lambda}$ and $L:=z_{\Lambda}$ then

\halfline

\begin{maplegroup}
\begin{mapleinput}
\mapleinline{active}{1d}{z := Z+L : y := Y+K : }{\[\]}
\end{mapleinput}
\end{maplegroup}
\begin{maplegroup}
\begin{mapleinput}
\mapleinline{active}{21}{F[Lambda] := collect(F[Lambda], [Y, Z], 'distributed'); }{\[\]}
\end{mapleinput}
\mapleresult
\begin{maplelatex}
\mapleinline{inert}{2d}{}
{\[\displaystyle F\Lambda \, := \,L{Y}^{2}+Z{Y}^{2}+2\,KYZ+ \left( \lambda+2\,KL \right) Y+{Z}^{5}+ \left( 5\,L+\lambda_{{4}} \right) {Z}^{4} +\]}
\mapleinline{inert}{2d}{}
{\[\displaystyle\  \left( 10\,{L}^{2}+\lambda_{{3}}+4\,\lambda_{{4}}L \right) {Z}^{3}
+ \left( 6\,\lambda_{{4}}{L}^{2}+\lambda_{{2}}+3\,\lambda_{{3}}L+10\,{L}^{3} \right) {Z}^{2} +\]}
\mapleinline{inert}{2d}{}
{\[\displaystyle\  \left( 4\,\lambda_{{4}}{L}^{3}+5\,{L}^{4}+{K}^{2}+\lambda_{{1}}+3\,\lambda_{{3}}{L}^{2}
+2\,\lambda_{{2}}L \right) Z+\]}
\mapleinline{inert}{2d}{}
{\[\displaystyle\ {K}^{2}L+\lambda_{{0}}+\lambda\,K+\lambda_{{4}}{L}^{4}+{L}^{5}+\lambda_{{2}}{L}^{2}
+\lambda_{{1}}L+\lambda_{{3}}{L}^{3}\]}
\end{maplelatex}
\end{maplegroup}

\halfline

\noindent Compare $F\Lambda$ in the latter output with (\ref{fLambda-Dn}).
As already done for $A_n$ singularities in (\ref{lamda-condizione:A_n-2}), let us call $v_k$ the coefficients of singular deformations $F\Lambda$ of $FD_6$. Precisely, let us type:

\halfline

\begin{maplegroup}
\begin{mapleinput}
\mapleinline{active}{1d}{SD := [Y*Z, Y^2, seq(Z^i, i = 2 .. 4)]; }{\[\]}
\end{mapleinput}
\mapleresult
\begin{maplelatex}
\mapleinline{inert}{2d}{SD := [Y*Z, Y^2, Z^2, Z^3, Z^4]}{\[\displaystyle {\it SD}\, := \,[YZ,{Y}^{2},{Z}^{2},{Z}^{3},{Z}^{4}]\]}
\end{maplelatex}
\end{maplegroup}
\begin{maplegroup}
\begin{mapleinput}
\mapleinline{active}{1d}{s := nops(SD):}{\[\]}
\end{mapleinput}
\end{maplegroup}
\begin{maplegroup}
\begin{mapleinput}
\mapleinline{active}{1d}{FL := Y^2*Z+Z^5+sum(v[i]*SD[i+1], i = 0 .. s-1); }{\[\]}
\end{mapleinput}
\mapleresult
\begin{maplelatex}
\mapleinline{inert}{2d}{FL := Y^2*Z+Z^5+v[0]*Y*Z+v[1]*Y^2+v[2]*Z^2+v[3]*Z^3+v[4]*Z^4}{\[\displaystyle {\it FL}\, := \,{Y}^{2}Z+{Z}^{5}+v_{{0}}YZ+v_{{1}}{Y}^{2}+v_{{2}}{Z}^{2}+v_{{3}}{Z}^{3}+v_{{4}}{Z}^{4}\]}
\end{maplelatex}
\end{maplegroup}
\begin{maplegroup}
\begin{mapleinput}
\mapleinline{active}{1d}{MGB := MilnorGroebnerBasis(FL, \{Y, Z\});}{\[\]}
\end{mapleinput}
\mapleresult
\begin{maplelatex}
\mapleinline{inert}{2d}{}
{\[\displaystyle {\it MGB}\, := \, \left\{ 2\,YZ+v_{{0}}Z+2\,v_{{1}}Y,{Y}^{2}+5\,{Z}^{4}+v_{{0}}Y+2\,v_{{2}}Z+3\,v_{{3}}{Z}^{2}+4\,v_{{4}}{Z}^{3},\right.\]}
\mapleinline{inert}{2d}{}
{\[\displaystyle\quad \left. 4\,v_{{2}}YZ-3\,v_{{0}}v_{{3}}{Z}^{2}-4\,v_{{0}}v_{{4}}{Z}^{3}-5\,v_{{0}}{Z}^{4}-v_{{0}}{Y}^{2}+ \left( -{v_{{0}}}^{2}+4\,v_{{2}}v_{{1}} \right) Y \right\} \]}
\end{maplelatex}
\end{maplegroup}
\begin{maplegroup}
\begin{mapleinput}
\mapleinline{active}{1d}{for i to nops(MGB) do LeadingTerm(MGB[i],tdeg_min(Y,Z))end do;}{\[\]}
\end{mapleinput}
\mapleresult
\begin{maplelatex}
\mapleinline{inert}{2d}{}{\[\displaystyle v_{{0}},\,Z\]}
\end{maplelatex}
\mapleresult
\begin{maplelatex}
\mapleinline{inert}{2d}{}{\[\displaystyle 2\,v_{{2}},\,Z\]}
\end{maplelatex}
\mapleresult
\begin{maplelatex}
\mapleinline{inert}{2d}{-v[0]^2+4*v[2]*v[1], Y}{\[\displaystyle -{v_{{0}}}^{2}+4\,v_{{2}}v_{{1}},\,Y\]}
\end{maplelatex}
\end{maplegroup}

\halfline

\noindent The previous output gives the couples (\emph{leading coefficient}, \emph{leading monomial})
of the three generators in the standard basis $MGB$ of the jacobian ideal $J_{FL}$ w.r.t. the l.m.o. \texttt{tdeg\_min(Y,Z)} \footnote{It is defined as the \emph{opposite} of the reverse lexicographic g.m.o. on variables $Y,Z$ i.e. it is a graduated and lexicographic l.m.o. for which $Y<Z$.}:
\begin{itemize}
  \item \emph{the possible relations on the coefficients $v_k$ making a Milnor number growing up have then to be chosen among the three leading coefficients}.
\end{itemize}
On the contrary of the $A_n$ case, since both $v_0$ and $v_2$ multiply the same leading monomial $Z$, they have to be annihilated together to get a Milnor number increasing, describing a codimension 2 algebraic subset of $\mathcal{L}$. Actually it turns out to be an algebraic subset of the codimension 1 algebraic subset of $\mathcal{L}$ obtained by annihilating the third leading coefficient
\begin{equation}\label{A2 in Dn}
    \mathcal{W}_2:=\left\{4v_1v_2-v_0^2=0\right\}
\end{equation}
Notice that this is the same codimension 1 algebraic subset of $\mathcal{L}$ obtained by imposing the vanishing of the Hessian determinant $\det(\Hess(0))=0$, in fact

\halfline

\begin{maplegroup}
\begin{mapleinput}
\mapleinline{active}{1d}{with(VectorCalculus): with(LinearAlgebra):
H:=unapply(Hessian(FL, [Y, Z]), [Y, Z]):}{\[\]}
\end{mapleinput}
\end{maplegroup}
\begin{maplegroup}
\begin{mapleinput}
\mapleinline{active}{1d}{print(Hess(0)=H(0, 0),det(Hess(0))=Determinant(H(0,0)));}{\[\]}
\end{mapleinput}
\mapleresult
\begin{maplelatex}
\mapleinline{inert}{2d}{Matrix(
{\[\displaystyle  \Hess(0)=\left[ \begin {array}{cc} 2\,v_{{1}}&v_{{0}}\\
\noalign{\medskip}v_{{0}}&2\,v_{{2}}\end {array} \right] ,\,\ det(\Hess(0))=-{v_{{0}}}^{2}+4\,v_{{2}}v_{{1}}\]}
\end{maplelatex}
\end{maplegroup}

\halfline

\noindent Moreover if at least one of the three coefficients
$v_0,v_1,v_2$ does not vanish then $\crk(0)\leq 1$ and we are
considering deformations of $0\in f^{-1}(0)$ to simple $A_m$
singularities, by Theorem \ref{classificazione}(2). Introducing
relation (\ref{A2 in Dn}) in the deformation $FL$ means to impose
that
\[
    v_0YZ +v_1Y^2 + v_2Z^2 = (w_1Y + w_2Z)^2
\]
where $w_i^2=v_i$. Then

\halfline

\begin{maplegroup}
\begin{mapleinput}
\mapleinline{active}{1d}{F2:=Y^2*Z+Z^5+(w[1]*Y+w[2]*Z)^2+sum(v[j]*SD[j+1],j=3..s-1);}{\[\]}
\end{mapleinput}
\mapleresult
\begin{maplelatex}
\mapleinline{inert}{2d}{F2 := Y^2*Z+Z^5+(w[1]*Y+w[2]*Z)^2+v[3]*Z^3+v[4]*Z^4}{\[\displaystyle {\it F2}\, := \,{Y}^{2}Z+{Z}^{5}+ \left( w_{{1}}Y+w_{{2}}Z \right) ^{2}+v_{{3}}{Z}^{3}+v_{{4}}{Z}^{4}\]}
\end{maplelatex}
\end{maplegroup}
\begin{maplegroup}
\begin{mapleinput}
\mapleinline{active}{21d}{print(mu=MilnorNumber(F2,{Y,Z}),tau=TYURINANumber(F2,{Y,Z}));}{\[\]}
\end{mapleinput}
\mapleresult
\begin{maplelatex}
\mapleinline{inert}{2d}{mu = 2, tau = 2}{\[\displaystyle \mu=2,\,\tau=2\]}
\end{maplelatex}
\end{maplegroup}

\halfline

\noindent implying, by Theorem \ref{classificazione}(2), that
$\mathcal{W}_2\subset\mathcal{L}$ parameterizes small deformations
of $0\in f^{-1}(0)$ to a simple $A_2$ singularity or equivalently
that $0\in f_{\Lambda}^{-1}(0)$ is a simple $A_2$ singular point
for generic $\Lambda\in \mathcal{W}_2$. As above, consider now the
following standard basis

\halfline
\begin{maplegroup}
\begin{mapleinput}
\mapleinline{active}{1d}{MGB2 := MilnorGroebnerBasis(F2, {Y, Z});}{\[\]}
\end{mapleinput}
\mapleresult
\begin{maplelatex}
\mapleinline{inert}{2d}{}
{\[\displaystyle {\it MGB2}\, := \, \left\{ 2\,YZ+2\, \left( w_{{1}}Y+w_{{2}}Z \right) w_{{1}},{Y}^{2}+5\,{Z}^{4}+2\, \left( w_{{1}}Y+w_{{2}}Z \right) w_{{2}}+3\,v_{{3}}{Z}^{2}+\right.\]}
\mapleinline{inert}{2d}{}
{\[\displaystyle\ 4\,v_{{4}}{Z}^{3},6\,w_{{1}}w_{{2}}v_{{3}}{Z}^{2}Y+ \left( -6\,v_{{3}}{w_{{1}}}^{2}-4\,{w_{{2}}}^{2}
 \right) Z{Y}^{2}\\
\mbox{}+ \left( -6\,v_{{3}}{w_{{1}}}^{4}-6\,{w_{{1}}}^{2}{w_{{2}}}^{2} \right) {Y}^{2}+\]}
\mapleinline{inert}{2d}{}
{\[\displaystyle\ \left. -10\,{w_{{1}}}^{2}{w_{{2}}}^{2}{Z}^{4}-8\,{w_{{1}}}^{2}{w_{{2}}}^{2}v_{{4}}{Z}^{3} \right\} \]}
\end{maplelatex}
\end{maplegroup}

\halfline

\noindent whose leading terms are listed as follows

\halfline

\begin{maplegroup}
\begin{mapleinput}
\mapleinline{active}{1d}{for i to nops(MGB2) do LeadingTerm(MGB2[i],tdeg_min(Y,Z))end do;}{\[\]}
\end{mapleinput}
\mapleresult
\begin{maplelatex}
\mapleinline{inert}{2d}{2*w[1]*w[2], Z}{\[\displaystyle 2\,w_{{1}}w_{{2}},\,Z\]}
\end{maplelatex}
\mapleresult
\begin{maplelatex}
\mapleinline{inert}{2d}{2*w[2]^2, Z}{\[\displaystyle 2\,{w_{{2}}}^{2},\,Z\]}
\end{maplelatex}
\mapleresult
\begin{maplelatex}
\mapleinline{inert}{2d}{-6*v[3]*w[1]^4-6*w[1]^2*w[2]^2, Y^2}{\[\displaystyle -6\,v_{{3}}{w_{{1}}}^{4}-6\,{w_{{1}}}^{2}{w_{{2}}}^{2},\,{Y}^{2}\]}
\end{maplelatex}
\end{maplegroup}

\halfline

\noindent Observe that the vanishing of $w_2$, hence $v_2$, does
not change the Milnor number since one get $Y^2$ as leading monomial
from the last generator. Then, as before, the interesting
relation is given by the third leading coefficient, precisely:
$v_1(v_1v_3+v_2)$. Set
\[
    \mathcal{V}_1:=\left\{v_1=0\right\}\quad,\quad\mathcal{W}_3:=\left\{v_1v_3+v_2=0\right\}
\]
and observe that:

\halfline

\begin{maplegroup}
\begin{mapleinput}
\mapleinline{active}{1d}{w[1] := 0;}{\[\]}
\end{mapleinput}
\begin{mapleinput}
\mapleinline{active}{1d}{print(mu=MilnorNumber(F2,{Y,Z}),tau=TyurinaNumber(F2,{Y,Z}));}{\[\]}
\end{mapleinput}
\mapleresult
\begin{maplelatex}
\mapleinline{inert}{2d}{}
{\[\displaystyle \mu=3,\,\tau=3\]}
\end{maplelatex}
\end{maplegroup}

\halfline

\noindent On the other hand, by setting $w_2=iw_1w_3$ with $w_i^2=v_i$, one gets \footnote{The reader may check that setting $w_2=-iw_1w_3$ leads to the same results.}

\halfline

\begin{maplegroup}
\begin{mapleinput}
\mapleinline{active}{1d}{F3 := Z*Y^2+Z^5+v[1]*(Y+I*w[3]*Z)^2+w[3]^2*Z^3+v[4]*Z^4;}{\[\]}
\end{mapleinput}
\mapleresult
\begin{maplelatex}
\mapleinline{inert}{2d}{F3 := Z*Y^2+Z^5+v[1]*(Y+I*w[3]*Z)^2+w[3]^2*Z^3+v[4]*Z^4}{\[\displaystyle {\it F3}\, := \,Z{Y}^{2}+{Z}^{5}+v_{{1}} \left( Y+iw_{{3}}Z \right) ^{2}+{w_{{3}}}^{2}{Z}^{3}+v_{{4}}{Z}^{4}\]}
\end{maplelatex}
\end{maplegroup}
\begin{maplegroup}
\begin{mapleinput}
\mapleinline{active}{1d}{print(mu=MilnorNumber(F3,{Y,Z}),tau=TyurinaNumber(F3,{Y,Z})); }{\[\]}
\end{mapleinput}
\mapleresult
\begin{maplelatex}
\mapleinline{inert}{2d}{mu = 3, tau = 3}{\[\displaystyle \mu=3,\,\tau=3\]}
\end{maplelatex}
\end{maplegroup}

\halfline

\noindent Once more Theorem \ref{classificazione}(2) gives that
\begin{itemize}
    \item \emph{$0\in f^{-1}_{\Lambda}(0)$ is a simple $A_3$ singularity for
a generic $\Lambda$ in}
\begin{equation}\label{Dn-3}
    \mathcal{W}_2\cap(\mathcal{V}_1\cup\mathcal{W}_3)=(\mathcal{V}_0\cap\mathcal{V}_1)\cup
    (\mathcal{W}_2\cap\mathcal{W}_3)\subset\mathcal{W}_2\subset\mathcal{L}
\end{equation}
\end{itemize}
where $\mathcal{V}_k:=\{v_k=0\}$. Let us go on by considering the
first connected component in (\ref{Dn-3}) and by applying the same
argument as above, precisely

\halfline

\begin{maplegroup}
\begin{mapleinput}
\mapleinline{active}{1d}{v[1] := 0: }{\[\]}
\end{mapleinput}
\end{maplegroup}
\begin{maplegroup}
\begin{mapleinput}
\mapleinline{active}{1d}{MGB3:=MilnorGroebnerBasis(F2,\{Y,Z\}); }{\[\]}
\end{mapleinput}
\mapleresult
\begin{maplelatex}
\mapleinline{inert}{2d}{}
{\[\displaystyle {\it MGB3}\, := \, \left\{ 2\,{Y}^{3}+10\,{Z}^{4}Y+8\,Yv_{{4}}{Z}^{3},{Y}^{2}+5\,{Z}^{4}+2\,{w_{{2}}}^{2}Z+3\,v_{{3}}{Z}^{2}+4\,v_{{4}}{Z}^{3}\\
\mbox{},2\,YZ \right\} \]}
\end{maplelatex}
\end{maplegroup}
\begin{maplegroup}
\begin{mapleinput}
\mapleinline{active}{1d}{for i to nops(MGB3) do LeadingTerm(MGB3[i],tdeg_min(Y,Z))end do; }{\[\]}
\end{mapleinput}
\mapleresult
\begin{maplelatex}
\mapleinline{inert}{2d}{2*w[2]^2, Z}{\[\displaystyle
2\,{w_{{2}}}^{2},\,Z\]}
\end{maplelatex}
\mapleresult
\begin{maplelatex}
\mapleinline{inert}{2d}{2, Y*Z}{\[\displaystyle 2,\,YZ\]}
\end{maplelatex}
\mapleresult
\begin{maplelatex}
\mapleinline{inert}{2d}{2, Y^3}{\[\displaystyle 2,\,{Y}^{3}\]}
\end{maplelatex}
\end{maplegroup}

\halfline

\noindent Therefore the only interesting relation is given by the first leading coefficient, precisely $v_2=0$:

\halfline

\begin{maplegroup}
\begin{mapleinput}
\mapleinline{active}{1d}{w[2] := 0; }{\[\]}
\end{mapleinput}
\end{maplegroup}
\begin{maplegroup}
\begin{mapleinput}
\mapleinline{active}{1d}{print(mu=MilnorNumber(F2,\{Y,Z\}),tau=TyurinaNumber(F2,\{Y,Z\}));}{\[\]}
\end{mapleinput}
\mapleresult
\begin{maplelatex}
\mapleinline{inert}{2d}{mu = 4, tau = 4}{\[\displaystyle
\mu=4,\,\tau=4\]}
\end{maplelatex}
\end{maplegroup}

\halfline

\noindent Since now $\crk(0)=2$, Theorem \ref{classificazione}(3)
gives that:
\begin{itemize}
    \item \emph{$0\in f_{\Lambda}^{-1}(0)$ is a simple $D_4$ singularity for a generic point
    $\Lambda\in\mathcal{V}_0\cap\mathcal{V}_1\cap\mathcal{V}_2$.}
\end{itemize}
Let us now consider the second connected component
$\mathcal{W}_2\cap\mathcal{W}_3$ in (\ref{Dn-3}). As before

\halfline

\begin{maplegroup}
\begin{mapleinput}
\mapleinline{active}{1d}{MGB3b := MilnorGroebnerBasis(F3,{Y,Z});}{\[\]}
\end{mapleinput}
\mapleresult
\begin{maplelatex}
\mapleinline{inert}{2d}{} {\[\displaystyle {\it MGB3b}\, := \,
\left\{2\,YZ+2\,v_{{1}} \left( Y+iw_{{3}}Z
\right),{Y}^{2}+5\,{Z}^{4}+2\,iv_{{1}} \left( Y+iw_{{3}}Z \right)
w_{{3}}+3\,{w_{{3}}}^{2}{Z}^{2}+\right.\]}
\mapleinline{inert}{2d}{} {\[\displaystyle\quad
4\,v_{{4}}{Z}^{3},10\,i{w_{{3}}}^{3}{v_{{1}}}^{2}{Z}^{4}+ \left(
-8\,iw_{{3}}v_{{1}}v_{{4}}-6\,i{w_{{3}}}^{3} \right)
{Y}^{2}{Z}^{2} +\left( 8\,{w_{{3}}}^{2}+8\,v_{{1}}v_{{4}} \right)
{Y}^{3}Z+\]} \mapleinline{inert}{2d}{} {\[\displaystyle\quad
\left. -8\, {w_{{3}}}^{2}v_{{1}}v_{{4}}{Z}^{3}Y+ \left(
8\,{v_{{1}}}^{2}v_{{4}}+8\,{w_{{3}}}^{2}v_{{1}} \right)
{Y}^{3}\right\} \]}
\end{maplelatex}
\end{maplegroup}
\begin{maplegroup}
\begin{mapleinput}
\mapleinline{active}{1d}{for i to nops(MGB3b) do LeadingTerm(MGB3b[i],tdeg_min(Y,Z))end do;}{\[\]}
\end{mapleinput}
\mapleresult
\begin{maplelatex}
\mapleinline{inert}{2d}{(2*I)*v[1]*w[3], Z}{\[\displaystyle
2\,iv_{{1}}w_{{3}},\,Z\]}
\end{maplelatex}
\mapleresult
\begin{maplelatex}
\mapleinline{inert}{2d}{-2*w[3]^2*v[1], Z}{\[\displaystyle
-2\,{w_{{3}}}^{2}v_{{1}},\,Z\]}
\end{maplelatex}
\mapleresult
\begin{maplelatex}
\mapleinline{inert}{2d}{8*v[1]^2*v[4]+8*w[3]^2*v[1],
Y^3}{\[\displaystyle
8\,{v_{{1}}}^{2}v_{{4}}+8\,{w_{{3}}}^{2}v_{{1}},\,{Y}^{3}\]}
\end{maplelatex}
\end{maplegroup}

\halfline

\noindent The only interesting relation is given by the last
leading coefficient, precisely $v_1(v_1v_4+v_3)$. Set then
\[
    \mathcal{W}_4:=\left\{v_1v_4+v_3=0\right\}
\]
and introduce the new relation in the algebraic set
$\mathcal{W}_2\cap\mathcal{W}_3$ by intersecting it with
$\mathcal{V}_1$ and $\mathcal{W}_4$ to get the following
codimension 1 algebraic subset
\begin{equation}\label{Dn-4}
    \left(\mathcal{V}_1\cap\mathcal{W}_2\cap\mathcal{W}_3\right)
    \cup\left(\mathcal{W}_2\cap\mathcal{W}_3\cap\mathcal{W}_4\right)\ .
\end{equation}
First of all observe that the first connected component in
(\ref{Dn-4}) reduces to the previous case since
$\mathcal{V}_1\cap\mathcal{W}_2\cap\mathcal{W}_3 =
    \mathcal{V}_0\cap\mathcal{V}_1\cap\mathcal{V}_2$.
Let us then consider the second connected component in
(\ref{Dn-4}) by introducing the relation $w_3=iw_1w_4$, with
$w_i^2=v_i$, in $F3$. Then

\halfline

\begin{maplegroup}
\begin{mapleinput}
\mapleinline{active}{1d}{F4:=Z*Y^2+Z^5+(w[1]*Y-w[1]^2*w[4]*Z)^2-w[1]^2*w[4]^2*Z^3+
\quad w[4]^2*Z^4;}{\[\]}
\end{mapleinput}
\mapleresult
\begin{maplelatex}
\mapleinline{inert}{2d}{F4 := Z*Y^2+Z^5+(w[1]*Y-w[1]^2*w[4]*Z)^2-w[1]^2*w[4]^2*Z^3+w[4]^2*Z^4}
{\[\displaystyle {\it F4}\, := \,Z{Y}^{2}+{Z}^{5}+ \left( w_{{1}}Y-{w_{{1}}}^{2}w_{{4}}Z \right) ^{2}
-{w_{{1}}}^{2}{w_{{4}}}^{2}{Z}^{3}\\
\mbox{}+{w_{{4}}}^{2}{Z}^{4}\]}
\end{maplelatex}
\end{maplegroup}
\begin{maplegroup}
\begin{mapleinput}
\mapleinline{active}{1d}{print(mu=MilnorNumber(F4,\{Y,Z\}),tau=TYURINANumber(F4,\{Y,Z\}));}{\[\]}
\end{mapleinput}
\mapleresult
\begin{maplelatex}
\mapleinline{inert}{2d}{mu = 4, tau = 4}{\[\displaystyle
\mu=4,\,\tau=4\]}
\end{maplelatex}
\end{maplegroup}

\halfline

\noindent Theorem \ref{classificazione}(2) then gives that
\begin{itemize}
    \item \emph{$0\in f^{-1}_{\Lambda}(0)$ is a simple
$A_4$ singularity for a generic
$\Lambda\in\mathcal{W}_2\cap\mathcal{W}_3\cap\mathcal{W}_4$}.
\end{itemize}
At the moment we have gotten the following chain of codimension 1
algebraic subsets
\[
    \xymatrix{\mathcal{L}&\ar@{_{(}->}[l]\ \mathcal{W}_2&\ar@{_{(}->}[l]\
                (\mathcal{V}_0\cap\mathcal{V}_1)\cup(\mathcal{W}_2\cap\mathcal{W}_3)&
                \ar@{_{(}->}[l]\ \mathcal{W}_2\cap\mathcal{W}_3\cap\mathcal{W}_4\\
                &&\ar@{_{(}->}[u]\ \mathcal{V}_0\cap\mathcal{V}_1\cap\mathcal{V}_2}
\]
representing the adjacency diagram
\[
    \xymatrix{A_1&\ar[l] A_2&\ar[l] A_3&\ar[l] A_4\\
              &&D_4\ar[u] }
\]
Go on by considering the standard basis of
$J_{f_\Lambda}$ when $\Lambda\in\mathcal{V}_0\cap\mathcal{V}_1\cap\mathcal{V}_2$:

\halfline

\begin{maplegroup}
\begin{mapleinput}
\mapleinline{active}{1d}{v[0]:=0 : v[1]:=0 : v[2]:=0 :}{\[\]}
\end{mapleinput}
\end{maplegroup}
\begin{maplegroup}
\begin{mapleinput}
\mapleinline{active}{1d}{MGB4 := MilnorGroebnerBasis(FL,\{Y,Z\})}{\[ \]}
\end{mapleinput}
\mapleresult
\begin{maplelatex}
\mapleinline{inert}{2d}{}
{\[\displaystyle {\it MGB4}\, := \,
\left\{
-2\,{Y}^{3}-10\,{Z}^{4}Y-8\,Yv_{{4}}{Z}^{3},2\,YZ,{Y}^{2}+5\,{Z}^{4}+3\,v_{{3}}{Z}^{2}+4\,v_{{4}}{Z}^{3}
\right\} \]}
\end{maplelatex}
\end{maplegroup}
\begin{maplegroup}
\begin{mapleinput}
\mapleinline{active}{1d}{for i to nops(MGB4) do LeadingTerm(MGB4[i],tdeg_min(Y,Z))end do;}{\[\]}
\end{mapleinput}
\mapleresult
\begin{maplelatex}
\mapleinline{inert}{2d}{2, Y*Z}{\[\displaystyle 2,\,YZ\]}
\end{maplelatex}
\mapleresult
\begin{maplelatex}
\mapleinline{inert}{2d}{-2, Y^3}{\[\displaystyle -2,\,{Y}^{3}\]}
\end{maplelatex}
\mapleresult
\begin{maplelatex}
\mapleinline{inert}{2d}{3*v[3], Z^2}{\[\displaystyle
3\,v_{{3}},\,{Z}^{2}\]}
\end{maplelatex}
\end{maplegroup}

\halfline

\noindent and then impose the further condition given by the
latter leading coefficient, namely $v_3=0$:

\halfline

\begin{maplegroup}
\begin{mapleinput}
\mapleinline{active}{1d}{v[3] := 0 :}{\[\]}
\end{mapleinput}
\end{maplegroup}
\begin{maplegroup}
\begin{mapleinput}
\mapleinline{active}{1d}{print(mu=MilnorNumber(FL,\{Y,Z\}),tau=TyurinaNumber(FL,\{Y,Z\}));}{\[\]}
\end{mapleinput}
\mapleresult
\begin{maplelatex}
\mapleinline{inert}{2d}{mu = 5, tau = 5}{\[\displaystyle
\mu=5,\,\tau=5\]}
\end{maplelatex}
\end{maplegroup}

\halfline

\noindent Theorem \ref{classificazione}(3) then gives that
\begin{itemize}
    \item \emph{$0\in f^{-1}_{\Lambda}$ is a simple $D_5$ singularity for generic
    $\Lambda\in\bigcap_{k=0}^3\mathcal{V}_k$}.
\end{itemize}
On the other hand, for
$\Lambda\in\mathcal{W}_2\cap\mathcal{W}_3\cap\mathcal{W}_4$ set:

\halfline

\begin{maplegroup}
\begin{mapleinput}
\mapleinline{active}{1d}{MGB4b := MilnorGroebnerBasis(F4, \{Y, Z\}); }{\[\]}
\end{mapleinput}
\mapleresult
\begin{maplelatex}
\mapleinline{inert}{2d}{}
{\[\displaystyle {\it MGB4b}\, := \,
\left\{ 2\,YZ+2\, \left( w_{{1}}Y-{w_{{1}}}^{2}w_{{4}}Z \right)
w_{{1}}, {Y}^{2}+5\,{Z}^{4}-2\, \left(
w_{{1}}Y-{w_{{1}}}^{2}w_{{4}}Z \right) {w_{{1}}}^{2}
w_{{4}}+\right.\]}
\mapleinline{inert}{2d}{}
{\[\displaystyle\ -3\,{w_{{1}}}^{2}{w_{{4}}}^{2}{Z}^{2}+4\,{w_{{4}}}^{2}{Z}^{3},
\left( -8\,{w_{{1}}}^{2}{w_{{4}}}^{4}-10\,{w_{{1}}}^{4}{w_{{4}}}^{2} \right) {Z}^{3}{Y}^{2}+
 \left( -10\,{w_{{4}}}^{2}-10\,{w_{{1}}}^{2} \right) Z{Y}^{4}+\]}
\mapleinline{inert}{2d}{}
{\[\displaystyle\ \left.-10\,{w_{{1}}}^{5}{w_{{4}}}^{3}{Z}^{4}Y\\
\mbox{}+ \left( -10\,{w_{{1}}}^{3}w_{{4}}-10\,w_{{1}}{w_{{4}}}^{3}
\right) {Y}^{3}{Z}^{2}+ \left(
-10\,{w_{{1}}}^{2}{w_{{4}}}^{2}-10\,{w_{{1}}}^{4} \right)
{Y}^{4}\right\} \]}
\end{maplelatex}
\end{maplegroup}
\begin{maplegroup}
\begin{mapleinput}
\mapleinline{active}{1d}{for i to nops(MGB4b) do LeadingTerm(MGB4b[i],tdeg_min(Y,Z))end do;}{\[\]}
\end{mapleinput}
\mapleresult
\begin{maplelatex}
\mapleinline{inert}{2d}{}{\[\displaystyle
-2\,{w_{{1}}}^{3}w_{{4}},\,Z\]}
\end{maplelatex}
\mapleresult
\begin{maplelatex}
\mapleinline{inert}{2d}{}{\[\displaystyle
2\,{w_{{1}}}^{4}{w_{{4}}}^{2},\,Z\]}
\end{maplelatex}
\mapleresult
\begin{maplelatex}
\mapleinline{inert}{2d}{}{\[\displaystyle
-10\,{w_{{1}}}^{2}{w_{{4}}}^{2}-10\,{w_{{1}}}^{4},\,{Y}^{4}\]}
\end{maplelatex}
\end{maplegroup}

\halfline

\noindent As usual the only interesting relation is the last one,
precisely $v_1(v_1+v_4)$. Set $\mathcal{W}_5:=\{v_1+v_4=0\}$ and intersect the algebraic subset defined by the latter relation with the second component in (\ref{Dn-4}) to get
\begin{eqnarray}\label{Dn-5}
\nonumber
  \left(\mathcal{V}_1\cup\mathcal{W}_5\right)\cap\left(\mathcal{W}_2\cap\mathcal{W}_3\cap\mathcal{W}_4\right) &=&
  \left(\mathcal{V}_1\cap\mathcal{W}_2\cap\mathcal{W}_3\cap\mathcal{W}_4\right)\cup
\left(\mathcal{W}_2\cap\mathcal{W}_3\cap\mathcal{W}_4\cap\mathcal{W}_5\right) \\
   &=&\left(\mathcal{V}_0\cap\mathcal{V}_1\cap\mathcal{V}_2\cap\mathcal{V}_3\right)\cup
   \left(\mathcal{W}_2\cap\mathcal{W}_3\cap\mathcal{W}_4\cap\mathcal{W}_5\right)
\end{eqnarray}
The first connected component in the right term gives the previous
case just considered, then look at the second component by
introducing the relation $w_4=iw_1$ in $F4$:

\halfline

\begin{maplegroup}
\begin{mapleinput}
\mapleinline{active}{1d}{F5:=Z*Y^2+Z^5+v[1]*(Y-I*v[1]*Z)^2+v[1]^2*Z^3-v[1]*Z^4;}{\[\]}
\end{mapleinput}
\mapleresult
\begin{maplelatex}
\mapleinline{inert}{2d}{F5 :=
Z*Y^2+Z^5+v[1]*(Y-I*v[1]*Z)^2+v[1]^2*Z^3-v[1]*Z^4}{\[\displaystyle
{\it F5}\, := \,Z{Y}^{2}+{Z}^{5}+v_{{1}} \left( Y-iv_{{1}}Z
\right) ^{2}+{v_{{1}}}^{2}{Z}^{3}-v_{{1}}{Z}^{4}\]}
\end{maplelatex}
\end{maplegroup}
\begin{maplegroup}
\begin{mapleinput}
\mapleinline{active}{1d}{print(mu=MilnorNumber(F5,\{Z,Y\}),tau=TyurinaNumber(F5,\{Z,Y\})); }{\[\]}
\end{mapleinput}
\mapleresult
\begin{maplelatex}
\mapleinline{inert}{2d}{mu = 5, tau = 5}{\[\displaystyle
\mu=5,\,\tau=5\]}
\end{maplelatex}
\end{maplegroup}

\halfline

\noindent and Theorem \ref{classificazione}(2) allows to conclude
that
\begin{itemize}
    \item \emph{$0\in f^{-1}_{\Lambda}$ is a simple $A_5$ singularity for generic
    $\Lambda\in\bigcap_{k=2}^5\mathcal{W}_k$}.
\end{itemize}
Since the algebraic subsets in (\ref{Dn-5}) parameterizes
1-parameter deformations, the further and last step is clearly the
trivial deformation given by $0\in T^1$, completing the inclusions
diagram as follows
\begin{equation}\label{D6-inclusioni}
    \xymatrix{\mathcal{L}&\ar@{_{(}->}[l]\ \mathcal{W}_2&\ar@{_{(}->}[l]\
                \mathcal{V}_0^1\cup\mathcal{W}_2^3&
                \ar@{_{(}->}[l]\ \mathcal{W}_2^4&\ar@{_{(}->}[l]\ \mathcal{W}_2^5\\
                &&\ar@{_{(}->}[u]\ \mathcal{V}_0^2&\ar@{_{(}->}[u]\ar@{_{(}->}[l]\
                \mathcal{V}_0^3&\ar@{_{(}->}[u]\ar@{_{(}->}[l]\ \{0\}}
\end{equation}
where $\mathcal{V}_a^b:=\bigcap_{k=a}^b\mathcal{V}_k$ and
$\mathcal{W}_a^b:=\bigcap_{k=a}^b\mathcal{W}_k$. Diagram
(\ref{D6-inclusioni}) gives a stratification, by nested
algebraic subsets, of the Kuranishi space of a
simple $D_6$ singular point, representing the following adjacency
diagram
\begin{equation}\label{D6-adiacenza}
    \xymatrix{A_1&\ar[l]A_2&\ar[l]A_3&\ar[l]A_4&\ar[l]A_5\\
                &&\ar[u]D_4&\ar[u]\ar[l]D_5&\ar[u]\ar[l]D_6}
\end{equation}
and geometrically represented by Figure 2, as already observed in Remark \ref{D6asFig2}.

\noindent The recursive structure is now sufficiently clear to pass at the
general step for $n\geq 4$. Recalling $f_{\Lambda}$ as written in
(\ref{fLambda-Dn}) define
\begin{eqnarray}\label{vk-Dn}
\nonumber
    v_0&:=&2y_{\Lambda}  \\
    v_{1}&:=&z_{\Lambda}\\
\nonumber
    v_k&:=&{n-1\choose k} z_{\Lambda}^{n-1-k} +
    \sum_{i=k}^{n-2}{i\choose k} \lambda_i
    z_{\Lambda}^{i-k} \ ,\ k=2,\ldots,n-2
\end{eqnarray}
and consider the associated codimension 1 sub--schemes of
$\mathcal{L}$
\begin{eqnarray}\label{VW-Dn}
  \mathcal{V}_k &:=& \left\{v_k = 0\right\}\ ,\quad 0\leq k\leq n-2 \\
\nonumber
  \mathcal{W}_k &:=& \left\{\begin{array}{cc}
    \left\{4v_1v_2-v_0^2=0\right\} & \text{for $k=2$ (vanishing of $\det(\Hess)$)} \\
    \left\{v_1v_k+v_{k-1}=0\right\} & \text{for $3\leq k\leq n-2$} \\
    \left\{v_1+v_{n-2}=0\right\} & \text{for $k=n-1$} \\
  \end{array}
  \right.
\end{eqnarray}
Then diagrams (\ref{D6-inclusioni}) and (\ref{D6-adiacenza})
generalizes to give diagrams (\ref{Dn-inclusioni}) and
(\ref{Dn-adiacenza}) in the statement, respectively. Let us
conclude by pointing out that in diagram (\ref{D6-inclusioni}),
$\mathcal{V}_0^3=\mathcal{V}_0^2\cap\mathcal{W}_2^4$. This fact
generalizes to diagram (\ref{Dn-inclusioni}) giving the stated
relations (\ref{intersezioni}) between spaces of such a
stratification.
\end{proof}

\subsection{Simple singularities of $E_6$ type.}

\begin{theorem}\label{E6}
Let $T^1$ be the Kuranishi space of a simple
$N$--dimensional singular point $0\in f^{-1}(0)$ with
\[
    f(x_1,\ldots,x_{N+1}) = \sum_{i=1}^{N-1} x_i^2\ +\
    x_N^3+ x_{N+1}^4
\]
The subset of $T^1$ parameterizing small deformations of $0\in f^{-1}(0)$ to a simple node is the union of 6 hypersurfaces. Moreover, calling $\mathcal{L}$ any of those hypersurfaces, there exists a stratification of nested
algebraic subsets giving rise to the following sequence of inclusions and c.i.p. squares
\begin{equation}\label{E6-inclusioni}
    \xymatrix{\mathcal{L}&\ar@{_{(}->}[l]\ \mathcal{W}_2&\ar@{_{(}->}[l]\
                \mathcal{W}_2^3&
                \ar@{_{(}->}[l]\ \widetilde{\mathcal{W}}_2^4&\ar@{_{(}->}[l]\ \mathcal{W}\cap\left(\bigcap_{k=0}^2\mathcal{V}_{2k}\right)\\
                &&\ar@{_{(}->}[u]\ \mathcal{V}_0^2&\ar@{_{(}->}[u]\ar@{_{(}->}[l]\
                \mathcal{V}\cap\mathcal{V}_0^2&&\\
                &&&\ar@{_{(}->}[ruu]\ar@{_{(}->}[u]\ \{0\}}
\end{equation}
verifying the Arnol'd's adjacency diagram
\begin{equation}\label{E6-adiacenza}
    \xymatrix{A_1&\ar[l]A_2&\ar[l]A_3&\ar[l]A_4&\ar[l]A_5\\
                &&\ar[u]D_4&\ar[u]\ar[l]D_5&&\\
                &&&\ar[ruu]\ar[u]E_6}
\end{equation}
where
\begin{itemize}
    \item $\mathcal{L}$ is the hypersurface of $T^1$ defined by equation
(\ref{lamda-condizione:E_6}), keeping in mind
(\ref{E_6,z,y}),
    \item $\mathcal{V}_0^m:=\bigcap_{k=0}^m \mathcal{V}_k$
and $\mathcal{W}_2^m:=\bigcap_{k=2}^m \mathcal{W}_k$ where
$\mathcal{V}_k, \mathcal{W}_k$ are hypersurfaces of $\mathcal{L}$
defined by equations (\ref{vk-E6}), (\ref{A_2 in E6}) and (\ref{A3
in E6}),
    \item $\mathcal{V}$ is a hypersurfaces of $\mathcal{L}$ defined
    by equation (\ref{D5 in E6}),
    \item $\widetilde{\mathcal{W}}_2^4$ is a codimension 3 Zariski closed
    subset of $\mathcal{L}$ defined in  (\ref{A_4bis in E6}),
    \item $\mathcal{W}$ is a hypersurfaces of $\mathcal{L}$ defined
    by equation (\ref{A5 in E6}).
\end{itemize}
\end{theorem}

\begin{remark}
Analogously to what observed in remark \ref{D6asFig2}, the fact that every square in diagram (\ref{E6-inclusioni}) is c.i.p. means that
\begin{equation}\label{intersezioni in E6}
    \widetilde{\mathcal{W}}_2^4\cap\mathcal{V}_0^2=\mathcal{V}\cap\mathcal{V}_0^2\quad,\quad
    \left(\mathcal{W}\cap\left(\bigcap_{k=0}^2\mathcal{V}_{2k}\right)\right)\cap\left(\mathcal{V}\cap\mathcal{V}_0^2\right)=\{0\}
\end{equation}
In particular Figure 2 represents geometrically the stratification of $\mathcal{W}_2^3\subset T^1$, by setting
\[
    F=\mathcal{W}_2^3\ ,\ E=\mathcal{V}_0^2\ ,\ D=\widetilde{\mathcal{W}}_2^4\ ,\ C=\mathcal{V}\cap\mathcal{V}_0^2\ ,\ B=\mathcal{W}\cap\left(\bigcap_{k=0}^2\mathcal{V}_{2k}\right)\ ,\ A=\{0\}\ .
\]
\end{remark}

\begin{proof} Following the outline \ref{outline}.

\vskip 5pt\noindent (1) By the Morse Splitting Lemma \ref{Morse}, our problem can be reduced to the case $N=1$ with $f(y,z)=  y^3 +
z^4$. To get an explicit basis of the Kuranishi space $T^1$ type

\halfline

\begin{maplegroup}
\begin{mapleinput}
\mapleinline{active}{1d}{F:= y\symbol{94}3 + z\symbol{94}4; }{}
\end{mapleinput}
\mapleresult
\begin{maplelatex}
\mapleinline{inert}{2d}{F := y^3+z^4}{\[\displaystyle F\, :=
\,{y}^{3}+{z}^{4}\]}
\end{maplelatex}
\end{maplegroup}
\begin{maplegroup}
\begin{mapleinput}
\mapleinline{active}{1d}{TyB := TyurinaBasis(F); }{\[\]}
\end{mapleinput}
\mapleresult
\begin{maplelatex}
\mapleinline{inert}{2d}{TyB := [y*z^2, z^2, y*z, z, y,
1]}{\[\displaystyle {\it TyB}\, :=
\,[y{z}^{2},{z}^{2},yz,z,y,1]\]}
\end{maplelatex}
\end{maplegroup}

\halfline

\noindent Therefore
\[
T^1 \cong  \langle 1,y,z,yz,z^2,yz^2\rangle_{\C}\ .
\]
Given $\Lambda = (\lambda_0, \lambda_1,\ldots,\lambda_5)\in T^1$,
the associated deformation of $U_0$ is
\begin{eqnarray*}
   U_{\Lambda} &=& \{f_{\Lambda}(y,z)=0\}\quad\text{where}  \\
   f_{\Lambda}(y,z) &:=& f(y,z) + \lambda_0 + \lambda_1 y
    + \lambda_2 z + \lambda_3 yz + \lambda_4 z^2 + \lambda_5 yz^2
\end{eqnarray*}
A solution of the jacobian system of partials is then given by a solution
$p_{\Lambda}=(y_{\Lambda},z_{\Lambda})$ of the following polynomial system in $\C[\lambda][y,z]$
\begin{equation}\label{E_6,z,y}
   \left\{\begin{array}{c}
      3y^2 + \lambda_1 + \lambda_3z + \lambda_5z^2 = 0 \\
      4z^3 + 2 \lambda_4z +\lambda_2 +
      y(\lambda_3 + 2\lambda_5 z)= 0\\
    \end{array}\right.
\end{equation}
giving precisely 6 critical points for $f_{\Lambda}$.

\vskip 5pt\noindent (2) Imposing that one of those critical points, say $p_{\Lambda}$, is actually a singular point of $U_{\Lambda}$ means to require that
\begin{eqnarray}\label{lamda-condizione:E_6}
    &p_{\Lambda}\in U_{\Lambda}\ \Longleftrightarrow \ \Lambda\in
    \mathcal{L}\quad\text{where}&\\
\nonumber
    &    \mathcal{L}:=\{ y_{\Lambda}^3 + z_{\Lambda}^4 + \lambda_0 + \lambda_1y_{\Lambda}
    + \lambda_2z_{\Lambda} + \lambda_3y_{\Lambda}z_{\Lambda} + \lambda_4z_{\Lambda}^2 +
    \lambda_5y_{\Lambda}z_{\Lambda}^2 =
    0\}&\subset T^1
\end{eqnarray}
which is \emph{one of the 6 hypersurfaces of $T^1$ parameterizing small
deformations of $0\in U_0$ to nodes.} After translating $y\mapsto
y+y_{\Lambda},z\mapsto z+z_{\Lambda}$, we get
\begin{eqnarray*}
  f_{\Lambda}(y+y_{\Lambda},z+z_{\Lambda})=& \\
  f(y,z)\quad  +
    &\left(y_{\Lambda}^3 + z_{\Lambda}^4 + \lambda_0 + \lambda_1y_{\Lambda}
    + \lambda_2z_{\Lambda} + \lambda_3y_{\Lambda}z_{\Lambda} + \lambda_4z_{\Lambda}^2 +
    \lambda_5y_{\Lambda}z_{\Lambda}^2\right) \\
    + &\left(3y_{\Lambda}^2 + \lambda_1 + \lambda_3z_{\Lambda} +
    \lambda_5z_{\Lambda}^2\right)y \\
    + &\left(4z_{\Lambda}^3 + 2 \lambda_4z_{\Lambda} +\lambda_2 +
      y_{\Lambda}(\lambda_3 + 2\lambda_5 z_{\Lambda})\right)z \\
    + &\left(\lambda_3 + 2\lambda_5z_{\Lambda}\right)y z + 3y_{\Lambda}\ y^2 +
    \left(6z_{\Lambda}^2 + \lambda_4 + \lambda_5y_{\Lambda}\right)z^2 \\
    + &4 z_{\Lambda}\ z^3 + \lambda_5\ yz^2
\end{eqnarray*}
as may be verified by setting $K:=y_{\Lambda}, L:=z_{\Lambda}$ and typing:

\halfline

\begin{maplegroup}
\begin{mapleinput}
\mapleinline{active}{1d}{F[Lambda]:= F+sum(lambda[i]*TyB[T-i],i = 0 .. T-1)}
{\[F\Lambda \, := \,F+\sum
_{i=0}^{T-1}\lambda_{{i}}{\it TyB}_{{T-i}}\]}
\end{mapleinput}
\mapleresult
\begin{maplelatex}
\mapleinline{inert}{2d}{}
{\[\displaystyle
F\Lambda \, :=
\,{y}^{3}+{z}^{4}+\lambda_{{0}}+\lambda_{{1}}y+\lambda_{{2}}z+\lambda_{{3}}yz+\lambda_{{4}}{z}^{2}+\lambda_{{5}}y{z}^{2}\]}
\end{maplelatex}
\end{maplegroup}
\begin{maplegroup}
\begin{mapleinput}
\mapleinline{active}{1d}{z:=Z+L : y:=Y+K :}{\[\]}
\end{mapleinput}
\end{maplegroup}
\begin{maplegroup}
\begin{mapleinput}
\mapleinline{active}{1d}{F[Lambda] := collect(F[Lambda], [Y, Z], 'distributed'); }{\[\]}
\end{mapleinput}
\mapleresult
\begin{maplelatex}
\mapleinline{inert}{2d}{}
{\[\displaystyle F\Lambda \, := \,{Y}^{3}+3\,K{Y}^{2}+\lambda_{{5}}Y{Z}^{2}+ \left( \lambda_{{3}}+2\,
\lambda_{{5}}L \right) YZ+ \left( \lambda_{{3}}L+\lambda_{{1}}+3\,{K}^{2}+\lambda_{{5}}{L}^{2} \right) Y+\]}
\mapleinline{inert}{2d}{}
{\[\displaystyle\  +{Z}^{4}+4\,L{Z}^{3}+ \left( \lambda_{{4}}+\lambda_{{5}}K+6\,{L}^{2} \right) {Z}^{2}+
\left( 2\,\lambda_{{4}}L+\lambda_{{3}}K+\lambda_{{2}}+2\,\lambda_{{5}}KL+4\,{L}^{3} \right) Z+\]}
\mapleinline{inert}{2d}{}
{\[\displaystyle \ +{K}^{3}+\lambda_{{0}}+\lambda_{{1}}K+\lambda_{{4}}{L}^{2}+\lambda_{{2}}L+{L}^{4}+
\lambda_{{3}}KL+\lambda_{{5}}K{L}^{2}\]}
\end{maplelatex}
\end{maplegroup}

\halfline

\vskip 5pt\noindent (3) Define:
\begin{equation}\label{vk-E6}
  v_0 = \lambda_3 + 2\lambda_5z_{\Lambda} \quad ,\quad
  v_1 =  3y_{\Lambda}\quad ,\quad v_2 =   6z_{\Lambda}^2 + \lambda_4 + \lambda_5y_{\Lambda}
  \quad ,\quad  v_3 = \lambda_5\quad ,\quad
  v_4 = 4z_{\Lambda}
\end{equation}
and $\mathcal{V}_k:=\{v_k=0\}$. Then type:

\halfline

\begin{maplegroup}
\begin{mapleinput}
\mapleinline{active}{1d}{SD := [Y*Z, Y^2, Z^2, Y*Z^2, Z^3]:}{\[\]}
\end{mapleinput}
\end{maplegroup}
\begin{maplegroup}
\begin{mapleinput}
\mapleinline{active}{1d}{s := nops(SD):}{\[\]}
\end{mapleinput}
\end{maplegroup}
\begin{maplegroup}
\begin{mapleinput}
\mapleinline{active}{1d}{FL := Y^3+Z^4+sum(v[i]*SD[i+1], i = 0 .. s-1);}{\[\]}
\end{mapleinput}
\mapleresult
\begin{maplelatex}
\mapleinline{inert}{2d}{}
{\[\displaystyle {\it FL}\, :=
\,{Y}^{3}+{Z}^{4}+v_{{0}}YZ+v_{{1}}{Y}^{2}+v_{{2}}{Z}^{2}+v_{{3}}Y{Z}^{2}+v_{{4}}{Z}^{3}\]}
\end{maplelatex}
\end{maplegroup}
\begin{maplegroup}
\begin{mapleinput}
\mapleinline{active}{1d}{MGB := MilnorGroebnerBasis(FL, \{Y,Z\});}{\[\]}
\end{mapleinput}
\mapleresult
\begin{maplelatex}
\mapleinline{inert}{2d}{} {\[\displaystyle {\it MGB}\, := \,
\left\{ -6\,v_{{2}}{Y}^{2}+ \left(
-2\,v_{{2}}v_{{3}}+3\,v_{{0}}v_{{4}} \right) {Z}^{2}
+2\,v_{{0}}v_{{3}}YZ+ \left( {v_{{0}}}^{2}-4\,v_{{2}}v_{{1}}
\right) Y+4\,v_{{0}}{Z}^{3},\right.\]} \mapleinline{inert}{2d}{}
{\[\displaystyle \quad \left.
3\,{Y}^{2}+v_{{0}}Z+2\,v_{{1}}Y+v_{{3}}{Z}^{2},4\,{Z}^{3}+v_{{0}}Y+2\,v_{{2}}Z+2\,v_{{3}}YZ+3\,
v_{{4}}{Z}^{2} \right\} \]}
\end{maplelatex}
\end{maplegroup}
\begin{maplegroup}
\begin{mapleinput}
\mapleinline{active}{12d}{for i to nops(MGB) do LeadingTerm(MGB[i],tdeg_min(Y,Z))end do;}{\[\]}
\end{mapleinput}
\mapleresult
\begin{maplelatex}
\mapleinline{inert}{2d}{v[0]^2-4*v[2]*v[1], Y}{\[\displaystyle
{v_{{0}}}^{2}-4\,v_{{2}}v_{{1}},\,Y\]}
\end{maplelatex}
\mapleresult
\begin{maplelatex}
\mapleinline{inert}{2d}{2*v[2], Z}{\[\displaystyle
2\,v_{{2}},\,Z\]}
\end{maplelatex}
\mapleresult
\begin{maplelatex}
\mapleinline{inert}{2d}{v[0], Z}{\[\displaystyle v_{{0}},\,Z\]}
\end{maplelatex}
\end{maplegroup}

\halfline

\noindent The only interesting relation is given by the first
leading coefficient giving the vanishing of
$\det(\Hess(0))=4v_1v_2-v_0^2$. Introduce it in $FL$ by setting
$v_k=w_k^2$ and considering

\halfline

\begin{maplegroup}
\begin{mapleinput}
\mapleinline{active}{1d}{F2:=Y^3+Z^4+(w[1]*Y+w[2]*Z)^2+sum(v[j]*SD[j+1], j = 3..s-1)}
{\[{\it
F2}\, := \,{Y}^{3}+{Z}^{4}+ \left( w_{{1}}Y+w_{{2}}Z \right)
^{2}+\sum _{j=3}^{s-1}v_{{j}}{\it SD}_{{j+1}}\]}
\end{mapleinput}
\mapleresult
\begin{maplelatex}
\mapleinline{inert}{2d}{F2 :=
Y^3+Z^4+(w[1]*Y+w[2]*Z)^2+v[3]*Y*Z^2+v[4]*Z^3}{\[\displaystyle
{\it F2}\, := \,{Y}^{3}+{Z}^{4}+ \left( w_{{1}}Y+w_{{2}}Z \right)
^{2}+v_{{3}}Y{Z}^{2}+v_{{4}}{Z}^{3}\]}
\end{maplelatex}
\end{maplegroup}
\begin{maplegroup}
\begin{mapleinput}
\mapleinline{active}{1d}{print(mu=MilnorNumber(F2,\{Y,Z\}),tau=TYURINANumber(F2,\{Y,Z\}));}{\[\]}
\end{mapleinput}
\mapleresult
\begin{maplelatex}
\mapleinline{inert}{2d}{mu = 2, tau = 2}{\[\displaystyle
\mu=2,\,\tau=2\]}
\end{maplelatex}
\end{maplegroup}

\halfline

\noindent Then Theorem \ref{classificazione}(2) ensures that
\begin{itemize}
    \item \emph{$0\in f_{\Lambda}^{-1}(0)$ is a simple $A_2$ singularity for $\Lambda$ generic in}
    \begin{equation}\label{A_2 in E6}
    \mathcal{W}_2:=\left\{4v_1v_2-v_0^2=0\right\}\ .
    \end{equation}
\end{itemize}
Consider the standard basis:

\halfline

\begin{maplegroup}
\begin{mapleinput}
\mapleinline{active}{1d}{MGB2 := MilnorGroebnerBasis(F2, \{Y,Z\}); }{\[\]}
\end{mapleinput}
\mapleresult
\begin{maplelatex}
\mapleinline{inert}{2d}{}
{\[\displaystyle {\it MGB2}\, := \, \left\{ 4\,{Z}^{3}+2\, \left( w_{{1}}Y+w_{{2}}Z \right) w_{{2}}+2\,v_{{3}}YZ
+3\,v_{{4}}{Z}^{2}, \left( 9\,{w_{{1}}}^{2}v_{{4}}-9\,w_{{1}}w_{{2}}v_{{3}} \right) {Y}^{3}\right.\]}
\mapleinline{inert}{2d}{}
{\[\displaystyle \quad +\left( -6\,{w_{{2}}}^{3}w_{{1}}+6\,{w_{{1}}}^{4}v_{{4}}-6\,{w_{{1}}}^{3}w_{{2}}v_{{3}} \right) {Y}^{2}+
 \left( -3\,w_{{1}}w_{{2}}{v_{{3}}}^{2}+3\,{w_{{1}}}^{2}v_{{4}}v_{{3}} \right)Y{Z}^{2}\]}
\mapleinline{inert}{2d}{}
{\[\displaystyle \quad + \left( -9\,w_{{1}}w_{{2}}v_{{4}}+3\,{w_{{2}}}^{2}v_{{3}} \right) {Y}^{2}Z+
\left( 8\,{w_{{1}}}^{2}{w_{{2}}}^{2}-3\,v_{{4}}v_{{3}}w_{{1}}w_{{2}}+{w_{{2}}}^{2}{v_{{3}}}^{2}\\
\mbox{} \right) {Z}^{3},\]}
\mapleinline{inert}{2d}{}
{\[\displaystyle \quad 3\,{Y}^{2}+2\, \left( w_{{1}}Y+w_{{2}}Z \right) w_{{1}}\\
\mbox{}+v_{{3}}{Z}^{2} \left.\right\} \]}
\end{maplelatex}
\end{maplegroup}
\begin{maplegroup}
\begin{mapleinput}
\mapleinline{active}{1d}{for i to nops(MGB2) do LeadingTerm(MGB2[i],tdeg_min(Y,Z))end do; }{\[\]}
\end{mapleinput}
\mapleresult
\begin{maplelatex}
\mapleinline{inert}{2d}{2*w[2]^2, Z}{\[\displaystyle
2\,{w_{{2}}}^{2},\,Z\]}
\end{maplelatex}
\mapleresult
\begin{maplelatex}
\mapleinline{inert}{2d}{2*w[1]*w[2], Z}{\[\displaystyle
2\,w_{{1}}w_{{2}},\,Z\]}
\end{maplelatex}
\mapleresult
\begin{maplelatex}
\mapleinline{inert}{2d}{-6*w[2]^3*w[1]+6*w[1]^4*v[4]-6*w[1]^3*w[2]*v[3],
Y^2}{\[\displaystyle
-6\,{w_{{2}}}^{3}w_{{1}}+6\,{w_{{1}}}^{4}v_{{4}}-6\,{w_{{1}}}^{3}w_{{2}}v_{{3}},\,{Y}^{2}\]}
\end{maplelatex}
\end{maplegroup}

\halfline

\noindent Relations coming from the first and the second leading
coefficients do not increase the Milnor number since:

\halfline

\begin{maplegroup}
\begin{mapleinput}
\mapleinline{active}{1d}{w[1] := 0:}{\[\]}
\end{mapleinput}
\end{maplegroup}
\begin{maplegroup}
\begin{mapleinput}
\mapleinline{active}{1d}{print(mu=MilnorNumber(F2,\{Z,Y\}),tau=TyurinaNumber(F2,\{Z,Y\})); }{\[\]}
\end{mapleinput}
\mapleresult
\begin{maplelatex}
\mapleinline{inert}{2d}{mu = 2, tau = 2}{\[\displaystyle
\mu=2,\,\tau=2\]}
\end{maplelatex}
\end{maplegroup}
\begin{maplegroup}
\begin{mapleinput}
\mapleinline{active}{1d}{unassign('w[1]'): w[2] := 0 :}{\[\]}
\end{mapleinput}
\end{maplegroup}
\begin{maplegroup}
\begin{mapleinput}
\mapleinline{active}{1d}{print(mu=MilnorNumber(F2,\{Z,Y\}), tau=TYURINANumber(F2,\{Z,Y\}))}
{\[{\mu={\it MilnorNumber} \left( {\it F2}, \left\{ Z,Y \right\}  \right) ,\tau={\it TYURINANumber}\\
\mbox{} \left( {\it F2}, \left\{ Z,Y \right\}  \right) }\]}
\end{mapleinput}
\mapleresult
\begin{maplelatex}
\mapleinline{inert}{2d}{mu = 2, tau = 2}{\[\displaystyle
\mu=2,\,\tau=2\]}
\end{maplelatex}
\end{maplegroup}

\halfline

\noindent Then consider the third leading coefficient, giving
$v_4={w_2\over w_1}v_3+\left({w_2\over w_1}\right)^3$, since we
can assume $w_1\neq 0$. Then set

\halfline

\begin{maplegroup}
\begin{mapleinput}
\mapleinline{active}{1d}{w[2] := u*w[1] : v[4] := u*(v[3]+u^2) : F2 ;}{\[\]}
\end{mapleinput}
\mapleresult
\begin{maplelatex}
\mapleinline{inert}{2d}{Y^3+Z^4+(w[1]*Y+u*w[1]*Z)^2+v[3]*Y*Z^2+u*(v[3]+u^2)*Z^3}{\[\displaystyle
{Y}^{3}+{Z}^{4}+ \left( w_{{1}}Y+uw_{{1}}Z \right)
^{2}+v_{{3}}Y{Z}^{2}+u \left( v_{{3}}+{u}^{2} \right) {Z}^{3}\]}
\end{maplelatex}
\end{maplegroup}
\begin{maplegroup}
\begin{mapleinput}
\mapleinline{active}{1d}{F3 := Y^3+Z^4+v[1]*(Y+u*Z)^2+v[3]*Y*Z^2+u*(v[3]+u^2)*Z^3}{\[{\it F3}\,
:= \,{Y}^{3}+{Z}^{4}+v_{{1}} \left( Y+uZ \right)
^{2}+v_{{3}}Y{Z}^{2}+u \left( v_{{3}}+{u}^{2} \right) {Z}^{3}\]}
\end{mapleinput}
\mapleresult
\begin{maplelatex}
\mapleinline{inert}{2d}{F3 :=
Y^3+Z^4+v[1]*(Y+u*Z)^2+v[3]*Y*Z^2+u*(v[3]+u^2)*Z^3}{\[\displaystyle
{\it F3}\, := \,{Y}^{3}+{Z}^{4}+v_{{1}} \left( Y+uZ \right)
^{2}+v_{{3}}Y{Z}^{2}+u \left( v_{{3}}+{u}^{2} \right) {Z}^{3}\]}
\end{maplelatex}
\end{maplegroup}
\begin{maplegroup}
\begin{mapleinput}
\mapleinline{active}{1d}{print(mu=MilnorNumber(F3, \{Z,Y\}),tau=TYURINANumber(F3,\{Z,Y\}))}
{\[{\mu={\it MilnorNumber} \left( {\it F3}, \left\{ Z,Y \right\}  \right) ,\tau={\it TYURINANumber}\\
\mbox{} \left( {\it F3}, \left\{ Z,Y \right\}  \right) }\]}
\end{mapleinput}
\mapleresult
\begin{maplelatex}
\mapleinline{inert}{2d}{mu = 3, tau = 3}{\[\displaystyle
\mu=3,\,\tau=3\]}
\end{maplelatex}
\end{maplegroup}

\halfline

\noindent and Theorem \ref{classificazione}(2) gives that \footnote{The reader may check that choosing $v_4=-u(v_3+u^2)$ leads to the same conclusion.}
\begin{itemize}
    \item \emph{$0\in f_{\Lambda}^{-1}(0)$ is a singularity of type $A_3$
    for any generic $\Lambda\in\mathcal{W}_2\cap\mathcal{W}_3$ where}
    \begin{equation}\label{A3 in E6}
    \mathcal{W}_3:=\left\{v_1^3v_4^2-v_2(v_1v_3+v_2)^2=0\right\}\ .
    \end{equation}
\end{itemize}
The latter equation is obtained by observing that, after
eliminating $w_1\neq 0$, the third leading coefficient in $MGB2$
gives
\[
    w_1^3v_4=\pm w_2(w_1^2v_3+w_2^2)
\]
(where the sign depends on the choice of the square roots
$w_i^2=v_i$) whose square is the equation of $\mathcal{W}_3$ in
(\ref{A3 in E6}). Let us now consider the standard basis:

\halfline

\begin{maplegroup}
\begin{mapleinput}
\mapleinline{active}{1d}{MGB3 := MilnorGroebnerBasis(F3,\{Y,Z\});}{\[\]}
\end{mapleinput}
\mapleresult
\begin{maplelatex}
\mapleinline{inert}{2d}{}
 {\[\displaystyle {\it MGB3}\, := \,
\left\{  \left(
16\,{v_{{1}}}^{2}-24\,v_{{1}}{u}^{2}v_{{3}}-4\,{v_{{3}}}^{2}v_{{1}}
-36\,v_{{1}}{u}^{4} \right) {Y}^{3}+ \left(
18\,{u}^{3}v_{{3}}-24\,v_{{1}}u+\right.\right.\]}
\mapleinline{inert}{2d}{} {\[\displaystyle\
\left.6\,u{v_{{3}}}^{2} \right) Z{Y}^{3}+ \left(
-3\,{u}^{4}{v_{{3}}}^{2}+8\,v_{{1}}{u}^{2}v_{{3}}-2\,{u}^{2}{v_{{3}}}^{3}
\right) {Z}^{4}+ \left(
-36\,{u}^{2}v_{{3}}-27\,{u}^{4}-6\,{v_{{3}}}^{2}+\right.\]}
\mapleinline{inert}{2d}{} {\[\displaystyle\ \left.24\,v_{{1}}
\right) {Y}^{4} +\left(
8\,v_{{1}}v_{{3}}-18\,{u}^{2}{v_{{3}}}^{2}+24\,v_{{1}}{u}^{2}-18\,{u}^{4}v_{{3}}-2\,
{v_{{3}}}^{3} \right) {Z}^{2}{Y}^{2}+ \]}
\mapleinline{inert}{2d}{} {\[\displaystyle\ \left(
-8\,v_{{1}}uv_{{3}}+6\,{u}^{3}{v_{{3}}}^{2}+2\,u{v_{{3}}}^{3}
\right) {Z}^{3}Y, 3\,{Y}^{2}+2\,v_{{1}} \left( Y+uZ \right)
+v_{{3}}{Z}^{2},\]} \mapleinline{inert}{2d}{} {\[\displaystyle\
4\,{Z}^{3}+2\,v_{{1}} \left( Y+uZ \right) u+2\,v_{{3}}YZ+3\,u
\left( v_{{3}}+{u}^{2} \right) {Z}^{2} \left.\right\} \]}
\end{maplelatex}
\end{maplegroup}
\begin{maplegroup}
\begin{mapleinput}
\mapleinline{active}{1d}{for i to nops(MGB3) do LeadingTerm(MGB3[i],tdeg_min(Y,Z))end do; }{\[\]}
\end{mapleinput}
\mapleresult
\begin{maplelatex}
\mapleinline{inert}{2d}{2*v[1]*u, Z}{\[\displaystyle
2\,v_{{1}}u,\,Z\]}
\end{maplelatex}
\mapleresult
\begin{maplelatex}
\mapleinline{inert}{2d}{2*v[1]*u^2, Z}{\[\displaystyle
2\,v_{{1}}{u}^{2},\,Z\]}
\end{maplelatex}
\mapleresult
\begin{maplelatex}
\mapleinline{inert}{2d}{16*v[1]^2-24*v[1]*u^2*v[3]-4*v[3]^2*v[1]-36*v[1]*u^4,
Y^3}{\[\displaystyle
16\,{v_{{1}}}^{2}-24\,v_{{1}}{u}^{2}v_{{3}}-4\,{v_{{3}}}^{2}v_{{1}}-36\,v_{{1}}{u}^{4},\,{Y}^{3}\]}
\end{maplelatex}
\end{maplegroup}

\halfline

\noindent The last leading coefficient can be rewritten as
$v_1\left[4v_1-(v_3+3u^2)^2\right]$, giving all the new
interesting relations, precisely $v_1=0$ and
$v_1=1/4(v_3+3u^2)^2$. Since:

\halfline

\begin{maplegroup}
\begin{mapleinput}
\mapleinline{active}{1d}{v[1] := 0:}{\[\]}
\end{mapleinput}
\end{maplegroup}
\begin{maplegroup}
\begin{mapleinput}
\mapleinline{active}{1d}{print(mu=MilnorNumber(F3,\{Z,Y\}),tau=TyurinaNumber(F3,\{Z,Y\}))}
{\[{\mu={\it MilnorNumber} \left( {\it F3}, \left\{ Z,Y \right\}  \right) ,\tau={\it TyurinaNumber}\\
\mbox{} \left( {\it F3}, \left\{ Z,Y \right\}  \right) }\]}
\end{mapleinput}
\mapleresult
\begin{maplelatex}
\mapleinline{inert}{2d}{mu = 4, tau = 4}{\[\displaystyle
\mu=4,\,\tau=4\]}
\end{maplelatex}
\end{maplegroup}
\begin{maplegroup}
\begin{mapleinput}
\mapleinline{active}{1d}{unassign('v[1]') : v[1] := (1/4)*(3*u^2+v[3])^2 :}{\[\]}
\end{mapleinput}
\end{maplegroup}
\begin{maplegroup}
\begin{mapleinput}
\mapleinline{active}{1d}{print(mu=MilnorNumber(F3,\{Z,Y\}),tau=TyurinaNumber(F3,\{Z,Y\}))}
{\[{\mu={\it MilnorNumber} \left( {\it F3}, \left\{ Z,Y \right\}  \right) ,\tau={\it TyurinaNumber}\\
\mbox{} \left( {\it F3}, \left\{ Z,Y \right\}  \right) }\]}
\end{mapleinput}
\mapleresult
\begin{maplelatex}
\mapleinline{inert}{2d}{mu = 4, tau = 4}{\[\displaystyle
\mu=4,\,\tau=4\]}
\end{maplelatex}
\end{maplegroup}

\halfline

\noindent the Milnor number is then increasing for $\Lambda$ in
\begin{equation}\label{E6 mu=4}
    \left(\mathcal{V}_1\cap\mathcal{W}_2\cap\mathcal{W}_3\right)\cup
    \left(\mathcal{W}_2\cap\mathcal{W}_3\cap\mathcal{W}_4\right)=
    \left(\bigcap_{k=0}^2\mathcal{V}_k\right)\cup\left(\bigcap_{k=2}^4\mathcal{W}_k\right)
\end{equation}
where, recalling that $u^2=v_2/v_1$,
\begin{equation}\label{A_4 in E6}
    \mathcal{W}_4:=\left\{4v_1^3-(v_1v_3+3v_2)^2=0\right\}\ .
\end{equation}
As usual by now, call
$\mathcal{V}_a^b:=\bigcap_{k=a}^b\mathcal{V}_k$ and
$\mathcal{W}_a^b:=\bigcap_{k=a}^b\mathcal{W}_k$. Then (\ref{E6
mu=4}) rewrites as $\mathcal{V}_0^2\cup\mathcal{W}_2^4$. Notice
that
\[
    \mathcal{V}_1^2\subseteq\mathcal{W}_3^4\quad\Longrightarrow\quad
    \mathcal{V}_0^2\subseteq\mathcal{W}_2^4\ .
\]
Define the following codimension 1 Zariski closed subset \footnote{\label{nota} Although (\ref{A_4bis in E6}) is a topological definition, $\widetilde{\mathcal{W}}_2^4$ is clearly an algebraic subset of $\mathcal{L}$. Unfortunately it is not a complete intersection but the reader may obtain the (long) list of generators of its defining ideal by employing the MAPLE command \texttt{EliminationIdeal} in the \texttt{PolynomialIdeals} package. Not all the so listed generators are effectively necessary. The reader may compare them with those obtained by employing the MAPLE command \texttt{eliminate}, being careful with possible multiplicity and reducibility of generators.} of
$\mathcal{W}_2^3$
\begin{equation}\label{A_4bis in E6}
    \widetilde{\mathcal{W}}_2^4:=\overline{\mathcal{W}_2^4\setminus\mathcal{V}_0^2}\
    .
\end{equation}
Then Theorem \ref{classificazione}(2) implies that
\begin{itemize}
    \item \emph{$0\in f_{\Lambda}^{-1}(0)$ is a singularity of type $A_4$ for any
    generic $\Lambda\in\widetilde{\mathcal{W}}_2^4$}
\end{itemize}
and Theorem \ref{classificazione}(3) ensures that
\begin{itemize}
    \item \emph{$0\in f_{\Lambda}^{-1}(0)$ is a singularity of type $D_4$ for any
    generic $\Lambda\in\mathcal{V}_0^2$}\ .
\end{itemize}
Let us firstly consider the latter case:

\halfline

\begin{maplegroup}
\begin{mapleinput}
\mapleinline{active}{1d}{v[0]:=0 : v[1]:=0 : v[2]:=0 : FL ;}{\[\]}
\end{mapleinput}
\mapleresult
\begin{maplelatex}
\mapleinline{inert}{2d}{Y^3+Z^4+v[3]*Y*Z^2+v[4]*Z^3}{\[\displaystyle {Y}^{3}+{Z}^{4}+v_{{3}}Y{Z}^{2}+v_{{4}}{Z}^{3}\]}
\end{maplelatex}
\end{maplegroup}
\begin{maplegroup}
\begin{mapleinput}
\mapleinline{active}{1d}{print(mu=MilnorNumber(FL,\{Y,Z\}),tau=TyurinaNumber(FL,\{Y,Z\}))}{\[{\mu={\it MilnorNumber} \left( {\it FL}, \left\{ Y,Z \right\}  \right) ,\tau={\it TyurinaNumber}\\
\mbox{} \left( {\it FL}, \left\{ Y,Z \right\}  \right) }\]}
\end{mapleinput}
\mapleresult
\begin{maplelatex}
\mapleinline{inert}{2d}{mu = 4, tau = 4}{\[\displaystyle \mu=4,\,\tau=4\]}
\end{maplelatex}
\end{maplegroup}
\begin{maplegroup}
\begin{mapleinput}
\mapleinline{active}{1d}{MGB4 := MilnorGroebnerBasis(FL,\{Y,Z\})}{\[{\it MGB4}\, := \,{\it MilnorGroebnerBasis} \left( {\it FL}, \left\{ Y,Z \right\} ,{\it tdeg\_min}\\
\mbox{} \left( Y,Z \right)  \right) \]}
\end{mapleinput}
\mapleresult
\begin{maplelatex}
\mapleinline{inert}{2d}{}
{\[\displaystyle {\it MGB4}\, := \, \left\{ 3\,{Y}^{2}+v_{{3}}{Z}^{2},4\,{Z}^{3}+2\,v_{{3}}YZ+3\,v_{{4}}{Z}^{2},9\,v_{{4}}{Y}^{2}-4\,v_{{3}}{Z}^{3}\\
\mbox{}-2\,{v_{{3}}}^{2}YZ,\right.\]}
\mapleinline{inert}{2d}{}
{\[\displaystyle\quad \left.-36\,v_{{4}}Yv_{{3}}{Z}^{3}-8\,{v_{{3}}}^{3}{Z}^{4}+ \left( 12\,{v_{{3}}}^{3}+81\,{v_{{4}}}^{2} \right) {Y}^{3}\\
\mbox{} \right\} \]}
\end{maplelatex}
\end{maplegroup}
\begin{maplegroup}
\begin{mapleinput}
\mapleinline{active}{1d}{for i to nops(MGB4) do LeadingTerm(MGB4[i],tdeg_min(Y,Z))end do; }{\[\]}
\end{mapleinput}
\mapleresult
\begin{maplelatex}
\mapleinline{inert}{2d}{v[3], Z^2}{\[\displaystyle v_{{3}},\,{Z}^{2}\]}
\end{maplelatex}
\mapleresult
\begin{maplelatex}
\mapleinline{inert}{2d}{3*v[4], Z^2}{\[\displaystyle 3\,v_{{4}},\,{Z}^{2}\]}
\end{maplelatex}
\mapleresult
\begin{maplelatex}
\mapleinline{inert}{2d}{-2*v[3]^2, Y*Z}{\[\displaystyle -2\,{v_{{3}}}^{2},\,YZ\]}
\end{maplelatex}
\mapleresult
\begin{maplelatex}
\mapleinline{inert}{2d}{12*v[3]^3+81*v[4]^2, Y^3}{\[\displaystyle 12\,{v_{{3}}}^{3}+81\,{v_{{4}}}^{2},\,{Y}^{3}\]}
\end{maplelatex}
\end{maplegroup}

\halfline

\noindent The only interesting relation is given by the last leading coefficient since:

\halfline

\begin{maplegroup}
\begin{mapleinput}
\mapleinline{active}{1d}{v[4] := 0 :}{\[\]}
\end{mapleinput}
\end{maplegroup}
\begin{maplegroup}
\begin{mapleinput}
\mapleinline{active}{1d}{print(mu=MilnorNumber(FL,\{Y,Z\}),tau=TyurinaNumber(FL,\{Y,Z\}))}
{\[{\mu={\it MilnorNumber} \left( {\it FL}, \left\{ Y,Z \right\}  \right) ,\tau={\it TyurinaNumber}\\
\mbox{} \left( {\it FL}, \left\{ Y,Z \right\}  \right) }\]}
\end{mapleinput}
\mapleresult
\begin{maplelatex}
\mapleinline{inert}{2d}{mu = 4, tau = 4}{\[\displaystyle \mu=4,\,\tau=4\]}
\end{maplelatex}
\end{maplegroup}
\begin{maplegroup}
\begin{mapleinput}
\mapleinline{active}{1d}{unassign('v[4]') : v[3]:=0 :}{\[\]}
\end{mapleinput}
\begin{mapleinput}
\mapleinline{active}{1d}{print(mu=MilnorNumber(FL,\{Y,Z\}),tau=TyurinaNumber(FL,\{Y,Z\}))}
{\[{\mu={\it MilnorNumber} \left( {\it FL}, \left\{ Y,Z \right\}  \right) ,\tau={\it TyurinaNumber}\\
\mbox{} \left( {\it FL}, \left\{ Y,Z \right\}  \right) }\]}
\end{mapleinput}
\mapleresult
\begin{maplelatex}
\mapleinline{inert}{2d}{mu = 4, tau = 4}{\[\displaystyle \mu=4,\,\tau=4\]}
\end{maplelatex}
\end{maplegroup}
\begin{maplegroup}
\begin{mapleinput}
\mapleinline{active}{1d}{unassign('v[3]') : v[3]:= -3*a^2 : v[4] := -2*a^3 : FL}{\[\]}
\end{mapleinput}
\end{maplegroup}
\begin{maplegroup}
\mapleresult
\begin{maplelatex}
\mapleinline{inert}{2d}{}{\[\displaystyle {Y}^{3}+{Z}^{4}-3{a}^{2}Y{Z}^{2}-2{a}^{3}{Z}^{3}\]}
\end{maplelatex}
\end{maplegroup}
\begin{maplegroup}
\begin{mapleinput}
\mapleinline{active}{2d}{print(mu = MilnorNumber(FL, {Y, Z}), tau = TyurinaNumber(FL, {Y, Z}))}{\[{\mu={\it MilnorNumber} \left( {\it FL}, \left\{ Y,Z \right\}  \right) ,\tau={\it TyurinaNumber}\\
\mbox{} \left( {\it FL}, \left\{ Y,Z \right\}  \right) }\]}
\end{mapleinput}
\mapleresult
\begin{maplelatex}
\mapleinline{inert}{2d}{mu = 5, tau = 5}{\[\displaystyle \mu=5,\,\tau=5\]}
\end{maplelatex}
\end{maplegroup}

\halfline

\noindent Therefore if
\begin{equation}\label{D5 in E6}
    \mathcal{V}:=\left\{4v_3^3+27v_4^2=0\right\}
\end{equation}
then Theorem \ref{classificazione}(3) allows to conclude that
\begin{itemize}
  \item \emph{$0\in f_{\Lambda}^{-1}(0)$ is a singularity of type $D_5$ for $\Lambda$
  generic in $\mathcal{V}\cap\mathcal{V}_0^2$}.
\end{itemize}
Let us now consider $\Lambda\in\widetilde{\mathcal{W}}_2^4$. Then

\halfline

\begin{maplegroup}
\begin{mapleinput}
\mapleinline{active}{1d}{F3;}{\[\]}
\end{mapleinput}
\mapleresult
\begin{maplelatex}
\mapleinline{inert}{2d}{Y^3+Z^4+(1/4)*(3*u^2+v[3])^2*(Y+u*Z)^2+v[3]*Y*Z^2+u*(v[3]+u^2)*Z^3}{\[\displaystyle {Y}^{3}+{Z}^{4}+1/4\, \left( 3\,{u}^{2}+v_{{3}} \right) ^{2} \left( Y+uZ \right) ^{2}+v_{{3}}Y{Z}^{2}+u \left( v_{{3}}+{u}^{2} \right) {Z}^{3}\]}
\end{maplelatex}
\end{maplegroup}
\begin{maplegroup}
\begin{mapleinput}
\mapleinline{active}{1d}{MGB4b:=MilnorGroebnerBasis(F3,\{Y,Z\}):}{\[{\it MGB4b}\, := \,{\it MilnorGroebnerBasis} \left( {\it F3}, \left\{ Y,Z \right\}  \right) \]}
\end{mapleinput}
\begin{mapleinput}
\mapleinline{active}{1d}{for i to nops(MGB4b) do
LeadingTerm(MGB4b[i],tdeg_min(Y,Z))end do;}{\[\]}
\end{mapleinput}
\mapleresult
\begin{maplelatex}
\mapleinline{inert}{2d}{(1/2)*(3*u^2+v[3])^2*u, Z}{\[\displaystyle 1/2\, \left( 3\,{u}^{2}+v_{{3}} \right) ^{2}u,\,Z\]}
\end{maplelatex}
\mapleresult
\begin{maplelatex}
\mapleinline{inert}{2d}{(1/2)*(3*u^2+v[3])^2*u^2, Z}{\[\displaystyle 1/2\, \left( 3\,{u}^{2}+v_{{3}} \right) ^{2}{u}^{2},\,Z\]}
\end{maplelatex}
\mapleresult
\begin{maplelatex}
\mapleinline{inert}{2d}{}{\[\displaystyle {\frac {12709329141645}{256}}\,{u}^{2}{v_{{3}}}^{2}+{\frac {38127987424935}{256}}\,{u}^{6}\\
\mbox{}+{\frac {38127987424935}{256}}\,{u}^{4}v_{{3}}+\]}
\mapleinline{inert}{2d}{}{\[\displaystyle\quad {\frac {1412147682405}{256}}\,{v_{{3}}}^{3},\,{Y}^{4}\]}
\end{maplelatex}
\end{maplegroup}
\begin{maplegroup}
\begin{mapleinput}
\mapleinline{active}{1d}{CC:=[seq(LeadingCoefficient(MGB4b[i],tdeg_min(Y,Z)),
i =1..nops(MGB4b))]:}{\[{\it CC}\, := \,[{\it seq} \left( {\it LeadingCoefficient} \left( {\it MGB4b}\\
\mbox{}_{{i}},{\it tdeg\_min} \left( Y,Z \right)  \right) ,i={1\ldots {\it nops} \left( {\it MGB4b}\\
\mbox{} \right) } \right) ]\]}
\end{mapleinput}
\begin{mapleinput}
\mapleinline{active}{1d}{factor(CC[3]);}{\[\]}
\end{mapleinput}
\mapleresult
\begin{maplelatex}
\mapleinline{inert}{2d}{(1412147682405/256)*(3*u^2+v[3])^3}{\[\displaystyle {\frac {1412147682405}{256}}\, \left( 3\,{u}^{2}+v_{{3}} \right) ^{3}\]}
\end{maplelatex}
\end{maplegroup}

\halfline
\noindent This last factorization clarifies that the third
leading coefficient does not give any further relations w.r.t.
those given by the previous leading coefficients, namely
$v_2(3v_2+v_1v_3)=0$. Observe that:

\halfline

\begin{maplegroup}
\begin{mapleinput}
\mapleinline{active}{1d}{u := 0 : F3}{\[\]}
\end{mapleinput}
\mapleresult
\begin{maplelatex}
\mapleinline{inert}{2d}{Y^3+Z^4+(1/4)*v[3]^2*Y^2+v[3]*Y*Z^2}{\[\displaystyle {Y}^{3}+{Z}^{4}+1/4\,{v_{{3}}}^{2}{Y}^{2}+v_{{3}}Y{Z}^{2}\]}
\end{maplelatex}
\end{maplegroup}
\begin{maplegroup}
\begin{mapleinput}
\mapleinline{active}{1d}{print(mu=MilnorNumber(F3,\{Y,Z\}),tau=TyurinaNumber(F3,\{Y,Z\}))}{\[{\mu={\it MilnorNumber} \left( {\it F3}, \left\{ Y,Z \right\}  \right) ,\tau={\it TyurinaNumber}\\
\mbox{} \left( {\it F3}, \left\{ Y,Z \right\}  \right) }\]}
\end{mapleinput}
\mapleresult
\begin{maplelatex}
\mapleinline{inert}{2d}{mu = 5, tau = 5}{\[\displaystyle \mu=5,\,\tau=5\]}
\end{maplelatex}
\end{maplegroup}

\halfline

\noindent implying, by Theorem \ref{classificazione}(2), that
\begin{itemize}
  \item \emph{$0\in f_{\Lambda}^{-1}(0)$ is a singularity of type $A_5$ for a generic $\Lambda$ in
      $$
      \widetilde{\mathcal{W}}^4_2\cap\mathcal{V}_2=\mathcal{W}\cap\left(\bigcap_{k=0}^2\mathcal{V}_{2k}\right)
      $$
where}
\begin{equation}\label{A5 in E6}
    \mathcal{W}:=\{v_3^2-4v_1=0\}\ .
\end{equation}
\end{itemize}
On the other hand, notice that

\halfline

\begin{maplegroup}
\begin{mapleinput}
\mapleinline{active}{1d}{unassign('u') : v[3] := -3*u^2 : F3}{\[\]}
\end{mapleinput}
\mapleresult
\begin{maplelatex}
\mapleinline{inert}{2d}{Y^3+Z^4-3*u^2*Y*Z^2-2*u^3*Z^3}{\[\displaystyle {Y}^{3}+{Z}^{4}-3\,{u}^{2}Y{Z}^{2}-2\,{u}^{3}{Z}^{3}\]}
\end{maplelatex}
\end{maplegroup}

\halfline

\noindent This means that, setting $W_5:=\left\{v_1v_3+3v_2=0\right\}$, then $\widetilde{\mathcal{W}}_2^4\cap\mathcal{W}_5=\mathcal{V}\cap\mathcal{V}_0^2$ and we are reduced to consider the $D_5$ case previously analyzed.
Since
$(\mathcal{V}\cap\mathcal{V}_0^2)\cup\left(\mathcal{W}\cap\left(\bigcap_{k=0}^2\mathcal{V}_{2k}\right)\right)$
parameterizes 1-parameter deformations, the further and last step
is clearly the trivial deformation given by $0\in T^1$. This fact ends up
the inclusions diagram (\ref{E6-inclusioni}), giving a
stratification, by nested algebraic subsets, of the
Kuranishi space of a simple $E_6$ singular point and representing
the adjacency diagram (\ref{E6-adiacenza}). At last, to prove
(\ref{intersezioni in E6}), let us observe that
\[
    \widetilde{\mathcal{W}}_2^4\cap\mathcal{V}_0^2=\mathcal{V}\cap\mathcal{V}_0^2
\]
since the common solutions of equations
\begin{equation*}
    v_0=v_1=v_2=0\ ,\ 4v_1v_2=v_0^2\ ,\ v_1u^2=v_2\ ,\
    v_4=u(v_3+u^2)\ ,\ 4v_1=(v_3+3u^2)^2
\end{equation*}
are given by $(v_0,v_1,v_2,v_3,v_4)=(0,0,0,-3u^2,-2u^3)$. Moreover
\[
    \left(\mathcal{W}\cap\left(\bigcap_{k=0}^2\mathcal{V}_{2k}\right)\right)\cap\left(\mathcal{V}\cap\mathcal{V}_0^2\right)
    =\mathcal{W}\cap\left(\bigcap_{k=0}^2\mathcal{V}_{2k}\right)\cap\mathcal{V} =\{0\}
\]
since $0$ is the only common solution of equations $v_0=v_2=v_4=0$, $4v_1=v_3^2$, $4v_3^3+27v_4^2=0$.
\end{proof}

\subsection{Simple singularities of $E_7$ type.}

\begin{theorem}\label{E7}
Let $T^1$ be the Kuranishi space of a simple
$N$--dimensional singular point $0\in f^{-1}(0)$ with
\[
    f(x_1,\ldots,x_{N+1}) = \sum_{i=1}^{N-1} x_i^2\ +\
    x_N^3+ x_Nx_{N+1}^3
\]
The subset of $T^1$ parameterizing small deformations of $0\in f^{-1}(0)$ to a simple node is the union of $n$ hypersurfaces. Moreover, calling $\mathcal{L}$ any of those hypersurfaces, there exists a stratification of nested algebraic subsets giving rise to the following sequence of inclusions, c.i.p. squares and a hinged union of c.i.p. squares
\begin{equation}\label{E7-inclusioni}
    \xymatrix{\mathcal{L}&\ar@{_{(}->}[l]\ \mathcal{W}_2&\ar@{_{(}->}[l]\
                \mathcal{W}_2^3&
                \ar@{_{(}->}[l]\ \widetilde{\mathcal{W}}_2^4&\ar@{_{(}->}[l]\
                \widetilde{\mathcal{W}}_2^5\cup\widetilde{\mathcal{W}'}_2^5&\ar@{_{(}->}[l]\
                \widetilde{\mathcal{W}'}_2^6\\
                &&\ar@{_{(}->}[u]\ \mathcal{V}_0^2&\ar@{_{(}->}[u]\ar@{_{(}->}[l]\
                \mathcal{V}\cap\mathcal{V}_0^2&\ar@{_{(}->}[u]\ar@{_{(}->}[l]\
                \mathcal{V}'\cap\mathcal{V}_0^2&\\
                &&&\ar@{_{(}->}[ruu]\ar@{_{(}->}[u]\ \mathcal{V}\cap\mathcal{V}_0^3&
                \ar@{_{(}->}[u]\ar@{_{(}->}[l]\ar@{_{(}->}[ruu]\ \{0\}&}
\end{equation}
verifying the Arnol'd's adjacency diagram
\begin{equation}\label{E7-adiacenza}
    \xymatrix{A_1&\ar[l]A_2&\ar[l]A_3&\ar[l]A_4&\ar[l]A_5&\ar[l]A_6\\
                &&\ar[u]D_4&\ar[u]\ar[l]D_5&\ar[u]\ar[l]D_6&\\
                &&&\ar[ruu]\ar[u]E_6&\ar[ruu]\ar[u]\ar[l]E_7&}
\end{equation}
where
\begin{itemize}
    \item $\mathcal{L}$ is the hypersurface of $T^1$ defined by equation
(\ref{lamda-condizione:E_7}), keeping in mind
(\ref{E_7,z,y}),
    \item $\mathcal{V}_0^m:=\bigcap_{k=0}^m \mathcal{V}_k$
and $\mathcal{W}_2^m:=\bigcap_{k=2}^m \mathcal{W}_k$ where
$\mathcal{V}_k, \mathcal{W}_k$ are hypersurfaces of $\mathcal{L}$
defined by equations (\ref{vk-E7}), (\ref{A2-3 in E7}) ,(\ref{A4
in E7}) and (\ref{A5 in E7}),
    \item $\mathcal{V}$ and $\mathcal{V}'$ are a hypersurfaces and
    a codimension 2 complete intersection in $\mathcal{L}$ defined
    by equations (\ref{D5 in E7}) and (\ref{D6 in E7}),
    respectively,
    \item $\widetilde{\mathcal{W}}_2^k$ and $\widetilde{\mathcal{W}'}_2^k$ are Zariski closed subsets of
    $\mathcal{L}$ defined by (\ref{A4bis in E7}), (\ref{A5 in
    E7}), (\ref{A5bis in E7}) and (\ref{A6 in E7}).
\end{itemize}
In particular complete intersection properties in diagram (\ref{E7-inclusioni}) are summarized by the following relations:
\begin{eqnarray}\label{intersezioni in E7}
\nonumber
    &\widetilde{\mathcal{W}}_2^4\cap\mathcal{V}_0^2=\mathcal{V}\cap\mathcal{V}_0^2\
    ,\
    \widetilde{\mathcal{W}}_2^5\cap\left(\mathcal{V}\cap\mathcal{V}_0^2\right)=
    \mathcal{V}'\cap\mathcal{V}_0^2\ ,\ \widetilde{\mathcal{W}'}_2^5\cap
    \left(\mathcal{V}\cap\mathcal{V}_0^2\right)= \mathcal{V}\cap\mathcal{V}_0^3& \\
    &\left(\mathcal{V}'\cap\mathcal{V}_0^2\right)
    \cap\left(\mathcal{V}\cap\mathcal{V}_0^3\right)=\{0\}=
    \widetilde{\mathcal{W}'}_2^6\cap\left(\mathcal{V}'\cap\mathcal{V}_0^2\right)\ .&
\end{eqnarray}
\end{theorem}

\begin{proof} Following the outline \ref{outline}.

\vskip 5pt\noindent (1) By the Morse Splitting Lemma \ref{Morse}, our problem can be reduced to the case $N=1$ with $f(y,z)=  y^3
+ yz^3$. To get an explicit basis of the Kuranishi space $T^1$
type

\halfline

 \begin{maplegroup}
\begin{mapleinput}
\mapleinline{active}{1d}{F:= y\symbol{94}3 + y*z\symbol{94}3: }{}
\end{mapleinput}
\end{maplegroup}
\begin{maplegroup}
\begin{mapleinput}
\mapleinline{active}{1d}{TyB := TyurinaBasis(F);}{\[\]}
\end{mapleinput}
\mapleresult
\begin{maplelatex}
\mapleinline{inert}{2d}{TyB := [z^4, z^3, z^2, y*z, z, y,
1]}{\[\displaystyle {\it TyB}\, :=
\,[{z}^{4},{z}^{3},{z}^{2},yz,z,y,1]\]}
\end{maplelatex}
\end{maplegroup}

\halfline

\noindent Therefore $T^1 \cong \langle
1,y,z,yz,z^2,z^3,z^4\rangle_{\C}$ and given $\Lambda = (\lambda_0,
\lambda_1,\ldots,\lambda_6)\in T^1$, the associated deformation of
$U_0$ is
\begin{eqnarray*}
   U_{\Lambda} &=& \{f_{\Lambda}(y,z)=0\}\quad\text{where}  \\
   f_{\Lambda}(y,z) &:=& f(y,z) + \lambda_0 + \lambda_1 y
    + \lambda_2 z + \lambda_3 yz + \lambda_4 z^2 + \lambda_5 z^3 +
    \lambda_6 z^4\ .
\end{eqnarray*}
A solution of the jacobian system of partial derivates is then
given by a solution $p_{\Lambda}=(y_{\Lambda},z_{\Lambda})$ of the following polynomial system in $\C[\lambda][y,z]$
\begin{equation}\label{E_7,z,y}
    \left\{\begin{array}{c}
      3y^2 + z^3 + \lambda_1 + \lambda_3z \equiv 0 \\
      3yz^2 +\lambda_2 + \lambda_3y + 2 \lambda_4z +
      3\lambda_5 z^2 + 4\lambda_6z^3= 0\\
    \end{array}\right.
\end{equation}
giving precisely 7 critical points for $f_{\Lambda}$.

\vskip 5pt\noindent (2) Imposing that one of those critical points, say $p_{\Lambda}$, is actually a singular point of $U_{\Lambda}$ means to require that
\begin{eqnarray}\label{lamda-condizione:E_7}
    &p_{\Lambda}\in U_{\Lambda}\ \Longleftrightarrow \ \Lambda\in
    \mathcal{L}\subset T^1\quad\text{where}&\\
\nonumber
    &    \mathcal{L}:=\{y_{\Lambda}^3 + y_{\Lambda}z_{\Lambda}^3 + \lambda_0 + \lambda_1y_{\Lambda}
    + \lambda_2z_{\Lambda} + \lambda_3y_{\Lambda}z_{\Lambda} + \lambda_4z_{\Lambda}^2 +
    \lambda_5z_{\Lambda}^3 + \lambda_6z_{\Lambda}^4  =
    0\}&\ .
\end{eqnarray}
After translating $y\mapsto y+y_{\Lambda},z\mapsto z+z_{\Lambda}$,
we get

\halfline

\noindent $f_{\Lambda}(y+y_{\Lambda},z+z_{\Lambda})=$
\begin{eqnarray*}
    f(y,z) &+&
    \left(y_{\Lambda}^3 + y_{\Lambda}z_{\Lambda}^3 + \lambda_0 + \lambda_1y_{\Lambda}
    + \lambda_2z_{\Lambda} + \lambda_3y_{\Lambda}z_{\Lambda} + \lambda_4z_{\Lambda}^2 +
    \lambda_5z_{\Lambda}^3 + \lambda_6z_{\Lambda}^4\right) \\
    &+& \left(3y_{\Lambda}^2 + z_{\Lambda}^3 + \lambda_1 + \lambda_3z_{\Lambda}\right)y \\
    &+& \left(3y_{\Lambda}z_{\Lambda}^2 +\lambda_2 + \lambda_3y_{\Lambda} + 2 \lambda_4z_{\Lambda} +
      3\lambda_5 z_{\Lambda}^2 + 4\lambda_6z_{\Lambda}^3\right)z \\
    &+& \left(3z_{\Lambda}^2 + \lambda_3\right)y z + 3y_{\Lambda}\ y^2 +
    \left(3y_{\Lambda}z_{\Lambda} + \lambda_4 + 3\lambda_5z_{\Lambda} + 6\lambda_6z_{\Lambda}^2\right)z^2 \\
    &+& 3 z_{\Lambda}\ y z^2 + \left(\lambda_5 + 4\lambda_6z_{\Lambda} + y_{\Lambda}\right)
    z^3 + \lambda_6\ z^4
\end{eqnarray*}
as can be checked by setting $K:=y_{\Lambda},L:=z_{\Lambda}$ and
typing:

\halfline

\begin{maplegroup}
\begin{mapleinput}
\mapleinline{active}{1d}{T := nops(TyB):}{\[\]}
\end{mapleinput}
\end{maplegroup}
\begin{maplegroup}
\begin{mapleinput}
\mapleinline{active}{1d}{F[Lambda] := F+sum(lambda[i]*TyB[T-i], i= 0 .. T-1)}
{\[F\Lambda \, := \,F+\sum
_{i=0}^{T-1}\lambda_{{i}}{\it TyB}_{{T-i}}\]}
\end{mapleinput}
\mapleresult
\begin{maplelatex}
\mapleinline{inert}{2d}{`F&Lambda;` := y^3+y*z^3+lambda[0]+lambda[1]*y+lambda[2]*z+lambda[3]*y*z+lambda[4]*z^2+lambda[5]*z^3+lambda[6]*z^4}{\[\displaystyle F\Lambda \, := \,{y}^{3}+y{z}^{3}+\lambda_{{0}}+\lambda_{{1}}y+\lambda_{{2}}z+\lambda_{{3}}yz+\lambda_{{4}}{z}^{2}+\lambda_{{5}}{z}^{3}\\
\mbox{}+\lambda_{{6}}{z}^{4}\]}
\end{maplelatex}
\end{maplegroup}
\begin{maplegroup}
\begin{mapleinput}
\mapleinline{active}{1d}{z := Z+L : y := Y+K :}{\[\]}
\end{mapleinput}
\end{maplegroup}
\begin{maplegroup}
\begin{mapleinput}
\mapleinline{active}{1d}{F[Lambda] := collect(F[Lambda], [Y, Z],
'distributed');}{\[\]}
\end{mapleinput}
\mapleresult
\begin{maplelatex}
\mapleinline{inert}{2d}{}
{\[\displaystyle F\Lambda \, := \,{Y}^{3}+3\,K{Y}^{2}+Y{Z}^{3}+3\,LY{Z}^{2}+ \left( \lambda_{{3}}+3\,{L}^{2} \right) YZ+
\]}
\mapleinline{inert}{2d}{}
{\[\displaystyle\ \left( \lambda_{{3}}L+3\,{K}^{2}+\lambda_{{1}}+{L}^{3} \right) Y+\lambda_{{6}}{Z}^{4}+ \left( \lambda_{{5}}+K+4\,\lambda_{{6}}L \right) {Z}^{3}+
\]}
\mapleinline{inert}{2d}{} {\[\displaystyle\ \left(
\lambda_{{4}}+3\,KL+3\,\lambda_{{5}}L+6\,\lambda_{{6}}{L}^{2}
\right) {Z}^{2}+\]}
\mapleinline{inert}{2d}{}
{\[\displaystyle\ \left(3\,K{L}^{2}+\lambda_{{2}}+4\,\lambda_{{6}}{L}^{3}+2\,\lambda_{{4}}L+\lambda_{{3}}K+3\,
\lambda_{{5}}{L}^{2} \right) Z+\]}
\mapleinline{inert}{2d}{}
{\[\displaystyle\ {K}^{3}+\lambda_{{6}}{L}^{4}+\lambda_{{0}}+\lambda_{{4}}{L}^{2}+K{L}^{3}+\lambda_{{3}}KL\\
\mbox{}+\lambda_{{1}}K+\lambda_{{2}}L+\lambda_{{5}}{L}^{3}\]}
\end{maplelatex}
\end{maplegroup}

\halfline

\vskip 5pt\noindent (3) Define:
\begin{eqnarray}\label{vk-E7}
\nonumber
   &v_0=3z_{\Lambda}^2 + \lambda_3\quad,\quad v_1=3y_{\Lambda}\quad, \quad \nu_2=3y_{\Lambda}z_{\Lambda} + \lambda_4
  + 3\lambda_5z_{\Lambda} + 6\lambda_6z_{\Lambda}^2\ , &  \\
   &v_3=3z_{\Lambda}\quad,\quad v_4=\lambda_5 + 4\lambda_6z_{\Lambda} +y_{\Lambda}\quad,
   \quad v_5=\lambda_6\ , &
\end{eqnarray}
and $\mathcal{V}_k:=\{v_k=0\}$. Then the proof goes on exactly as
in the $E_6$ case until $D_4$ singularities: precisely setting
\begin{eqnarray}\label{A2-3 in E7}
  \mathcal{W}_2 &=& \left\{4v_1v_2-v_0^2=0\right\} \\
\nonumber
  \mathcal{W}_3 &=& \left\{v_1^3v_4^2-v_2(v_1v_3+v_2)^2=0\right\}
\end{eqnarray}
we have the inclusions's chain
$$T^1\supset\mathcal{L}\supset\mathcal{W}_2\supset\mathcal{W}_2^3
\supset\mathcal{V}_1\cap\mathcal{W}_2^3=\mathcal{V}_0^2$$ of
subsets parameterizing small deformations whose generic fibre is
either smooth or admits a singularity of type $A_1$, $A_2$, $A_3$
or $D_4$, respectively: this fact can be checked as follows
\footnote{In the following computation we do not employ the
default l.m.o. but \texttt{plex\_min(Z,Y)}, which is the l.m.o.
defined as the opposite of the pure lexicographic g.m.o. with
$Y<Z$. In fact this last l.m.o. turned out to be considerably
more efficient than \texttt{tdeg\_min(Y,Z)}.}

\halfline

\begin{maplegroup}
\begin{mapleinput}
\mapleinline{active}{1d}{SD := [Y*Z, Y^2, Z^2, Y*Z^2, Z^3, Z^4]:
}{\[\]}
\end{mapleinput}
\end{maplegroup}
\begin{maplegroup}
\begin{mapleinput}
\mapleinline{active}{1d}{s := nops(SD):}{\[\]}
\end{mapleinput}
\end{maplegroup}
\begin{maplegroup}
\begin{mapleinput}
\mapleinline{active}{2d}{FL := Y^3+Y*Z^3+sum(v[i]*SD[i+1], i = 0
.. s-1); }{\[\]}
\end{mapleinput}
\mapleresult
\begin{maplelatex}
\mapleinline{inert}{2d}{FL :=
Y^3+Y*Z^3+v[0]*Y*Z+v[1]*Y^2+v[2]*Z^2+v[3]*Y*Z^2+v[4]*Z^3+v[5]*Z^4}{\[\displaystyle
{\it FL}\, :=
\,{Y}^{3}+Y{Z}^{3}+v_{{0}}YZ+v_{{1}}{Y}^{2}+v_{{2}}{Z}^{2}+v_{{3}}Y{Z}^{2}+v_{{4}}{Z}^{3}+v_{{5}}{Z}^{4}\]}
\end{maplelatex}
\end{maplegroup}
\begin{maplegroup}
\begin{mapleinput}
\mapleinline{active}{2d}{}{\[\]}
\end{mapleinput}
\end{maplegroup}
\begin{maplegroup}
\begin{mapleinput}
\mapleinline{active}{1d}{MGB:=MilnorGroebnerBasis(FL,\{Y,Z\},plex_min(Z,Y))}
{\[{\it MGB}\, := \,{\it MilnorGroebnerBasis} \left( {\it FL}, \left\{ Y,Z \right\} ,{\it plex\_min}\\
\mbox{} \left( Z,Y \right)  \right) \]}
\end{mapleinput}
\mapleresult
\begin{maplelatex}
\mapleinline{inert}{2d}{}
{\[\displaystyle {\it MGB}\, := \, \left\{ 3\,{Y}^{2}+{Z}^{3}+v_{{0}}Z+2\,v_{{1}}Y+v_{{3}}{Z}^{2},6\,v_{{3}}{Y}^{2}Z+
\left( 4\,v_{{1}}v_{{2}}-{v_{{0}}}^{2} \right) Z +\right.\]}
\mapleinline{inert}{2d}{}
{\[\displaystyle\  \left( 8\,v_{{1}}v_{{5}}-v_{{0}} \right) {Z}^{3}+ \left( -v_{{0}}v_{{3}}+6\,v_{{1}}v_{{4}} \right)
{Z}^{2}+9\,{Y}^{2}{Z}^{2}+12\,Yv_{{5}}{Z}^{3}+ \left( 9\,v_{{4}}+6\,v_{{1}} \right) Y{Z}^{2}\]}
\mapleinline{inert}{2d}{}
{\[\displaystyle\ + \left( 6\,v_{{2}}+4\,v_{{1}}v_{{3}} \right) YZ\\
\mbox{},3\,Y{Z}^{2}+v_{{0}}Y+2\,v_{{2}}Z+2\,v_{{3}}YZ+3\,v_{{4}}{Z}^{2}+4\,v_{{5}}{Z}^{3}
\left.\right\} \]}
\end{maplelatex}
\end{maplegroup}
\begin{maplegroup}
\begin{mapleinput}
\mapleinline{active}{1d}{for i to nops(MGB) do LeadingTerm(MGB[i],plex_min(Z,Y))end do;}{\[\]}
\end{mapleinput}
\mapleresult
\begin{maplelatex}
\mapleinline{inert}{2d}{v[0], Y}{\[\displaystyle v_{{0}},\,Y\]}
\end{maplelatex}
\mapleresult
\begin{maplelatex}
\mapleinline{inert}{2d}{2*v[1], Y}{\[\displaystyle
2\,v_{{1}},\,Y\]}
\end{maplelatex}
\mapleresult
\begin{maplelatex}
\mapleinline{inert}{2d}{4*v[1]*v[2]-v[0]^2, Z}{\[\displaystyle
4\,v_{{1}}v_{{2}}-{v_{{0}}}^{2},\,Z\]}
\end{maplelatex}
\end{maplegroup}
\begin{maplegroup}
\begin{mapleinput}
\mapleinline{active}{1d}{v[0] := 0 :}{\[\]}
\end{mapleinput}
\end{maplegroup}
\begin{maplegroup}
\begin{mapleinput}
\mapleinline{active}{1d}{print(mu=MilnorNumber(FL,\{Y,Z\},plex_min(Z,Y));}{\[\]}
\end{mapleinput}
\mapleresult
\begin{maplelatex}
\mapleinline{inert}{2d}{}{\[\displaystyle \mu = 1\]}
\end{maplelatex}
\end{maplegroup}
\begin{maplegroup}
\begin{mapleinput}
\mapleinline{active}{1d}{unassign('v[0]') : v[2]:= 0 :}{\[\]}
\end{mapleinput}
\end{maplegroup}
\begin{maplegroup}
\begin{mapleinput}
\mapleinline{active}{1d}{print(mu=MilnorNumber(FL,\{Y,Z\},plex_min(Z,Y));}{\[\]}
\end{mapleinput}
\mapleresult
\begin{maplelatex}
\mapleinline{inert}{2d}{}
{\[\displaystyle \mu=1\]}
\end{maplelatex}
\end{maplegroup}
\begin{maplegroup}
\begin{mapleinput}
\mapleinline{active}{1d}{unassign('v[2]'):}{\[\]}
\end{mapleinput}
\end{maplegroup}
\begin{maplegroup}
\begin{mapleinput}
\mapleinline{active}{1d}{F2:=Y^3+Y*Z^3+(w[1]*Y+w[2]*Z)^2+sum(v[j]*SD[j+1],j=3..s-1)}
{\[{\it F2}\, := \,{Y}^{3}+Y{Z}^{3}+ \left( w_{{1}}Y+w_{{2}}Z
\right) ^{2}+\sum _{j=3}^{s-1}v_{{j}}{\it SD}_{{j+1}}\]}
\end{mapleinput}
\mapleresult
\begin{maplelatex}
\mapleinline{inert}{2d}{F2 :=
Y^3+Y*Z^3+(w[1]*Y+w[2]*Z)^2+v[3]*Y*Z^2+v[4]*Z^3+v[5]*Z^4}{\[\displaystyle
{\it F2}\, := \,{Y}^{3}+Y{Z}^{3}+ \left( w_{{1}}Y+w_{{2}}Z \right)
^{2}+v_{{3}}Y{Z}^{2}+v_{{4}}{Z}^{3}+v_{{5}}{Z}^{4}\]}
\end{maplelatex}
\end{maplegroup}
\begin{maplegroup}
\begin{mapleinput}
\mapleinline{active}{1d}{print(mu=MilnorNumber(F2,\{Y,Z\},plex_min(Z,Y)))}{\[\mu={\it MilnorNumber} \left( {\it F2}, \left\{ Y,Z \right\} ,{\it plex\_min}\\
\mbox{} \left( Z,Y \right)  \right) \]}
\end{mapleinput}
\mapleresult
\begin{maplelatex}
\mapleinline{inert}{2d}{mu = 2}{\[\displaystyle \mu=2\]}
\end{maplelatex}
\end{maplegroup}
\begin{maplegroup}
\begin{mapleinput}
\mapleinline{active}{1d}{MGB2:=MilnorGroebnerBasis(F2,\{Y,Z\},plex_min(Z,Y));}{\[\]}
\end{mapleinput}
\mapleresult
\begin{maplelatex}
\mapleinline{inert}{2d}{}
{\[\displaystyle {\it MGB2}\, := \, \left\{  \left( 2\,w_{{1}}w_{{2}}v_{{3}}+12\,{w_{{2}}}^{2}v_{{5}}
\right) {Z}^{4}+ \right.\]}
\mapleinline{inert}{2d}{}
{\[\displaystyle\ \left( 2\,w_{{1}}{v_{{3}}}^{2}w_{{2}}+2\,{w_{{2}}}^{2}{w_{{1}}}^{2}+9\,{w_{{2}}}^{2}v_{{4}}-8\,
w_{{2}}{w_{{1}}}^{3}v_{{5}} \right) {Z}^{3}+ \]}
\mapleinline{inert}{2d}{}
{\[\displaystyle\ \left( 6\,{w_{{2}}}^{2}{w_{{1}}}^{2}v_{{3}}+6\,{w_{{2}}}^{4}-6\,w_{{2}}{w_{{1}}}^{3}v_{{4}} \right) {Z}^{2}-9\,w_{{2}}w_{{1}}{Y}^{2}{Z}^{2}\\
\mbox{}+ \left( 9\,{w_{{2}}}^{2}-12\,w_{{2}}w_{{1}}v_{{5}} \right) Y{Z}^{3}\]}
\mapleinline{inert}{2d}{}
{\[\displaystyle\ + \left( -9\,w_{{2}}w_{{1}}v_{{4}}-6\,w_{{2}}{w_{{1}}}^{3}+6\,{w_{{2}}}^{2}v_{{3}} \right) Y{Z}^{2},
3\,Y{Z}^{2}+2\, \left( w_{{1}}Y+w_{{2}}Z \right) w_{{2}}+2\,v_{{3}}YZ\]}
\mapleinline{inert}{2d}{}
{\[\displaystyle\ +3\,v_{{4}}{Z}^{2}+4\,v_{{5}}{Z}^{3},3\,{Y}^{2}+{Z}^{3}+2\, \left( w_{{1}}Y+w_{{2}}Z \right) w_{{1}}\\
\mbox{}+v_{{3}}{Z}^{2} \left.\right\} \]}
\end{maplelatex}
\end{maplegroup}
\begin{maplegroup}
\begin{mapleinput}
\mapleinline{active}{1d}{for i to nops(MGB2) do LeadingTerm(MGB2[i],plex_min(Z,Y))end do;}{\[\]}
\end{mapleinput}
\mapleresult
\begin{maplelatex}
\mapleinline{inert}{2d}{6*w[2]^2*w[1]^2*v[3]+6*w[2]^4-6*w[2]*w[1]^3*v[4],
Z^2}{\[\displaystyle
6\,{w_{{2}}}^{2}{w_{{1}}}^{2}v_{{3}}+6\,{w_{{2}}}^{4}-6\,w_{{2}}{w_{{1}}}^{3}v_{{4}},\,{Z}^{2}\]}
\end{maplelatex}
\mapleresult
\begin{maplelatex}
\mapleinline{inert}{2d}{2*w[2]*w[1], Y}{\[\displaystyle
2\,w_{{2}}w_{{1}},\,Y\]}
\end{maplelatex}
\mapleresult
\begin{maplelatex}
\mapleinline{inert}{2d}{2*w[1]^2, Y}{\[\displaystyle
2\,{w_{{1}}}^{2},\,Y\]}
\end{maplelatex}
\end{maplegroup}
\begin{maplegroup}
\begin{mapleinput}
\mapleinline{active}{2d}{w[1] := 0; -1}{\[\]}
\end{mapleinput}
\end{maplegroup}
\begin{maplegroup}
\begin{mapleinput}
\mapleinline{active}{1d}{print(mu=MilnorNumber(F2,\{Z,Y\},plex_min(Z,Y)),
tau=TYURINANumber(F2,\{Z,Y\},[plex_min(Z,Y),tdeg(Z,Y)]));}{\[\]}
\end{mapleinput}
\mapleresult
\begin{maplelatex}
\mapleinline{inert}{2d}{mu = 2, tau = 2}{\[\displaystyle
\mu=2,\,\tau=2\]}
\end{maplelatex}
\end{maplegroup}
\begin{maplegroup}
\begin{mapleinput}
\mapleinline{active}{1d}{unassign('w[1]') : w[2]:=0 :}{\[\]}
\end{mapleinput}
\end{maplegroup}
\begin{maplegroup}
\begin{mapleinput}
\mapleinline{active}{1d}{print(mu=MilnorNumber(F2,\{Y,Z\},plex_min(Z,Y)),
tau=TyurinaNumber(F2,\{Y,Z\},plex_min(Z,Y)))}
{\[{\mu={\it MilnorNumber} \left( {\it F2}, \left\{ Y,Z \right\} ,{\it plex\_min}\\
\mbox{} \left( Z,Y \right)  \right) ,\tau={\it TyurinaNumber} \left( {\it F2}, \left\{ Y,Z \right\} ,{\it plex\_min}\\
\mbox{} \left( Z,Y \right)  \right) }\]}
\end{mapleinput}
\mapleresult
\begin{maplelatex}
\mapleinline{inert}{2d}{mu = 2, tau = 2}{\[\displaystyle
\mu=2,\,\tau=2\]}
\end{maplelatex}
\end{maplegroup}
\begin{maplegroup}
\begin{mapleinput}
\mapleinline{active}{1d}{unassign('w[2]') : w[2] := u*w[1] : v[4]:= u*(v[3]+u^2) :}{\[\]}
\end{mapleinput}
\end{maplegroup}
\begin{maplegroup}
\begin{mapleinput}
\mapleinline{active}{1d}{F2 ;}{\[{\it F2}\]}
\end{mapleinput}
\mapleresult
\begin{maplelatex}
\mapleinline{inert}{2d}{Y^3+Y*Z^3+(w[1]*Y+u*w[1]*Z)^2+v[3]*Y*Z^2+u*(v[3]+u^2)*Z^3+v[5]*Z^4}{\[\displaystyle
{Y}^{3}+Y{Z}^{3}+ \left( w_{{1}}Y+uw_{{1}}Z \right)
^{2}+v_{{3}}Y{Z}^{2}+u \left( v_{{3}}+{u}^{2} \right)
{Z}^{3}+v_{{5}}{Z}^{4}\]}
\end{maplelatex}
\end{maplegroup}
\begin{maplegroup}
\begin{mapleinput}
\mapleinline{active}{1d}{F3:=Y^3+Y*Z^3+v[1]*(Y+u*Z)^2+v[3]*Y*Z^2+u*(v[3]+u^2)*Z^3
+v[5]*Z^4}
{\[{\it F3}\, := \,{Y}^{3}+Y{Z}^{3}+v_{{1}} \left( Y+uZ \right)
^{2}+v_{{3}}Y{Z}^{2}+u \left( v_{{3}}+{u}^{2} \right)
{Z}^{3}+v_{{5}}{Z}^{4}\]}
\end{mapleinput}
\mapleresult
\begin{maplelatex}
\mapleinline{inert}{2d}{F3 :=
Y^3+Y*Z^3+v[1]*(Y+u*Z)^2+v[3]*Y*Z^2+u*(v[3]+u^2)*Z^3+v[5]*Z^4}{\[\displaystyle
{\it F3}\, := \,{Y}^{3}+Y{Z}^{3}+v_{{1}} \left( Y+uZ \right)
^{2}+v_{{3}}Y{Z}^{2}+u \left( v_{{3}}+{u}^{2} \right)
{Z}^{3}+v_{{5}}{Z}^{4}\]}
\end{maplelatex}
\end{maplegroup}
\begin{maplegroup}
\begin{mapleinput}
\mapleinline{active}{1d}{print(mu=MilnorNumber(F3,\{Y,Z\},plex_min(Z,Y)),
tau=TYURINANumber(F3,\{Y,Z\},[plex_min(Z,Y),tdeg(Z,Y)]))}
{\[{\mu={\it MilnorNumber} \left( {\it F3}, \left\{ Y,Z \right\} ,{\it plex\_min}\\
\mbox{} \left( Z,Y \right)  \right) ,\tau={\it TYURINANumber} \left( {\it F3}, \left\{ Y,Z \right\} ,[{\it plex\_min}\\
\mbox{} \left( Z,Y \right) ,{\it tdeg}\\
\mbox{} \left( Z,Y \right) ] \right) }\]}
\end{mapleinput}
\mapleresult
\begin{maplelatex}
\mapleinline{inert}{2d}{mu = 3, tau = 3}{\[\displaystyle
\mu=3,\,\tau=3\]}
\end{maplelatex}
\end{maplegroup}
\begin{maplegroup}
\begin{mapleinput}
\mapleinline{active}{1d}{MGB3:=MilnorGroebnerBasis(F3,\{Y,Z\},plex_min(Z,Y));}{\[\]}
\end{mapleinput}
\mapleresult
\begin{maplelatex}
\mapleinline{inert}{2d}{} {\[\displaystyle{\it MGB3}\, :=
\,\left\{ \left(
-27\,{u}^{6}-9\,{u}^{2}{v_{{3}}}^{2}-10\,v_{{1}}uv_{{3}}-24\,v_{{1}}{u}^{2}v_{{5}}-36\,
{u}^{4}v_{{3}} \right) {Z}^{4}+\right.\]}
\mapleinline{inert}{2d}{} {\[\displaystyle \  \left(
-18\,v_{{1}}{u}^{2}+24\,v_{{1}}uv_{{5}}-6\,u{v_{{3}}}^{2}-18\,{u}^{3}v_{{3}}
\right) Y{Z}^{3}+ \]} \mapleinline{inert}{2d}{}
{\[\displaystyle\  \left( -24\,v_{{1}}v_{{3}}{u}^{3}+16\,u{v_{{1}}}^{2}v_{{5}}-36\,{u}^{5}v_{{1}}-4\,v_{{1}}u{v_{{3}}}^{2}\\
\mbox{}-16\,{v_{{1}}}^{2}{u}^{2} \right) {Z}^{3}+\]}
\mapleinline{inert}{2d}{}
{\[\displaystyle\  \left( -9\,uv_{{3}}-27\,{u}^{3} \right) Y{Z}^{4}\\
\mbox{}+ \left(
-36\,{u}^{3}v_{{5}}-6\,v_{{1}}u-12\,uv_{{3}}v_{{5}} \right)
{Z}^{5},\]}
\mapleinline{inert}{2d}{}
{\[\displaystyle 3\,Y{Z}^{2}+2\,v_{{1}} \left( Y+uZ \right) u+2\,v_{{3}}YZ+3\,u \left( v_{{3}}+{u}^{2} \right) {Z}^{2}+4\,v_{{5}}{Z}^{3}\\
\mbox{} , \]} \mapleinline{inert}{2d}{} {\[\displaystyle
3\,{Y}^{2}+{Z}^{3}+2\,v_{{1}} \left( Y+uZ \right) +v_{{3}}{Z}^{2}
\left.\right\}\]}
\end{maplelatex}
\end{maplegroup}
\begin{maplegroup}
\begin{mapleinput}
\mapleinline{active}{1d}{for i to nops(MGB3) do LeadingTerm(MGB3[i],plex_min(Z,Y))end do; }{\[\]}
\end{mapleinput}
\mapleresult
\begin{maplelatex}
\mapleinline{inert}{2d}{-24*v[1]*v[3]*u^3+16*u*v[1]^2*v[5]-36*u^5*v[1]-4*v[1]*u*v[3]^2-16*v[1]^2*u^2, Z^3}{\[\displaystyle -24\,v_{{1}}v_{{3}}{u}^{3}+16\,u{v_{{1}}}^{2}v_{{5}}-36\,{u}^{5}v_{{1}}-4\,v_{{1}}u{v_{{3}}}^{2}\\
\mbox{}-16\,{v_{{1}}}^{2}{u}^{2},\,{Z}^{3}\]}
\end{maplelatex}
\mapleresult
\begin{maplelatex}
\mapleinline{inert}{2d}{2*v[1]*u, Y}{\[\displaystyle
2\,v_{{1}}u,\,Y\]}
\end{maplelatex}
\mapleresult
\begin{maplelatex}
\mapleinline{inert}{2d}{2*v[1], Y}{\[\displaystyle
2\,v_{{1}},\,Y\]}
\end{maplelatex}
\end{maplegroup}
\begin{maplegroup}
\begin{mapleinput}
\mapleinline{active}{1d}{u := 0 :}{\[\]}
\end{mapleinput}
\end{maplegroup}
\begin{maplegroup}
\begin{mapleinput}
\mapleinline{active}{1d}{print(mu=MilnorNumber(F3,\{Z,Y\},plex_min(Z,Y)),
tau = TyurinaNumber(F3,\{Z,Y\},plex_min(Z,Y)))}{\[{\mu={\it MilnorNumber} \left( {\it F3}, \left\{ Z,Y \right\} ,{\it plex\_min}\\
\mbox{} \left( Z,Y \right)  \right) ,\tau={\it TyurinaNumber} \left( {\it F3}, \left\{ Z,Y \right\} ,{\it plex\_min}\\
\mbox{} \left( Z,Y \right)  \right) }\]}
\end{mapleinput}
\mapleresult
\begin{maplelatex}
\mapleinline{inert}{2d}{mu = 3, tau = 3}{\[\displaystyle
\mu=3,\,\tau=3\]}
\end{maplelatex}
\end{maplegroup}
\begin{maplegroup}
\begin{mapleinput}
\mapleinline{active}{1d}{unassign('u') : v[1] := 0 :}{\[\]}
\end{mapleinput}
\end{maplegroup}
\begin{maplegroup}
\begin{mapleinput}
\mapleinline{active}{1d}{print(mu=MilnorNumber(F3,\{Y,Z\},plex_min(Z,Y)),
tau=TyurinaNumber(F3,\{Y,Z\},plex_min(Z,Y)))}
{\[{\mu={\it MilnorNumber} \left( {\it F3}, \left\{ Y,Z \right\} ,{\it plex\_min}\\
\mbox{} \left( Z,Y \right)  \right) ,\tau={\it TyurinaNumber} \left( {\it F3}, \left\{ Y,Z \right\} ,{\it plex\_min}\\
\mbox{} \left( Z,Y \right)  \right) }\]}
\end{mapleinput}
\mapleresult
\begin{maplelatex}
\mapleinline{inert}{2d}{mu = 4, tau = 4}{\[\displaystyle
\mu=4,\,\tau=4\]}
\end{maplelatex}
\end{maplegroup}

\halfline

\noindent Also the singularities' specialization in fibers
parameterized by $\mathcal{V}_0^2$ proceeds as in the case of
$E_6$ singularities. Precisely:

\halfline

\begin{maplegroup}
\begin{mapleinput}
\mapleinline{active}{1d}{v[0] := 0 : v[1] := 0 : v[2] := 0 : FL}{\[\]}
\end{mapleinput}
\mapleresult
\begin{maplelatex}
\mapleinline{inert}{2d}{Y^3+Y*Z^3+v[3]*Y*Z^2+v[4]*Z^3+v[5]*Z^4}{\[\displaystyle
{Y}^{3}+Y{Z}^{3}+v_{{3}}Y{Z}^{2}+v_{{4}}{Z}^{3}+v_{{5}}{Z}^{4}\]}
\end{maplelatex}
\end{maplegroup}
\begin{maplegroup}
\begin{mapleinput}
\mapleinline{active}{1d}{MGB4 := MilnorGroebnerBasis(FL,\{Z,Y\},plex_min(Z,Y))}
{\[{\it MGB4}\, := \,{\it MilnorGroebnerBasis} \left( {\it FL}, \left\{ Z,Y \right\} ,{\it plex\_min}\\
\mbox{} \left( Z,Y \right)  \right) \]}
\end{mapleinput}
\mapleresult
\begin{maplelatex}
\mapleinline{inert}{2d}{} {\[\displaystyle {\it MGB4}\, := \,
\left\{ 3\,{Y}^{2}+{Z}^{3}+v_{{3}}{Z}^{2}, \left(
10\,{v_{{3}}}^{2}+36\, v_{{4}}v_{{5}} \right) {Z}^{4}+ \left(
-24\,v_{{3}}v_{{5}}+27\,v_{{4}} \right) Y{Z}^{3}\right.\]}
\mapleinline{inert}{2d}{}
{\[\displaystyle\ + \left( 4\,{v_{{3}}}^{3}+27\,{v_{{4}}}^{2} \right) {Z}^{3}\\
\mbox{}+6\,v_{{3}}{Z}^{5},3\,Y{Z}^{2}+2\,v_{{3}}YZ+3\,v_{{4}}{Z}^{2}+4\,v_{{5}}{Z}^{3}
\left.\right\} \]}
\end{maplelatex}
\end{maplegroup}
\begin{maplegroup}
\begin{mapleinput}
\mapleinline{active}{1d}{for i to nops(MGB4) do LeadingTerm(MGB4[i],plex_min(Z,Y))end do;}{\[\]}
\end{mapleinput}
\mapleresult
\begin{maplelatex}
\mapleinline{inert}{2d}{4*v[3]^3+27*v[4]^2, Z^3}{\[\displaystyle
4\,{v_{{3}}}^{3}+27\,{v_{{4}}}^{2},\,{Z}^{3}\]}
\end{maplelatex}
\mapleresult
\begin{maplelatex}
\mapleinline{inert}{2d}{3, Y^2}{\[\displaystyle 3,\,{Y}^{2}\]}
\end{maplelatex}
\mapleresult
\begin{maplelatex}
\mapleinline{inert}{2d}{2*v[3], Y*Z}{\[\displaystyle
2\,v_{{3}},\,YZ\]}
\end{maplelatex}
\end{maplegroup}
\begin{maplegroup}
\begin{mapleinput}
\mapleinline{active}{1d}{v[3] := 0 :}{\[\]}
\end{mapleinput}
\end{maplegroup}
\begin{maplegroup}
\begin{mapleinput}
\mapleinline{active}{1d}{print(mu=MilnorNumber(FL,\{Z,Y\},plex_min(Z,Y)),
tau=TyurinaNumber(FL,\{Z,Y\},plex_min(Z,Y)))}
{\[{\mu={\it MilnorNumber} \left( {\it FL}, \left\{ Z,Y \right\} ,{\it plex\_min}\\
\mbox{} \left( Z,Y \right)  \right) ,\tau={\it TyurinaNumber} \left( {\it FL}, \left\{ Z,Y \right\} ,{\it plex\_min}\\
\mbox{} \left( Z,Y \right)  \right) }\]}
\end{mapleinput}
\mapleresult
\begin{maplelatex}
\mapleinline{inert}{2d}{mu = 4, tau = 4}{\[\displaystyle
\mu=4,\,\tau=4\]}
\end{maplelatex}
\end{maplegroup}
\begin{maplegroup}
\begin{mapleinput}
\mapleinline{active}{1d}{unassign('v[3]') : v[3] := -3*a^2 : v[4] := -2*a^3 : FL}{\[\]}
\end{mapleinput}
\mapleresult
\begin{maplelatex}
\mapleinline{inert}{2d}{Y^3+Y*Z^3-3*a^2*Y*Z^2-2*a^3*Z^3+v[5]*Z^4}{\[\displaystyle
{Y}^{3}+Y{Z}^{3}-3\,{a}^{2}Y{Z}^{2}-2\,{a}^{3}{Z}^{3}+v_{{5}}{Z}^{4}\]}
\end{maplelatex}
\end{maplegroup}
\begin{maplegroup}
\begin{mapleinput}
\mapleinline{active}{1d}{print(mu=MilnorNumber(FL,\{Z,Y\},plex_min(Z,Y)),
tau=TyurinaNumber(FL,\{Z,Y\},plex_min(Z,Y)))}
{\[{\mu={\it MilnorNumber} \left( {\it FL}, \left\{ Z,Y \right\} ,{\it plex\_min}\\
\mbox{} \left( Z,Y \right)  \right) ,\tau={\it TyurinaNumber} \left( {\it FL}, \left\{ Z,Y \right\} ,{\it plex\_min}\\
\mbox{} \left( Z,Y \right)  \right) }\]}
\end{mapleinput}
\mapleresult
\begin{maplelatex}
\mapleinline{inert}{2d}{mu = 5, tau = 5}{\[\displaystyle
\mu=5,\,\tau=5\]}
\end{maplelatex}
\end{maplegroup}

\halfline

\noindent Therefore if
\begin{equation}\label{D5 in E7}
    \mathcal{V}:=\left\{4v_3^3+27v_4^2=0\right\}
\end{equation}
then Theorem \ref{classificazione}(3) allows to conclude that
\begin{itemize}
  \item \emph{$0\in f_{\Lambda}^{-1}(0)$ is a singularity of type $D_5$ for $\Lambda$
  generic in $\mathcal{V}\cap\mathcal{V}_0^2$}.
\end{itemize}
Moreover:

\halfline

\begin{maplegroup}
\begin{mapleinput}
\mapleinline{active}{1d}{MGB5 := MilnorGroebnerBasis(FL,\{Z,Y\},plex_min(Z,Y));}{\[\]}
\end{mapleinput}
\mapleresult
\begin{maplelatex}
\mapleinline{inert}{2d}{} {\[\displaystyle {\it MGB5}\, := \,
\left\{  \left( 48\,{a}^{4}-48\,{a}^{3}v_{{5}} \right) {Z}^{4}+
\left( 12\, v_{{5}}-9\,a \right) {Z}^{4}Y+ \left(
16\,{v_{{5}}}^{2}-12\,av_{{5}}-6\,{a}^{2} \right)
{Z}^{5},\right.\]} \mapleinline{inert}{2d}{}
{\[\displaystyle\ 3\,Y{Z}^{2}-6\,{a}^{2}YZ-6\,{a}^{3}{Z}^{2}+4\,v_{{5}}{Z}^{3}\\
\mbox{},3\,{Y}^{2}+{Z}^{3}-3\,{a}^{2}{Z}^{2} \left.\right\} \]}
\end{maplelatex}
\end{maplegroup}
\begin{maplegroup}
\begin{mapleinput}
\mapleinline{active}{1d}{for i to nops(MGB5) do LeadingTerm(MGB5[i],plex_min(Z,Y))end do;}{\[\]}
\end{mapleinput}
\mapleresult
\begin{maplelatex}
\mapleinline{inert}{2d}{48*a^4-48*a^3*v[5], Z^4}{\[\displaystyle
48\,{a}^{4}-48\,{a}^{3}v_{{5}},\,{Z}^{4}\]}
\end{maplelatex}
\mapleresult
\begin{maplelatex}
\mapleinline{inert}{2d}{-6*a^2, Y*Z}{\[\displaystyle
-6\,{a}^{2},\,YZ\]}
\end{maplelatex}
\mapleresult
\begin{maplelatex}
\mapleinline{inert}{2d}{3, Y^2}{\[\displaystyle 3,\,{Y}^{2}\]}
\end{maplelatex}
\end{maplegroup}
\begin{maplegroup}
\begin{mapleinput}
\mapleinline{active}{1d}{a := 0 : FL}{\[\]}
\end{mapleinput}
\mapleresult
\begin{maplelatex}
\mapleinline{inert}{2d}{Y^3+Y*Z^3+v[5]*Z^4}{\[\displaystyle
{Y}^{3}+Y{Z}^{3}+v_{{5}}{Z}^{4}\]}
\end{maplelatex}
\end{maplegroup}
\begin{maplegroup}
\begin{mapleinput}
\mapleinline{active}{1d}{print(mu=MilnorNumber(FL,\{Z,Y\},plex_min(Z,Y)),
tau=TyurinaNumber(FL,\{Z,Y\},plex_min(Z,Y)))}
{\[{\mu={\it MilnorNumber} \left( {\it FL}, \left\{ Z,Y \right\} ,{\it plex\_min}\\
\mbox{} \left( Z,Y \right)  \right) ,\tau={\it TyurinaNumber} \left( {\it FL}, \left\{ Z,Y \right\} ,{\it plex\_min}\\
\mbox{} \left( Z,Y \right)  \right) }\]}
\end{mapleinput}
\mapleresult
\begin{maplelatex}
\mapleinline{inert}{2d}{mu = 6, tau = 6}{\[\displaystyle
\mu=6,\,\tau=6\]}
\end{maplelatex}
\end{maplegroup}
\begin{maplegroup}
\begin{mapleinput}
\mapleinline{active}{1d}{unassign('a') : v[5] := a : FL}{\[\]}
\end{mapleinput}
\mapleresult
\begin{maplelatex}
\mapleinline{inert}{2d}{Y^3+Y*Z^3-3*a^2*Y*Z^2-2*a^3*Z^3+a*Z^4}{\[\displaystyle
{Y}^{3}+Y{Z}^{3}-3\,{a}^{2}Y{Z}^{2}-2\,{a}^{3}{Z}^{3}+a{Z}^{4}\]}
\end{maplelatex}
\end{maplegroup}
\begin{maplegroup}
\begin{mapleinput}
\mapleinline{active}{1d}{print(mu=MilnorNumber(FL,\{Z,Y\},plex_min(Z,Y)),
tau=TyurinaNumber(FL,\{Z,Y\},plex_min(Z,Y)))}
{\[{\mu={\it MilnorNumber} \left( {\it FL}, \left\{ Z,Y \right\} ,{\it plex\_min}\\
\mbox{} \left( Z,Y \right)  \right) ,\tau={\it TyurinaNumber} \left( {\it FL}, \left\{ Z,Y \right\} ,{\it plex\_min}\\
\mbox{} \left( Z,Y \right)  \right) }\]}
\end{mapleinput}
\mapleresult
\begin{maplelatex}
\mapleinline{inert}{2d}{mu = 6, tau = 6}{\[\displaystyle
\mu=6,\,\tau=6\]}
\end{maplelatex}
\end{maplegroup}

\halfline

\noindent Then Theorem \ref{classificazione}(4) gives that
\begin{itemize}
    \item \emph{$0\in f_{\Lambda}^{-1}(0)$ is a singularity of type $E_6$ for $\Lambda$
  generic in $\mathcal{V}\cap\mathcal{V}_0^3=\mathcal{V}_0^4$.}
\end{itemize}
On the other hand, by defining the following codimension 2 subset
of $\mathcal{L}$
\begin{equation}\label{D6 in E7}
    \mathcal{V}'=\left\{v_3+3v_5^2=v_4+2v_5^3=0\right\}\ ,
\end{equation}
Theorem \ref{classificazione}(3) gives that
\begin{itemize}
    \item \emph{$0\in f_{\Lambda}^{-1}(0)$ is a singularity of type $D_6$ for $\Lambda$
  generic in $\mathcal{V}'\cap\mathcal{V}_0^2$}
\end{itemize}
A further specialization here gives then the trivial deformation
since
$(\mathcal{V}\cap\mathcal{V}_0^3)\cap(\mathcal{V}'\cap\mathcal{V}_0^2)=\mathcal{V}'\cap\mathcal{V}_0^3=\{0\}$.

\noindent Let us the come back to consider $\mathcal{W}_2^3$ and
look at the first leading coefficient in $MGB3$, giving the relation
\[
    4v_1(v_5-u)-(v_3+3u^2)^2=0\ .
\]
Then set $v_5:=t^2+u$ and $v_3:=\pm 2w_1t -3u^2$ and type \footnote{The reader may check that choosing $v_3:=- 2w_1t -3u^2$ gives the same result.}

\halfline

\begin{maplegroup}
\begin{mapleinput}
\mapleinline{active}{1d}{v[3]:=2*t*w[1]-3*u^2 : v[5]:=t^2+u : F2}{\[\]}
\end{mapleinput}
\mapleresult
\begin{maplelatex}
\mapleinline{inert}{2d}{} {\[\displaystyle {Y}^{3}+Y{Z}^{3}+
\left( w_{{1}}Y+uw_{{1}}Z \right) ^{2}+ \]}
\mapleinline{inert}{2d}{} {\[\displaystyle \left(
2\,tw_{{1}}-3\,{u}^{2} \right) Y{Z}^{2}+
u \left( 2\,tw_{{1}}-2\,{u}^{2} \right) {Z}^{3}\\
\mbox{}+ \left( {t}^{2}+u \right) {Z}^{4}\]}
\end{maplelatex}
\end{maplegroup}
\begin{maplegroup}
\begin{mapleinput}
\mapleinline{active}{1d}{print(mu=MilnorNumber(F2,\{Y,Z\},plex_min(Z,Y)),
tau = TYURINANumber(F2,\{Y,Z\},[plex_min(Z,Y),tdeg(Z,Y)]))}
{\[{\mu={\it MilnorNumber} \left( {\it F2}, \left\{ Y,Z \right\} ,{\it plex\_min}\\
\mbox{} \left( Z,Y \right)  \right) ,\tau={\it TYURINANumber} \left( {\it F2}, \left\{ Y,Z \right\} ,[{\it plex\_min}\\
\mbox{} \left( Z,Y \right) ,{\it tdeg}\\
\mbox{} \left( Z,Y \right) ] \right) }\]}
\end{mapleinput}
\mapleresult
\begin{maplelatex}
\mapleinline{inert}{2d}{mu = 4, tau = 4}{\[\displaystyle
\mu=4,\,\tau=4\]}
\end{maplelatex}
\end{maplegroup}

\halfline

\noindent Define
\begin{equation}\label{A4 in E7}
    \mathcal{W}_4:=\left\{16v_1^5v_2-[(v_1v_2+3v_2)^2-4v_1^3v_5]^2=0\right\}
\end{equation}
whose equation is obtained by eliminating $u$ and $t$ from the following set of equations, parameterizing $\mathcal{W}_3^4\setminus\mathcal{V}_1^2$,
\begin{equation}\label{parametrizzazioneE7}
    v_2=u^2v_1\ ,\ v_4 = u(v_3+u^2)\ ,\ v_5=u+t^2\ , v_1=w_1^2\ ,\  v_3=\pm 2w_1t-3u^2\ .
\end{equation}
Observe that $\mathcal{V}_1^2\subseteq\mathcal{W}_3^4\
\Rightarrow\ \mathcal{V}_0^2\subseteq\mathcal{W}_2^4$ and define
the following codimension 1 Zariski closed subset \footnote{The same considerations explained by footnote \ref{nota} are still holding, here.} of
$\mathcal{W}_2^3$
\begin{equation}\label{A4bis in E7}
    \widetilde{\mathcal{W}}_2^4:=\overline{\mathcal{W}_2^4\setminus\mathcal{V}_0^2}\
    .
\end{equation}
Then Theorem \ref{classificazione}(2) implies that
\begin{itemize}
    \item \emph{$0\in f_{\Lambda}^{-1}(0)$ is a singularity of type $A_4$ for any
    generic $\Lambda\in\widetilde{\mathcal{W}}_2^4$}.
\end{itemize}
Moreover:

\halfline

\begin{maplegroup}
\begin{mapleinput}
\mapleinline{active}{1d}{MGB4b := MilnorGroebnerBasis(F2,\{Y,Z\},plex_min(Z,Y));}{\[\]}
\end{mapleinput}
\mapleresult
\begin{maplelatex}
\mapleinline{inert}{2d}{} {\[\displaystyle {\it MGB4b}\, := \,
\left\{ 3\,{Y}^{2}+{Z}^{3}+2\, \left( w_{{1}}Y+uw_{{1}}Z \right)
w_{{1}}+ \left( 2\,tw_{{1}}-3\,{u}^{2} \right) {Z}^{2},
\right.\]}\mapleinline{inert}{2d}{}
{\[\displaystyle 3\,Y{Z}^{2}+2\, \left( w_{{1}}Y+uw_{{1}}Z \right) uw_{{1}}\\
\mbox{}+2\, \left( 2\,tw_{{1}}-3\,{u}^{2} \right) YZ+3\,u \left(
2\,tw_{{1}}-2\,{u}^{2} \right) {Z}^{2}+ \]}
\mapleinline{inert}{2d}{} {\[\displaystyle 4\, \left( {t}^{2}+u
\right) {Z}^{3},\left( -9\,w_{{1}}-54\,ut \right) Y{Z}^{5}+ \left(
-9\,uw_{{1}}-12\,w_{{1}}{t}^{2}-72\,u{t}^{3}-54\,{u}^{2}t \right)
{Z}^{6}+ \]} \mapleinline{inert}{2d}{} {\[\displaystyle
\left(-20\,t{w_{{1}}}^{4}-80\,{w_{{1}}}^{3}u{t}^{2}-60\,{u}^{2}{t}^{3}{w_{{1}}}^{2}
\right) {Z}^{4}+ \]} \mapleinline{inert}{2d}{} {\[\displaystyle
\left(
-30\,t{w_{{1}}}^{2}-102\,w_{{1}}u{t}^{2}+90\,{u}^{3}t-72\,{u}^{2}{t}^{3}+15\,
w_{{1}}{u}^{2} \right) Y{Z}^{4}+\]} \mapleinline{inert}{2d}{}
{\[\displaystyle \left(
-24\,{w_{{1}}}^{2}{t}^{3}+90\,{u}^{4}t-38\,{w_{{1}}}^{2}ut+72\,{u}^{3}{t}^{3}+15\,
{u}^{3}w_{{1}}-60\,w_{{1}}{u}^{2}{t}^{2}-24\,w_{{1}}u{t}^{4}-6\,{w_{{1}}}^{3}
\right) {Z}^{5}\left.\right\}\]}
\end{maplelatex}
\end{maplegroup}
\begin{maplegroup}
\begin{mapleinput}
\mapleinline{active}{1d}{for i to nops(MGB4b) do LeadingTerm(MGB4b[i],plex_min(Z,Y))end do; }{\[\]}
\end{mapleinput}
\mapleresult
\begin{maplelatex}
\mapleinline{inert}{2d}{2*w[1]^2, Y}{\[\displaystyle
2\,{w_{{1}}}^{2},\,Y\]}
\end{maplelatex}
\mapleresult
\begin{maplelatex}
\mapleinline{inert}{2d}{2*u*w[1]^2, Y}{\[\displaystyle
2\,u{w_{{1}}}^{2},\,Y\]}
\end{maplelatex}
\mapleresult
\begin{maplelatex}
\mapleinline{inert}{2d}{-20*t*w[1]^4-80*w[1]^3*u*t^2-60*u^2*t^3*w[1]^2,
Z^4}{\[\displaystyle
-20\,t{w_{{1}}}^{4}-80\,{w_{{1}}}^{3}u{t}^{2}-60\,{u}^{2}{t}^{3}{w_{{1}}}^{2},\,{Z}^{4}\]}
\end{maplelatex}
\end{maplegroup}
\begin{maplegroup}
\begin{mapleinput}
\mapleinline{active}{1d}{u := 0 : F2}{\[\]}
\end{mapleinput}
\mapleresult
\begin{maplelatex}
\mapleinline{inert}{2d}{Y^3+Y*Z^3+w[1]^2*Y^2+2*t*w[1]*Y*Z^2+t^2*Z^4}{\[\displaystyle
{Y}^{3}+Y{Z}^{3}+{w_{{1}}}^{2}{Y}^{2}+2\,tw_{{1}}Y{Z}^{2}+{t}^{2}{Z}^{4}\]}
\end{maplelatex}
\end{maplegroup}
\begin{maplegroup}
\begin{mapleinput}
\mapleinline{active}{1d}{print(mu=MilnorNumber(F2,\{Z,Y\},plex_min(Z,Y)),
tau=TYURINANumber(F2,\{Z,Y\},[plex_min(Z,Y), tdeg(Z,Y)]))}{\[{\mu={\it MilnorNumber} \left( {\it F2}, \left\{ Z,Y \right\} ,{\it plex\_min}\\
\mbox{} \left( Z,Y \right)  \right) ,\tau={\it TYURINANumber} \left( {\it F2}, \left\{ Z,Y \right\} ,[{\it plex\_min}\\
\mbox{} \left( Z,Y \right) ,{\it tdeg}\\
\mbox{} \left( Z,Y \right) ] \right) }\]}
\end{mapleinput}
\mapleresult
\begin{maplelatex}
\mapleinline{inert}{2d}{mu = 4, tau = 4}{\[\displaystyle
\mu=4,\,\tau=4\]}
\end{maplelatex}
\end{maplegroup}
\begin{maplegroup}
\begin{mapleinput}
\mapleinline{active}{1d}{unassign('u') : w[1] := 0 : F2}{\[\]}
\end{mapleinput}
\mapleresult
\begin{maplelatex}
\mapleinline{inert}{2d}{Y^3+Y*Z^3-3*u^2*Y*Z^2-2*u^3*Z^3+(t^2+u)*Z^4}{\[\displaystyle
{Y}^{3}+Y{Z}^{3}-3\,{u}^{2}Y{Z}^{2}-2\,{u}^{3}{Z}^{3}+ \left(
{t}^{2}+u \right) {Z}^{4}\]}
\end{maplelatex}
\end{maplegroup}
\begin{maplegroup}
\begin{mapleinput}
\mapleinline{active}{1d}{print(mu=MilnorNumber(F2,\{Z,Y\},plex_min(Z,Y)),
tau=TyurinaNumber(F2,\{Z,Y\},plex_min(Z,Y)))}
{\[{\mu={\it MilnorNumber} \left( {\it F2}, \left\{ Z,Y \right\} ,{\it plex\_min}\\
\mbox{} \left( Z,Y \right)  \right) ,\tau={\it TyurinaNumber} \left( {\it F2}, \left\{ Z,Y \right\} ,{\it plex\_min}\\
\mbox{} \left( Z,Y \right)  \right) }\]}
\end{mapleinput}
\mapleresult
\begin{maplelatex}
\mapleinline{inert}{2d}{mu = 5, tau = 5}{\[\displaystyle
\mu=5,\,\tau=5\]}
\end{maplelatex}
\end{maplegroup}
\begin{maplegroup}
\begin{mapleinput}
\mapleinline{active}{1d}{unassign('w[1]') : t:=0 : F2}{\[\]}
\end{mapleinput}
\mapleresult
\begin{maplelatex}
\mapleinline{inert}{2d}{Y^3+Y*Z^3+(w[1]*Y+u*w[1]*Z)^2-3*u^2*Y*Z^2-2*u^3*Z^3+Z^4*u}{\[\displaystyle
{Y}^{3}+Y{Z}^{3}+ \left( w_{{1}}Y+uw_{{1}}Z \right)
^{2}-3\,{u}^{2}Y{Z}^{2}-2\,{u}^{3}{Z}^{3}+{Z}^{4}u\]}
\end{maplelatex}
\end{maplegroup}
\begin{maplegroup}
\begin{mapleinput}
\mapleinline{active}{1d}{print(mu=MilnorNumber(F2,\{Z,Y\},plex_min(Z,Y)),
tau=TYURINANumber(F2,\{Z,Y\},[plex_min(Z,Y),tdeg(Z,Y)]))}
{\[{\mu={\it MilnorNumber} \left( {\it F2}, \left\{ Z,Y \right\} ,{\it plex\_min}\\
\mbox{} \left( Z,Y \right)  \right) ,\tau={\it TYURINANumber} \left( {\it F2}, \left\{ Z,Y \right\} ,[{\it plex\_min}\\
\mbox{} \left( Z,Y \right) ,{\it tdeg}\\
\mbox{} \left( Z,Y \right) ] \right) }\]}
\end{mapleinput}
\mapleresult
\begin{maplelatex}
\mapleinline{inert}{2d}{mu = 5, tau = 5}{\[\displaystyle
\mu=5,\,\tau=5\]}
\end{maplelatex}
\end{maplegroup}

\halfline

\noindent Then
$\widetilde{\mathcal{W}}_2^4\cap\mathcal{V}_1=\mathcal{V}\cap\mathcal{V}_0^2$
reduces
to the already considered case of generic $D_5$
singularities. On the other hand setting $t=0$ means defining
\begin{equation}\label{A5 in E7}
    \mathcal{W}_5:=\left\{v_1v_5^2-v_2=0\right\}\quad ,\quad
    \widetilde{\mathcal{W}}_2^5:=\overline{\mathcal{W}_2^5\setminus\mathcal{V}_0^2}\
    ,
\end{equation}
and observing that, by Theorem \ref{classificazione}(2),
\begin{itemize}
    \item \emph{$0\in f_{\Lambda}^{-1}(0)$ is a singularity of type $A_5$
    for    generic
    $\Lambda\in\widetilde{\mathcal{W}}_2^4\cap\mathcal{W}_5=\widetilde{\mathcal{W}}_2^5$}.
\end{itemize}
Consider then

\halfline

\begin{maplegroup}
\begin{mapleinput}
\mapleinline{active}{1d}{MGB5b := MilnorGroebnerBasis(F2,\{Z,Y\},plex_min(Z,Y));}{\[\]}
\end{mapleinput}
\mapleresult
\begin{maplelatex}
\mapleinline{inert}{2d}{} {\[\displaystyle {\it MGB5b}\, := \,
\left\{3\,{Y}^{2}+{Z}^{3}+2\, \left( w_{{1}}Y+uw_{{1}}Z \right)
w_{{1}}-3\,{u}^{2}{Z}^{2} ,\right.\]} \mapleinline{inert}{2d}{}
{\[\displaystyle\
-12\,{w_{{1}}}^{4}{Z}^{5}+12\,{w_{{1}}}^{2}{Z}^{6}u-
18\,{w_{{1}}}^{2}{Z}^{5}Y,\]} \mapleinline{inert}{2d}{}
{\[\displaystyle \ 3\,Y{Z}^{2}+2\, \left( w_{{1}}Y+uw_{{1}}Z
\right) uw_{{1}}-6\,YZ{u}^{2}-6\,{u}^{3}{Z}^{2}+4\,{Z}^{3}u
\left.\right\}
\]}
\end{maplelatex}
\end{maplegroup}
\begin{maplegroup}
\begin{mapleinput}
\mapleinline{active}{1d}{for i to nops(MGB5b) do LeadingTerm(MGB5b[i],plex_min(Z,Y))end do; }{\[\]}
\end{mapleinput}
\mapleresult
\begin{maplelatex}
\mapleinline{inert}{2d}{2*w[1]^2, Y}{\[\displaystyle
2\,{w_{{1}}}^{2},\,Y\]}
\end{maplelatex}
\mapleresult
\begin{maplelatex}
\mapleinline{inert}{2d}{-12*w[1]^4, Z^5}{\[\displaystyle
-12\,{w_{{1}}}^{4},\,{Z}^{5}\]}
\end{maplelatex}
\mapleresult
\begin{maplelatex}
\mapleinline{inert}{2d}{2*u*w[1]^2, Y}{\[\displaystyle
2\,u{w_{{1}}}^{2},\,Y\]}
\end{maplelatex}
\end{maplegroup}
\begin{maplegroup}
\begin{mapleinput}
\mapleinline{active}{1d}{u := 0 : F2}{\[\]}
\end{mapleinput}
\mapleresult
\begin{maplelatex}
\mapleinline{inert}{2d}{Y^3+Y*Z^3+w[1]^2*Y^2}{\[\displaystyle
{Y}^{3}+Y{Z}^{3}+{w_{{1}}}^{2}{Y}^{2}\]}
\end{maplelatex}
\end{maplegroup}
\begin{maplegroup}
\begin{mapleinput}
\mapleinline{active}{1d}{print(mu=MilnorNumber(F2,\{Z,Y\},plex_min(Z,Y)),
tau=TyurinaNumber(F2,\{Z,Y\},plex_min(Z,Y)))}
{\[{\mu={\it MilnorNumber} \left( {\it F2}, \left\{ Z,Y \right\} ,{\it plex\_min}\\
\mbox{} \left( Z,Y \right)  \right) ,\tau={\it TyurinaNumber} \left( {\it F2}, \left\{ Z,Y \right\} ,{\it plex\_min}\\
\mbox{} \left( Z,Y \right)  \right) }\]}
\end{mapleinput}
\mapleresult
\begin{maplelatex}
\mapleinline{inert}{2d}{mu = 5, tau = 5}{\[\displaystyle
\mu=5,\,\tau=5\]}
\end{maplelatex}
\end{maplegroup}
\begin{maplegroup}
\begin{mapleinput}
\mapleinline{active}{1d}{unassign('u') : w[1]:=0 : F2}{\[\]}
\end{mapleinput}
\mapleresult
\begin{maplelatex}
\mapleinline{inert}{2d}{Y^3+Y*Z^3-3*u^2*Y*Z^2-2*u^3*Z^3+Z^4*u}{\[\displaystyle
{Y}^{3}+Y{Z}^{3}-3\,{u}^{2}Y{Z}^{2}-2\,{u}^{3}{Z}^{3}+{Z}^{4}u\]}
\end{maplelatex}
\end{maplegroup}
\begin{maplegroup}
\begin{mapleinput}
\mapleinline{active}{1d}{print(mu=MilnorNumber(F2,\{Z,Y\},plex_min(Z,Y)),
tau=TyurinaNumber(F2,\{Z,Y\},plex_min(Z,Y)))}
{\[{\mu={\it MilnorNumber} \left( {\it F2}, \left\{ Z,Y \right\} ,{\it plex\_min}\\
\mbox{} \left( Z,Y \right)  \right) ,\tau={\it TyurinaNumber} \left( {\it F2}, \left\{ Z,Y \right\} ,{\it plex\_min}\\
\mbox{} \left( Z,Y \right)  \right) }\]}
\end{mapleinput}
\mapleresult
\begin{maplelatex}
\mapleinline{inert}{2d}{mu = 6, tau = 6}{\[\displaystyle
\mu=6,\,\tau=6\]}
\end{maplelatex}
\end{maplegroup}

\halfline

\noindent which means that
$\mathcal{V}_1\cap\widetilde{\mathcal{W}}_2^5 =
\mathcal{V}'\cap\mathcal{V}_0^2$ reduces to the already considered
case of generic $D_6$ singularities.

\noindent Let us then come back to $\widetilde{\mathcal{W}}_2^4$
and consider the last leading coefficient in the standard basis
$MGB4b$ giving the further relations $w_1+ut=0$ and $w_1+3ut=0$,
not yet analyzed. Then

\halfline

\begin{maplegroup}
\begin{mapleinput}
\mapleinline{active}{1d}{F2}{\[{\it F2}\]}
\end{mapleinput}
\mapleresult
\begin{maplelatex}
\mapleinline{inert}{2d}{} {\[\displaystyle {Y}^{3}+Y{Z}^{3}+
\left( w_{{1}}Y+uw_{{1}}Z \right) ^{2}+ \left(
2\,tw_{{1}}-3\,{u}^{2} \right) Y{Z}^{2}\]}
\mapleinline{inert}{2d}{}
{\[\displaystyle\quad +u \left( 2\,tw_{{1}}-2\,{u}^{2} \right) {Z}^{3}\\
\mbox{}+ \left( {t}^{2}+u \right) {Z}^{4}\]}
\end{maplelatex}
\end{maplegroup}
\begin{maplegroup}
\begin{mapleinput}
\mapleinline{active}{1d}{w[1] := -u*t : F2}{\[\]}
\end{mapleinput}
\mapleresult
\begin{maplelatex}
\mapleinline{inert}{2d}{} {\[\displaystyle {Y}^{3}+Y{Z}^{3}+
\left( -tuY-{u}^{2}tZ \right) ^{2}+ \left( -2\,{t}^{2}u-3\,{u}^{2}
\right) Y{Z}^{2}\]}
\mapleinline{inert}{2d}{}
{\[\displaystyle\quad +u \left( -2\,{t}^{2}u-2\,{u}^{2} \right) {Z}^{3}\\
\mbox{}+ \left( {t}^{2}+u \right) {Z}^{4}\]}
\end{maplelatex}
\end{maplegroup}
\begin{maplegroup}
\begin{mapleinput}
\mapleinline{active}{1d}{print(mu=MilnorNumber(F2,\{Z,Y\},plex_min(Z,Y)),
tau=TYURINANumber(F2,\{Z,Y\},[plex_min(Z,Y),tdeg(Z,Y)]))}
{\[{\mu={\it MilnorNumber} \left( {\it F2}, \left\{ Z,Y \right\} ,{\it plex\_min}\\
\mbox{} \left( Z,Y \right)  \right) ,\tau={\it TYURINANumber} \left( {\it F2}, \left\{ Z,Y \right\} ,[{\it plex\_min}\\
\mbox{} \left( Z,Y \right) ,{\it tdeg}\\
\mbox{} \left( Z,Y \right) ] \right) }\]}
\end{mapleinput}
\mapleresult
\begin{maplelatex}
\mapleinline{inert}{2d}{mu = 4, tau = 4}{\[\displaystyle
\mu=4,\,\tau=4\]}
\end{maplelatex}
\end{maplegroup}
\begin{maplegroup}
\begin{mapleinput}
\mapleinline{active}{12d}{w[1] := -3*u*t : F2}{\[\]}
\end{mapleinput}
\mapleresult
\begin{maplelatex}
\mapleinline{inert}{2d}{} {\[\displaystyle {Y}^{3}+Y{Z}^{3}+
\left( -3\,tuY-3\,{u}^{2}tZ \right) ^{2}+ \left(
-6\,{t}^{2}u-3\,{u}^{2} \right) Y{Z}^{2}\]}
\mapleinline{inert}{2d}{}
{\[\displaystyle\quad +u \left( -6\,{t}^{2}u-2\,{u}^{2} \right) {Z}^{3}\\
\mbox{}+ \left( {t}^{2}+u \right) {Z}^{4}\]}
\end{maplelatex}
\end{maplegroup}
\begin{maplegroup}
\begin{mapleinput}
\mapleinline{active}{1d}{print(mu=MilnorNumber(F2,\{Z,Y\},plex_min(Z,Y)),
tau=TYURINANumber(F2,\{Z,Y\},[plex_min(Z,Y),tdeg(Z,Y)]))}
{\[{\mu={\it MilnorNumber} \left( {\it F2}, \left\{ Z,Y \right\} ,{\it plex\_min}\\
\mbox{} \left( Z,Y \right)  \right) ,\tau={\it TYURINANumber} \left( {\it F2}, \left\{ Z,Y \right\} ,[{\it plex\_min}\\
\mbox{} \left( Z,Y \right) ,{\it tdeg}\\
\mbox{} \left( Z,Y \right) ] \right) }\]}
\end{mapleinput}
\mapleresult
\begin{maplelatex}
\mapleinline{inert}{2d}{mu = 5, tau = 5}{\[\displaystyle
\mu=5,\,\tau=5\]}
\end{maplelatex}
\end{maplegroup}

\halfline

\noindent and, by eliminating $u,t$ and $w_1$ from equations $w_1=-3ut$ and (\ref{parametrizzazioneE7}), define
\begin{equation}\label{A5bis in E7}
    \mathcal{W}_5':=\left\{v_1(v_1^2-9v_2v_5)^2-81v_2^3=0\right\}\quad,
    \quad\widetilde{\mathcal{W}'}_2^5:=\widetilde{\mathcal{W}}_2^4\cap\mathcal{W}_5'\ .
\end{equation}
Hence Theorem \ref{classificazione}(2) gives that
\begin{itemize}
    \item \emph{$0\in f_{\Lambda}^{-1}(0)$ is a singularity of type $A_5$
    for    generic
    $\Lambda\in\widetilde{\mathcal{W}'}_2^5$}.
\end{itemize}
Consider the associated standard basis:

\halfline

\begin{maplegroup}
\begin{mapleinput}
\mapleinline{active}{1d}{MGB5c:=MilnorGroebnerBasis(F2,\{Z,Y\},plex_min(Z,Y))}
{\[{\it MGB5c}\, := \,{\it MilnorGroebnerBasis} \left( {\it F2}, \left\{ Z,Y \right\} ,{\it plex\_min}\\
\mbox{} \left( Z,Y \right)  \right) \]}
\end{mapleinput}
\mapleresult
\begin{maplelatex}
\mapleinline{inert}{2d}{}
{\[\displaystyle {\it MGB5c}\, := \, \left\{  3\,Y{Z}^{2}-6\, \left( -3\,tuY-3\,{u}^{2}tZ \right) {u}^{2}t\\
\mbox{}+2\, \left( -6\,{t}^{2}u-3\,{u}^{2} \right) YZ\right.\]}
\mapleinline{inert}{2d}{}
{\[\displaystyle +3\,u \left(
-6\,{t}^{2}u-2\,{u}^{2} \right) {Z}^{2}+4\, \left( {t}^{2}+u
\right) {Z}^{3},\]}
\mapleinline{inert}{2d}{}
{\[\displaystyle 3\,{Y}^{2}\\
\mbox{}+{Z}^{3}-6\, \left( -3\,tuY-3\,{u}^{2}tZ \right) tu+ \left(
-6\,{t}^{2}u-3\,{u}^{2} \right) {Z}^{2},\]}
\mapleinline{inert}{2d}{}
{\[\displaystyle \left( 48\,{t}^{4}-24\,{t}^{2}u \right) {Z}^{5}Y+ \left( 192\,{t}^{6}{u}^{2}-144\,{u}^{3}{t}^{4} \right) {Z}^{5}\\
\mbox{}+ \left( -16\,{t}^{4}u+16\,{u}^{2}{t}^{2}+64\,{t}^{6}
\right) {Z}^{6} \left.\right\} \]}
\end{maplelatex}
\end{maplegroup}
\begin{maplegroup}
\begin{mapleinput}
\mapleinline{active}{1d}{for i to nops(MGB5c) do LeadingTerm(MGB5c[i],plex_min(Z,Y))end do;}{\[\]}
\end{mapleinput}
\mapleresult
\begin{maplelatex}
\mapleinline{inert}{2d}{18*u^3*t^2, Y}{\[\displaystyle
18\,{u}^{3}{t}^{2},\,Y\]}
\end{maplelatex}
\mapleresult
\begin{maplelatex}
\mapleinline{inert}{2d}{18*u^2*t^2, Y}{\[\displaystyle
18\,{u}^{2}{t}^{2},\,Y\]}
\end{maplelatex}
\mapleresult
\begin{maplelatex}
\mapleinline{inert}{2d}{192*t^6*u^2-144*u^3*t^4,
Z^5}{\[\displaystyle
192\,{t}^{6}{u}^{2}-144\,{u}^{3}{t}^{4},\,{Z}^{5}\]}
\end{maplelatex}
\end{maplegroup}
\begin{maplegroup}
\begin{mapleinput}
\mapleinline{active}{1d}{u := 0 : F2}{\[\]}
\end{mapleinput}
\mapleresult
\begin{maplelatex}
\mapleinline{inert}{2d}{Y^3+Y*Z^3+t^2*Z^4}{\[\displaystyle
{Y}^{3}+Y{Z}^{3}+{t}^{2}{Z}^{4}\]}
\end{maplelatex}
\end{maplegroup}
\begin{maplegroup}
\begin{mapleinput}
\mapleinline{active}{1d}{print(mu=MilnorNumber(F2,\{Z,Y\},plex_min(Z,Y)),
tau=TYURINANumber(F2,\{Z,Y\},[plex_min(Z,Y),tdeg(Z,Y)]))}
{\[{\mu={\it MilnorNumber} \left( {\it F2}, \left\{ Z,Y \right\} ,{\it plex\_min}\\
\mbox{} \left( Z,Y \right)  \right) ,\tau={\it TYURINANumber} \left( {\it F2}, \left\{ Z,Y \right\} ,[{\it plex\_min}\\
\mbox{} \left( Z,Y \right) ,{\it tdeg}\\
\mbox{} \left( Z,Y \right) ] \right) }\]}
\end{mapleinput}
\mapleresult
\begin{maplelatex}
\mapleinline{inert}{2d}{mu = 6, tau = 6}{\[\displaystyle
\mu=6,\,\tau=6\]}
\end{maplelatex}
\end{maplegroup}
\begin{maplegroup}
\begin{mapleinput}
\mapleinline{active}{1d}{unassign('u') : t := 0 : F2}{\[\]}
\end{mapleinput}
\mapleresult
\begin{maplelatex}
\mapleinline{inert}{2d}{Y^3+Y*Z^3-3*u^2*Y*Z^2-2*u^3*Z^3+Z^4*u}{\[\displaystyle
{Y}^{3}+Y{Z}^{3}-3\,{u}^{2}Y{Z}^{2}-2\,{u}^{3}{Z}^{3}+{Z}^{4}u\]}
\end{maplelatex}
\end{maplegroup}

\halfline

\noindent This gives precisely the same situation of generic $E_6$
and $D_6$ singularities previously considered, meaning that
\begin{itemize}
    \item \emph{$0\in f_{\Lambda}^{-1}(0)$ is a singularity of type $E_6$
    for    generic
    $\Lambda\in\mathcal{V}_2\cap\widetilde{\mathcal{W}'}_2^5=\mathcal{V}_0^4$},
    \item \emph{$0\in f_{\Lambda}^{-1}(0)$ is a singularity of type $D_6$
    for    generic
    $\Lambda\in\mathcal{W}_5\cap\widetilde{\mathcal{W}'}_2^5=\widetilde{\mathcal{W}}_2^5\cap\mathcal{W}'_5=
    \mathcal{V}'\cap\mathcal{V}_0^2$}.
\end{itemize}
Moreover the last leading coefficient in $MGB5c$ gives the further
relation

\halfline

\begin{maplegroup}
\begin{mapleinput}
\mapleinline{active}{1d}{u := (4/3)*t^2 : F2}{\[\]}
\end{mapleinput}
\mapleresult
\begin{maplelatex}
\mapleinline{inert}{2d}{Y^3+Y*Z^3+(-4*t^3*Y-(16/3)*t^5*Z)^2-(40/3)*t^4*Y*Z^2-(416/27)*t^6*Z^3+(7/3)*t^2*Z^4}{\[\displaystyle {Y}^{3}+Y{Z}^{3}+ \left( -4\,{t}^{3}Y-16/3\,{t}^{5}Z \right) ^{2}-{\frac {40}{3}}\,{t}^{4}Y{Z}^{2}-{\frac {416}{27}}\,{t}^{6}{Z}^{3}\\
\mbox{}+7/3\,{t}^{2}{Z}^{4}\]}
\end{maplelatex}
\end{maplegroup}
\begin{maplegroup}
\begin{mapleinput}
\mapleinline{active}{1d}{print(mu=MilnorNumber(F2,\{Z,Y\},plex_min(Z,Y)),
tau=TYURINANumber(F2,\{Z,Y\},[plex_min(Z,Y),tdeg(Z,Y)]))}
{\[{\mu={\it MilnorNumber} \left( {\it F2}, \left\{ Z,Y \right\} ,{\it plex\_min}\\
\mbox{} \left( Z,Y \right)  \right) ,\tau={\it TYURINANumber} \left( {\it F2}, \left\{ Z,Y \right\} ,[{\it plex\_min}\\
\mbox{} \left( Z,Y \right) ,{\it tdeg}\\
\mbox{} \left( Z,Y \right) ] \right) }\]}
\end{mapleinput}
\mapleresult
\begin{maplelatex}
\mapleinline{inert}{2d}{mu = 6, tau = 6}{\[\displaystyle
\mu=6,\,\tau=6\]}
\end{maplelatex}
\end{maplegroup}

\halfline

\noindent meaning that
\begin{itemize}
    \item \emph{$0\in f_{\Lambda}^{-1}(0)$ is a singularity of type $A_6$
    for    generic
    $\Lambda\in\widetilde{\mathcal{W}'}_2^6$}
\end{itemize}
where
\begin{equation}\label{A6 in E7}
    \widetilde{\mathcal{W}'}_2^6:=\widetilde{\mathcal{W}'}_2^5\cap\mathcal{W}_6\quad
    \text{and}\quad\mathcal{W}_6:=\left\{16v_1^5-729v_2^3=0\right\}
\end{equation}
The reader can then easily check that (\ref{intersezioni in E7})
holds.
\end{proof}

\subsection{Simple singularities of $E_8$ type.}

\begin{theorem}\label{E8}
Let $T^1$ be the Kuranishi space of a simple
$N$--dimensional singular point $0\in f^{-1}(0)$ with
\[
    f(x_1,\ldots,x_{N+1}) = \sum_{i=1}^{N-1} x_i^2\ +\
    x_N^3+ x_{N+1}^5
\]
The subset of $T^1$ parameterizing small deformations of $0\in f^{-1}(0)$ to a simple node is the union of $n$ hypersurfaces. Moreover, calling $\mathcal{L}$ any of those hypersurfaces, there exists a stratification of nested
algebraic subsets giving rise to the following sequence of inclusions, c.i.p. squares and reducible c.i.p squares
\begin{equation}\label{E8-inclusioni}
    \xymatrix{\mathcal{L}&\ar@{_{(}->}[l]\ \mathcal{W}_2&\ar@{_{(}->}[l]\
                \mathcal{W}_2^3&
                \ar@{_{(}->}[l]\ \widetilde{\mathcal{W}}_2^4&\ar@{_{(}->}[l]\
                \widetilde{\mathcal{W}}_2^5&\ar@{_{(}->}[l]\
                \widetilde{\mathcal{W}}_2^6&\ar@{_{(}->}[l]\
                \widetilde{\mathcal{W}}_2^7\\
                &&\ar@{_{(}->}[u]\ \mathcal{V}_0^2&\ar@{_{(}->}[u]\ar@{_{(}->}[l]\
                \mathcal{V}\cap\mathcal{V}_0^2&\ar@{_{(}->}[u]\ar@{_{(}->}[l]\
                \mathcal{V}'\cap\mathcal{V}_0^2&\ar@{_{(}->}[u]\ar@{_{(}->}[l]\
                \mathcal{V}''\cap\mathcal{V}_0^2&\\
                &&&\ar@{_{(}->}[ruu]\ar@{_{(}->}[u]\ \mathcal{V}_0^4&
                \ar@{_{(}->}[u]\ar@{_{(}->}[l]\ar@{_{(}->}[ruu]\ \mathcal{V}_0^4\cap\mathcal{V}_6&
                \ar@{_{(}->}[u]\ar@{_{(}->}[l]\ar@{_{(}->}[ruu]\ \{0\}}
\end{equation}
verifying the Arnol'd's adjacency diagram
\begin{equation}\label{E8-adiacenza}
    \xymatrix{A_1&\ar[l]A_2&\ar[l]A_3&\ar[l]A_4&\ar[l]A_5&\ar[l]A_6&\ar[l]A_7\\
                &&\ar[u]D_4&\ar[u]\ar[l]D_5&\ar[u]\ar[l]D_6&\ar[u]\ar[l]D_7&\\
                &&&\ar[ruu]\ar[u]E_6&\ar[ruu]\ar[u]\ar[l]E_7&\ar[ruu]\ar[u]\ar[l]E_8&}
\end{equation}
where
\begin{itemize}
    \item $\mathcal{L}$ is the hypersurface of $T^1$ defined by equation
(\ref{lamda-condizione:E_8}), keeping in mind
(\ref{E_8,z,y}),
    \item $\mathcal{V}_0^m:=\bigcap_{k=0}^m \mathcal{V}_k$
and $\mathcal{W}_2^m:=\bigcap_{k=2}^m \mathcal{W}_k$ where
$\mathcal{V}_k, \mathcal{W}_k$ are hypersurfaces of $\mathcal{L}$
defined by equations (\ref{vk-E8}), (\ref{A3 in E8}) ,(\ref{A4 in E8}), (\ref{A5 in E8bis}), (\ref{A6 in E8bis})
and (\ref{A7 in E8}),
    \item $\mathcal{V}, \mathcal{V}'$ and $\mathcal{V}''$ are a hypersurface, a codimension 2 and
    a codimension 3 complete intersections in $\mathcal{L}$, respectively, defined
    by the latter equation in (\ref{A3 in E8}),  by (\ref{D6 in E8}) and by (\ref{D6 in E8 bis}),
    respectively,
    \item $\widetilde{\mathcal{W}}_2^k$ are Zariski closed subsets of
    $\mathcal{L}$ defined by (\ref{A4bis in E8}), (\ref{A5 in E8}), (\ref{A6 in E8}) and (\ref{A7 in E8}).
\end{itemize}
In particular, complete intersection properties in diagram (\ref{E8-inclusioni}) are summarized by the following relations:
\begin{eqnarray}\label{intersezioni in E8}
    \widetilde{\mathcal{W}}_2^4\cap\mathcal{V}_0^2&=&\mathcal{V}\cap\mathcal{V}_0^2\ ,\\
\nonumber
    \widetilde{\mathcal{W}}_2^5\cap\left(\mathcal{V}\cap\mathcal{V}_0^2\right)&=&
    (\mathcal{V}'\cap\mathcal{V}_0^2)\cup\mathcal{V}_0^4\ ,\\
\nonumber
    \widetilde{\mathcal{W}}_2^6\cap\left(\mathcal{V}'\cap\mathcal{V}_0^2\right)&=& (\mathcal{V}''\cap\mathcal{V}_0^2)\cup(\mathcal{V}_0^4\cap\mathcal{V}_6)\ , \\
\nonumber
    \left(\mathcal{V}''\cap\mathcal{V}_0^2\right)
    \cap\left(\mathcal{V}_0^4\cap\mathcal{V}_6\right)&=&
    \widetilde{\mathcal{W}}_2^7\cap\left(\mathcal{V}''\cap\mathcal{V}_0^2\right)\ =\ \{0\}\ .
\end{eqnarray}
\end{theorem}

\begin{proof} Following the outline \ref{outline}.

\vskip 5pt\noindent (1) By the Morse Splitting Lemma \ref{Morse}, our
problem can be reduced to the case $N=1$ with $f(y,z)= y^3 +
z^5$. To get an explicit basis of the Kuranishi space $T^1$ type:

\halfline

\begin{maplegroup}
\begin{mapleinput}
\mapleinline{active}{1d}{F:= y\symbol{94}3 + z\symbol{94}5; }{}
\end{mapleinput}
\end{maplegroup}
\begin{maplegroup}
\begin{mapleinput}
\mapleinline{active}{1d}{TyB := TyurinaBasis(F);}{\[\]}
\end{mapleinput}
\mapleresult
\begin{maplelatex}
\mapleinline{inert}{2d}{TyB := [y*z^3, z^3, y*z^2, z^2, y*z, z, y,
1]} {\[\displaystyle {\it TyB}\, :=
\,[y{z}^{3},{z}^{3},y{z}^{2},{z}^{2},yz,z,y,1]\]}
\end{maplelatex}
\end{maplegroup}

\halfline

\noindent Therefore
\[
T^1 \cong \langle 1,y,z,yz,z^2,yz^2,z^3,yz^3\rangle_{\C}
\]
and given $\Lambda = (\lambda_0, \lambda_1,\ldots,\lambda_7)\in
T^1$, the associated deformation of $U_0$ is
\begin{eqnarray*}
   U_{\Lambda} &=& \{f_{\Lambda}(y,z)=0\}\quad\text{where}\quad f_{\Lambda}(y,z) :=  \\
   &=& f(y,z) + \lambda_0 + \lambda_1 y
    + \lambda_2 z + \lambda_3 yz + \lambda_4 z^2 + \lambda_5 yz^2
    + \lambda_6 z^3 + \lambda_7 yz^3
\end{eqnarray*}
A solution of the jacobian system of partial derivatives is then
given by a solution $p_{\Lambda}=(y_{\Lambda},z_{\Lambda})$ of the following polynomial system in $\C[\lambda][y,z]$
\begin{equation}\label{E_8,z,y}
   \left\{\begin{array}{c}
      3y^2 + \lambda_1 + \lambda_3z + \lambda_5z^2
      + \lambda_7z^3 = 0 \\
      5z^4 +\lambda_2 + 2 \lambda_4z + 3\lambda_6z^2
      + y(\lambda_3 + 2\lambda_5 z
      + 3\lambda_7z^2)= 0\\
    \end{array}\right.
\end{equation}
giving 8 critical points for $f_{\Lambda}$.

\vskip 5pt\noindent (2) Imposing that one of those critical points, say $p_{\Lambda}$, is actually a singular point of $U_{\Lambda}$ means to require that
\begin{eqnarray}\label{lamda-condizione:E_8}
    p_{\Lambda}&\in& U_{\Lambda}\ \Longleftrightarrow \ \Lambda\in
    \mathcal{L}\subset T^1\quad\text{where}\\
\nonumber
    \mathcal{L}&:=&\{ y_{\Lambda}^3 + z_{\Lambda}^5 + \lambda_0 + \lambda_1y_{\Lambda}
    + \lambda_2z_{\Lambda} + \lambda_3y_{\Lambda}z_{\Lambda} \\
    &&\quad + \lambda_4z_{\Lambda}^2 +
    \lambda_5y_{\Lambda}z_{\Lambda}^2 + \lambda_6 z_{\Lambda}^3
    + \lambda_7 y_{\Lambda}z_{\Lambda}^3=
    0\}\ .
\end{eqnarray}
After translating $y\mapsto y+y_{\Lambda},z\mapsto z+z_{\Lambda}$,
we get
\begin{eqnarray*}
  &&f_{\Lambda}(y+y_{\Lambda},z+z_{\Lambda})= f(y,z) + \\
    &+&\left(y_{\Lambda}^3 + z_{\Lambda}^5 + \lambda_0 + \lambda_1y_{\Lambda}
    + \lambda_2z_{\Lambda} + \lambda_3y_{\Lambda}z_{\Lambda}
    + \lambda_4z_{\Lambda}^2 +
    \lambda_5y_{\Lambda}z_{\Lambda}^2 + \lambda_6 z_{\Lambda}^3
    + \lambda_7 y_{\Lambda}z_{\Lambda}^3\right) \\
    &+& \left(3y_{\Lambda}^2 + \lambda_1 + \lambda_3z_{\Lambda} + \lambda_5z_{\Lambda}^2
      + \lambda_7z_{\Lambda}^3\right)y \\
    &+& \left(5z_{\Lambda}^4 +\lambda_2 + 2 \lambda_4z_{\Lambda} + 3\lambda_6z_{\Lambda}^2
      + y_{\Lambda}(\lambda_3 + 2\lambda_5 z_{\Lambda}
      + 3\lambda_7z_{\Lambda}^2)\right)z \\
    &+& \left(\lambda_3 + 2\lambda_5z_{\Lambda} + 3\lambda_7z_{\Lambda}^2\right)y z + 3y_{\Lambda}\ y^2 +
    \left(10z_{\Lambda}^3 + \lambda_4 + \lambda_5y_{\Lambda} + 3\lambda_6z_{\Lambda}
    + 3\lambda_7y_{\Lambda}z_{\Lambda}\right)z^2 \\
    &+& \left(\lambda_5 + 3\lambda_7z_{\Lambda}\right)y z^2
    + \left(10z_{\Lambda}^2 + \lambda_6 + \lambda_7y_{\Lambda}\right) z^3
    + \lambda_7\ yz^3 + 5 z_{\Lambda}\ z^4
\end{eqnarray*}
as can be checked by setting $K:=y_{\Lambda},L:=z_{\Lambda}$ and
typing:

\halfline

\begin{maplegroup}
\begin{mapleinput}
\mapleinline{active}{1d}{T := nops(TyB):}{\[\]}
\end{mapleinput}
\end{maplegroup}
\begin{maplegroup}
\begin{mapleinput}
\mapleinline{active}{1d}{F[Lambda] := F+sum(lambda[i]*TyB[T-i],i=0 .. T-1)}
 {\[F\Lambda \, := \,F+\sum
_{i=0}^{T-1}\lambda_{{i}}{\it TyB}_{{T-i}}\]}
\end{mapleinput}
\mapleresult
\begin{maplelatex}
\mapleinline{inert}{2d}{`F&Lambda;` := y^3+z^5+lambda[0]+lambda[1]*y+lambda[2]*z+lambda[3]*y*z+lambda[4]*z^2+lambda[5]*y*z^2+lambda[6]*z^3+lambda[7]*y*z^3}{\[\displaystyle F\Lambda \, := \,{y}^{3}+{z}^{5}+\lambda_{{0}}+\lambda_{{1}}y+\lambda_{{2}}z+\lambda_{{3}}yz+\lambda_{{4}}{z}^{2}+\lambda_{{5}}y{z}^{2}\\
\mbox{}+\lambda_{{6}}{z}^{3}+\lambda_{{7}}y{z}^{3}\]}
\end{maplelatex}
\end{maplegroup}
\begin{maplegroup}
\begin{mapleinput}
\mapleinline{active}{1d}{z := Z+L: y := Y+K:}{\[\]}
\end{mapleinput}
\begin{mapleinput}
\mapleinline{active}{1d}{F[Lambda] := collect(F[Lambda], [Y, Z],'distributed'):}{\[\]}
\end{mapleinput}
\mapleresult
\begin{maplelatex}
\mapleinline{inert}{2d}{}
{\[\displaystyle F\Lambda \, := \,{Y}^{3}+3\,K{Y}^{2}+\lambda_{{7}}Y{Z}^{3}+ \left( \lambda_{{5}}+3\,\lambda_{{7}}L \right)
 Y{Z}^{2}+ \left( 3\,\lambda_{{7}}{L}^{2}+\lambda_{{3}}+2\,\lambda_{{5}}L \right) YZ\]}
\mapleinline{inert}{2d}{}
{\[\displaystyle +\left(\lambda_{{7}}{L}^{3}+\lambda_{{1}}+\lambda_{{5}}{L}^{2}+
\lambda_{{3}}L+3\,{K}^{2}\right) Y+{Z}^{5}+5\,L{Z}^{4}+ \left(
\lambda_{{6}}+\lambda_{{7}}K+10\,{L}^{2} \right) {Z}^{3}\]}
\mapleinline{inert}{2d}{}
{\[\displaystyle +\left(3\,\lambda_{{7}}KL+\lambda_{{4}}+3\,\lambda_{{6}}L+10\,{L}^{3}+\lambda_{{5}}K
\right) {Z}^{2}\]}
\mapleinline{inert}{2d}{}
{\[\displaystyle + \left( \lambda_{{2}}+3\,\lambda_{{7}}K{L}^{2}+5\,{L}^{4}+3\,\lambda_{{6}}{L}^{2}+\lambda_{{3}}K+2\,
\lambda_{{5}}KL+2\,\lambda_{{4}}L \right) Z\]}
\mapleinline{inert}{2d}{}
{\[\displaystyle +{K}^{3}+\lambda_{{6}}{L}^{3}+\lambda_{{1}}K+\lambda_{{2}}L+{L}^{5}+\lambda_{{7}}K{L}^{3}+\lambda_{{0}}+\lambda_{{5}}K{L}^{2}+\lambda_{{4}}{L}^{2}\\
\mbox{}+\lambda_{{3}}KL\]}
\end{maplelatex}
\end{maplegroup}

\halfline

\vskip 5pt\noindent (3) Define:
\begin{eqnarray}\label{vk-E8}
\nonumber
   &v_0=\lambda_3 + 2\lambda_5z_{\Lambda} + 3\lambda_7z_{\Lambda}^2\quad,
   \quad v_1=3y_{\Lambda}&\\
   &v_2=10z_{\Lambda}^3 + \lambda_4 + \lambda_5y_{\Lambda} + 3\lambda_6z_{\Lambda}
    + 3\lambda_7y_{\Lambda}z_{\Lambda}\quad,\quad v_3=\lambda_5 + 3\lambda_7z_{\Lambda}&\\
\nonumber
   &v_4=10z_{\Lambda}^2 + \lambda_6 +
    \lambda_7y_{\Lambda}\quad ,\quad v_5=\lambda_7\quad,\quad v_6=5z_{\Lambda}&
\end{eqnarray}
and $\mathcal{V}_k:=\{v_k=0\}$. Then the proof goes on exactly as
in the $E_7$ case until $D_5$ singularities: precisely by setting
\begin{eqnarray}\label{A3 in E8}
\nonumber
    \mathcal{W}_2 &:=& \left\{4v_1v_2-v_0^2=0\right\} \\
    \mathcal{W}_3 &:=& \left\{v_1^3v_4^2-v_2(v_1v_3+v_2)^2=0\right\}\\
\nonumber
    \mathcal{V} &:=&\left\{4v_3^3+27v_4^2=0\right\}
\end{eqnarray}
we have the inclusions' chain
$$T^1\supset\mathcal{L}\supset\mathcal{W}_2\supset\mathcal{W}_2^3
\supset\mathcal{V}_1\cap\mathcal{W}_2^3=\mathcal{V}_0^2\supset\mathcal{V}\cap
\mathcal{V}_0^2$$ of subsets parameterizing small deformations
whose generic fibre is either smooth or admits a singularity of
type $A_1$, $A_2$, $A_3$ , $D_4$ and $D_5$, respectively, as can be checked by
typing

 \halfline

\begin{maplegroup}
\begin{mapleinput}
\mapleinline{active}{1d}{SD := [Y*Z, Y^2, Z^2, Y*Z^2, Z^3, Y*Z^3, Z^4]:}{\[\]}
\end{mapleinput}
\begin{mapleinput}
\mapleinline{active}{1d}{s := nops(SD): FL := Y^3+Z^5+sum(v[i]*SD[i+1], i = 0 .. s-1)}{\[\]}
\end{mapleinput}
\mapleresult
\begin{maplelatex}
\mapleinline{inert}{2d}{FL := Y^3+Z^5+v[0]*Y*Z+v[1]*Y^2+v[2]*Z^2+v[3]*Y*Z^2+v[4]*Z^3+v[5]*Y*Z^3+v[6]*Z^4}
{\[\displaystyle {\it FL}\, := \,{Y}^{3}+{Z}^{5}+v_{{0}}YZ+v_{{1}}{Y}^{2}+v_{{2}}{Z}^{2}+v_{{3}}Y{Z}^{2}
+v_{{4}}{Z}^{3}+v_{{5}}Y{Z}^{3}\\
\mbox{}+v_{{6}}{Z}^{4}\]}
\end{maplelatex}
\begin{mapleinput}
\mapleinline{active}{1d}{MGB := MilnorGroebnerBasis(FL,\{Y,Z\},plex_min(Z,Y)):}{\[\]}
\end{mapleinput}
\end{maplegroup}
\begin{maplegroup}
\begin{mapleinput}
\mapleinline{active}{1d}{for i to nops(MGB) do LeadingTerm(MGB[i],plex_min(Z, Y))end do}{\[\]}
\end{mapleinput}
\mapleresult
\begin{maplelatex}
\mapleinline{inert}{2d}{v[0], Y}{\[\displaystyle v_{{0}},\,Y\]}
\end{maplelatex}
\mapleresult
\begin{maplelatex}
\mapleinline{inert}{2d}{2*v[1], Y}{\[\displaystyle
2\,v_{{1}},\,Y\]}
\end{maplelatex}
\mapleresult
\begin{maplelatex}
\mapleinline{inert}{2d}{-4*v[1]^2*v[0]^2+16*v[1]^3*v[2],
Z}{\[\displaystyle
-4\,{v_{{1}}}^{2}{v_{{0}}}^{2}+16\,{v_{{1}}}^{3}v_{{2}},\,Z\]}
\end{maplelatex}
\end{maplegroup}
\begin{maplegroup}
\begin{mapleinput}
\mapleinline{active}{1d}{F2 := Y^3+Z^5+(w[1]*Y+w[2]*Z)^2+sum(v[j]*SD[j+1], j = 3 .. s-1)}{\[{\it
F2}\, := \,{Y}^{3}+{Z}^{5}+ \left( w_{{1}}Y+w_{{2}}Z \right)
^{2}+\sum _{j=3}^{s-1}v_{{j}}{\it SD}_{{j+1}}\]}
\end{mapleinput}
\mapleresult
\begin{maplelatex}
\mapleinline{inert}{2d}{F2 :=
Y^3+Z^5+(w[1]*Y+w[2]*Z)^2+v[3]*Y*Z^2+v[4]*Z^3+v[5]*Y*Z^3+v[6]*Z^4}{\[\displaystyle
{\it F2}\, := \,{Y}^{3}+{Z}^{5}+ \left( w_{{1}}Y+w_{{2}}Z \right)
^{2}+v_{{3}}Y{Z}^{2}+v_{{4}}{Z}^{3}+v_{{5}}Y{Z}^{3}+v_{{6}}{Z}^{4}\]}
\end{maplelatex}
\end{maplegroup}
\begin{maplegroup}
\begin{mapleinput}
\mapleinline{active}{1d}{print(mu = MilnorNumber(F2,\{Y,Z\}, plex_min(Z, Y)))}
{\[\mu={\it MilnorNumber} \left( {\it F2}, \left\{ Y,Z \right\} ,{\it plex\_min}\\
\mbox{} \left( Z,Y \right)  \right) \]}
\end{mapleinput}
\mapleresult
\begin{maplelatex}
\mapleinline{inert}{2d}{mu = 2}{\[\displaystyle \mu=2\]}
\end{maplelatex}
\end{maplegroup}
\begin{maplegroup}
\begin{mapleinput}
\mapleinline{active}{1d}{MGB2 := MilnorGroebnerBasis(F2,\{Y,Z\},plex_min(Z,Y)):}{\[\]}
\end{mapleinput}
\end{maplegroup}
\begin{maplegroup}
\begin{mapleinput}
\mapleinline{active}{1d}{for i to nops(MGB2) do LeadingTerm(MGB2[i],plex_min(Z, Y))end do}{\[\]}
\end{mapleinput}
\mapleresult
\begin{maplelatex}
\mapleinline{inert}{2d}{2*w[2]*w[1], Y}{\[\displaystyle
2\,w_{{2}}w_{{1}},\,Y\]}
\end{maplelatex}
\mapleresult
\begin{maplelatex}
\mapleinline{inert}{2d}{2*w[1]^2, Y}{\[\displaystyle
2\,{w_{{1}}}^{2},\,Y\]}
\end{maplelatex}
\mapleresult
\begin{maplelatex}
\mapleinline{inert}{2d}{-12*w[1]^3*w[2]^3-12*w[1]^5*w[2]*v[3]+12*w[1]^6*v[4],
Z^2}{\[\displaystyle
-12\,{w_{{1}}}^{3}{w_{{2}}}^{3}-12\,{w_{{1}}}^{5}w_{{2}}v_{{3}}+12\,{w_{{1}}}^{6}v_{{4}},\,{Z}^{2}\]}
\end{maplelatex}
\end{maplegroup}
\begin{maplegroup}
\begin{mapleinput}
\mapleinline{active}{1d}{w[2] := u*w[1] : v[4]:= u*(v[3]+u^2): F2}{\[\]}
\end{mapleinput}
\mapleresult
\begin{maplelatex}
\mapleinline{inert}{2d}{Y^3+Z^5+(w[1]*Y+u*w[1]*Z)^2+v[3]*Y*Z^2+u*(v[3]+u^2)*Z^3+v[5]*Y*Z^3+v[6]*Z^4}
{\[\displaystyle {Y}^{3}+{Z}^{5}+ \left( w_{{1}}Y+uw_{{1}}Z \right) ^{2}+v_{{3}}Y{Z}^{2}+u \left( v_{{3}}+{u}^{2} \right) {Z}^{3}+v_{{5}}Y{Z}^{3}\\
\mbox{}+v_{{6}}{Z}^{4}\]}
\end{maplelatex}
\end{maplegroup}
\begin{maplegroup}
\begin{mapleinput}
\mapleinline{active}{1d}{M1 := Milnor(F2,\{Y,Z\},plex_min(Z,Y)): print(mu = M1[4]);}{\[\]}
\end{mapleinput}
\mapleresult
\begin{maplelatex}
\mapleinline{inert}{2d}{mu = 3}{\[\displaystyle \mu=3\]}
\end{maplelatex}
\end{maplegroup}
\begin{maplegroup}
\begin{mapleinput}
\mapleinline{active}{1d}{MGB3 := M1[1] : }{\[\]}
\end{mapleinput}
\begin{mapleinput}
\mapleinline{active}{1d}{for i to nops(MGB3) do LeadingTerm(MGB3[i], plex_min(Z, Y))end do}{\[\]}
\end{mapleinput}
\mapleresult
\begin{maplelatex}
\mapleinline{inert}{2d}{2*u*w[1]^2, Y}{\[\displaystyle
2\,u{w_{{1}}}^{2},\,Y\]}
\end{maplelatex}
\mapleresult
\begin{maplelatex}
\mapleinline{inert}{2d}{-8*w[1]^4*v[3]^2+32*w[1]^6*v[6]-32*w[1]^6*u*v[5]-48*w[1]^4*u^2*v[3]-72*w[1]^4*u^4, Z^3}{\[\displaystyle -8\,{w_{{1}}}^{4}{v_{{3}}}^{2}+32\,{w_{{1}}}^{6}v_{{6}}-32\,{w_{{1}}}^{6}uv_{{5}}-48\,{w_{{1}}}^{4}{u}^{2}v_{{3}}\\
\mbox{}-72\,{w_{{1}}}^{4}{u}^{4},\,{Z}^{3}\]}
\end{maplelatex}
\mapleresult
\begin{maplelatex}
\mapleinline{inert}{2d}{2*w[1]^2, Y}{\[\displaystyle
2\,{w_{{1}}}^{2},\,Y\]}
\end{maplelatex}
\end{maplegroup}
\begin{maplegroup}
\begin{mapleinput}
\mapleinline{active}{1d}{w[1] := 0:}{\[\]}
\end{mapleinput}
\end{maplegroup}
\begin{maplegroup}
\begin{mapleinput}
\mapleinline{active}{1d}{print(mu = MilnorNumber(F2,\{Y,Z\}, plex_min(Z,Y)),
tau = TYURINANumber(F2,\{Y,Z\}, [plex_min(Z, Y), tdeg(Z, Y)]))}
{\[{\mu={\it MilnorNumber} \left( {\it F2}, \left\{ Y,Z \right\} ,{\it plex\_min}\\
\mbox{} \left( Z,Y \right)  \right) ,\tau={\it TYURINANumber} \left( {\it F2}, \left\{ Y,Z \right\} ,[{\it plex\_min}\\
\mbox{} \left( Z,Y \right) ,{\it tdeg}\\
\mbox{} \left( Z,Y \right) ] \right) }\]}
\end{mapleinput}
\mapleresult
\begin{maplelatex}
\mapleinline{inert}{2d}{mu = 4, tau = 4}{\[\displaystyle
\mu=4,\,\tau=4\]}
\end{maplelatex}
\end{maplegroup}
\begin{maplegroup}
\begin{mapleinput}
\mapleinline{active}{1d}{unassign('w[1]'):unassign('v[4]'):v[0]:=0:v[1]:=0:v[2]:=0:FL}{\[\]}
\end{mapleinput}
\mapleresult
\begin{maplelatex}
\mapleinline{inert}{2d}{Y^3+Z^5+v[3]*Y*Z^2+v[4]*Z^3+v[5]*Y*Z^3+v[6]*Z^4}{\[\displaystyle
{Y}^{3}+{Z}^{5}+v_{{3}}Y{Z}^{2}+v_{{4}}{Z}^{3}+v_{{5}}Y{Z}^{3}+v_{{6}}{Z}^{4}\]}
\end{maplelatex}
\end{maplegroup}
\begin{maplegroup}
\begin{mapleinput}
\mapleinline{active}{1d}{MGB4 := MilnorGroebnerBasis(FL,\{Y,Z\},plex_min(Z,Y)):}{\[ \]}
\end{mapleinput}
\begin{mapleinput}
\mapleinline{active}{1d}{for i to nops(MGB4) do LeadingTerm(MGB4[i],plex_min(Z,Y))end do}{\[\]}
\end{mapleinput}
\mapleresult
\begin{maplelatex}
\mapleinline{inert}{2d}{2*v[3], Y*Z}{\[\displaystyle
2\,v_{{3}},\,YZ\]}
\end{maplelatex}
\mapleresult
\begin{maplelatex}
\mapleinline{inert}{2d}{-27*v[4]^2-4*v[3]^3, Z^3}{\[\displaystyle
-27\,{v_{{4}}}^{2}-4\,{v_{{3}}}^{3},\,{Z}^{3}\]}
\end{maplelatex}
\mapleresult
\begin{maplelatex}
\mapleinline{inert}{2d}{3, Y^2}{\[\displaystyle 3,\,{Y}^{2}\]}
\end{maplelatex}
\end{maplegroup}
\begin{maplegroup}
\begin{mapleinput}
\mapleinline{active}{1d}{v[3]:=-3*a^2 : v[4]:=-2*a^3 : FL}{\[\]}
\end{mapleinput}
\mapleresult
\begin{maplelatex}
\mapleinline{inert}{2d}{Y^3+Z^5-3*a^2*Y*Z^2-2*a^3*Z^3+v[5]*Y*Z^3+v[6]*Z^4}{\[\displaystyle
{Y}^{3}+{Z}^{5}-3\,{a}^{2}Y{Z}^{2}-2\,{a}^{3}{Z}^{3}+v_{{5}}Y{Z}^{3}+v_{{6}}{Z}^{4}\]}
\end{maplelatex}
\end{maplegroup}
\begin{maplegroup}
\begin{mapleinput}
\mapleinline{active}{1d}{print(mu = MilnorNumber(FL,\{Y,Z\},plex_min(Z,Y)),
tau = TyurinaNumber(FL,\{Y,Z\},plex_min(Z,Y)))}
{\[{\mu={\it MilnorNumber} \left( {\it FL}, \left\{ Y,Z \right\} ,{\it plex\_min}\\
\mbox{} \left( Z,Y \right)  \right) ,\tau={\it TyurinaNumber} \left( {\it FL}, \left\{ Y,Z \right\} ,{\it plex\_min}\\
\mbox{} \left( Z,Y \right)  \right) }\]}
\end{mapleinput}
\mapleresult
\begin{maplelatex}
\mapleinline{inert}{2d}{mu = 5, tau = 5}{\[\displaystyle
\mu=5,\,\tau=5\]}
\end{maplelatex}
\end{maplegroup}
\begin{maplegroup}
\begin{mapleinput}
\mapleinline{active}{1d}{v[4]:=2*a^3 : FL}{\[\]}
\end{mapleinput}
\mapleresult
\begin{maplelatex}
\mapleinline{inert}{2d}{Y^3+Z^5-3*a^2*Y*Z^2-2*a^3*Z^3+v[5]*Y*Z^3+v[6]*Z^4}{\[\displaystyle
{Y}^{3}+{Z}^{5}-3\,{a}^{2}Y{Z}^{2}+2\,{a}^{3}{Z}^{3}+v_{{5}}Y{Z}^{3}+v_{{6}}{Z}^{4}\]}
\end{maplelatex}
\end{maplegroup}
\begin{maplegroup}
\begin{mapleinput}
\mapleinline{active}{1d}{print(mu = MilnorNumber(FL,\{Y,Z\},plex_min(Z,Y)),
tau = TyurinaNumber(FL,\{Y,Z\},plex_min(Z,Y)))}
{\[{\mu={\it MilnorNumber} \left( {\it FL}, \left\{ Y,Z \right\} ,{\it plex\_min}\\
\mbox{} \left( Z,Y \right)  \right) ,\tau={\it TyurinaNumber} \left( {\it FL}, \left\{ Y,Z \right\} ,{\it plex\_min}\\
\mbox{} \left( Z,Y \right)  \right) }\]}
\end{mapleinput}
\mapleresult
\begin{maplelatex}
\mapleinline{inert}{2d}{mu = 5, tau = 5}{\[\displaystyle
\mu=5,\,\tau=5\]}
\end{maplelatex}
\end{maplegroup}

\halfline

\noindent Let us now go on by considering:

\halfline

\begin{maplegroup}
\begin{mapleinput}
\mapleinline{active}{1d}{MGB5 := MilnorGroebnerBasis(FL, \{Y,Z\},plex_min(Z,Y)):}{\[\]}
\end{mapleinput}
\begin{mapleinput}
\mapleinline{active}{21d}{for i to nops(MGB5) do LeadingTerm(MGB5[i],plex_min(Z,Y))end do}{\[\]}
\end{mapleinput}
\mapleresult
\begin{maplelatex}
\mapleinline{inert}{2d}{3, Y^2}{\[\displaystyle 3,\,{Y}^{2}\]}
\end{maplelatex}
\mapleresult
\begin{maplelatex}
\mapleinline{inert}{2d}{-48*a^4*v[5]+48*a^3*v[6],
Z^4}{\[\displaystyle
-48\,{a}^{4}v_{{5}}+48\,{a}^{3}v_{{6}},\,{Z}^{4}\]}
\end{maplelatex}
\mapleresult
\begin{maplelatex}
\mapleinline{inert}{2d}{-6*a^2, Y*Z}{\[\displaystyle
-6\,{a}^{2},\,YZ\]}
\end{maplelatex}
\end{maplegroup}
\begin{maplegroup}
\begin{mapleinput}
\mapleinline{active}{1d}{a := 0 : FL}{\[\]}
\end{mapleinput}
\mapleresult
\begin{maplelatex}
\mapleinline{inert}{2d}{Y^3+Z^5+v[5]*Y*Z^3+v[6]*Z^4}{\[\displaystyle
{Y}^{3}+{Z}^{5}+v_{{5}}Y{Z}^{3}+v_{{6}}{Z}^{4}\]}
\end{maplelatex}
\end{maplegroup}
\begin{maplegroup}
\begin{mapleinput}
\mapleinline{active}{1d}{print(mu = MilnorNumber(FL,\{Y,Z\},plex_min(Z,Y)),
tau = TyurinaNumber(FL,\{Y,Z\},plex_min(Z,Y)))}
{\[{\mu={\it MilnorNumber} \left( {\it FL}, \left\{ Y,Z \right\} ,{\it plex\_min}\\
\mbox{} \left( Z,Y \right)  \right) ,\tau={\it TyurinaNumber} \left( {\it FL}, \left\{ Y,Z \right\} ,{\it plex\_min}\\
\mbox{} \left( Z,Y \right)  \right) }\]}
\end{mapleinput}
\mapleresult
\begin{maplelatex}
\mapleinline{inert}{2d}{mu = 6, tau = 6}{\[\displaystyle
\mu=6,\,\tau=6\]}
\end{maplelatex}
\end{maplegroup}
\begin{maplegroup}
\begin{mapleinput}
\mapleinline{active}{1d}{MGB6 := MilnorGroebnerBasis(FL,\{Y,Z\},plex_min(Z,Y)):}{\[\]}
\end{mapleinput}
\begin{mapleinput}
\mapleinline{active}{1d}{for i to nops(MGB6) do LeadingTerm(MGB6[i],plex_min(Z, Y))end do}{\[\]}
\end{mapleinput}
\mapleresult
\begin{maplelatex}
\mapleinline{inert}{2d}{3, Y^2}{\[\displaystyle 3,\,{Y}^{2}\]}
\end{maplelatex}
\mapleresult
\begin{maplelatex}
\mapleinline{inert}{2d}{3*v[5], Y*Z^2}{\[\displaystyle
3\,v_{{5}},\,Y{Z}^{2}\]}
\end{maplelatex}
\mapleresult
\begin{maplelatex}
\mapleinline{inert}{2d}{16*v[6]^2, Z^4}{\[\displaystyle
16\,{v_{{6}}}^{2},\,{Z}^{4}\]}
\end{maplelatex}
\end{maplegroup}
\begin{maplegroup}
\begin{mapleinput}
\mapleinline{active}{1d}{v[5] := 0 : FL}{\[\]}
\end{mapleinput}
\mapleresult
\begin{maplelatex}
\mapleinline{inert}{2d}{Y^3+Z^5+v[6]*Z^4}{\[\displaystyle
{Y}^{3}+{Z}^{5}+v_{{6}}{Z}^{4}\]}
\end{maplelatex}
\end{maplegroup}
\begin{maplegroup}
\begin{mapleinput}
\mapleinline{active}{1d}{print(mu = MilnorNumber(FL,\{Y,Z\},plex_min(Z,Y)),
tau = TyurinaNumber(FL,\{Y,Z\},plex_min(Z, Y)))}
{\[{\mu={\it MilnorNumber} \left( {\it FL}, \left\{ Y,Z \right\} ,{\it plex\_min}\\
\mbox{} \left( Z,Y \right)  \right) ,\tau={\it TyurinaNumber} \left( {\it FL}, \left\{ Y,Z \right\} ,{\it plex\_min}\\
\mbox{} \left( Z,Y \right)  \right) }\]}
\end{mapleinput}
\mapleresult
\begin{maplelatex}
\mapleinline{inert}{2d}{mu = 6, tau = 6}{\[\displaystyle
\mu=6,\,\tau=6\]}
\end{maplelatex}
\end{maplegroup}
\begin{maplegroup}
\begin{mapleinput}
\mapleinline{active}{1d}{unassign('v[5]') : v[6]:=0 : FL}{\[\]}
\end{mapleinput}
\mapleresult
\begin{maplelatex}
\mapleinline{inert}{2d}{Y^3+Z^5+v[5]*Y*Z^3}{\[\displaystyle
{Y}^{3}+{Z}^{5}+v_{{5}}Y{Z}^{3}\]}
\end{maplelatex}
\end{maplegroup}
\begin{maplegroup}
\begin{mapleinput}
\mapleinline{active}{1d}{print(mu = MilnorNumber(FL,\{Y,Z\},plex_min(Z,Y)),
tau = TyurinaNumber(FL,\{Y,Z\},plex_min(Z,Y)))}{\[{\mu={\it MilnorNumber} \left( {\it FL}, \left\{ Y,Z \right\} ,{\it plex\_min}\\
\mbox{} \left( Z,Y \right)  \right) ,\tau={\it TyurinaNumber} \left( {\it FL}, \left\{ Y,Z \right\} ,{\it plex\_min}\\
\mbox{} \left( Z,Y \right)  \right) }\]}
\end{mapleinput}
\mapleresult
\begin{maplelatex}
\mapleinline{inert}{2d}{mu = 7, tau = 7}{\[\displaystyle
\mu=7,\,\tau=7\]}
\end{maplelatex}
\end{maplegroup}

\halfline

\noindent meaning that
\begin{itemize}
    \item \emph{$0\in f_{\Lambda}^{-1}(0)$ is a singularity of type $E_6$ for $\Lambda$ generic
    in $\mathcal{V}_0^4$,}
    \item \emph{$0\in f_{\Lambda}^{-1}(0)$ is a singularity of type $E_7$ for $\Lambda$ generic
    in $\mathcal{V}_0^4\cap \mathcal{V}_6$.}
\end{itemize}
On the other hand, by considering the second leading coefficient in MGB5, define the following codimension 2 subset of
$\mathcal{L}$
\begin{equation}\label{D6 in E8}
    \mathcal{V}':=\left\{v_3v_5^2+3v_6^2=v_4v_5^3+2v_6^3=0\right\}\subset\mathcal{V}
\end{equation}
and type

\halfline

\begin{maplegroup}
\begin{mapleinput}
\mapleinline{active}{1d}{unassign('v[6]') : unassign('a') : FL}{\[\]}
\end{mapleinput}
\mapleresult
\begin{maplelatex}
\mapleinline{inert}{2d}{Y^3+Z^5-3*a^2*Y*Z^2-2*a^3*Z^3+v[5]*Y*Z^3+v[6]*Z^4}{\[\displaystyle
{Y}^{3}+{Z}^{5}-3\,{a}^{2}Y{Z}^{2}-2\,{a}^{3}{Z}^{3}+v_{{5}}Y{Z}^{3}+v_{{6}}{Z}^{4}\]}
\end{maplelatex}
\end{maplegroup}
\begin{maplegroup}
\begin{mapleinput}
\mapleinline{active}{1d}{v[6]:=a*v[5] : FL}{\[\]}
\end{mapleinput}
\mapleresult
\begin{maplelatex}
\mapleinline{inert}{2d}{Y^3+Z^5-3*a^2*Y*Z^2-2*a^3*Z^3+v[5]*Y*Z^3+a*v[5]*Z^4}{\[\displaystyle
{Y}^{3}+{Z}^{5}-3\,{a}^{2}Y{Z}^{2}-2\,{a}^{3}{Z}^{3}+v_{{5}}Y{Z}^{3}+av_{{5}}{Z}^{4}\]}
\end{maplelatex}
\end{maplegroup}
\begin{maplegroup}
\begin{mapleinput}
\mapleinline{active}{1d}{print(mu = MilnorNumber(FL,\{Y,Z\},plex_min(Z,Y)),
tau = TyurinaNumber(FL,\{Y,Z\},plex_min(Z,Y)))}{\[{\mu={\it MilnorNumber} \left( {\it FL}, \left\{ Y,Z \right\} ,{\it plex\_min}\\
\mbox{} \left( Z,Y \right)  \right) ,\tau={\it TyurinaNumber} \left( {\it FL}, \left\{ Y,Z \right\} ,{\it plex\_min}\\
\mbox{} \left( Z,Y \right)  \right) }\]}
\end{mapleinput}
\mapleresult
\begin{maplelatex}
\mapleinline{inert}{2d}{mu = 6, tau = 6}{\[\displaystyle
\mu=6,\,\tau=6\]}
\end{maplelatex}
\end{maplegroup}
\begin{maplegroup}
\begin{mapleinput}
\mapleinline{active}{1d}{MGB6b := MilnorGroebnerBasis(FL,\{Y,Z\},plex_min(Z,Y)):}{\[\]}
\end{mapleinput}
\begin{mapleinput}
\mapleinline{active}{1d}{for i to nops(MGB6b) do LeadingTerm(MGB6b[i],plex_min(Z,Y))end do}{\[\]}
\end{mapleinput}
\mapleresult
\begin{maplelatex}
\mapleinline{inert}{2d}{3, Y^2}{\[\displaystyle 3,\,{Y}^{2}\]}
\end{maplelatex}
\mapleresult
\begin{maplelatex}
\mapleinline{inert}{2d}{-6*a^2, Y*Z}{\[\displaystyle
-6\,{a}^{2},\,YZ\]}
\end{maplelatex}
\mapleresult
\begin{maplelatex}
\mapleinline{inert}{2d}{60*a^3*v[5]+5*a^2*v[5]^3,
Z^5}{\[\displaystyle
60\,{a}^{3}v_{{5}}+5\,{a}^{2}{v_{{5}}}^{3},\,{Z}^{5}\]}
\end{maplelatex}
\end{maplegroup}
\begin{maplegroup}
\begin{mapleinput}
\mapleinline{active}{1d}{a:=0 : FL}{\[\]}
\end{mapleinput}
\mapleresult
\begin{maplelatex}
\mapleinline{inert}{2d}{Y^3+Z^5+v[5]*Y*Z^3}{\[\displaystyle
{Y}^{3}+{Z}^{5}+v_{{5}}Y{Z}^{3}\]}
\end{maplelatex}
\end{maplegroup}
\begin{maplegroup}
\begin{mapleinput}
\mapleinline{active}{1d}{unassign('a') : a:=-(1/12)*v[5]^2 : FL}{\[\]}
\end{mapleinput}
\mapleresult
\begin{maplelatex}
\mapleinline{inert}{2d}{Y^3+Z^5-(1/48)*v[5]^4*Y*Z^2+(1/864)*v[5]^6*Z^3+v[5]*Y*Z^3-(1/12)*v[5]^3*Z^4}{\[\displaystyle {Y}^{3}+{Z}^{5}-1/48\,{v_{{5}}}^{4}Y{Z}^{2}+{\frac {1}{864}}\,{v_{{5}}}^{6}{Z}^{3}+v_{{5}}Y{Z}^{3}\\
\mbox{}-1/12\,{v_{{5}}}^{3}{Z}^{4}\]}
\end{maplelatex}
\end{maplegroup}
\begin{maplegroup}
\begin{mapleinput}
\mapleinline{active}{1d}{print(mu = MilnorNumber(FL,\{Y,Z\},plex_min(Z,Y)),
tau = TyurinaNumber(FL,\{Y,Z\},plex_min(Z,Y)))}{\[{\mu={\it MilnorNumber} \left( {\it FL}, \left\{ Y,Z \right\} ,{\it plex\_min}\\
\mbox{} \left( Z,Y \right)  \right) ,\tau={\it TyurinaNumber} \left( {\it FL}, \left\{ Y,Z \right\} ,{\it plex\_min}\\
\mbox{} \left( Z,Y \right)  \right) }\]}
\end{mapleinput}
\mapleresult
\begin{maplelatex}
\mapleinline{inert}{2d}{mu = 7, tau = 7}{\[\displaystyle
\mu=7,\,\tau=7\]}
\end{maplelatex}
\end{maplegroup}

\halfline

\noindent Consider the latter leading coefficient in MGB6b and define \footnote{The reader may check that choosing
    $v_4=2a^3\ ,\ v_6=-av_5\ ,\ 12a=v_5^2$
leads to same results on Milnor and Tyurina numbers and to the same condition defining $\mathcal{V}''$ in (\ref{D6 in E8 bis}).}
the following codimension 3
subset of $\mathcal{L}$
\begin{equation}\label{D6 in E8 bis}
    \mathcal{V}'':=\mathcal{V}'\cap \left\{12v_6+v_5^3=0\right\}\ ,
\end{equation}
Then, by Theorem \ref{classificazione}(3), we get that
\begin{itemize}
    \item \emph{$0\in f_{\Lambda}^{-1}(0)$ is a singularity of type $D_6$ for $\Lambda$ generic
    in $\mathcal{V}'\cap\mathcal{V}_0^2$,}
    \item \emph{$0\in f_{\Lambda}^{-1}(0)$ is a singularity of type $D_7$ for $\Lambda$ generic
    in $\mathcal{V}''\cap\mathcal{V}_0^2$,}
\end{itemize}
while
$\mathcal{V}'\cap\mathcal{V}_0^3=\mathcal{V}_0^4\cap\mathcal{V}_6$
and we get the $E_7$ singularities already discussed. A further
specialization here gives the trivial deformation.

\noindent Let us then come back to consider $\mathcal{W}_2^3$ and
look at the second leading coefficient in MGB3. Define
\begin{equation}\label{A4 in E8}
    \mathcal{W}_4:=\left\{16v_1^5v_2v_5^2-[(v_1v_3+3v_2)^2-4v_1^3v_6]^2=0\right\}
\end{equation}
whose equation is obtained by eliminating $u$ and $t$ from the following set of equations, parameterizing $\mathcal{W}_3^4\setminus\mathcal{V}_1^2$,
\begin{equation}\label{parametrizzazione}
    v_2-u^2v_1=0\ ,\ v_4 - u(v_3+u^2)=0\ ,\ 4v_1(v_6-uv_5)- (v_3+3u^2)^2=0\ .
\end{equation}
Observe that $\mathcal{V}_1^2\subseteq\mathcal{W}_3^4\
\Rightarrow\ \mathcal{V}_0^2\subseteq\mathcal{W}_2^4$ and define
the following codimension 1 Zariski closed subset \footnote{The same considerations explained by footnote \ref{nota} are still holding, here.} of
$\mathcal{W}_2^3$
\begin{equation}\label{A4bis in E8}
    \widetilde{\mathcal{W}}_2^4:=\overline{\mathcal{W}_2^4\setminus\mathcal{V}_0^2}\
    .
\end{equation}
Then \footnote{The reader may check that the choice $v_3 = - 2tw_1- 3u^2$ leads to the same conclusion.}

\halfline

\begin{maplegroup}
\begin{mapleinput}
\mapleinline{active}{1d}{v[6]:=t^2+u*v[5] : v[3]:=2*t*w[1]-3*u^2 : F2}{\[{\it F2}\]}
\end{mapleinput}
\mapleresult
\begin{maplelatex}
\mapleinline{inert}{2d}{} {\[\displaystyle {Y}^{3}+{Z}^{5}+ \left(
w_{{1}}Y+uw_{{1}}Z \right) ^{2}+ \left( 2\,tw_{{1}}-3\,{u}^{2}
\right) Y{Z}^{2}+\]} \mapleinline{inert}{2d}{}
{\[\displaystyle \quad u \left( 2\,tw_{{1}}-2\,{u}^{2} \right) {Z}^{3}\\
\mbox{}+v_{{5}}Y{Z}^{3}+ \left( {t}^{2}+uv_{{5}} \right)
{Z}^{4}\]}
\end{maplelatex}
\end{maplegroup}
\begin{maplegroup}
\begin{mapleinput}
\mapleinline{active}{1d}{print(mu = MilnorNumber(F2,\{Y,Z\},plex_min(Z,Y)),
tau = TYURINANumber(F2,\{Y,Z\},[plex_min(Z,Y), tdeg(Z,Y)]))}{\[{\mu={\it MilnorNumber} \left( {\it F2}, \left\{ Y,Z \right\} ,{\it plex\_min}\\
\mbox{} \left( Z,Y \right)  \right) ,\tau={\it TYURINANumber} \left( {\it F2}, \left\{ Y,Z \right\} ,[{\it plex\_min}\\
\mbox{} \left( Z,Y \right) ,{\it tdeg}\\
\mbox{} \left( Z,Y \right) ] \right) }\]}
\end{mapleinput}
\mapleresult
\begin{maplelatex}
\mapleinline{inert}{2d}{mu = 4, tau = 4}{\[\displaystyle
\mu=4,\,\tau=4\]}
\end{maplelatex}
\end{maplegroup}

\halfline

\noindent and Theorem \ref{classificazione}(2) allows to conclude
that
\begin{itemize}
    \item \emph{$0\in f_{\Lambda}^{-1}(0)$ is a singularity of type $A_4$ for $\Lambda$ generic
    in $\widetilde{\mathcal{W}}_2^4$.}
\end{itemize}
Then type:

\halfline

\begin{maplegroup}
\begin{mapleinput}
\mapleinline{active}{1d}{MGB4b := MilnorGroebnerBasis(F2,\{Y,Z\},plex_min(Z,Y)):}{\[\]}
\end{mapleinput}
\begin{mapleinput}
\mapleinline{active}{1d}{for i to nops(MGB4b) do LeadingTerm(MGB4b[i],plex_min(Z,Y))end do}{\[\]}
\end{mapleinput}
\mapleresult
\begin{maplelatex}
\mapleinline{inert}{2d}{120*t^3*w[1]^4*u+40*t^2*w[1]^5*v[5]-40*t*w[1]^6,
Z^4}{\[\displaystyle
120\,{t}^{3}{w_{{1}}}^{4}u+40\,{t}^{2}{w_{{1}}}^{5}v_{{5}}-40\,t{w_{{1}}}^{6},\,{Z}^{4}\]}
\end{maplelatex}
\mapleresult
\begin{maplelatex}
\mapleinline{inert}{2d}{2*u*w[1]^2, Y}{\[\displaystyle
2\,u{w_{{1}}}^{2},\,Y\]}
\end{maplelatex}
\mapleresult
\begin{maplelatex}
\mapleinline{inert}{2d}{2*w[1]^2, Y}{\[\displaystyle
2\,{w_{{1}}}^{2},\,Y\]}
\end{maplelatex}
\end{maplegroup}
\begin{maplegroup}
\begin{mapleinput}
\mapleinline{active}{1d}{w[1]:=0 : F2}{\[\]}
\end{mapleinput}
\mapleresult
\begin{maplelatex}
\mapleinline{inert}{2d}{Y^3+Z^5-3*u^2*Y*Z^2-2*u^3*Z^3+v[5]*Y*Z^3+(t^2+u*v[5])*Z^4}{\[\displaystyle
{Y}^{3}+{Z}^{5}-3\,{u}^{2}Y{Z}^{2}-2\,{u}^{3}{Z}^{3}+v_{{5}}Y{Z}^{3}+
\left( {t}^{2}+uv_{{5}} \right) {Z}^{4}\]}
\end{maplelatex}
\end{maplegroup}
\begin{maplegroup}
\begin{mapleinput}
\mapleinline{active}{1d}{print(mu = MilnorNumber(F2,\{Y,Z\},plex_min(Z,Y)),
tau = TyurinaNumber(F2,\{Y,Z\},plex_min(Z,Y)))}{\[{\mu={\it MilnorNumber} \left( {\it F2}, \left\{ Y,Z \right\} ,{\it plex\_min}\\
\mbox{} \left( Z,Y \right)  \right) ,\tau={\it TyurinaNumber} \left( {\it F2}, \left\{ Y,Z \right\} ,{\it plex\_min}\\
\mbox{} \left( Z,Y \right)  \right) }\]}
\end{mapleinput}
\mapleresult
\begin{maplelatex}
\mapleinline{inert}{2d}{mu = 5, tau = 5}{\[\displaystyle
\mu=5,\,\tau=5\]}
\end{maplelatex}

\halfline

\noindent This means that $\mathcal{V}_1\cap
\widetilde{\mathcal{W}}_2^4=\mathcal{V}\cap \mathcal{V}_0^2$
obtaining a generic $D_5$ singularity, as already described above.
Go on:

\halfline

\end{maplegroup}
\begin{maplegroup}
\begin{mapleinput}
\mapleinline{active}{1d}{unassign('w[1]') : t:=0 : F2}{\[\]}
\end{mapleinput}
\mapleresult
\begin{maplelatex}
\mapleinline{inert}{2d}{Y^3+Z^5+(w[1]*Y+u*w[1]*Z)^2-3*u^2*Y*Z^2-2*u^3*Z^3+v[5]*Y*Z^3+u*v[5]*Z^4}{\[\displaystyle
{Y}^{3}+{Z}^{5}+ \left( w_{{1}}Y+uw_{{1}}Z \right)
^{2}-3\,{u}^{2}Y{Z}^{2}-2\,{u}^{3}{Z}^{3}+v_{{5}}Y{Z}^{3}+uv_{{5}}{Z}^{4}\]}
\end{maplelatex}
\end{maplegroup}
\begin{maplegroup}
\begin{mapleinput}
\mapleinline{active}{1d}{print(mu = MilnorNumber(F2,\{Y,Z\},plex_min(Z,Y)),
tau = TYURINANumber(F2,\{Y,Z\},[plex_min(Z,Y),tdeg(Z,Y)]))}{\[{\mu={\it MilnorNumber} \left( {\it F2}, \left\{ Y,Z \right\} ,{\it plex\_min}\\
\mbox{} \left( Z,Y \right)  \right) ,\tau={\it TYURINANumber} \left( {\it F2}, \left\{ Y,Z \right\} ,[{\it plex\_min}\\
\mbox{} \left( Z,Y \right) ,{\it tdeg}\\
\mbox{} \left( Z,Y \right) ] \right) }\]}
\end{mapleinput}
\mapleresult
\begin{maplelatex}
\mapleinline{inert}{2d}{mu = 4, tau = 4}{\[\displaystyle
\mu=4,\,\tau=4\]}
\end{maplelatex}
\end{maplegroup}
\begin{maplegroup}
\begin{mapleinput}
\mapleinline{active}{1d}{unassign('t') : u:=0 : F2}{\[\]}
\end{mapleinput}
\end{maplegroup}
\begin{maplegroup}
\begin{mapleinput}
\mapleinline{active}{1d}{print(mu = MilnorNumber(F2,\{Y,Z\},plex_min(Z,Y)),
tau = TYURINANumber(F2,\{Y,Z\},[plex_min(Z,Y),tdeg(Z,Y)]))}{\[{\mu={\it MilnorNumber} \left( {\it F2}, \left\{ Y,Z \right\} ,{\it plex\_min}\\
\mbox{} \left( Z,Y \right)  \right) ,\tau={\it TYURINANumber} \left( {\it F2}, \left\{ Y,Z \right\} ,[{\it plex\_min}\\
\mbox{} \left( Z,Y \right) ,{\it tdeg}\\
\mbox{} \left( Z,Y \right) ] \right) }\]}
\end{mapleinput}
\mapleresult
\begin{maplelatex}
\mapleinline{inert}{2d}{mu = 4, tau = 4}{\[\displaystyle
\mu=4,\,\tau=4\]}
\end{maplelatex}
\end{maplegroup}

\halfline

\noindent Then we have to consider the relation assigned by the first
leading coefficient in MGB4b, giving the further condition
\[
    v_1-w_1v_5 t - 3ut^2=0\ ,
\]
which added to (\ref{parametrizzazione}) parameterizes the following (reducible) codimension 4 Zariski closed subset of $\mathcal{L}$
\begin{equation}\label{A5 in E8}
    \widetilde{\mathcal{W}}_2^5:=\overline{\mathcal{W}_2^5\setminus\mathcal{V}_0^2}
\end{equation}
where \footnote{The following equation of $\mathcal{W}_5$ is obtained by carefully applying Maple's commands \texttt{eliminate} and \texttt{EliminationIdeal} in the \texttt{PolynomialIdeals} package.}
\begin{eqnarray}\label{A5 in E8bis}
 \nonumber
    &\mathcal{W}_5:=\left\{[4v_1^5v_5^2(v_1v_3+3v_2)^2+ 16 v_1^7(v_1^2+v_1v_3v_5+3v_2v_5)-9v_2(v_1v_3+3v_2)^4]\cdot \right.  &\\
    &\quad\quad\left. [4v_1^5v_5^2(v_1v_3+3v_2)^2+ 16 v_1^7(v_1^2-v_1v_3v_5-3v_2v_5)-9v_2(v_1v_3+3v_2)^4]=0\right\}&
\end{eqnarray}
Since

\halfline

\begin{maplegroup}
\begin{mapleinput}
\mapleinline{active}{1d}{u:=1/12*(b^2-v[5]^2) : w[1]:=(1/2)*t*(v[5]-b) : F2}{\[\]}
\end{mapleinput}
\mapleresult
\begin{maplelatex}
\mapleinline{inert}{2d}{}
{\[\displaystyle {Y}^{3}+{Z}^{5}+ \left( 1/2\,t \left( v_{{5}}-b \right) Y+1/2\, \left( 1/12\,{b}^{2}-1/12\,{v_{{5}}}^{2} \right) t \left( v_{{5}}-b \right) Z\right) ^{2}+\]}
\mapleinline{inert}{2d}{}
{\[\displaystyle \left( {t}^{2} \left( v_{{5}}-b \right) -3\, \left( 1/12\,{b}^{2}-1/12\,{v_{{5}}}^{2} \right) ^{2} \right) Y{Z}^{2}+\]}
\mapleinline{inert}{2d}{}
{\[\displaystyle \left( 1/12\,{b}^{2}-1/12\,{v_{{5}}}^{2} \right)  \left( {t}^{2} \left( v_{{5}}-b \right) -2\, \left( 1/12\,{b}^{2}-1/12\,{v_{{5}}}^{2} \right) ^{2} \right) {Z}^{3}+v_{{5}}Y{Z}^{3}+\]}
\mapleinline{inert}{2d}{}
{\[\displaystyle \left( {t}^{2}+ \left( 1/12\,{b}^{2}-1/12\,{v_{{5}}}^{2} \right) v_{{5}} \right) {Z}^{4}\]}
\end{maplelatex}
\end{maplegroup}
\begin{maplegroup}
\begin{mapleinput}
\mapleinline{active}{1d}{M5 := Milnor(F2,\{Y,Z\},plex_min(Z, Y)) : print(mu = M5[4])}{\[\]}
\end{mapleinput}
\mapleresult
\begin{maplelatex}
\mapleinline{inert}{2d}{mu = 5}{\[\displaystyle \mu=5\]}
\end{maplelatex}
\end{maplegroup}

\halfline

\noindent then Theorem \ref{classificazione} (2) gives that \footnote{The reader may check that setting $w_1:=t(v_5+b)/2$ leads to the same conclusion.}
\begin{itemize}
  \item \emph{$0\in f_{\Lambda}^{-1}(0)$ is a singularity of type $A_5$ for $\Lambda$ generic
    in $\widetilde{\mathcal{W}}_2^5$.}
\end{itemize}
Furthermore:

\halfline

\begin{maplegroup}
\begin{mapleinput}
\mapleinline{active}{1d}{MGB5:=M5[1] : for i to nops(MGB5) do
factor(LeadingTerm(MGB5[i],plex_min(Z,Y)))end do}{\[\]}
\end{mapleinput}
\mapleresult
\begin{maplelatex}
\mapleinline{inert}{2d}{}{\[\displaystyle {\frac {1594323}{576460752303423488}}\,{t}^{6} \left( b-v_{{5}} \right) ^{3}\\
\mbox{} \left( {b}^{3}-{b}^{2}v_{{5}}-8\,{t}^{2} \right) ,\,Z^5 \]}
\end{maplelatex}
\mapleresult
\begin{maplelatex}
\mapleinline{inert}{2d}{(1/2)*t^2*(v[5]-b)^2*((1/12)*b^2-(1/12)*v[5]^2), Y}{\[\displaystyle 1/2\,{t}^{2} \left( v_{{5}}-b \right) ^{2} \left( 1/12\,{b}^{2}-1/12\,{v_{{5}}}^{2} \right) ,\,Y\]}
\end{maplelatex}
\mapleresult
\begin{maplelatex}
\mapleinline{inert}{2d}{(1/2)*t^2*(v[5]-b)^2, Y}{\[\displaystyle 1/2\,{t}^{2} \left( v_{{5}}-b \right) ^{2},\,Y\]}
\end{maplelatex}
\end{maplegroup}
\begin{maplegroup}
\begin{mapleinput}
\mapleinline{active}{1d}{t := 0 : F2}{\[\]}
\end{mapleinput}
\mapleresult
\begin{maplelatex}
\mapleinline{inert}{2d}{}
{\[\displaystyle {Y}^{3}+{Z}^{5}-3\, \left( 1/12\,{b}^{2}-1/12\,{v_{{5}}}^{2} \right) ^{2}Y{Z}^{2}-2\, \left( 1/12\,{b}^{2}-1/12\,{v_{{5}}}^{2} \right) ^{3}{Z}^{3}+\]}
\mapleinline{inert}{2d}{}
{\[\displaystyle \quad v_{{5}}Y{Z}^{3}+ \left( 1/12\,{b}^{2}-1/12\,{v_{{5}}}^{2} \right) v_{{5}}{Z}^{4}\]}
\end{maplelatex}
\end{maplegroup}
\begin{maplegroup}
\begin{mapleinput}
\mapleinline{active}{1d}{M5b := Milnor(F2,\{Y,Z\},plex_min(Z,Y)) : print(mu = M5b[4])}{\[\]}
\end{mapleinput}
\mapleresult
\begin{maplelatex}
\mapleinline{inert}{2d}{mu = 6}{\[\displaystyle \mu=6\]}
\end{maplelatex}
\end{maplegroup}
\begin{maplegroup}
\begin{mapleinput}
\mapleinline{active}{1d}{T5b := TYURINA(F2,\{Y,Z\},[plex_min(Z,Y), tdeg(Z,Y)]) : print(tau = T5b[4])}{\[\]}
\end{mapleinput}
\mapleresult
\begin{maplelatex}
\mapleinline{inert}{2d}{tau = 6}{\[\displaystyle \tau=6\]}
\end{maplelatex}
\end{maplegroup}

\halfline

\noindent Observe that, by identifying $a=(b^2-v_5^2)/12$, we are precisely dealing with deformations parameterized by the above considered algebraic closed subset $\mathcal{V}'\cap\mathcal{V}_0^2$ of $\mathcal{L}$, whose generic fibre has a singularity of type $D_6$. On the other hand

\halfline

\begin{maplegroup}
\begin{mapleinput}
\mapleinline{active}{1d}{unassign('t') : v[5]:= -b : F2}{\[\]}
\end{mapleinput}
\mapleresult
\begin{maplelatex}
\mapleinline{inert}{2d}{Y^3+Z^5+t^2*b^2*Y^2-2*t^2*b*Y*Z^2-b*Y*Z^3+t^2*Z^4}{\[\displaystyle {Y}^{3}+{Z}^{5}+{t}^{2}{b}^{2}{Y}^{2}-2\,{t}^{2}bY{Z}^{2}-bY{Z}^{3}+{t}^{2}{Z}^{4}\]}
\end{maplelatex}
\end{maplegroup}
\begin{maplegroup}
\begin{mapleinput}
\mapleinline{active}{1d}{M5c := Milnor(F2, {Y, Z}, plex_min(Z, Y)) : print(mu = M5c[4])}{\[\]}
\end{mapleinput}
\mapleresult
\begin{maplelatex}
\mapleinline{inert}{2d}{mu = 5}{\[\displaystyle \mu=5\]}
\end{maplelatex}
\end{maplegroup}
\begin{maplegroup}
\begin{mapleinput}
\mapleinline{active}{1d}{v[5] := b : F2}{\[\]}
\end{mapleinput}
\mapleresult
\begin{maplelatex}
\mapleinline{inert}{2d}{Y^3+Z^5+b*Y*Z^3+t^2*Z^4}{\[\displaystyle {Y}^{3}+{Z}^{5}+bY{Z}^{3}+{t}^{2}{Z}^{4}\]}
\end{maplelatex}
\end{maplegroup}

\halfline

\noindent still giving the above considered deformations parameterized by $\mathcal{V}_0^4$ whose generic fibre has a singularity of type $E_6$. Let us then consider the relation given by the last factor in the first leading coefficient of MGB5:

\halfline

\begin{maplegroup}
\begin{mapleinput}
\mapleinline{active}{1d}{t:=c*b : v[5]:=b-8*c^2 : F2}{\[\]}
\end{mapleinput}
\mapleresult
\begin{maplelatex}
\mapleinline{inert}{2d}{}
{\[\displaystyle {Y}^{3}+{Z}^{5}+ \left( -4\,{c}^{3}bY-4\, \left( 1/12\,{b}^{2}-1/12\, \left( b-8\,{c}^{2} \right) ^{2} \right) {c}^{3}bZ \right) ^{2} +\]}
\mapleinline{inert}{2d}{}
{\[\displaystyle \left( -8\,{c}^{4}{b}^{2}-3\, \left( 1/12\,{b}^{2}-1/12\, \left( b-8\,{c}^{2} \right) ^{2} \right) ^{2} \right) Y{Z}^{2}+ \]}
\mapleinline{inert}{2d}{}
{\[\displaystyle \left( 1/12\,{b}^{2}-1/12\, \left( b-8\,{c}^{2} \right) ^{2} \right)  \left( -8\,{c}^{4}{b}^{2}-2\, \left( 1/12\,{b}^{2}-1/12\, \left( b-8\,{c}^{2} \right) ^{2} \right) ^{2} \right) {Z}^{3}+ \]}
\mapleinline{inert}{2d}{}
{\[\displaystyle \left( b-8\,{c}^{2} \right) Y{Z}^{3}+ \left( {c}^{2}{b}^{2}+ \left( 1/12\,{b}^{2}-1/12\, \left( b-8\,{c}^{2} \right) ^{2} \right)  \left( b-8\,{c}^{2} \right)  \right) {Z}^{4}\]}
\end{maplelatex}
\end{maplegroup}
\begin{maplegroup}
\begin{mapleinput}
\mapleinline{active}{1d}{M5d := Milnor(F2, \{Y,Z\}, plex_min(Z,Y)) : print(mu = M5d[4])}{\[\]}
\end{mapleinput}
\mapleresult
\begin{maplelatex}
\mapleinline{inert}{2d}{mu = 6}{\[\displaystyle \mu=6\]}
\end{maplelatex}
\end{maplegroup}
\begin{maplegroup}
\begin{mapleinput}
\mapleinline{active}{1d}{T5d := TYURINA(F2, \{Y,Z\}, [plex_min(Z,Y), tdeg(Z,Y)]) :
print(tau = T5d[4])}{\[\]}
\end{mapleinput}
\mapleresult
\begin{maplelatex}
\mapleinline{inert}{2d}{tau = 6}{\[\displaystyle \tau=6\]}
\end{maplelatex}
\end{maplegroup}
\begin{maplegroup}
\begin{mapleinput}
\mapleinline{active}{1d}{MGB6a := M5d[1] : for i to nops(MGB6a) do
LeadingTerm(MGB6a[i], plex_min(Z, Y)) end do:}{\[\]}
\end{mapleinput}
\end{maplegroup}
\begin{maplegroup}
\begin{mapleinput}
\mapleinline{active}{1d}{for i to nops(MGB6a) do
factor(LeadingCoefficient(MGB6a[i], plex_min(Z, Y))) end do}{\[\]}
\end{mapleinput}
\mapleresult
\begin{maplelatex}
\mapleinline{inert}{2d}{32*c^6*b^2}{\[\displaystyle 32\,{c}^{6}{b}^{2}\]}
\end{maplelatex}
\mapleresult
\begin{maplelatex}
\mapleinline{inert}{2d}{-3670016*c^16*b^8*(b+8*c^2)*(b+16*c^2)^2}{\[\displaystyle -3670016\,{c}^{16}{b}^{8} \left( b+8\,{c}^{2} \right)  \left( b+16\,{c}^{2} \right) ^{2}\]}
\end{maplelatex}
\mapleresult
\begin{maplelatex}
\mapleinline{inert}{2d}{(128/3)*c^8*b^2*(-4*c^2+b)}{\[\displaystyle {\frac {128}{3}}\,{c}^{8}{b}^{2} \left( -4\,{c}^{2}+b \right) \]}
\end{maplelatex}
\end{maplegroup}
\begin{maplegroup}
\begin{mapleinput}
\mapleinline{active}{1d}{c := 0 : F2}{\[\]}
\end{mapleinput}
\mapleresult
\begin{maplelatex}
\mapleinline{inert}{2d}{Y^3+Z^5+b*Y*Z^3}{\[\displaystyle {Y}^{3}+{Z}^{5}+bY{Z}^{3}\]}
\end{maplelatex}
\end{maplegroup}
\begin{maplegroup}
\begin{mapleinput}
\mapleinline{active}{1d}{unassign('c') : b:=0 : F2}{\[\]}
\end{mapleinput}
\mapleresult
\begin{maplelatex}
\mapleinline{inert}{2d}{Y^3+Z^5-(256/3)*Y*Z^2*c^8+(8192/27)*Z^3*c^12-8*Y*Z^3*c^2+(128/3)*Z^4*c^6}{\[\displaystyle {Y}^{3}+{Z}^{5}-{\frac {256}{3}}\,Y{Z}^{2}{c}^{8}+{\frac {8192}{27}}\,{Z}^{3}{c}^{12}\\
\mbox{}-8\,Y{Z}^{3}{c}^{2}+{\frac {128}{3}}\,{Z}^{4}{c}^{6}\]}
\end{maplelatex}
\end{maplegroup}
\begin{maplegroup}
\begin{mapleinput}
\mapleinline{active}{1d}{M6a := Milnor(F2, \{Y,Z\}, plex_min(Z,Y)) : print(mu = M6a[4])}{\[\]}
\end{mapleinput}
\mapleresult
\begin{maplelatex}
\mapleinline{inert}{2d}{mu = 7}{\[\displaystyle \mu=7\]}
\end{maplelatex}
\end{maplegroup}
\begin{maplegroup}
\begin{mapleinput}
\mapleinline{active}{1d}{T6a := TYURINA(F2, \{Y,Z\}, [plex_min(Z,Y), tdeg(Z,Y)]) :
 print(tau = T6a[4]) :}{\[\]}
\end{mapleinput}
\mapleresult
\begin{maplelatex}
\mapleinline{inert}{2d}{tau = 7}{\[\displaystyle \tau=7\]}
\end{maplelatex}
\end{maplegroup}
\begin{maplegroup}
\begin{mapleinput}
\mapleinline{active}{1d}{b := 4*c^2 : F2}{\[\]}
\end{mapleinput}
\mapleresult
\begin{maplelatex}
\mapleinline{inert}{2d}{}{\[\displaystyle {Y}^{3}+{Z}^{5}+256\,{Y}^{2}{c}^{10}-128\,Y{Z}^{2}{c}^{8}-4\,Y{Z}^{3}{c}^{2}+16\,{Z}^{4}{c}^{6}\]}
\end{maplelatex}
\end{maplegroup}
\begin{maplegroup}
\begin{mapleinput}
\mapleinline{active}{1d}{M6b := Milnor(F2, \{Y,Z\}, plex_min(Z, Y)) : print(mu = M6b[4])}{\[\]}
\end{mapleinput}
\mapleresult
\begin{maplelatex}
\mapleinline{inert}{2d}{mu = 6}{\[\displaystyle \mu=6\]}
\end{maplelatex}
\end{maplegroup}
\begin{maplegroup}
\begin{mapleinput}
\mapleinline{active}{1d}{T6b := TYURINA(F2, \{Y,Z\}, [plex_min(Z, Y), tdeg(Z, Y)]) :
print(tau = T6b[4])}{\[\]}
\end{mapleinput}
\mapleresult
\begin{maplelatex}
\mapleinline{inert}{2d}{tau = 6}{\[\displaystyle \tau=6\]}
\end{maplelatex}
\end{maplegroup}
\begin{maplegroup}
\begin{mapleinput}
\mapleinline{active}{1d}{b := -8*c^2 : F2}{\[\]}
\end{mapleinput}
\mapleresult
\begin{maplelatex}
\mapleinline{inert}{2d}{}
{\[\displaystyle {Y}^{3}+{Z}^{5}+ \left( 32\,{c}^{5}Y-512\,{c}^{9}Z \right) ^{2}-1280\,Y{Z}^{2}{c}^{8}+16384\,{Z}^{3}{c}^{12}-16\,Y{Z}^{3}{c}^{2}+320\,{Z}^{4}{c}^{6}\]}
\end{maplelatex}
\end{maplegroup}
\begin{maplegroup}
\begin{mapleinput}
\mapleinline{active}{1d}{M6c := Milnor(F2, \{Y,Z\}, plex_min(Z, Y)) : print(mu = M6c[4])}{\[\]}
\end{mapleinput}
\mapleresult
\begin{maplelatex}
\mapleinline{inert}{2d}{mu = 7}{\[\displaystyle \mu=7\]}
\end{maplelatex}
\end{maplegroup}
\begin{maplegroup}
\begin{mapleinput}
\mapleinline{active}{1d}{T6c := TYURINA(F2, \{Y,Z\}, [plex_min(Z, Y), tdeg(Z, Y)]) : print(tau = T6c[4])}{\[\]}
\end{mapleinput}
\mapleresult
\begin{maplelatex}
\mapleinline{inert}{2d}{tau = 7}{\[\displaystyle \tau=7\]}
\end{maplelatex}
\end{maplegroup}
\begin{maplegroup}
\begin{mapleinput}
\mapleinline{active}{1d}{b := -16*c^2 : F2}{\[\]}
\end{mapleinput}
\mapleresult
\begin{maplelatex}
\mapleinline{inert}{2d}{}
{\[\displaystyle {Y}^{3}+{Z}^{5}+ \left( 64\,{c}^{5}Y-{\frac {5120}{3}}\,{c}^{9}Z \right) ^{2}-{\frac {12544}{3}}\,Y{Z}^{2}{c}^{8} +\]}
\mapleinline{inert}{2d}{} {\[\displaystyle{\frac {2498560}{27}}\,{Z}^{3}{c}^{12}-24\,Y{Z}^{3}{c}^{2}+896\,{Z}^{4}{c}^{6}\]}
\end{maplelatex}
\end{maplegroup}
\begin{maplegroup}
\begin{mapleinput}
\mapleinline{active}{1d}{M6d := Milnor(F2, \{Y,Z\}, plex_min(Z, Y)) : print(mu = M6d[4])}{\[\]}
\end{mapleinput}
\mapleresult
\begin{maplelatex}
\mapleinline{inert}{2d}{mu = 6}{\[\displaystyle \mu=6\]}
\end{maplelatex}
\end{maplegroup}
\begin{maplegroup}
\begin{mapleinput}
\mapleinline{active}{1d}{T6d := TYURINA(F2, \{Y,Z\}, [plex_min(Z, Y), tdeg(Z, Y)]) :
print(tau = T6d[4])}{\[\]}
\end{mapleinput}
\mapleresult
\begin{maplelatex}
\mapleinline{inert}{2d}{tau = 6}{\[\displaystyle \tau=6\]}
\end{maplelatex}
\end{maplegroup}

\halfline

\noindent This means that:
\begin{itemize}
  \item equation ${b}^{3}-{b}^{2}v_{{5}}-8{t}^{2} = 0$ defines the following codimension 5 (reducible) algebraic subset of $\mathcal{L}$
\begin{equation}\label{A6 in E8}
    \widetilde{\mathcal{W}}_2^6:=\mathcal{W}_6\cap\widetilde{\mathcal{W}}_2^5\ ,
\end{equation}
where \footnote{The following equation of $\mathcal{W}_6$ is obtained by carefully applying Maple's commands \texttt{eliminate} and \texttt{EliminationIdeal} in the \texttt{PolynomialIdeals} package.}
\begin{eqnarray}\label{A6 in E8bis}
\nonumber
    &\mathcal{W}_6:=\left\{[32v_1^9-2v_1^5v_5(v_1v_3+3v_2)(8v_1^2-3v_2v_5-v_1v_3v_5)+(v_1v_3+3v_2)^5]\cdot\right.&\\
    &\quad\quad\left.[32v_1^9+2v_1^5v_5(v_1v_3+3v_2)(8v_1^2+3v_2v_5+v_1v_3v_5)-(v_1v_3+3v_2)^5]=0\right\}&
\end{eqnarray}
\emph{whose generic point $\Lambda$ is such that $0\in f_{\Lambda}^{-1}(0)$ is an $A_6$ singularity, }
  \item setting $c=0$ means studying deformations parameterized by $\mathcal{V}_0^4\cap\mathcal{V}_6$, generically admitting an $E_7$ singularity, as already observed above,
  \item and setting $b=0$ means studying deformations parameterized by $\mathcal{V}_0^2\cap\mathcal{V}''$, generically admitting a $D_7$ singularity and already considered above, too,
  \item the further step is then to consider the relation $b+8c^2=0$ defining the following codimension 6 (reducible) algebraic subset of $\mathcal{L}$
\begin{equation}\label{A7 in E8}
    \widetilde{\mathcal{W}}_2^7:=\mathcal{W}_7\cap\widetilde{\mathcal{W}}_2^6\quad\text{with}\quad
    \mathcal{W}_7:=\left\{ 256v_2-v_1v_5^4=0\right\}
\end{equation}
\emph{whose generic point $\Lambda$ is such that $0\in f_{\Lambda}^{-1}(0)$ is an $A_7$ singularity.}
\end{itemize}
The check of relations (\ref{intersezioni in E8}) are then left to the reader.
\end{proof}

\subsection{A list of very special adjacencies}\label{speciali}

As a consequence of the analysis performed in the previous sections, we are now able to concretely write down some very special small 1-parameter deformations of a $A_n$ , $D_n$ , $E_6$ , $E_7$ or $E_8$, realizing adjacencies not directly mentioned in \cite{Arnol'd} and \cite{Arnol'd&c} (except for those in \ref{adj:Dn-An-1}). The 1-parameter deformations we are going to list in the following are obtained by last steps in proofs of Theorems \ref{An}, \ref{Dn}, \ref{E6}, \ref{E7} and \ref{E8}, giving precisely 1-parameter deformations, after some possible parameter's re-scaling.

\subsubsection{$A_{n-1}\longleftarrow D_n$}\label{adj:Dn-An-1}

Assume that $n\geq 4$ is either $n=2m+4$ when $n$ is even, or $n=2m+5$ when odd. Then consider the 1-parameter family $\mathcal{X}_t:=\{f_t(\mathbf{x})=0\}$ , $t\in \C$, where
\[
   f_0(\mathbf{x}):= \sum_{i=1}^{N-1} x_i^2\ +\
    x_N^2\ x_{N+1}+ x_{N+1}^{n-1}
\]
and either
\[
   f_t(\mathbf{x}):=  f_0(\mathbf{x})+t(x_N+it^m x_{N+1})^2+\sum_{k=3}^{2m+2}(-t)^{2m+3-k}x_{N+1}^k \quad\text{($n$ even)}
\]
or
\[
   f_t(\mathbf{x}):= f_0(\mathbf{x}) +t^2(x_N+t^{2m+1} x_{N+1})^2+\sum_{k=3}^{2m+3}(-t^2)^{2m+4-k}x_{N+1}^k\quad\text{($n$ odd)}.
\]
Then $\mathcal{X}_0$ is an isolated $D_n$ point and, for generic $t\in\C$, $\mathcal{X}_t$ admits the unique singular point $\mathbf{0}\in f_t^{-1}(0)$ which is of type $A_{n-1}$.

\subsubsection{$A_5\longleftarrow E_6$}

Consider the 1-parameter family $\mathcal{X}_t:=\{f_t(\mathbf{x})=0\}$ , $t\in \C$, with
\begin{eqnarray*}
\nonumber
  f_0(\mathbf{x}) &:=& \sum_{i=1}^{N-1} x_i^2\ +\
    x_N^3+ x_{N+1}^4\\
\nonumber
  f_t(\mathbf{x})&:=& f_0(\mathbf{x}) + t^2 x_N^2 + 2t x_N x_{N+1}^2 \ .
\end{eqnarray*}
Then $\mathcal{X}_0$ is an isolated $E_6$ singular point and, for generic $t\in\C$, $\mathcal{X}_t$ admits the unique singular point $\mathbf{0}\in f_t^{-1}(0)$ which is of type $A_{5}$.

\subsubsection{$D_5\longleftarrow E_6$}

Let $f_0(\mathbf{x})$ be as in the previous case and assume
\[
  f_t(\mathbf{x}):= f_0(\mathbf{x}) - 3t^2 x_Nx_{N+1}^2 - 2t^3 x_{N+1}^3 \ .
\]
Then, for generic $t\in\C$, $\mathcal{X}_t$ admits the unique singular point $\mathbf{0}\in f_t^{-1}(0)$ which is of type $D_{5}$.

\subsubsection{$A_6\longleftarrow E_7$}

Consider the 1-parameter family $\mathcal{X}_t:=\{f_t(\mathbf{x})=0\}$ , $t\in \C$, with
\[
  f_0(\mathbf{x}) := \sum_{i=1}^{N-1} x_i^2\ +\
    x_N^3+ x_Nx_{N+1}^3
\]
\[
  f_t(\mathbf{x}):= f_0(\mathbf{x}) +432t^3(x_N + 4t x_{N+1})^2-120t^2x_N x_{N+1}^2-416t^3x_{N+1}^3+7tx_{N+1}^4
\]
Then $\mathcal{X}_0$ is an isolated $E_7$ singular point and, for generic $t\in\C$, $\mathcal{X}_t$ admits the unique singular point $\mathbf{0}\in f_t^{-1}(0)$ which is of type $A_{6}$.

\subsubsection{$D_6\longleftarrow E_7$}

Let $f_0(\mathbf{x})$ be as in the previous case and assume
\[
  f_t(\mathbf{x}):= f_0(\mathbf{x}) -3t^2x_N x_{N+1}^2-2t^3x_{N+1}^3+tx_{N+1}^4 \ .
\]
Then, for generic $t\in\C$, $\mathcal{X}_t$ admits the unique singular point $\mathbf{0}\in f_t^{-1}(0)$ which is of type $D_{6}$.

\subsubsection{$A_7\longleftarrow E_8$}

Consider the 1-parameter family $\mathcal{X}_t:=\{f_t(\mathbf{x})=0\}$ , $t\in \C$, with
\[
  f_0(\mathbf{x}) := \sum_{i=1}^{N-1} x_i^2\ +\
    x_N^3+ x_{N+1}^5
\]
\[
  f_t(\mathbf{x}):= f_0(\mathbf{x}) +t^{5}(x_N-t^2x_{N+1})^2-5t^4x_N x_{N+1}^2+4t^{6}x_{N+1}^3-4tx_N x_{N+1}^3+5t^3 x_{N+1}^4
\]
Then $\mathcal{X}_0$ is an isolated $E_8$ singular point and, for generic $t\in\C$, $\mathcal{X}_t$ admits the unique singular point $\mathbf{0}\in f_t^{-1}(0)$ which is of type $A_{7}$.

\subsubsection{$D_7\longleftarrow E_8$}

Let $f_0(\mathbf{x})$ be as in the previous case and assume
\[
  f_t(\mathbf{x}):= f_0(\mathbf{x}) -27t^4x_N x_{N+1}^2+54t^{6}x_{N+1}^3-6tx_N x_{N+1}^3 +18t^3x_{N+1}^4\ .
\]
Then, for generic $t\in\C$, $\mathcal{X}_t$ admits the unique singular point $\mathbf{0}\in f_t^{-1}(0)$ which is of type $D_{7}$.

\end{document}